\documentclass[10pt,reqno]{amsbook}

\usepackage[T1]{fontenc}
\usepackage[utf8x]{inputenc}
\usepackage[english]{babel}
\usepackage{amsmath,amsfonts,amssymb,amsthm,bm,shuffle}
\usepackage[math]{anttor}

\usepackage[top=3.7cm,bottom=3.7cm,left=3.7cm,right=3.7cm]{geometry}

\usepackage{xcolor}

\definecolor{ColBlack}{RGB}{0,0,0} 
\definecolor{ColWhite}{RGB}{255,255,255} 
\definecolor{Col1}{RGB}{133,6,6} 
\definecolor{Col2}{RGB}{198,8,0} 
\definecolor{Col3}{RGB}{174,74,52} 
\definecolor{Col4}{RGB}{103,113,121} 
\definecolor{Col5}{RGB}{90,94,107} 
\definecolor{Col6}{RGB}{70,63,50} 

\newcommand{\ColA}[1]{\textcolor{Col1}{#1}}
\newcommand{\ColB}[1]{\textcolor{Col2}{#1}}

\newcommand{\ColD}[1]{\textcolor{Col4}{#1}}

\newcommand{\ColF}[1]{\textcolor{Col6}{#1}}

\usepackage[hyperindex=true,frenchlinks=true,colorlinks=true,
citecolor=Col2,linkcolor=Col3,urlcolor=Col4,linktocpage,
pagebackref=true]{hyperref}

\usepackage{tikz}
\usetikzlibrary{patterns}

\usepackage{dsfont}
\usepackage{stmaryrd}
\usepackage{cite}
\usepackage{xspace}
\usepackage{enumitem}
\usepackage{multirow}
\usepackage{multicol}
\usepackage{chngcntr}

\tikzstyle{Centering}=[{baseline={([yshift=-0.5ex]current
    bounding box.center)}}]

\tikzstyle{Node}=[circle,draw=Col1!80,fill=Col1!8,inner sep=1pt,
minimum size=2mm,thick,font=\scriptsize]
\tikzstyle{Edge}=[draw=Col2!80,cap=round,thick,rounded corners=2.5pt]
\tikzstyle{Leaf}=[rectangle,draw=ColBlack!70,fill=ColBlack!16,
inner sep=0pt,minimum size=1mm,thick]
\tikzstyle{NodeClear}=[Node,fill=ColWhite!100]
\tikzstyle{NodeColorA}=[Node,draw=Col4!80,fill=Col4!8]
\tikzstyle{NodeColorB}=[Node,draw=Col6!80,fill=Col6!8]
\tikzstyle{NodeColorC}=[Node,draw=Col2!90,fill=Col2!20]
\tikzstyle{NodeColorD}=[Node,draw=Col5!90,fill=Col5!20]
\tikzstyle{NodeST}=[font=\footnotesize]
\tikzstyle{EdgeLabel}=[midway,inner sep=1pt,fill=ColWhite!0,
font=\scriptsize]
\tikzstyle{LeafLabel}=[font=\scriptsize,node distance=2mm]
\tikzstyle{EdgeValue}=[regular polygon,regular polygon sides=6,
draw=Col2!80,fill=Col2!2,inner sep=.2pt,line width=.75pt,
minimum size=2.8mm,midway,text=ColBlack,font=\tiny]
\tikzstyle{Subtree}=[regular polygon,regular polygon sides=3,
draw=Col1!80,fill=Col1!20,thick,minimum size=5mm,font=\scriptsize]

\tikzstyle{PathNode}=[circle,draw=Col1!90,fill=Col1!30,thick,
inner sep=0pt,minimum size=2.0mm]
\tikzstyle{PathNodeColorA}=[PathNode,draw=Col4!80,fill=Col4!18]
\tikzstyle{PathNodeColorB}=[PathNode,draw=Col6!80,fill=Col6!48]
\tikzstyle{PathStep}=[color=Col1!60,thick]
\tikzstyle{PathStepColorA}=[color=Col4!70,thick]

\tikzstyle{Box}=[rectangle,rounded corners=2pt,draw=Col1!80,fill=Col1!8,
inner sep=1pt,minimum size=.25cm,thick,font=\scriptsize]
\tikzstyle{BoxClear}=[Box,fill=ColWhite!100]
\tikzstyle{BoxColorA}=[Box,draw=Col4!80,fill=Col4!8]
\tikzstyle{BoxColorB}=[Box,draw=Col6!80,fill=Col6!8]
\tikzstyle{Arc}=[Edge,->,draw=Col5!80]
\tikzstyle{Grid}=[color=ColBlack!30]
\tikzstyle{EdgeRew}=[->,Col2!80,cap=round,thick]

\tikzstyle{Injection}=[ColBlack!100,draw,>->]
\tikzstyle{Surjection}=[ColBlack!100,draw,->>]
\tikzstyle{Map}=[ColBlack!100,draw,->]

\tikzstyle{PosetVertex}=[circle,draw=Col1!80,fill=Col1!8,
inner sep=1pt,minimum size=2mm,thick,font=\scriptsize]

\tikzstyle{CliqueEdge}=[draw=Col2!90,thick]
\tikzstyle{CliqueEdgeColorA}=[CliqueEdge,draw=Col4!80,fill=Col4!8]
\tikzstyle{CliqueEmptyEdge}=[draw=Col4!90,thick,densely dashed]
\tikzstyle{CliqueLabel}=[midway,inner sep=1pt,fill=ColWhite!0,
font=\scriptsize]
\tikzstyle{CliquePoint}=[circle,inner sep=1pt,fill=Col2!25,
draw=Col2!70]
\tikzstyle{CliqueEdgeBlue}=[Col1!80,thick,draw,cap=round]
\tikzstyle{CliqueEdgeRed}=[Col2!80,thick,draw,cap=round,dotted]

\tikzstyle{Operator}=[rectangle,rounded corners,draw=Col1!100,
fill=Col1!10,minimum size=4mm,thick,font=\footnotesize]
\tikzstyle{OperatorColorA}=[Operator,draw=Col4!80,fill=Col4!8]
\tikzstyle{OperatorColorB}=[Operator,draw=Col6!80,fill=Col6!8]
\tikzstyle{OperatorColorC}=[Operator,draw=Col2!80,fill=Col2!8]

\tikzstyle{Circle}=[circle,draw=Col1!100,fill=Col1!5,inner sep=0,
minimum size=14mm,thick,font=\scriptsize]
\tikzstyle{CircleColorA}=[Circle,draw=Col4!80,fill=Col4!16]
\tikzstyle{CircleColorB}=[Circle,draw=Col6!80,fill=Col6!30]

\newcommand{\Hide}[1]{\textcolor{Col4}{\tt [hidden]}}

\newcommand{\Def}[1]{\textcolor{Col3}{\em #1}}
\newcommand{\OEIS}[1]{\href{http://oeis.org/#1}{{\bf #1}}}

\DeclareRobustCommand{\gobblefive}[5]{}
\newcommand*{\SkipTocEntry}{\addtocontents{toc}{\gobblefive}}

\let\SavedCaption=\caption
\renewcommand*{\caption}[2][\shortcaption]{%
    \def\shortcaption{#2}
    \SavedCaption[\; #1]{#2}}

\let\SavedParagraph=\paragraph
\renewcommand{\paragraph}[1]{%
    \SavedParagraph{\it #1}}

\newcommand*\ClearToLeftPage{%
    \clearpage
    \ifodd\value{page}
        \hbox{}
        \vspace*{\fill}
        \thispagestyle{empty}
        \newpage
    \fi
}

\newcommand{\N}{\mathbb{N}}
\newcommand{\Z}{\mathbb{Z}}

\newcommand{\K}{\mathbb{K}}

\newcommand{\Aca}{\mathcal{A}}
\newcommand{\Bca}{\mathcal{B}}

\newcommand{\Lca}{\mathcal{L}}
\newcommand{\Mca}{\mathcal{M}}
\newcommand{\Nca}{\mathcal{N}}
\newcommand{\Oca}{\mathcal{O}}
\newcommand{\Pca}{\mathcal{P}}
\newcommand{\Qca}{\mathcal{Q}}
\newcommand{\Rca}{\mathcal{R}}
\newcommand{\Sca}{\mathcal{S}}
\newcommand{\Tca}{\mathcal{T}}
\newcommand{\Vca}{\mathcal{V}}

\newcommand{\Afr}{\mathfrak{a}}

\newcommand{\Cfr}{\mathfrak{c}}
\newcommand{\Dfr}{\mathfrak{d}}
\newcommand{\Efr}{\mathfrak{e}}

\newcommand{\Mfr}{\mathfrak{m}}
\newcommand{\Sfr}{\mathfrak{s}}
\newcommand{\Rfr}{\mathfrak{r}}
\newcommand{\Pfr}{\mathfrak{p}}
\newcommand{\Qfr}{\mathfrak{q}}
\newcommand{\Tfr}{\mathfrak{t}}
\newcommand{\Ufr}{\mathfrak{u}}
\newcommand{\Vfr}{\mathfrak{v}}

\newcommand{\CFr}{\mathfrak{C}}
\newcommand{\GFr}{\mathfrak{G}}

\newcommand{\Gbb}{\mathbb{G}}
\newcommand{\Hbb}{\mathbb{H}}

\newcommand{\Sbb}{\mathbb{S}}

\newcommand{\Fbf}{\mathbf{f}}
\newcommand{\Gbf}{\mathbf{g}}

\newcommand{\Asf}{\mathsf{a}}
\newcommand{\Bsf}{\mathsf{b}}
\newcommand{\Csf}{\mathsf{c}}
\newcommand{\Dsf}{\mathsf{d}}
\newcommand{\Esf}{\mathsf{e}}

\newcommand{\LambdaB}{\bm{\lambda}}
\newcommand{\MuB}{\bm{\mu}}

\newcommand{\Unit}{\mathds{1}}
\newcommand{\UnitSeries}{\mathbf{1}}
\newcommand{\Counit}{\upsilon}
\newcommand{\OneElement}{\epsilon}
\newcommand{\AtomElement}{\bullet}
\newcommand{\GeneratingSeries}{\Gbb}
\newcommand{\HilbertSeries}{\Hbb}
\newcommand{\Leaf}{\perp}
\newcommand{\Node}{\bullet}
\newcommand{\Alg}{\Aca}

\newcommand{\PrefixOrder}{\preccurlyeq_{\mathrm{p}}}
\newcommand{\NormalForms}{\Nca}
\newcommand{\RelationSpaceRewriteRule}{\Rca}
\newcommand{\RelationSpace}{\Rca}
\newcommand{\GeneratingSet}{\GFr}
\newcommand{\Comparable}{\mathbf{Com}}

\newcommand{\Congr}{\equiv}

\newcommand{\Identity}{\mathrm{Id}}

\DeclareMathOperator{\Ord}{\preccurlyeq}
\newcommand{\OrdStrict}{\prec}

\newcommand{\BinRel}{\,\mathfrak{R}\,}

\newcommand{\EqRel}{\equiv}
\newcommand{\Height}{\mathrm{ht}}
\newcommand{\CoveringRel}{\lessdot}
\newcommand{\SymmetricAction}{\ocircle}
\newcommand{\CyclicAction}{\circlearrowright}
\newcommand{\SetMultisets}{\Mca}
\newcommand{\SetSets}{\Sca}
\newcommand{\SetTuples}{\Tca}
\newcommand{\IndexSeries}{\mathbf{I}}

\DeclareMathOperator{\Casting}{\mathbf{Cast}}
\DeclareMathOperator{\List}{\mathbf{List}}
\DeclareMathOperator{\Multiset}{\mathbf{MSet}}
\DeclareMathOperator{\Set}{\mathbf{Set}}
\DeclareMathOperator{\Coloration}{\mathbf{Col}}
\DeclareMathOperator{\Regularization}{\mathbf{Reg}}
\DeclareMathOperator{\Cycle}{\mathbf{Cyc}}

\DeclareMathOperator{\Support}{\mathrm{Supp}}
\DeclareMathOperator{\Charac}{\mathrm{ch}}
\newcommand{\Angle}[1]{\left\langle#1\right\rangle}
\newcommand{\AAngle}[1]{%
    \left\langle\left\langle#1\right\rangle\right\rangle}
\newcommand{\BBrack}[1]{%
    \pmb{\left[\vphantom{#1}\right.}%
    #1%
    \pmb{\left.\vphantom{#1}\right]}}
\newcommand{\Par}[1]{\left(#1\right)}
\DeclareMathOperator{\ExtProd}{\cdot}
\DeclareMathOperator{\Suspension}{\mathbf{Sus}}
\DeclareMathOperator{\Augmentation}{\mathbf{Aug}}
\DeclareMathOperator{\Length}{\ell}
\DeclareMathOperator{\LDup}{\ll}
\DeclareMathOperator{\RDup}{\gg}
\DeclareMathOperator{\Tensor}{\mathsf{T}}
\newcommand{\Symmetric}{\mathsf{S}}
\newcommand{\Exterior}{\mathsf{E}}
\newcommand{\Product}{\star}
\newcommand{\Coproduct}{\Delta}
\newcommand{\Dual}{\star}
\DeclareMathOperator{\Rew}{\to}
\DeclareMathOperator{\RewRT}{\overset{\ast}{\Rew}}
\DeclareMathOperator{\RewContext}{\Rightarrow}

\DeclareMathOperator{\Std}{\mathrm{std}}
\DeclareMathOperator{\Cmp}{\mathrm{cmp}}
\DeclareMathOperator{\Des}{\mathrm{Des}}

\newcommand{\Conc}{\cdot}
\DeclareMathOperator{\ProductTriangleRight}{\triangleright}

\newcommand{\Narayana}{\mathrm{nar}}

\newcommand{\PreLieProduct}{\curvearrowleft}
\DeclareMathOperator{\Graft}{\curlywedge}
\DeclareMathOperator{\Action}{\bullet}
\DeclareMathOperator{\TreeLanguage}{\Nca}
\DeclareMathOperator{\Arity}{\mathrm{ari}}
\DeclareMathOperator{\TamariInvariant}{\mathrm{tam}}

\newcommand{\Corolla}{\mathrm{c}}
\newcommand{\Eval}{\mathrm{ev}}
\newcommand{\In}{\mathrm{in}}
\newcommand{\Out}{\mathrm{out}}
\newcommand{\LDias}{\dashv}
\newcommand{\RDias}{\vdash}
\newcommand{\LDendr}{\prec}
\newcommand{\RDendr}{\succ}
\newcommand{\LNCT}{\leftharpoonup}
\newcommand{\RNCT}{\rightharpoonup}

\newcommand{\Statistics}{\mathrm{s}}
\newcommand{\Biproduct}{\square}
\newcommand{\Convolution}{*}
\newcommand{\Compo}{\odot}
\newcommand{\Hadamard}{\boxtimes}
\newcommand{\Index}{\mathrm{ind}}

\DeclareMathOperator{\RelBin}{\Rca}
\DeclareMathOperator{\CompoContext}{\circledcirc}

\DeclareMathOperator{\Basic}{\bullet}
\DeclareMathOperator{\PreLieGrafting}{\hookleftarrow}
\DeclareMathOperator{\Labels}{\Lca}
\newcommand{\Image}{\mathrm{Im}}
\newcommand{\ProHoriz}{\ast}
\newcommand{\ProVerti}{\circ}

\newcommand{\SymmetricGroup}{\mathfrak{S}}
\newcommand{\Arcs}{\mathcal{A}}
\newcommand{\Diagonals}{\mathcal{D}}
\newcommand{\Edges}{\mathcal{E}}
\newcommand{\Bool}{\mathbb{B}}

\newcommand{\ColComp}{\mathfrak{Com}} 
\newcommand{\ColIP}{\mathfrak{Par}} 
\newcommand{\ColBT}{\mathfrak{BT}} 
\newcommand{\ColRT}{\mathfrak{RT}} 
\newcommand{\ColSRT}{\mathfrak{SRT}} 
\newcommand{\ColPRT}{\mathfrak{PRT}} 
\newcommand{\ColST}{\mathfrak{ST}} 
\newcommand{\ColCST}{\mathfrak{CST}} 
\newcommand{\ColLad}{\mathfrak{Lad}} 
\newcommand{\ColCor}{\mathfrak{Cor}} 
\newcommand{\ColAry}{\mathfrak{Ary}} 
\newcommand{\ColSch}{\mathfrak{Sch}} 
\newcommand{\ColBSch}{\mathfrak{B}\ColSch} 
\newcommand{\ColCC}{\mathfrak{CC}} 
\newcommand{\ColNCC}{\mathfrak{N}\ColCC} 
\newcommand{\ColBNC}{\mathfrak{BNC}} 
\newcommand{\ColGravCC}{\mathfrak{G}\ColCC} 
\newcommand{\ColNCT}{\mathfrak{NCT}} 
\newcommand{\ColMat}{\mathfrak{Mat}} 
\newcommand{\ColMap}{\mathfrak{Map}} 
\newcommand{\ColNDMap}{\mathfrak{ND}\ColMap} 
\newcommand{\ColInj}{\mathfrak{I}\ColMap} 
\newcommand{\ColSur}{\mathfrak{S}\ColMap} 
\newcommand{\ColBij}{\mathfrak{B}\ColMap} 
\newcommand{\ColPath}{\mathfrak{Path}} 
\newcommand{\ColPathSch}{\mathfrak{Path}_{\mathfrak{Sch}}} 
\newcommand{\ColPathDyck}{\mathfrak{Path}_{\mathfrak{Dyck}}} 
\newcommand{\ColPathMotz}{\mathfrak{Path}_{\mathfrak{Motz}}} 
\newcommand{\ColPathFib}{\mathfrak{Path}_{\mathfrak{Fib}}} 
\newcommand{\ColM}{\mathfrak{M}}

\newcommand{\FreeOperad}{\mathbf{FO}}
\newcommand{\FreeColoredOperad}{\mathbf{FCO}}
\newcommand{\As}{\mathbf{As}}
\newcommand{\Com}{\mathbf{Com}}
\newcommand{\Dup}{\mathbf{Dup}}
\newcommand{\Dendr}{\mathbf{Dendr}}
\newcommand{\PreLie}{\mathbf{PLie}}
\newcommand{\Per}{\mathbf{Per}}
\newcommand{\Mag}{\mathbf{Mag}}
\newcommand{\BS}{\mathbf{BS}}
\newcommand{\Dias}{\mathbf{Dias}}
\newcommand{\Trias}{\mathbf{Trias}}
\newcommand{\NAP}{\mathbf{NAP}}
\newcommand{\BNC}{\mathbf{BNC}}
\newcommand{\NCT}{\mathbf{NCT}}
\newcommand{\Mould}{\mathbf{Mould}}

\newcommand{\Grav}{\mathbf{Grav}}
\newcommand{\T}{\mathbf{T}}
\newcommand{\CliMonoid}{\mathbf{C}}
\newcommand{\CliMagma}{\mathbf{C}_{\mathrm{u}}}
\newcommand{\NCliMonoid}{\mathbf{N}\CliMonoid}
\newcommand{\Bud}{\mathbf{Bud}}
\newcommand{\Motz}{\mathbf{Motz}}

\newcommand{\PAs}{\mathbf{PAs}}
\newcommand{\Mat}{\mathbf{Mat}}
\newcommand{\BRel}{\mathbf{BRel}}

\newcommand{\Sym}{\mathbf{NCSym}}

\newcommand{\FQSym}{\mathbf{FQSym}}
\newcommand{\FSym}{\mathbf{FSym}}
\newcommand{\PBT}{\mathbf{PBT}}
\newcommand{\Bell}{\mathbf{Bell}}
\newcommand{\Baxter}{\mathbf{Baxter}}
\newcommand{\WQSym}{\mathbf{WQSym}}
\newcommand{\PQSym}{\mathbf{PQSym}}
\newcommand{\Camb}{\mathbf{Camb}}

\newcommand{\BasisB}{\mathsf{B}}

\newcommand{\BasisF}{\mathsf{F}}

\newcommand{\BasisR}{\mathsf{R}}
\newcommand{\BasisS}{\mathsf{S}}

\newcommand{\LeafPic}{%
    \raisebox{.39em}{%
    \begin{tikzpicture}[yscale=.28,Centering]
        \node[Leaf](0)at(0,0){};
        \node(r)at(0,1){};
        \draw[Edge](r)--(0);
    \end{tikzpicture}}}

\newcommand{\NodePic}{%
    \raisebox{.39em}{%
    \begin{tikzpicture}[yscale=.25,Centering]
        \node[Node,minimum size=1.5mm](0)at(0,0){};
        \node(r)at(0,1){};
        \draw[Edge](r)--(0);
    \end{tikzpicture}}}

\newcommand{\EdgePic}{%
    \begin{tikzpicture}[xscale=.25,Centering]
        \draw[Edge](0,0)--(1,0);
    \end{tikzpicture}}

\newcommand{\BinaryNode}{
    \raisebox{.25em}{%
    \begin{tikzpicture}[xscale=.35,yscale=.26,Centering]
        \node[Leaf](0)at(0.50,-.75){};
        \node[Leaf](2)at(1.5,-.75){};
        \node[Node](1)at(1.00,0.00){};
        \draw[Edge](0)--(1);
        \draw[Edge](2)--(1);
        \node(r)at(1.00,1.15){};
        \draw[Edge](r)--(1);
    \end{tikzpicture}}}

\newcommand{\BinaryTreeLeft}{
    \begin{tikzpicture}[scale=.17,Centering]
        \node[Leaf](0)at(0.00,-3.33){};
        \node[Leaf](2)at(2.00,-3.33){};
        \node[Leaf](4)at(4.00,-1.67){};
        \node[Node](1)at(1.00,-1.67){};
        \node[Node](3)at(3.00,0.00){};
        \draw[Edge](0)--(1);
        \draw[Edge](1)--(3);
        \draw[Edge](2)--(1);
        \draw[Edge](4)--(3);
        \node(r)at(3.00,1.75){};
        \draw[Edge](r)--(3);
    \end{tikzpicture}}

\newcommand{\BinaryTreeRight}{
    \begin{tikzpicture}[scale=.17,Centering]
        \node[Leaf](0)at(0.00,-1.67){};
        \node[Leaf](2)at(2.00,-3.33){};
        \node[Leaf](4)at(4.00,-3.33){};
        \node[Node](1)at(1.00,0.00){};
        \node[Node](3)at(3.00,-1.67){};
        \draw[Edge](0)--(1);
        \draw[Edge](2)--(3);
        \draw[Edge](3)--(1);
        \draw[Edge](4)--(3);
        \node(r)at(1.00,1.75){};
        \draw[Edge](r)--(1);
    \end{tikzpicture}}

\newcommand{\BinaryRoot}[2]{%
    \begin{tikzpicture}[xscale=.35,yscale=.3,Centering]
        \node(0)at(-.50,-1.50){\begin{math}#1\end{math}};
        \node(2)at(2.50,-1.50){\begin{math}#2\end{math}};
        \node[Node](1)at(1.00,.50){};
        \draw[Edge](0)--(1);
        \draw[Edge](2)--(1);
        \node(r)at(1.00,1.75){};
        \draw[Edge](r)--(1);
    \end{tikzpicture}}

\newcommand{\LeafPicST}{%
    \begin{tikzpicture}[yscale=.35,Centering]
        \node(0)at(0,0){};
        \node(r)at(0,1){};
        \draw[Edge](r)--(0);
    \end{tikzpicture}}

\newcommand{\LeafPicCST}[1]{%
    \begin{tikzpicture}[yscale=.6,Centering]
        \node(0)at(0,0){};
        \node(r)at(0,1){};
        \draw[Edge](r)
            edge[]node[EdgeLabel]{\begin{math}#1\end{math}}(0);
    \end{tikzpicture}}

\newcommand{\UnitClique}{
\begin{tikzpicture}[scale=.6,Centering]
    \node[CliquePoint](1)at(0,0){};
    \node[CliquePoint](2)at(.75,0){};
    \draw[CliqueEmptyEdge](1)edge[]node[CliqueLabel]{}(2);
\end{tikzpicture}}

\newcommand{\CliqueOne}[1]{
\begin{tikzpicture}[scale=1.2,Centering]
    \node[CliquePoint](1)at(0,0){};
    \node[CliquePoint](2)at(.75,0){};
    \draw[CliqueEdge](1)edge[]node[CliqueLabel]
        {\begin{math}#1\end{math}}(2);
\end{tikzpicture}}

\newcommand{\UnitSolidClique}{
\begin{tikzpicture}[scale=.6,Centering]
    \node[CliquePoint](1)at(0,0){};
    \node[CliquePoint](2)at(.75,0){};
    \draw[CliqueEdgeBlue](1)--(2);
\end{tikzpicture}}

\newcommand{\TriangleRight}{
\begin{tikzpicture}[scale=.3,Centering]
    \node[CliquePoint](1)at(0,1){};
    \node[CliquePoint](2)at(0.87,-0.5){};
    \node[CliquePoint](3)at(-0.87,-0.5){};
    \draw[CliqueEdgeBlue](1)--(2);
    \draw[CliqueEmptyEdge](1)--(3);
    \draw[CliqueEdgeBlue](2)--(3);
\end{tikzpicture}}

\newcommand{\TriangleLeft}{
\begin{tikzpicture}[scale=.3,Centering]
    \node[CliquePoint](1)at(0,1){};
    \node[CliquePoint](2)at(0.87,-0.5){};
    \node[CliquePoint](3)at(-0.87,-0.5){};
    \draw[CliqueEmptyEdge](1)--(2);
    \draw[CliqueEdgeBlue](1)--(3);
    \draw[CliqueEdgeBlue](2)--(3);
\end{tikzpicture}}

\newcommand{\TriangleEmpty}{
\begin{tikzpicture}[scale=.3,Centering]
    \node[CliquePoint](1)at(0,1){};
    \node[CliquePoint](2)at(0.87,-0.5){};
    \node[CliquePoint](3)at(-0.87,-0.5){};
    \draw[CliqueEmptyEdge](1)--(2);
    \draw[CliqueEmptyEdge](1)--(3);
    \draw[CliqueEmptyEdge](2)--(3);
\end{tikzpicture}}

\newcommand{\UnitPath}{
\begin{tikzpicture}[scale=.28,Centering]
    \draw[Grid] (0,0) grid (0,0);
    \node[PathNode,minimum size=1.5mm](0)at(0,0){};
\end{tikzpicture}}

\newcommand{\PathUp}{
\begin{tikzpicture}[scale=.18,Centering]
    \draw[Grid](0,0) grid (1,1);
    \node[PathNode,minimum size=1.5mm](Dyck1)at(0,0){};
    \node[PathNode,minimum size=1.5mm](Dyck2)at(1,1){};
    \draw[PathStep](Dyck1)--(Dyck2);
\end{tikzpicture}}

\newcommand{\PathDown}{
\begin{tikzpicture}[scale=.18,Centering]
    \draw[Grid](0,0) grid (1,1);
    \node[PathNode,minimum size=1.5mm](Dyck1)at(0,1){};
    \node[PathNode,minimum size=1.5mm](Dyck2)at(1,0){};
    \draw[PathStep](Dyck1)--(Dyck2);
\end{tikzpicture}}

\newcommand{\PathStable}{
\begin{tikzpicture}[scale=.24,Centering]
    \draw[Grid] (0,0) grid (1,0);
    \node[PathNode,minimum size=1.5mm](0)at(0,0){};
    \node[PathNode,minimum size=1.5mm](1)at(1,0){};
    \draw[PathStep](0)--(1);
\end{tikzpicture}}

\newcommand{\PathStableStable}{
\begin{tikzpicture}[scale=.2,Centering]
    \draw[Grid] (0,0) grid (2,0);
    \node[PathNode,minimum size=1.5mm](0)at(0,0){};
    \node[PathNode,minimum size=1.5mm](1)at(1,0){};
    \node[PathNode,minimum size=1.5mm](2)at(2,0){};
    \draw[PathStep](0)--(1);
    \draw[PathStep](1)--(2);
\end{tikzpicture}}

\newcommand{\PathPeak}{
\begin{tikzpicture}[scale=.2,Centering]
    \draw[Grid] (0,0) grid (2,1);
    \node[PathNode,minimum size=1.5mm](0)at(0,0){};
    \node[PathNode,minimum size=1.5mm](1)at(1,1){};
    \node[PathNode,minimum size=1.5mm](2)at(2,0){};
    \draw[PathStep](0)--(1);
    \draw[PathStep](1)--(2);
\end{tikzpicture}}

\newcommand{\PathStableOp}{
\begin{tikzpicture}[scale=.24,Centering]
    \draw[Grid] (0,0) grid (1,0);
    \node[PathNodeColorB,minimum size=1.5mm](0)at(0,0){};
    \node[PathNodeColorB,minimum size=1.5mm](1)at(1,0){};
    \draw[PathStepColorA](0)--(1);
\end{tikzpicture}}

\newcommand{\PathPeakOp}{
\begin{tikzpicture}[scale=.2,Centering]
    \draw[Grid] (0,0) grid (2,1);
    \node[PathNodeColorB,minimum size=1.5mm](0)at(0,0){};
    \node[PathNodeColorB,minimum size=1.5mm](1)at(1,1){};
    \node[PathNodeColorB,minimum size=1.5mm](2)at(2,0){};
    \draw[PathStepColorA](0)--(1);
    \draw[PathStepColorA](1)--(2);
\end{tikzpicture}}

\newcommand{\Matrix}[1]{
\left[
    \begin{smallmatrix}
        #1
    \end{smallmatrix}
\right]}

\newcommand{\Keywords}{
    Combinatorics;
    Algebraic combinatorics;
    Computer science;
    Tree;
    Formal power series;
    Rewrite system;
    Operad;
    Bialgebra;
    Hopf bialgebra;
    Pre-Lie algebra;
    Dendriform algebra.
}

\newcommand{\Subjclass}{
    05-00, 
    05C05, 
    05E15, 
    16T05, 
    18D50. 
}

\newcommand{\Abstract}{
    \normalsize
    Operads are algebraic devices offering a formalization of the
    concept of operations with several inputs and one output. Such
    operations can be naturally composed to form bigger and more complex
    ones. Coming historically from algebraic topology, operads intervene
    now as important objects in computer science and in combinatorics.
    The theory of operads, together with the algebraic setting and the
    tools accompanying it, promises advances in these two areas. On the
    one hand, operads provide a useful abstraction of formal
    expressions, and also, provide connections with the theory of
    rewrite systems. On the other hand, a lot of operads involving
    combinatorial objects highlight some of their properties and allow
    to discover new ones.
    \smallbreak

    This book presents the theory of nonsymmetric operads under a
    combinatorial point of view. It portrays the main elements of this
    theory and the links it maintains with several areas of computer
    science and combinatorics. A lot of examples of operads appearing in
    combinatorics are studied and some constructions relating operads
    with known algebraic structures are presented. The modern treatment
    of operads consisting in considering the space of formal power
    series associated with an operad is developed. Enrichments of
    nonsymmetric operads as colored, cyclic, and symmetric operads are
    reviewed.
    \smallbreak

    This text is addressed to any computer scientist or combinatorist
    who looks a complete and a modern description of the theory of
    nonsymmetric operads. Evenly, this book is intended to an audience
    of algebraists who are looking for an original point of view fitting
    in the context of combinatorics.
}

\title{Nonsymmetric operads in combinatorics}
\subjclass[2010]{\Subjclass}
\date{\today}
\author{Samuele Giraudo}
\address{\scriptsize Université Paris-Est, LIGM (UMR $8049$), CNRS,
    ENPC, ESIEE Paris, UPEM, F-$77454$, Marne-la-Vallée, France}
\email{samuele.giraudo@u-pem.fr}

\allowdisplaybreaks

\numberwithin{equation}{subsection}

\counterwithin{figure}{chapter}
\counterwithin{table}{chapter}

\makeatletter
\def\l@part{\@tocline{-1}{20pt}{0pc}{5pc}{\Large \bf}}
\def\l@chapter{\@tocline{0}{10pt}{0pc}{5pc}{\bf}}
\def\l@section{\@tocline{1}{3pt}{1pc}{5pc}{}}
\def\l@subsection{\@tocline{2}{2pt}{2pc}{5pc}{}}
\makeatother

\renewcommand{\leq}{\leqslant}
\renewcommand{\geq}{\geqslant}

\newtheorem{Theorem}{Theorem}[subsection]
\newtheorem{Proposition}[Theorem]{Proposition}
\newtheorem{Lemma}[Theorem]{Lemma}

\begin{document}

\begin{abstract}
    \Abstract
\end{abstract}

\frontmatter

\maketitle

\bigbreak

\begin{center}
    \today
\end{center}
\bigbreak

\tableofcontents

\mainmatter


\chapter*{Introduction}

\SkipTocEntry\section*{Combinatorics, algebra, and programming}
Let us start this text by considering two fields of mathematics:
combinatorics and algebra, and a field of computer science: programming.
These topics have far more in common than it appears at first glance and
we shall explain why.
\medbreak

\SkipTocEntry\subsection*{Combinatorics}
One of the most common activities in combinatorics consists in studying
families of combinatorial objects such as permutations, words, trees, or
graphs. One of the most common activities about such families is to
provide formulas to enumerate them according to a decent size notion. A
large number of tools exists for this purpose such as manipulation of
generating series, use of bijections, and use of decompositions rules to
break objects into elementary pieces. This last tool is crucial and
relies on the fact that if one knows a way to {\em compose} objects to
build a bigger one, then without much effort one knows a way to
{\em decompose} them. Let us, for instance, consider the celebrated set 
of Motzkin paths. Recall that a Motzkin path is a walk in the quarter 
plane connecting the origin to a point on the $x$-axis by means of steps
$\PathUp$, $\PathStable$, and $\PathDown$. Let us denote by
$\PathStableOp\Par{\Ufr_1, \Ufr_2}$ the binary composition rule
consisting in substituting $\Ufr_1$ (resp. $\Ufr_2$) with the first
(resp. second) point of $\PathStable$, and by
$\PathPeakOp\Par{\Ufr_1, \Ufr_2, \Ufr_3}$ the ternary composition rule
consisting in substituting $\Ufr_1$ (resp. $\Ufr_2$, $\Ufr_3$) with the
first (resp. second, third) point of $\PathPeak$. By using these two
operations, one can, for instance, provide the decomposition
\begin{equation}
    \begin{tikzpicture}[scale=.3,Centering]
        \draw[Grid] (0,0) grid (8,2);
        \node[PathNode](0)at(0,0){};
        \node[PathNode](1)at(1,1){};
        \node[PathNode](2)at(2,2){};
        \node[PathNode](3)at(3,1){};
        \node[PathNode](4)at(4,0){};
        \node[PathNode](5)at(5,0){};
        \node[PathNode](6)at(6,1){};
        \node[PathNode](7)at(7,1){};
        \node[PathNode](8)at(8,0){};
        \draw[PathStep](0)--(1);
        \draw[PathStep](1)--(2);
        \draw[PathStep](2)--(3);
        \draw[PathStep](3)--(4);
        \draw[PathStep](4)--(5);
        \draw[PathStep](5)--(6);
        \draw[PathStep](6)--(7);
        \draw[PathStep](7)--(8);
    \end{tikzpicture}
    \enspace = \enspace
    \PathStableOp\Par{
    \PathPeakOp\Par{\UnitPath, \PathPeak, \UnitPath},
    \PathPeakOp\Par{\UnitPath, \PathStable, \UnitPath}}
\end{equation}
of the Motzkin path appearing as left member. An easy result says that a
nonempty Motzkin path decomposes in a unique way either as a step
$\PathStable$ followed by a Motzkin path, or as a step $\PathUp$
followed by a Motzkin path, a step $\PathDown$, and another Motzkin
path. This leads to the expression
\begin{equation} \label{equ:expression_motzkin_paths_set}
    \ColM \enspace = \enspace
    \{\UnitPath\}
    \enspace \sqcup \enspace
    \left\{
        \PathStableOp\Par{\UnitPath, \Ufr_2} : \Ufr_2 \in \ColM
    \right\}
    \enspace \sqcup \enspace
    \left\{
        \PathPeakOp\Par{\UnitPath, \Ufr_2, \Ufr_3} :
        \Ufr_2, \Ufr_3 \in \ColM
    \right\}
\end{equation}
describing the set $\ColM$ of all Motzkin paths.
From~\eqref{equ:expression_motzkin_paths_set}, one harvests
combinatorial data as a formula to count Motzkin paths of a given size
(with respect to an adequate notion of size), bijections between Motzkin
paths and other families of combinatorial objects admitting similar
decompositions (like Motzkin trees or Motzkin configurations of chords
in polygons), or even the definition of generalizations of these objects
(by adding colors on the steps, or by adding other possible sorts of
steps).
\medbreak

\SkipTocEntry\subsection*{Algebra}
On another side, in algebra, the study of algebraic structures
consisting in a set of operations and the relations they satisfy is very
habitual. To specify such an algebraic structure, one provides the
nontrivial relations satisfied by its operations. For instance, any
binary operation $\Product$ satisfying for any inputs $x_1$, $x_2$, and
$x_3$ the relation
\begin{equation} \label{equ:example_relation_pre_lie}
    \Par{x_1 \Product x_2} \Product x_3
    - x_1 \Product \Par{x_2 \Product x_3}
    =
    \Par{x_1 \Product x_3} \Product x_2
    - x_1 \Product \Par{x_3 \Product x_2}
\end{equation}
specifies the class of the so-called pre-Lie algebras. Classical
examples include among others monoids, groups, lattices, associative
algebras, commutative algebras, and dendriform algebras. An obvious but
important fact is that expressions and relations can be {\em composed}:
the variables in a relation refer to any term of the algebraic
structure. For instance, in~\eqref{equ:example_relation_pre_lie}, all
the occurrences of $x_2$ can be replaced by a term
${\bf x_4 \Product x_5}$, leading to the relation
\begin{equation} \label{equ:example_relation_pre_lie_consequence}
    \Par{x_1 \Product \Par{\bf x_4 \Product x_5}} \Product x_3
    - x_1 \Product \Par{\Par{\bf x_4 \Product x_5} \Product x_3}
    =
    \Par{x_1 \Product x_3} \Product \Par{\bf x_4 \Product x_5}
    - x_1 \Product \Par{x_3 \Product \Par{\bf x_4 \Product x_5}}
\end{equation}
that still holds. Furthermore, by seeing the binary operation $\Product$
as a device
\begin{equation}

\end{equation}
having two inputs $x_1$ and $x_2$ and one output $x_1 \Product x_2$,
Relation~\eqref{equ:example_relation_pre_lie} reformulates as
\begin{equation}
    %
\,.
\end{equation}
Therefore, one observes that algebraic expressions translates as syntax
trees, and that composition and substitution of expressions translates
as grafting of trees.
\medbreak

\SkipTocEntry\subsection*{Programming}
In programming, and more precisely in the context of functional
programming, one encounters expressions of the form
\begin{equation} \label{equ:example_composition_functions}
    \mathtt{f(g1(v1, v2), g2(v3), v4)}
\end{equation}
where $\mathtt{f}$, $\mathtt{g1}$, and $\mathtt{g2}$ are function names,
and $\mathtt{v1}$, $\mathtt{v2}$, $\mathtt{v3}$, and $\mathtt{v4}$ are
value names. This function call can be considered as a program, and the
returned value is the value computed by the execution of the program.
When our programming language satisfies referential transparency (that
is, any expression can be substituted by its value without changing the
overall computation), the call~\eqref{equ:example_composition_functions}
is equivalent to both the calls
\begin{equation}
    \mathtt{f(\mathbf{w1}, g2(v3), v4)}
    \qquad \mbox{and} \qquad
    \mathtt{f(g1(v1, v2), \mathbf{w2}, v4)}
\end{equation}
where $\mathtt{\mathbf{w1}}$ is the value of $\mathtt{g1(v1, v2)}$ and
$\mathtt{\mathbf{w2}}$ is the value of $\mathtt{g2(v3)}$. As a
consequence, this leads to the fact that the order of evaluation of the
sub-expressions $\mathtt{g1(v1, v2)}$ and $\mathtt{g2(v3)}$ does not
influence the computation of~\eqref{equ:example_composition_functions}.
Furthermore, by seeing the functions $\mathtt{f1}$, $\mathtt{g1}$, and
$\mathtt{g2}$ as black boxes
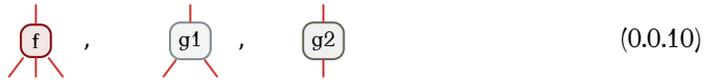
\begin{equation}

\end{equation}
where inputs are depicted below and outputs above the boxes, the
program~\eqref{equ:example_composition_functions} is made of
{\em compositions} of such black boxes where inputs are completed with
values (that are arguments of the function calls). Under this formalism,
this program and its evaluation are depicted as
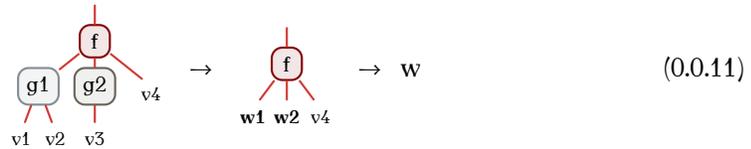
\begin{equation}
    %
    \enspace \Rew \enspace
    \mathtt{w}
\end{equation}
where $\mathtt{w}$ is the value computed by the program. In this
context, the evaluation of a program can be carried out by using rewrite
rules on such composition diagrams.
\medbreak

\SkipTocEntry\section*{Nonsymmetric operads as a meeting point}
The three previous examples highlight the importance of the notion of
composition in combinatorics, algebra, and programming. This is
precisely the common point between these three fields we want to
emphasize. Let us develop this concept.
\medbreak

\SkipTocEntry\subsection*{Coherent compositions}
The compositions of combinatorial objects, of syntax trees of algebraic
expressions, and of function calls have to be coherent. Indeed, to lead
to interesting and substantial consequences, they have to mimic the
usual functional composition. Composing two objects $x$ and $y$ consists
in choosing a {\em substitution sector} $a$ of $x$ and in replacing $a$
it by a copy of~$y$. This composition is denoted by $x \circ_a y$ and is
schematically represented, in the following arbitrary example, as
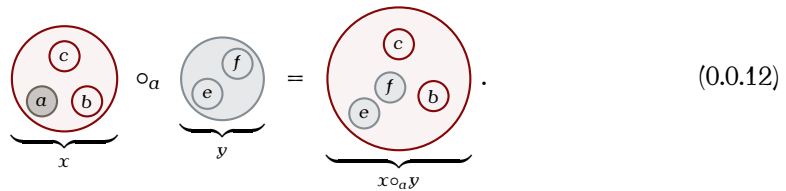
\begin{equation}
    \underbrace{
}_{x \circ_a y}\,.
\end{equation}
To be coherent, this composition has to satisfy the two relations
\begin{equation} \label{equ:relation_composition_1}
    \Par{
    %
}\,.
\end{equation}
First relation can be thought as a horizontal compatibility, and the
second as a vertical one.
\medbreak

To come back on the previous examples:
\begin{enumerate}[fullwidth,label={\bf (\arabic*)}]
    \item In a Motzkin path $\Ufr$, the substitution sectors are its
    points, specified by their positions~$i$. The composition
    $\Ufr \circ_i \Vfr$ of two Motzkin paths $\Ufr$ and $\Vfr$ consists
    in replacing the $i$th point of $\Ufr$ by a copy of~$\Vfr$.
    \medbreak

    \item In a syntax tree $\Tfr$ of an algebraic expression, the
    substitution sectors are its leaves, specified by their labels
    $x_i$. The composition $\Tfr \circ_{x_i} \Sfr$ of two syntax trees
    $\Tfr$ and $\Sfr$ consists in grafting the root of~$\Sfr$ onto the
    leaf $x_i$ of~$\Tfr$.
    \medbreak

    \item In a function $\mathtt{f}$ of a functional programming
    language, the substitution sectors are its parameters $\mathtt{x}_i$
    provided that $\mathtt{f}$ admits the prototype
    \begin{math}
        \mathtt{f}(\mathtt{x}_1, \dots,
        \mathtt{x}_{i - 1}, \mathtt{x}_i, \mathtt{x}_{i + 1},
        \dots, \mathtt{x}_n).
    \end{math}
    The composition $\mathtt{f} \circ_{\mathtt{x}_i} \mathtt{g}$ where
    $\mathtt{f}$ is the function just considered and $\mathtt{g}$ is a
    function admitting the prototype
    \begin{math}
        \mathtt{g}(\mathtt{y}_1, \dots, \mathtt{y}_m)
    \end{math}
    is the function with prototype
    \begin{math}
        \mathtt{f}(\mathtt{x}_1, \dots,
        \mathtt{x}_{i - 1},
        \mathtt{y}_1, \dots, \mathtt{y}_m,
        \mathtt{x}_{i + 1},
        \dots, \mathtt{x}_n)
    \end{math}
    defined as the function calling $\mathtt{f}$ wherein its $i$th
    argument substituting $\mathtt{x}_i$ is set to the value returned
    by~$\mathtt{g}$.
\end{enumerate}
It is easy to see that Relations~\eqref{equ:relation_composition_1}
and~\eqref{equ:relation_composition_2} are satisfied in these cases.
\medbreak

\SkipTocEntry\subsection*{Nonsymmetric operads}
Nonsymmetric operads are precisely algebraic structures furnishing an
abstraction of this concept of generalized compositions. Intuitively, a
nonsymmetric operad $\Oca$ is a set (or a vector space) equipped with
a map $|-|$ associating a positive integer with each of its elements,
and with composition maps $\circ_i$. Each element $x$ of $\Oca$ is
seen as an operator of arity $|x|$ and the composition maps of $\Oca$
satisfy the coherence relations stated above. Nonsymmetric operads mimic
and generalize in this way the usual functional composition for any
families of objects. There exist, for instance, nonsymmetric operads on
words, permutations, binary trees, Schröder trees, configurations of
chords, and paths.
\medbreak

There are a lot of reasons motivating the study of nonsymmetric operads.
Here follow the main ones:
\begin{enumerate}[fullwidth,label={\bf (\Alph*)}]
    \item Endowing a set of combinatorial objects with the structure
    of a nonsymmetric operad provides an algebraic framework for 
    studying it. This framework can potentially stress some of the 
    properties of the combinatorial objects, as for instance, 
    enumerative results and the discovery of hidden symmetries.
    \medbreak

    \item Continuing this last point, with any nonsymmetric operad
    $\Oca$ defined in the category of sets, it is associated a space of
    formal power series $\K \AAngle{\Oca}$ on the elements of $\Oca$.
    This notion of series on nonsymmetric operads generalizes the usual
    one. Moreover, the extension of the composition maps of $\Oca$ on
    $\K \AAngle{\Oca}$ leads to generalizations of the multiplication 
    and the composition products of series. All this provides 
    alternative ways to obtain expressions for the generating series of 
    families of combinatorial objects.
    \medbreak

    \item Nonsymmetric operads admit close connections with the
    combinatorics of planar rooted trees. This is due to the fact that
    free nonsymmetric operads can be constructed on sets of such trees.
    Moreover, given a nonsymmetric operad $\Oca$, a classical question
    consists in exhibiting a presentation by generators and relations of
    $\Oca$. Since any nonsymmetric operad can be described as a quotient
    of a free operad, the computation of a presentation is based upon
    manipulation of trees. In this context, tools coming from rewrite
    systems and some of their properties like termination and confluence
    intervene.
    \medbreak

    \item There are a lot of generalizations of nonsymmetric operads,
    increasing their fields of applications. For instance, in a
    symmetric operad $\Oca$, all symmetric groups $\SymmetricGroup(n)$,
    $n \in \N$, act on the sets of elements of arity $n$ of $\Oca$ by
    permuting the inputs of their elements. Besides, in a colored operad
    $\Oca$, the inputs and outputs have a color and the composition of
    two elements is defined only if the colors of the involved input and
    output match. In fact, our previous example about functional
    programming and composition of functions lies in this context of
    colors operads when the language is typed: the colors play the role
    of data types.
    \medbreak

    \item Given a nonsymmetric operad $\Oca$, there is a notion of
    algebras over $\Oca$. More precisely, an algebra over $\Oca$ is a
    vector space $\Alg$ wherein elements of $\Oca$ behave as operations
    on $\Alg$ by respecting the arities and the composition maps of the
    nonsymmetric operad. For instance, any algebra over the associative
    operad is an associative algebra, and any algebra over the pre-Lie
    operad is a pre-Lie algebra.
    \medbreak

    \item Related to the previous point, all the algebras over $\Oca$
    form a category of algebras. In this way, by studying $\Oca$, one
    obtains general results on all the algebras over $\Oca$. For
    instance, it is possible to show that the sum of the two usual
    generators of the dendriform operad behaves as an associative
    operation. This implies that the sum of the two operations of any
    dendriform algebra is associative. Besides, operads furnish here a
    framework for operation calculus.
    \medbreak

    \item Again related to the two previous points, nonsymmetric operads
    lead to the discovery of link between different categories of
    algebras. Indeed, if $\phi : \Oca_1 \to \Oca_2$ is a morphism
    between two nonsymmetric operads $\Oca_1$ and $\Oca_2$, one can
    construct from $\phi$ a functor from the category of algebras over
    $\Oca_2$ to the category of algebras over $\Oca_1$. For instance,
    there is a morphism from the two-associative operad to the duplicial
    operad leading to a functor from the category of duplicial algebras
    (algebras endowed with two binary associative operations satisfying
    one extra relation) to the category of two-associative algebras
    (algebras endowed with two binary associative operations).
    \medbreak

    \item Given a nonsymmetric operad $\Oca$ satisfying some precise
    properties, it is possible to compute a presentation by generators
    and relations of its so-called Koszul dual $\Oca^!$. This duality
    is an extension of the Koszul duality for associative algebras and
    establishes connections between nonsymmetric operads at first sight
    very different. For instance, the dendriform operad and the
    diassociative operad are Koszul dual one of the other. Moreover,
    linked to this notion of Koszul duality, there is a notion of
    Koszulity for nonsymmetric operads. This property, defined
    originally algebraically, admits equivalent reformulations in terms
    of properties of rewrite systems on trees associated with the
    presentation of the nonsymmetric operads. Moreover, given a Koszul
    nonsymmetric operad $\Oca$ admitting an Hilbert series, the
    alternating version of the Hilbert series of $\Oca^!$ and the one of
    $\Oca$ are inverse one of the other for series composition. This
    property admits for instance applications for enumerative prospects.
\end{enumerate}
\medbreak

\SkipTocEntry\section*{Construction of the book}
We now give some practical information about this text.
\medbreak

\SkipTocEntry\subsection*{Point of view}
All the algebraic structures considered here are linear spans
$\K \Angle{C}$ of some sets of combinatorial objects $C$. For this
reason, these spaces admit always an explicit basis $C$. Moreover, to
handle families of combinatorial objects, we introduce the notion of
collections. A collection is a set of combinatorial objects presented as
a disjoint union of subsets of objects satisfying a same property. For
instance, a graded collection is a set of combinatorial objects defined
as a disjoint union of sets of objects of a same size. We shall consider
different sorts of collections as colored, multigraded, or symmetric
collections, suited to the study of some particular algebraic
structures.
\medbreak

\SkipTocEntry\subsection*{Main purposes and audience}
The aim of this text is to give a presentation of the theory of
nonsymmetric operads under the context of algebraic combinatorics. We
orient our exposition of properties of operads toward enumerative and
combinatorial directions. This book is also intended to be an
introduction to algebraic combinatorics: the first chapters deal with
combinatorics and combinatorics of trees, while next one deals with
general properties of algebraic structures on combinatorial families.
For instance, Hopf bialgebras, forming a vast and rich topic in
algebraic combinatorics, are studied here.
\medbreak

We decide to not overload the text with bibliographic and historical
references. For this reason, each chapter ends with a section containing
bibliographic material. Besides, to present as much results as possible
about combinatorics, algebraic combinatorics, and operads, we decide to
not mention their proofs. For a vast majority of them, they can be
treated as not so hard exercises.
\medbreak

This book is addressed to any computer scientist or combinatorist who is
aiming to establish a first contact with the theory of operads. Evenly,
this book is intended to an audience of algebraists who are looking for
an original point of view fitting in the context of combinatorics.
\medbreak

\SkipTocEntry\subsection*{Structure}
The book contains five chapters. Each chapter depends on the previous
ones. They are organized as follows.
\begin{enumerate}[fullwidth,label={\bf (\arabic*)}]
    \item Chapter~\ref{chap:collections} presents general notions
    about combinatorics and collections. It also provides definitions
    about collections endowed with a poset structure. These structures
    intervene in changes of bases of algebraic structures. It presents
    finally collections endowed with rewrite rules. These collections
    intervene as tools to establish presentations of nonsymmetric
    operads.
    \medbreak

    \item Chapter~\ref{chap:trees} is devoted to general treelike
    structures. It presents syntax trees, that are sorts of trees
    appearing in the study of nonsymmetric operads. Rewrite systems on
    syntax trees are considered and tools to prove their termination or
    confluence are provided.
    \medbreak

    \item Chapter~\ref{chap:algebra} concerns algebraic structures
    defined on the linear span of collections. These vector spaces are
    called polynomial spaces. It presents also the notion of biproducts
    on polynomial spaces and of types of bialgebras. Classical types of
    bialgebras appearing in combinatorics are given: associative,
    dendriform, and pre-Lie algebras, and Hopf bialgebras.
    \medbreak

    \item Chapter~\ref{chap:operads} presents nonsymmetric operads and
    related notions. It exposes the notions of algebras over operads,
    free operads, presentations by generators and relations,
    Koszul duality, and Koszulity. Several examples of operads
    appearing in algebraic combinatorics are reviewed.
    \medbreak

    \item Chapter~\ref{chap:generalizations} contains generalizations,
    applications, and apertures of the theory of nonsymmetric
    operads. It reviews three topics in this vein. First, formal
    power series on nonsymmetric operads are considered and applications
    to enumeration are provided. Next, enrichments of nonsymmetric
    operads are discussed: colored operads, cyclic operads, and
    symmetric operads. Finally, it provides an overview of product
    categories, a generalization of operads wherein elements can have
    several outputs.
\end{enumerate}
\medbreak


\chapter{Enriched collections} \label{chap:collections}
This preliminary chapter contains general notions about combinatorics
used in the rest of the book. We introduce the notion of collections
of combinatorial objects and then the notions of posets and rewrite
systems, which are seen as collections endowed with some extra
structure.
\medbreak

\section{Collections} \label{sec:collections}
A collection is a set of combinatorial objects partitioned into subsets
of objects sharing a same property. Operations over collections are then
introduced and some classical examples of collections are provided.
\medbreak

\subsection{Structured collections} \label{subsec:collections}
There are many kinds of collections such as graded, multigraded,
colored, cyclic, and symmetric. Definitions about them are provided
here.
\medbreak

\subsubsection{Elementary definitions}
\label{subsubsec:collections_definition}
Let $I$ be a nonempty set called \Def{index set}. An
\Def{$I$-collection} is a set $C$ expressible as a disjoint union
\begin{equation} \label{equ:collection}
    C = \bigsqcup_{i \in I} C(i)
\end{equation}
where all $C(i)$, $i \in I$, are (possibly infinite) sets. All the
elements of $C$ (resp. $C(i)$ for an $i \in I$) are called
\Def{objects} (resp. \Def{$i$-objects}) of $C$. If $x$ is an $i$-object
of $C$, we say that the \Def{index} $\Index(x)$ of $x$ is $i$. When for
all $i \in I$, all $C(i)$ are finite sets, $C$ is \Def{combinatorial}.
Besides, $C$ is \Def{finite} if $C$ is finite as a set. The
\Def{empty $I$-collection} is the set $\emptyset$. When $I$ is a
singleton, $C$ is \Def{simple}. Any set can thus be seen as a simple
collection and conversely.
\medbreak

A \Def{relation} on $C$ is a binary relation $\RelBin$ on $C$ such that
for any objects $x$ and $y$ of $C$ satisfying $x \RelBin y$,
$\Index(x) = \Index(y)$. Let $C_1$ and $C_2$ be two $I$-collections. A
map $\phi : C_1 \to C_2$ is an \Def{$I$-collection morphism} if, for
all $x \in C_1$, $\Index(x) = \Index(\phi(x))$. We express by
$C_1 \simeq C_2$ the fact that there exists an isomorphism between
$C_1$ and $C_2$. Besides, if for all $i \in I$,
$C_1(i) \subseteq C_2(i)$, $C_1$ is a \Def{subcollection} of $C_2$. For
any $i \in I$, we can regard each $C(i)$ as a subcollection of $C$
consisting only in all its $i$-objects. Moreover, for any subset $J$
of $I$, we denote by $C(J)$ the subcollection of $C$ consisting only
in all its $j$-objects for all $j \in J$.
\medbreak

Let us now consider particular $I$-collections for precise sets~$I$.
Table~\ref{tab:sorts_collections} contains an overview of the
properties that such collections can satisfy.
\begin{table}[ht]
    \centering
    \begin{tabular}{|c|c|c|} \hline
        \multicolumn{3}{|c|}{{\bf Collections}} \\
        \multicolumn{3}{|c|}{
            {\it Combinatorial} \qquad
            {\it Finite} \qquad
            {\it Simple} \qquad
            {\it With products}}
            \\ \hline
        \multicolumn{2}{|c|}{{\bf $k$-graded}}
            & \multirow{3}{*}{{\bf Colored}} \\ \cline{1-2}
        \multicolumn{2}{|c|}{{\bf $1$-graded}}
            & \multirow{3}{*}{{\it Monochrome} \quad
                {\it $k$-colored}} \\
        \multicolumn{2}{|c|}{{\it Connected} \quad {\it Augmented}
            \quad {\it Monatomic}} & \\ \cline{1-2}
        \hspace{1cm} {\bf Cyclic} \hspace{1cm}
            & {\bf Symmetric} & \\ \hline
    \end{tabular}
    \bigbreak

    \caption[Properties of $I$-collections.]
    {\footnotesize The most common sorts of $I$-collections (in bold)
    and the properties (in italic) they can satisfy. The inclusions
    relations between these sorts of collections read from bottom to
    top. For instance, cyclic collections are particular $1$-graded
    collections which are themselves particular $k$-graded collections
    which are themselves particular collections.}
    \label{tab:sorts_collections}
\end{table}
\medbreak

\subsubsection{Graded collections} \label{subsubsec:graded_collections}
An $\N$-collection is called a \Def{graded collection}. If $C$ is a
graded collection, for any object $x$ of $C$, the \Def{size} $|x|$ of
$x$ is the integer $\Index(x)$. The map $|-| : C \to \N$ is the
\Def{size function} of~$C$.
\medbreak

We say that $C$ is \Def{connected} if $C(0)$ is a singleton, and that
$C$ is \Def{augmented} if $C(0) = \emptyset$. Moreover, $C$ is
\Def{monatomic} if it is augmented and $C(1)$ is a singleton. We denote
by $\{\OneElement\}$ the graded collection such that $\OneElement$ is an
object satisfying $|\OneElement| = 0$. This collection is called the
\Def{unit collection}. Observe that $\{\OneElement\}$ is connected, and
that $C$ is connected if and only if there is a unique collection
morphism from $\{\OneElement\}$ to $C$. We denote by $\{\AtomElement\}$
the collection such that $\AtomElement$ is an \Def{atom}, that is an
object satisfying $|\AtomElement| = 1$. This collection is called the
\Def{neutral collection}. Observe that $\{\AtomElement\}$ is monatomic,
and that $C$ is monatomic if and only if $C$ is augmented and there is a
unique collection morphism from $\{\AtomElement\}$ to~$C$. When $C$ is a
combinatorial graded collection, the \Def{generating series} of $C$ is
the series
\begin{equation} \label{equ:generating_series_graded_collection}
    \GeneratingSeries_C(t) := \sum_{n \in \N} \# C(n) t^n,
\end{equation}
where $t$ is a formal parameter and $\# S$ denotes the cardinality of
any finite set $S$. This formal power series encodes the \Def{integer
sequence} of $C$, that is the sequence $\Par{\# C(n)}_{n \in \N}$.
Observe that if $C_1$ and $C_2$ are two combinatorial graded
collections, $C_1 \simeq C_2$ holds if and only if
$\GeneratingSeries_{C_1}(t) = \GeneratingSeries_{C_2}(t)$.
\medbreak

\subsubsection{Multigraded collections and statistics}
\label{subsubsec:multigraded_collections}
A \Def{$k$-graded collection} (also called \Def{multigraded collection})
is an $\N^k$-collection for an integer $k \geq 1$. To not overload the
notation, we denote by $C\Par{n_1, \dots, n_k}$ the subset
$C\Par{\Par{n_1, \dots, n_k}}$ of any $k$-graded collection $C$. Recall
that a \Def{statistics} on an $I$-collection $C$ is a map
$\Statistics : C \to \N$, associating a nonnegative integer value with
any object of $C$. Multigraded collections are useful to work with
objects endowed with many statistics. Indeed, if $x$ is an
$\Par{n_1, \dots, n_k}$-object of a $k$-graded collection $C$, one sets
$\Statistics_j(x) := n_j$ for each $1 \leq j \leq k$. This defines in
this way $k$ statistics $\Statistics_j : C \to \N$, $1 \leq j \leq k$.
\medbreak

\subsubsection{Colored collections}
\label{subsubsec:colored_collections}
Let $\CFr$ be a finite set, called \Def{set of colors}. A
\Def{$\CFr$-colored collection} $C$ is an $I$-collection such that
\begin{equation}
    I :=
    \left\{
        (a, u) : a \in \CFr
        \mbox{ and } u \in \CFr^\ell \mbox{ for an } \ell \in \N
    \right\}.
\end{equation}
In other terms, any object $x$ of $C$ has an index
$(a, u) \in \CFr \times \CFr^\ell$, $\ell \in \N$, called
\Def{$\CFr$-colored index}. Moreover, the \Def{output color} of $x$ is
$\Out(x) := a$, and the \Def{word of input colors} of $x$ is
$\In(x) := u$. The \Def{$j$th input color} of $x$ is the $j$th letter of
$\In(x)$, denoted by $\In_j(x)$. To not overload the notation, we denote
by $C(a, u)$ the subset $C((a, u))$ of $C$. We say that $C$ is
\Def{monochrome} if $\CFr$ is a singleton. For any nonnegative integer
$k$, a \Def{$k$-colored collection} is a $\CFr$-colored collection where
$\CFr$ is the set of integers $\{1, \dots, k\}$.
\medbreak

\subsubsection{Cyclic collections} \label{subsubsec:cyclic_collections}
Let $C$ be a graded collection endowed for all $n \in \N$ with maps
\begin{equation}
    \CyclicAction_n : C(n) \to C(n)
\end{equation}
such that each $n + 1$st functional power $\CyclicAction_n^{n + 1}$ is
the identity map on $C(n)$. Then, one says that $C$ is a
\Def{cyclic collection} and that the $\CyclicAction_n$, $n \in \N$, are
the \Def{cycle maps} of $C$. Observe that by setting for any $n \in \N$,
$\Action : \Z/_{(n + 1) \Z} \times C(n) \to C(n)$ as the map defined for
any $k \in \Z/_{(n + 1) \Z}$ and $x \in C(n)$ by
$k \Action x := \CyclicAction_n^k(x)$, $\Action$ is a left group action
of the cyclic group of order $n + 1$ on $C(n)$. The reason why we demand
that each $\CyclicAction_n$ is of order $n + 1$ (and not of order $n$)
will appear in the context of cyclic operads.
\medbreak

\subsubsection{Symmetric collections}
\label{subsubsec:symmetric_collections}
Let us first denote by $\SymmetricGroup$ the graded collection of all
the bijections on the set $\{1, \dots, n\}$, $n \in \N$, such that the
size of a bijection is the cardinality of its domain. Let $C$ be a
graded collection endowed, for all $n \in \N$ and
$\sigma \in \SymmetricGroup(n)$, with maps
\begin{equation}
    \SymmetricAction_\sigma : C(n) \to C(n)
\end{equation}
such that $\SymmetricAction_{\Identity_n}$ is the identity map on
$C(n)$, where $\Identity_n$ denotes the identity map of
$\SymmetricGroup(n)$, and
\begin{math}
    \SymmetricAction_{\sigma_1} \circ \SymmetricAction_{\sigma_2}
    = \SymmetricAction_{\sigma_2 \circ \sigma_1}
\end{math}
for any bijections $\sigma_1$ and $\sigma_2$ of $\SymmetricGroup(n)$.
Then, one says that $C$ is a \Def{symmetric collection} and that the
$\SymmetricAction_\sigma$, $\sigma \in \SymmetricGroup(n)$, are the
\Def{symmetric maps} of $C$. Observe that by setting for any $n \in \N$,
$\Action : \SymmetricGroup(n) \times C(n) \to C(n)$ as the map defined
for any $\sigma \in \SymmetricGroup(n)$ and $x \in C(n)$ by
$\sigma \Action x := \SymmetricAction_\sigma(x)$, $\Action$ is a left
group action of the symmetric group of order $n$ on~$C(n)$.
\medbreak

\subsubsection{Collections with products}
\label{subsubsec:collections_with_products}
Let $C$ be an $I$-collection. A \Def{product} on $C$ is a map
\begin{equation} \label{equ:product_on_collection}
   \Product : C\Par{J_1} \times \dots \times C\Par{J_p} \to C
\end{equation}
where $p \in \N$ and $J_1$, \dots, $J_p$ are nonempty subsets of $I$.
The \Def{arity} of $\Product$ is $p$ and the \Def{index domain} of
$\Product$ is the set $J_1 \times \dots \times J_p$. A sequence
$\Par{x_1, \dots, x_p}$ of objects of $C$ is a \Def{valid input} for
$\Product$ whenever $\Product\Par{x_1, \dots, x_p}$ is defined, that is,
\begin{math}
    \Par{\Index\Par{x_1}, \dots, \Index\Par{x_p}}
\end{math}
belongs to the index domain of $\Product$. When $C$ is endowed with a
set of such products, we say that $C$ is an $I$-collection
\Def{with products}. Such a product $\Product$ can be seen as an
operation taking $p$ elements of $C$ as input and outputting one element
of $C$. Let us now review some properties a product $\Product$ of the
form~\eqref{equ:product_on_collection} can satisfy.
\medbreak

First, when the index domain of $\Product$ is $I^p$, $\Product$ is
\Def{complete}. When $J_1 = \dots = J_p = \{i\}$ for a certain index $i$
of $I$, and, for any valid input $\Par{x_1, \dots, x_p}$ for $\Product$,
$\Index\Par{\Product\Par{x_1, \dots, x_p}} = i$, we say that $\Product$
is \Def{internal}. Besides, when there is a map
$\omega : J_1 \times \dots \times J_p \to I$ satisfying, for any valid
input $\Par{x_1, \dots, x_p}$ for $\Product$,
\begin{equation} \label{equ:concentrated_product_on_collection}
    \Product\Par{x_1, \dots, x_p}
    \in C\Par{\omega\Par{\Index\Par{x_1}, \dots, \Index\Par{x_p}}},
\end{equation}
we say that $\Product$ is \Def{$\omega$-concentrated} (or simply
\Def{concentrated} when it not useful to specify $\omega$). In intuitive
terms, this means that the index of the result of a product depends only
of the indexes of its operands. Finally, in the particular case where
$C$ is a graded collection, $\Product$ is \Def{graded} if $\Product$ is
$\omega$-concentrated for the map $\omega : \N^p \to \N$ defined by
$\omega\Par{n_1, \dots, n_p} := n_1 + \dots + n_p$. As a side remark,
the cycle maps (resp. symmetric maps) of a cyclic (resp. symmetric)
collection $C$ are unary internal products on~$C$.
\medbreak

\subsection{Operations over collections}
\label{subsec:operations_collections}
We list here the most important operations that take as input
$I$-collections and output new ones. Some of these are defined only on
graded collections. Table~\ref{tab:operations_collections} shows an
overview of some properties of these operations. Since combinatorial
graded collections have generating series, we provide expressions for
the generating series of the collection produced by the exposed
operations.
\medbreak

\begin{table}[ht]
    \centering
    \begin{tabular}{|c|c|c|c|} \hline
        {\bf Name} & {\bf Arity} & {\bf Inputs} & {\bf Output}
            \\ \hline \hline
        Sum & $2$ & $I$-coll. $C_1$ and $C_2$
            & $I$-coll. $C_1 + C_2$ \\
        Casting & $1$ & $I$-coll $C$ & $J$-coll. $\Casting^\omega(C)$ \\
        Cartesian product & $p \in \N$
            & $I_k$-coll. $C_k$, $1 \leq k \leq p$
            & $I_1 \times \dots \times I_p$-coll.
            $\BBrack{C_1, \dots, C_p}_\times$ \\ \hline
        Hadamard product & $p \in \N$
            & $I$-coll. $C_1$, \dots, $C_p$
            & $I$-coll. $\BBrack{C_1, \dots, C_p}_\Hadamard$
            \\ \hline
        List collection & $1$ & $I$-coll. $C$
            & $\SetTuples(I)$-coll. $\List(C)$ \\
        Multiset collection & $1$ & $I$-coll. $C$
            & $\SetMultisets(I)$-coll. $\Multiset(C)$ \\
        Set collection & $1$ & $I$-coll. $C$
            & $\SetSets(I)$-coll. $\Set(C)$ \\ \hline
        $\ell$-suspension & $1$ & graded coll. $C$
            & graded coll. $\Suspension_\ell(C)$ \\
        Augmentation & $1$ & graded coll. $C$
            & graded coll. $\Augmentation(C)$ \\ \hline
        Composition & $2$ & graded coll. $C_1$ and $C_2$
            & graded coll. $C_1 \Compo C_2$ \\ \hline
        $\CFr$-coloration & $1$ & graded coll. $C$
            & $\CFr$-colored coll. $\Coloration_\CFr(C)$ \\
        Cycle & $1$ & graded coll. $C$
            & Cyclic coll. $\Cycle(C)$ \\
        Regularization & $1$ & graded coll. $C$
            & symmetric coll. $\Regularization(C)$ \\ \hline
    \end{tabular}
    \bigbreak

    \caption[Operations over collections.]
    {\footnotesize Main properties of some operations over collections.
    Here, $I$, $J$, and $I_1$, \dots, $I_p$, $p \in \N$, are index sets,
    $\CFr$ is a set of colors, $\SetTuples(I)$ is the set of all the
    finite tuples of elements of $I$, $\SetMultisets(I)$ is the set of
    all the finite multisets of elements of $I$, $\SetSets(I)$ is the
    set of all the finite subsets of $I$, and $\ell$ is an integer.}
    \label{tab:operations_collections}
\end{table}
\medbreak

\subsubsection{Sum operation} \label{subsubsec:sum_collections}
Let $C_1$ and $C_2$ be two $I$-collections. The \Def{sum} of $C_1$ and
$C_2$ is the $I$-collection $C_1 + C_2$ such that, for all $i \in I$,
\begin{equation} \label{equ:sum_operation}
    \Par{C_1 + C_2}(i) := C_1(i) \sqcup C_2(i).
\end{equation}
In other words, each object of index $i$ of $C_1 + C_2$ is either an
object of index $i$ of $C_1$ or an object of index $i$ of $C_2$. Since
the sum operation~\eqref{equ:sum_operation} is defined through a
disjoint union, when the sets $C_1(i)$ and $C_2(i)$ are not disjoint,
there are in $\Par{C_1 + C_2}(i)$ two copies of each element belonging
to the intersection $C_1(i) \cap C_2(i)$, one coming from $C_1(i)$, the
other from $C_2(i)$. Moreover, observe that the sum operation admits the
empty $I$-collection $\emptyset$ as unit and that it is associative and
commutative. The iterated version of the operation $+$ shall be denoted
by $\bigsqcup$ in the sequel.
\medbreak

When $C_1$ and $C_2$ are combinatorial, $C_1 + C_2$ is combinatorial.
Moreover, when $C_1$ and $C_2$ are combinatorial and graded, the
generating series of $C_1 + C_2$ satisfies
\begin{equation} \label{equ:generating_series_sum}
    \GeneratingSeries_{C_1 + C_2}(t)
    = \GeneratingSeries_{C_1}(t) + \GeneratingSeries_{C_2}(t).
\end{equation}
\medbreak

\subsubsection{Casting operation} \label{subsubsec:casting_collections}
Let $C$ be an $I$-collection, $J$ be an index set, and
$\omega : I \to J$ be a map. The \Def{$\omega$-casting} of $C$ is the
$J$-collection $\Casting^\omega(C)$ defined for any $j \in J$ by
\begin{equation} \label{equ:casting_collections}
    \Par{\Casting^\omega(C)}(j) :=
    \bigsqcup_{\substack{
        i \in I \\
        \omega(i) = j
    }}
    C(i).
\end{equation}
In other words, each object of index $j \in J$ of $\Casting^\omega(C)$
comes from an object of index $i \in I$ of $C$ such that
$\omega(i) = j$. Observe also that the right member
of~\eqref{equ:casting_collections} is equal to~$C\Par{\omega^{-1}(j)}$.
\medbreak

When the codomain of $\omega$ is $\N$, $\Casting^\omega(C)$ is a graded
collection  called \Def{$\omega$-graduation} of $C$. Let us detail two
particular cases of $\omega$-graduations. When $C$ is a $k$-graded
collection, for any $1 \leq p \leq k$ we call \Def{$p$-graduation} of
$C$ the $\pi_p$-graduation of $C$ for the map $\pi_p : \N^k \to \N$
defined by
\begin{math}
    \pi_p\Par{\Par{n_1, \dots, n_k}} := n_p.
\end{math}
Besides, when $C$ is a $\CFr$-colored collection where $\CFr$ is a set
of colors, we call \Def{graduation} of $C$ the $\omega$-graduation of
$C$ for the map $\omega$ sending any $\CFr$-colored index $(a, u)$ to
the length of the tuple $u$.
\medbreak

Observe that when $C$ is combinatorial and each fiber $\omega^{-1}(j)$
is finite for any $j \in J$, $\Casting^\omega(C)$ is combinatorial.
\medbreak

\subsubsection{Cartesian product operation}
\label{subsubsec:cartesian_product_collections}
Let $p \in \N$, $I_1$, \dots, $I_p$ be index sets, and $C_1$ be an
$I_1$-collection, \dots, $C_p$ be an $I_p$-collection. The
\Def{Cartesian product} of $C_1$, \dots, $C_p$ is the
$I_1 \times \dots \times I_p$-collection
$\BBrack{C_1, \dots, C_p}_\times$ such that, for all
\begin{math}
    \Par{i_1, \dots, i_p} \in I_1 \times \dots \times I_p,
\end{math}
\begin{equation} \label{equ:cartesian_product_operation}
    \BBrack{C_1, \dots, C_p}_\times
    \Par{\Par{i_1, \dots, i_p}} :=
    C_1\Par{i_1} \times \dots \times C_p\Par{i_p}.
\end{equation}
In other words, each object of index $\Par{i_1, \dots, i_p}$ of
$\BBrack{C_1, \dots, C_p}_\times$ is a tuple $\Par{x_1, \dots, x_p}$
such that for any $1 \leq k \leq p$, each $x_k$ is an object of index
$i_k$ of $C_k$. As a special but important case, the Cartesian product
$\BBrack{}_\times$ of zero $I$-collections is the $I$-collection
containing the empty tuple (that is, the unique tuple of length~$0$).
The index of this object is the empty tuple on $I$ (that is, the unique
element of $I^0$). To not overload the notation, we denote by
\begin{math}
    \BBrack{C_1, \dots, C_p}_\times\Par{i_1, \dots, i_p}
\end{math}
the set
\begin{math}
    \BBrack{C_1, \dots, C_p}_\times \Par{\Par{i_1, \dots, i_p}}
\end{math}
for any index $\Par{i_1, \dots, i_p}$ of~$I_1 \times \dots \times I_p$.
\medbreak

Observe that when all the $C_1$, \dots, $C_p$ are combinatorial,
$\BBrack{C_1, \dots, C_p}_\times$ is combinatorial. Besides, when all
the $C_1$, \dots, $C_p$ are graded, $\BBrack{C_1, \dots, C_p}_\times$ is
a $p$-graded collection.
\medbreak

When $J$ is an index set and
\begin{math}
    \omega : I_1 \times \dots \times I_p \to J
\end{math}
is a map, the \Def{$\omega$-Cartesian product operation} of $C_1$,
\dots, $C_p$ is the $J$-collection
\begin{equation}
    \BBrack{C_1, \dots, C_p}_\times^\omega :=
    \Casting^\omega\Par{\BBrack{C_1, \dots, C_p}_\times}.
\end{equation}
By definition of the casting and the Cartesian product operations, one
has, for any $j \in J$,
\begin{equation}
    \BBrack{C_1, \dots, C_p}_\times^\omega(j)
    =
    \left\{
    \Par{x_1, \dots, x_p} \in C_1 \times \dots \times C_p :
    \omega\Par{\Index\Par{x_1}, \dots, \Index\Par{x_p}} = j
    \right\}.
\end{equation}
In other words, each object of index $j \in J$ of
$\BBrack{C_1, \dots, C_p}_\times^\omega$ is a tuple
$\Par{x_1, \dots, x_p}$ such that for any $1 \leq k \leq p$, each $x_k$
is an object of $C_k$, and the image by $\omega$ of the tuple formed by
the indexes of $x_1$, \dots, $x_p$ is~$j$.
\medbreak

Observe that when all the $C_1$, \dots, $C_p$ are combinatorial and each
fiber $\omega^{-1}(j)$ is finite for any $j \in J$,
$\BBrack{C_1, \dots, C_p}_\times^\omega$ is combinatorial. When all the
index sets $I_1$, \dots, $I_p$, and $J$ are equal to $\N$, we denote by
$+ : \N^p \to \N$ the map defined by
\begin{math}
    +\Par{\Par{n_1, \dots, n_p}} := n_1 + \dots + n_p.
\end{math}
When all the $C_1$, \dots, $C_p$ are combinatorial and graded,
$\BBrack{C_1, \dots, C_p}_\times^+$ is combinatorial and its generating
series satisfies
\begin{equation} \label{equ:generating_series_multiplication}
    \GeneratingSeries_{\BBrack{C_1, \dots, C_p}_\times^+}(t)
    = \prod_{1 \leq k \leq p} \GeneratingSeries_{C_k}(t).
\end{equation}
\medbreak

\subsubsection{Hadamard product operation}
\label{subsubsec:hadamard_product_collections}
Let $p \in \N$ and $C_1$, \dots, $C_p$ be $I$-collections. The
\Def{Hadamard product} of $C_1$, \dots, $C_p$ is the $I$-collection
$\BBrack{C_1, \dots, C_p}_\Hadamard$ such that, for all $i \in I$,
\begin{equation} \label{equ:hadamard_product_operation}
    \BBrack{C_1, \dots, C_p}_\Hadamard(i) :=
    \BBrack{C_1, \dots, C_p}_\times
    \Par{\underbrace{i, \dots, i}_{p \mbox{ \footnotesize terms}}}.
\end{equation}
In other words, each object of index $i$ of
$\BBrack{C_1, \dots, C_p}_\Hadamard$ is a tuple $\Par{x_1, \dots, x_p}$
such that for all $1 \leq k \leq p$, all the $x_k$ are objects of index
$i$ of $C_k$.
\medbreak

Observe that when all the $C_1$, \dots, $C_p$ are combinatorial,
$\BBrack{C_1, \dots, C_p}_\Hadamard$ is combinatorial. When all the
$C_1$, \dots, $C_p$ are combinatorial and graded,
$\BBrack{C_1, \dots, C_p}_\Hadamard$ is combinatorial and its generating
series satisfies
\begin{equation} \label{equ:generating_series_hadamard_product}
    \GeneratingSeries_{\BBrack{C_1, \dots, C_p}_\Hadamard}(t)
    = \sum_{n \in \N}
    \Par{\prod_{1 \leq k \leq p} \# C_k(n)} t^n.
\end{equation}
\medbreak

\subsubsection{List collection operation}
\label{subsubsec:list_collections}
Let $C$ be an $I$-collection. Let us denote by $\SetTuples(I)$ the index
set of all the finite tuples $\Par{i_1, \dots, i_p}$, $p \in \N$, of
elements of $I$. The \Def{list collection} of $C$ is the
$\SetTuples(I)$-collection $\List(C)$ such that
\begin{equation}
    \List(C) := \bigsqcup_{p \in \N}
    \BBrack{\underbrace{C, \dots, C}_{p \mbox{ \footnotesize terms}}}
    _\times
\end{equation}
In other words, each object of $\List(C)$ is a tuple
$\Par{x_1, \dots, x_p}$ of objects of $C$ and its index is the element
$\Par{\Index\Par{x_1}, \dots, \Index\Par{x_p}}$ of $\SetTuples(I)$.
Moreover, for any subset $S$ of $\N$, let $\List_S(C)$ be the
subcollection of $\List(C)$ restrained on tuples that have a length
in~$S$.
\medbreak

Observe that when $C$ is combinatorial, $\List(C)$ is combinatorial.
\medbreak

When $J$ is an index set and $\omega : \SetTuples(I) \to J$ is a map,
the \Def{$\omega$-list collection} of $C$ is the $J$-collection
\begin{equation}
    \List^\omega(C) := \Casting^\omega\Par{\List(C)}.
\end{equation}
By definition of the casting and the list collection operations, one
has, for any $j \in J$,
\begin{equation}
    \Par{\List^\omega(C)}(j) =
    \left\{
    \Par{x_1, \dots, x_p} \in C^p : p \in \N \mbox{ and }
    \omega\Par{\Par{\Index\Par{x_1}, \dots, \Index\Par{x_p}}} = j
    \right\}.
\end{equation}
In other words, each object of index $j \in J$ of $\List^\omega(C)$ is a
tuple $\Par{x_1, \dots, x_p}$, $p \in \N$, such that for all
$1 \leq k \leq p$, the $x_k$ are objects of $C$, and the image by
$\omega$ of the tuple formed by the indexes of $x_1$, \dots, $x_p$ is
$j$. Moreover, for any subset $S$ of $\N$, let $\List^\omega_S(C)$ be
the subcollection of $\List^\omega(C)$ restrained on tuples that have a
length in $S$.
\medbreak

Observe that when $C$ is combinatorial and each fiber $\omega^{-1}(j)$
is finite for any $j \in J$, $\List^\omega(C)$ is combinatorial. When
$I$ and $J$ are equal to $\N$, we denote by $+ : \SetTuples(\N) \to \N$
the map defined by $+\Par{\Par{n_1, \dots, n_p}} := n_1 + \dots + n_p$.
When $C$ is combinatorial and graded, $\List^+(C)$ is combinatorial if
and only if $C$ is augmented. In this case, its generating series
satisfies
\begin{equation} \label{equ:generating_series_list}
    \GeneratingSeries_{\List^+(C)}(t)
    = \frac{1}{1 - \GeneratingSeries_C(t)}.
\end{equation}
\medbreak

\subsubsection{Multiset collection operation}
\label{subsubsec:multiset_operation}
Let $C$ be an $I$-collection. Let us denote by $\SetMultisets(I)$ the
index set formed by all finite multisets $\lbag i_1, \dots, i_p \rbag$,
$p \in \N$, of elements of $I$. The \Def{multiset collection} of $C$ is
the $\SetMultisets(I)$-collection $\Multiset(C)$ such that, for all
$\lbag i_1, \dots, i_p \rbag \in \SetMultisets(I)$, $p \in \N$,
\begin{equation}
    \Par{\Multiset(C)}\Par{\lbag i_1, \dots, i_p \rbag} :=
    \Par{
    \bigcup_{\sigma \in \SymmetricGroup(p)}
    \Par{\List(C)}\Par{i_{\sigma(1)}, \dots, i_{\sigma(p)}}}/_\EqRel,
\end{equation}
where $\SymmetricGroup$ is the graded collection of bijections defined
in Section~\ref{subsubsec:symmetric_collections} and $\EqRel$ is the
equivalence relation on $\List(C)$ satisfying, for any
$\Par{x_1, \dots, x_p} \in \List(C)$, $p \in \N$, and any bijection
$\sigma \in \SymmetricGroup(p)$,
\begin{equation}
    \Par{x_1, \dots, x_p}
    \EqRel
    \Par{x_{\sigma(1)}, \dots, x_{\sigma(p)}}.
\end{equation}
In other words, each object of $\Multiset(C)$ is an
$\EqRel$-equivalence class of tuples of $\List(C)$ and such an
$\EqRel$-equivalence class $\left[\Par{x_1, \dots, x_p}\right]_\EqRel$
can be represented by the multiset $\lbag x_1, \dots, x_p \rbag$.
Therefore, the objects of $\Multiset(C)$ can be regarded as multisets of
objects of~$C$.
\medbreak

When $J$ is an index set and $\omega : \SetMultisets(I) \to J$ is a map,
the \Def{$\omega$-multiset collection} of $C$ is the $J$-collection
\begin{equation}
    \Multiset^\omega(C) :=
    \Casting^\omega(\Multiset(C)).
\end{equation}
By definition of the casting and the multiset collection operations, one
has, for any $j \in J$,
\begin{equation}
    \Par{\Multiset^\omega(C)}(j) =
    \left\{
    \lbag x_1, \dots, x_p \rbag :
    p \in \N \mbox{ and }
    \omega\Par{\lbag \Index\Par{x_1}, \dots, \Index\Par{x_p}\rbag}
    = j
    \right\}.
\end{equation}
In other words, each object of index $j \in J$ of $\Multiset^\omega(C)$
is a finite multiset $\lbag x_1, \dots, x_p \rbag$, $p \in \N$, such
that for all $1 \leq k \leq p$, the $x_k$ are objects of $C$, and the
image by $\omega$ of the multiset formed by the indexes of $x_1$, \dots,
$x_p$ is $j$.
\medbreak

Observe that when $C$ is combinatorial and each fiber $\omega^{-1}(j)$
is finite for any $j \in J$, $\Multiset^\omega(C)$ is combinatorial.
When $I$ and $J$ are equal to $\N$, we denote by
$+ : \SetMultisets(\N) \to \N$ the map defined by
\begin{math}
    +\Par{\lbag n_1, \dots, n_p \rbag} := n_1 + \dots + n_p.
\end{math}
When $C$ is combinatorial and graded, $\Multiset^+(C)$ is combinatorial
if and only if $C$ is augmented. In this case, its generating series
satisfies
\begin{equation} \label{equ:generating_series_multiset}
    \GeneratingSeries_{\Multiset^+(C)}(t)
    = \prod_{n \in \N \setminus \{0\}}
    \Par{\frac{1}{1 - t^n}}^{\# C(n)}.
\end{equation}
\medbreak

\subsubsection{Set collection operation}
\label{subsubsec:set_collections}
Let $C$ be an $I$-collection. Let us denote by $\SetSets(I)$ the index
set formed by all finite subsets of elements of $I$. The
\Def{set collection} of $C$ is the $\SetSets(I)$-collection $\Set(C)$
defined as the subcollection of $\Multiset(C)$ containing only the
multisets having all their elements with multiplicity~$1$. In this way,
the objects of $\Set(C)$ can be represented as finite sets of objects
of~$C$.
\medbreak

When $J$ is an index set and $\omega : \SetSets(I) \to J$ is a map, the
\Def{$\omega$-set collection} of $C$ is the $J$-collection
\begin{equation}
    \Set^\omega(C) := \Casting^\omega(\Set(C)).
\end{equation}
By definition of the casting and the set collection operations, one has,
for any $j \in J$,
\begin{equation}
    \Set^\omega(C)(j) =
    \left\{\{x_1, \dots, x_p\} \subseteq C :
    p \in \N \mbox{ and }
    \omega\Par{\left\{\Index\Par{x_1}, \dots, \Index\Par{x_p}\right\}}
    = j \right\}.
\end{equation}
In other words, each object of index $j \in J$ of $\Set^\omega(C)$ is a
finite set $\left\{x_1, \dots, x_p\right\}$, $p \in \N$, such that for
all $1 \leq k \leq p$, the $x_k$ are objects of $C$, and the image by
$\omega$ of the set formed by the indexes of $x_1$, \dots, $x_p$ is $j$.
\medbreak

Observe that when $C$ is combinatorial and each fiber $\omega^{-1}(j)$
is finite for any $j \in J$, $\Set^\omega(C)$ is combinatorial. When $I$
and $J$ are equal to $\N$, we denote by $+ : \SetSets(\N) \to \N$ the
map defined by
\begin{math}
    +\Par{\left\{x_1, \dots, x_p \right\}} := x_1 + \dots + x_p.
\end{math}
When $C$ is combinatorial and graded, $\Set^\omega(C)$ is combinatorial
(without requiring any additional condition contrariwise to the similar
cases for the list and multiset collection operations) and its
generating series satisfies
\begin{equation} \label{equ:generating_series_set}
    \GeneratingSeries_{\Set^+(C)}(t)
    = \prod_{n \in \N \setminus \{0\}}
    \Par{1 + t^n}^{\# C(n)}.
\end{equation}
\medbreak

\subsubsection{Suspension and augmentation operations}
\label{subsubsec:suspension_collections}
Let $C$ be a graded collection. For any $\ell \in \Z$, the
\Def{$\ell$-suspension} of $C$ is the graded collection
$\Suspension_\ell(C)$ such that, for all $n \in \N$,
\begin{equation} \label{equ:suspension_operation}
    \Par{\Suspension_\ell(C)}(n) :=
    \begin{cases}
        C(n - \ell) & \mbox{if } n - \ell \in \N, \\
        \emptyset & \mbox{otherwise}.
    \end{cases}
\end{equation}
Observe that $\Suspension_1\Par{\Suspension_{-1}(C)}$ is the
subcollection $C \setminus C(0)$ of $C$, that is the augmented
collection having the objects of $C$ without its objects of size $0$. We
call this collection the \Def{augmentation} of~$C$ and we denote it
by~$\Augmentation(C)$.
\medbreak

When $C$ is combinatorial, $\Suspension_\ell(C)$ and $\Augmentation(C)$
are combinatorial and their generating series satisfy, respectively,
\begin{equation}
    \GeneratingSeries_{\Suspension_\ell(C)}(t) =
    t^\ell \,
    \Par{\GeneratingSeries_{C}(t)
        - \GeneratingSeries_{C([0, -\ell - 1])}(t)}
\end{equation}
where $[0, -\ell - 1]$ is the set of the integers $n$ satisfying
$0 \leq n \leq -\ell - 1$, and
\begin{equation}
    \GeneratingSeries_{\Augmentation(C)}(t) =
    \GeneratingSeries_{C}(t) - \#C(0).
\end{equation}
\medbreak

\subsubsection{Composition operation}
\label{subsubsec:composition_collections}
Let $C_1$ and $C_2$ be two graded collections. The \Def{composition} of
$C_1$ and $C_2$ is the graded collection $C_1 \Compo C_2$ such that, for
all $n \in \N$,
\begin{equation} \label{equ:composition_operation}
    \Par{C_1 \Compo C_2}(n) :=
    \bigsqcup_{k \in \N}
    \BBrack{C_1(k),
    \Par{\List^+_{\{k\}}\Par{C_2}}(n)}_\times^\omega
\end{equation}
where $\omega : \N^2 \to \N$ is the map defined by
$\omega\Par{n_1, n_2} := n_2$. In other words, each object of size $n$
of $C_1 \Compo C_2$ is an ordered pair $\Par{x, \Par{y_1, \dots, y_k}}$,
$k \in \N$, where $x$ is an object of $C_1$ of size $k$, and
$\Par{y_1, \dots, y_k}$ is a tuple of objects of $C_2$ such that the sum
of the sizes of the $y_j$, $1 \leq j \leq k$ is $n$. Observe that, if
$C_3$ is a graded collection,
\begin{subequations}
\begin{equation} \label{equ:relation_unit_composition_collection}
    C_1 \Compo \{\AtomElement\} \simeq C_1 \simeq
    \{\AtomElement\} \Compo C_1,
\end{equation}
\begin{equation}
    \Par{C_1 \Compo C_2} \Compo C_3
    \simeq C_1 \Compo \Par{C_2 \Compo C_3}.
\end{equation}
\end{subequations}
\medbreak

When $C_1$ and $C_2$ are combinatorial and graded, $C_1 \Compo C_2$ is
combinatorial if and only if $C_2$ is augmented. In this case, its
generating series satisfies
\begin{equation} \label{equ:generating_series_composition}
    \GeneratingSeries_{C_1 \Compo C_2}(t) =
    \GeneratingSeries_{C_1}\Par{\GeneratingSeries_{C_2}(t)}.
\end{equation}
\medbreak

\subsubsection{Coloration operation}
\label{subsubsec:coloration_collections}
Let $C$ be a graded collection and $\CFr$ be a set of colors. The
\Def{$\CFr$-coloration} of $C$ is the $\CFr$-colored collection
$\Coloration_\CFr(C)$ defined, for all $\CFr$-colored indexes
$(a, u) \in \CFr \times \CFr^\ell$, $\ell \in \N$, by
\begin{equation}
    \Par{\Coloration_\CFr(C)}(a, u) :=
    \left\{(a, x, u) : x \in C(\ell)\right\}.
\end{equation}
In other words, each object of $\Coloration_\CFr(C)$ is built from an
object $x$ of $C$ by equipping it freely with an output color from
$\CFr$ and a word of input colors from $\CFr$ having the size of $x$ as
length.
\medbreak

When $C$ is combinatorial, the graduation
$\Casting^\omega\Par{\Coloration_\CFr(C)}$ of $\Coloration_\CFr(C)$ is
combinatorial if and only if the set of colors $\CFr$ is finite. In this
case, by setting $m := \# \CFr$, its generating series satisfies
\begin{equation} \label{equ:generating_series_coloration}
    \GeneratingSeries_{\Casting^\omega\Par{\Coloration_\CFr(C)}}(t)
    = \sum_{n \in \N} \# C(n) m^{n + 1} \, t^n.
\end{equation}
\medbreak

\subsubsection{Cycle operation} \label{subsubsec:cycle_collections}
Let $C$ be a graded collection. The \Def{cycle collection} of $C$ is
the graded collection $\Cycle(C)$ defined, for all $n \in \N$, by
\begin{equation}
    \Par{\Cycle(C)}(n) :=
    \left\{(x, k) : x \in C(n) \mbox{ and } 0 \leq k \leq n\right\}.
\end{equation}
In other words, each object of $\Cycle(C)$ is built from an object $x$
of $C$  by equipping it freely with a nonnegative integer nongreater
than its size.
\medbreak

Let us observe that by defining, for any $n \in \N$, the map
\begin{math}
    \CyclicAction_n : (\Cycle(C))(n) \to (\Cycle(C))(n)
\end{math}
by
\begin{math}
    \CyclicAction_n((x, k)) := (x, k + 1 \pmod{n + 1})
\end{math}
for any $(x, k) \in (\Cycle(C))(n)$, each $\CyclicAction_n$ is a cycle
map of $\Cycle(C)$. Therefore $\Cycle(C)$ is cyclic.
\medbreak

When $C$ is combinatorial, $\Cycle(C)$ is combinatorial and its
generating series satisfies
\begin{equation}
    \GeneratingSeries_{\Cycle(C)}(t) =
    \sum_{n \in \N}  \# C(n) (n + 1) \, t^n.
\end{equation}
\medbreak

\subsubsection{Regularization operation}
\label{subsubsec:regularization}
Let $C$ be a graded collection. The \Def{regularization} of $C$ is the
graded collection $\Regularization(C)$ defined by
\begin{equation}
    \Regularization(C) := \BBrack{C, \SymmetricGroup}_\Hadamard
\end{equation}
where $\SymmetricGroup$ is the graded collection of bijections defined
in Section~\ref{subsubsec:symmetric_collections}. In other words, each
object of $\Regularization(C)$ is built from an object $x$ of $C$ by
equipping it freely with a bijection of $\SymmetricGroup(n)$ where~$n$
is the size of~$x$.
\medbreak

Let us observe that by defining, for any
$\sigma \in \SymmetricGroup(n)$, $n \in \N$, the map
\begin{math}
    \SymmetricAction_\sigma :
    (\Regularization(C))(n) \to (\Regularization(C))(n)
\end{math}
by
\begin{math}
    \SymmetricAction_\sigma((x, \nu)) := \Par{x, \sigma^{-1} \circ \nu}
\end{math}
for any $(x, \nu) \in (\Regularization(C))(n)$, each
$\SymmetricAction_\sigma$ is a symmetric map of $\Regularization(C)$.
Therefore, $\Regularization(C)$ is symmetric.
\medbreak

When $C$ is combinatorial, $\Regularization(C)$ is combinatorial and
its generating series satisfies
\begin{equation}
    \GeneratingSeries_{\Regularization(C)}(t)
    = \sum_{n \in \N} \# C(n) n! \, t^n.
\end{equation}
\medbreak

\subsection{Examples}
We define, in some cases by using the operations of
Section~\ref{subsec:operations_collections}, some usual graded
combinatorial collections. At the same time, we set here our main
notations and definitions about their objects.
\medbreak

\subsubsection{Natural numbers} \label{subsubsec:natural_numbers}
We can regard the set $\N$ as the graded collection satisfying
$\N(n) := \{n\}$ for all $n \in \N$. Hence,
$\List^+(\{\AtomElement\}) \simeq \N$. Moreover, for any $\ell \in \N$,
let $\N_{\geq \ell}$ be the graded collection defined by
\begin{equation}
    \N_{\geq \ell} := \Suspension_\ell\Par{\Suspension_{-\ell}(\N)}.
\end{equation}
By definition of the suspension operation over graded collections,
$\N_{\geq \ell}$ is the set of all integers greater than or equal to
$\ell$. Observe that $\N_{\geq 1} = \Augmentation(\N)$. The generating
series of $\N_{\geq \ell}$ satisfies
\begin{equation} \label{equ:generating_series_natural_numbers}
    \GeneratingSeries_{\N_{\geq \ell}}(t) = \frac{t^\ell}{1 - t}
    =
    \sum_{n \in \N_{\geq \ell}} t^n.
\end{equation}
Observe also that the list collection operation over graded collections
can be expressed as a composition involving $\N$ since
\begin{equation}
    \List^+(C) \simeq \N \Compo C
\end{equation}
for any graded collection $C$. We shall consider in the sequel, for any
$x, z \in \N$, the subcollections
\begin{math}
    [x, z] := \left\{y \in \N : x \leq y \leq z\right\},
\end{math}
and $[x] := [1, x]$ of $\N$. These examples of graded collections are
among the simplest nontrivial ones.
\medbreak

\subsubsection{Words} \label{subsubsec:words}
Let $A$ be an \Def{alphabet}, that is a set whose elements are called
\Def{letters}. One can see $A$ as a graded collection wherein all
letters are atoms. In this case, we denote by $A^*$ the graded
collection $\List^+(A)$. By definition, the objects of $A^*$ are finite
sequences of elements of $A$. We call \Def{words} on $A$ these objects.
When $A$ is finite, $A^*$ is combinatorial and it follows
from~\eqref{equ:generating_series_list} that the generating series of
$A^*$ satisfies
\begin{equation}
    \GeneratingSeries_{A^*}(t) = \frac{1}{1 - mt}
    = \sum_{n \in \N} m^n \, t^n
\end{equation}
where $m := \# A$. If $u := \Par{a_1, \dots, a_n}$ is a word on $A$, it
follows from the definition of $A^*$ that the size $|u|$ of $u$ is $n$.
The unique word on $A$ of size $0$ is denoted by $\epsilon$ and is
called \Def{empty word}.
\medbreak

Let $u := \Par{a_1, \dots, a_n}$ be a word on $A$. The
\Def{$i$th letter} of $u$ is $a_i$ and is denoted by $u(i)$. For any
letter $\Bsf \in A$, the \Def{number of occurrences} $|u|_\Bsf$ of
$\Bsf$ in $u$ is the cardinality of the set
\begin{math}
    \left\{i \in \left[|u|\right] : u(i) = \Bsf\right\}.
\end{math}
When $A$ is endowed with a total order $\Ord$ and $u$ is nonempty,
$\max_{\Ord}(u)$ is the greatest letter appearing in $u$ with respect to
$\Ord$. Moreover, an \Def{inversion} of $u$ is a pair $(i, j)$ such that
$i < j$, $u(i) \ne u(j)$, and $u(j) \Ord u(i)$. Given two words $u$ and
$v$ on $A$, the \Def{concatenation} of $u$ and $v$ is the word
$u \Conc v$ containing from left to right the letters of $u$ and then
the ones of $v$. The concatenation $\Conc$ is a graded complete product
on $A^*$. If $u$ can be expressed as $u = u_1 \Conc u_2$ where
$u_1, u_2 \in A^*$, we say that $u_1$ is a \Def{prefix} of $u$ and this
property is denoted by $u \PrefixOrder v$. For any subset
$S := \left\{s_1 \leq \cdots \leq s_k\right\}$ of $[|u|]$, $u_{|S}$ is
the word $u\Par{s_1} \dots u\Par{s_k}$. Moreover, when $v$ is a word
such that there exists $S \subseteq [|u|]$ satisfying $v = u_{|S}$, $v$
is a \Def{subword} of~$u$.
\medbreak

A \Def{language} on $A$ is subcollection of $A^*$. A language $\Lca$ on
$A$ is \Def{prefix} if for all $u \in \Lca$ and $v \in A^*$,
$v \PrefixOrder u$ implies~$v \in \Lca$. We denote by $A^+$ the language
$\Augmentation(A^*)$ containing all nonempty words on $A$, and, for any
$n \in \N$, by $A^n$ the language $A^*(n)$.
\medbreak

\subsubsection{Integer compositions}
\label{subsubsec:integer_compositions}
By regarding the set $\N$ as a graded collection as explained in
Section~\ref{subsubsec:natural_numbers}, let $\ColComp$ be the
combinatorial graded collection $\List^+\Par{\N_{\geq 1}}$. It follows
from~\eqref{equ:generating_series_list}
and~\eqref{equ:generating_series_natural_numbers} that the generating
series of $\ColComp$ is
\begin{equation}
    \GeneratingSeries_{\ColComp}(t)
    = \frac{1 - t}{1 - 2t}
    = 1 + \sum_{n \in \N_{\geq 1}} 2^{n - 1} \, t^n.
\end{equation}
Hence, the integer sequence of $\ColComp$ begins by
\begin{equation}
    1, 1, 2, 4, 8, 16, 32, 64, 128
\end{equation}
and is Sequence~\OEIS{A011782} of~\cite{Slo}. By definition, the objects
of $\ColComp$ are finite sequences of positive numbers. We call
\Def{integer compositions} (or, for short, \Def{compositions}) these
objects. If $\LambdaB := \Par{\LambdaB_1, \dots, \LambdaB_k}$ is a
composition, it follows from the definition of $\ColComp$ that the size
$|\LambdaB|$ of $\LambdaB$ is $\LambdaB_1 + \dots + \LambdaB_k$. The
\Def{length} $\Length(\LambdaB)$ of $\LambdaB$ is $k$, and for any
$i \in [\Length(\LambdaB)]$, the \Def{$i$th part} of $\LambdaB$ is
$\LambdaB_i$. The unique composition of size $0$ is denoted by
$\epsilon$ and is called \Def{empty composition} (even if $\epsilon$ is
already used to express the empty word, this overloading of notation is
not a problem in practice).
\medbreak

The \Def{descents set} of $\LambdaB$ is the set
\begin{equation}
    \Des(\LambdaB) :=
    \left\{\LambdaB_1, \LambdaB_1 + \LambdaB_2, \dots,
    \LambdaB_1 + \LambdaB_2 + \dots + \LambdaB_{k - 1}\right\}.
\end{equation}
For instance, $\Des(4131) = \{4, 5, 8\}$. Moreover, for any word $u$
defined on an alphabet $A$ equipped with a total order $\Ord$, the
\Def{composition} $\Cmp(u)$ of $u$ is the composition of size $|u|$
defined by
\begin{equation}
    \Cmp(u) := \Par{\left|u_1\right|, \dots, \left|u_k\right|},
\end{equation}
where $u = u_1 \Conc \dots \Conc u_k$ is the factorization of $u$ in
longest nondecreasing factors (with respect to the order $\Ord$). For
instance, if
\begin{math}
    u := \ColA{a_2 a_2 a_3} \ColD{a_1 a_3} \ColA{a_2} \ColD{a_1 a_2}
\end{math}
is a word on the alphabet $A := \{a_1, a_2, a_3\}$ ordered by
$a_1 \Ord a_2 \Ord a_3$, $\Cmp(u) = 3212$. When $\# A \geq 2$, this map
$\Cmp$ is a surjective collection morphism from $A^*$ to~$\ColComp$.
\medbreak

Integer compositions are drawn as \Def{ribbon diagrams} in the following
way. For each part $\LambdaB_i$ of $\LambdaB$, we draw a horizontal line
of $\LambdaB_i$ boxes. These lines are organized so that the line for
the first part of $\LambdaB$ is the uppermost, and the first box of the
line of the part $\LambdaB_{i + 1}$ is glued below the last box of the
line of the part $\LambdaB_i$, for all $i \in [\Length(\LambdaB) - 1]$.
For instance, the ribbon diagram of the composition $4131$ is
\begin{equation}
    \begin{tikzpicture}[Centering]
        \node[Box](1)at(0,0){};
        \node[Box,right of=1,node distance=0.25cm](2){};
        \node[Box,right of=2,node distance=0.25cm](3){};
        \node[Box,right of=3,node distance=0.25cm](4){};
        \node[Box,below of=4,node distance=0.25cm](5){};
        \node[Box,below of=5,node distance=0.25cm](6){};
        \node[Box,right of=6,node distance=0.25cm](7){};
        \node[Box,right of=7,node distance=0.25cm](8){};
        \node[Box,below of=8,node distance=0.25cm](9){};
    \end{tikzpicture}\,.
\end{equation}
\medbreak

\subsubsection{Integer partitions} \label{subsubsec:integer_partitions}
Again by regarding the set $\N$ as a graded collection as considered in
Section~\ref{subsubsec:natural_numbers}, let $\ColIP$ be the graded
combinatorial collection $\Multiset^+\Par{\N_{\geq 1}}$. Since
$\# \N_{\geq 1}(n) = 1$ for all $n \geq 1$, it follows
from~\eqref{equ:generating_series_multiset} that the generating series
of $\ColIP$ is
\begin{equation}
    \GeneratingSeries_{\ColIP}(t)
    = \prod_{n \in \N_{\geq 1}} \frac{1}{1 - t^n}.
\end{equation}
Hence, the integer sequence of $\ColIP$ begins by
\begin{equation}
    1, 1, 2, 3, 5, 7, 11, 15, 22
\end{equation}
and is Sequence~\OEIS{A000041} of~\cite{Slo}. By definition, the objects
of $\ColIP$ are finite multisets of positive integers. We call
\Def{integer partitions} (or, for short, \Def{partitions}) these
objects. As a consequence of the definition of $\ColIP$, the size
$|\lambda|$ of any partition $\lambda$ is the sum of the integers
appearing in the multiset $\lambda$. Due to the definition of partitions
as multisets, we can present a partition as an ordered sequence of
positive integers with respect to any total order on $\N_{\geq 1}$. For
this reason, we denote any partition $\lambda$ by a nonincreasing
sequence $\Par{\lambda_1, \dots, \lambda_k}$ of positive integers
(that is, $\lambda_i \geq \lambda_{i + 1}$ for all $i \in [k - 1]$).
Under this convention, the \Def{length} $\Length(\lambda)$ of $\lambda$
is $k$, and for any $i \in [\Length(\lambda)]$, the \Def{i$th$ part}
of $\lambda$ is~$\lambda_i$.
\medbreak

\subsubsection{Permutations} \label{subsubsec:permutations}
A \Def{permutation} of \Def{size} $n$ is a bijection $\sigma$ from $[n]$
to $[n]$. The combinatorial graded collection of all permutations is,
in accordance with the notations of
Section~\ref{subsubsec:symmetric_collections}, denoted by
$\SymmetricGroup$. The generating series of $\SymmetricGroup$ is
\begin{equation}
    \GeneratingSeries_\SymmetricGroup(t)
    = \sum_{n \in \N} n! \, t^n.
\end{equation}
Hence, the integer sequence of $\SymmetricGroup$ begins by
\begin{equation}
    1, 1, 2, 6, 24, 120, 720, 5040, 40320
\end{equation}
and is Sequence~\OEIS{A000142} of~\cite{Slo}. Any permutation $\sigma$
of $\SymmetricGroup(n)$ is denoted as a word $\sigma(1) \dots \sigma(n)$
on $\N_{\geq 1}$. Under this convention, a permutation of size $n$ is a
word on the alphabet $[n]$ with exactly one occurrence of each letter
of~$[n]$. The composition operation $\circ$ of maps is a concentrated
product on $\SymmetricGroup$ and the valid inputs of $\circ$ are
the ordered pairs $\Par{\sigma_1, \sigma_2}$ such that
$|\sigma_1| = |\sigma_2|$.
\medbreak

A \Def{descent} of $\sigma \in \SymmetricGroup(n)$ is a position
$i \in [n - 1]$ such that $\sigma(i) > \sigma(i + 1)$. The set of all
descents of $\sigma$ is denoted by $\Des(\sigma)$. For any word $u$
defined on an alphabet $A$ equipped with a total order $\Ord$, the
\Def{standardized} $\Std(u)$ of $u$ is the permutation of size $|u|$
having the same inversions as the ones of $u$. In other terms $\Std(u)$
has its letters in the same relative order as those of $u$, with respect
to $\Ord$, where equal letters of $u$ are ordered from left to right as
the smallest to the greatest. For example, by considering the alphabet
$\N$ equipped with the natural order of integers,
$\Std(211241) = 412563$. This map $\Std$ is a surjective collection
morphism from $\N^*$ to~$\SymmetricGroup$.
\medbreak

\subsubsection{Binary trees} \label{subsubsec:binary_trees}
Let $\ColBT_\Node$ be the combinatorial graded collection
satisfying the relation
\begin{equation}
    \ColBT_\Node = \{\Leaf\} +
    \BBrack{\{\Node\},
    \BBrack{\ColBT_\Node, \ColBT_\Node}_\times^+}_\times^+,
\end{equation}
where $\Leaf$ is an object of size $0$ called \Def{leaf} and $\Node$ is
an atomic object called \Def{internal node}. We call \Def{binary tree}
each object of $\ColBT_\Node$. By definition, a binary tree $\Tfr$
is either the leaf $\Leaf$ or an ordered pair
$\Par{\Node, \Par{\Tfr_1, \Tfr_2}}$ where $\Tfr_1$ and $\Tfr_2$ are
binary trees. Observe that this description of binary trees is
recursive. For instance,
\begin{equation} \label{equ:examples_binary_trees}
    \Leaf, \quad
    (\Node, (\Leaf, \Leaf)), \quad
    (\Node, ((\Node, (\Leaf, \Leaf)), \Leaf)), \quad
    (\Node, (\Leaf, (\Node, (\Leaf, \Leaf)))), \quad
    (\Node, ((\Node, (\Leaf, \Leaf)), (\Node, (\Leaf, \Leaf)))),
\end{equation}
are binary trees.  If $\Tfr$ is a binary tree different from the leaf,
by definition, $\Tfr$ can be expressed as
$\Tfr = (\Node, (\Tfr_1, \Tfr_2))$ where $\Tfr_1$ and $\Tfr_2$ are two
binary trees. In this case, $\Tfr_1$ (resp. $\Tfr_2$) is the
\Def{left subtree} (resp. \Def{right subtree}) of~$\Tfr$. By drawing
each leaf by $\LeafPic$ and each binary tree with at least one internal
node by an internal node $\NodePic$ attached below it, from left to
right, to its left and right subtrees by means of edges $\EdgePic$, the
binary trees of~\eqref{equ:examples_binary_trees} are depicted by
\begin{equation}
    \LeafPic, \quad
    \begin{tikzpicture}[xscale=.2,yscale=.17,Centering]
        \node[Leaf](0)at(0.00,-1.50){};
        \node[Leaf](2)at(2.00,-1.50){};
        \node[Node](1)at(1.00,0.00){};
        \draw[Edge](0)--(1);
        \draw[Edge](2)--(1);
        \node(r)at(1.00,1.75){};
        \draw[Edge](r)--(1);
    \end{tikzpicture}\,, \quad
    \begin{tikzpicture}[xscale=.16,yscale=.14,Centering]
        \node[Leaf](0)at(0.00,-3.33){};
        \node[Leaf](2)at(2.00,-3.33){};
        \node[Leaf](4)at(4.00,-1.67){};
        \node[Node](1)at(1.00,-1.67){};
        \node[Node](3)at(3.00,0.00){};
        \draw[Edge](0)--(1);
        \draw[Edge](1)--(3);
        \draw[Edge](2)--(1);
        \draw[Edge](4)--(3);
        \node(r)at(3.00,2){};
        \draw[Edge](r)--(3);
    \end{tikzpicture}\,, \quad
    \begin{tikzpicture}[xscale=.16,yscale=.14,Centering]
        \node[Leaf](0)at(0.00,-1.67){};
        \node[Leaf](2)at(2.00,-3.33){};
        \node[Leaf](4)at(4.00,-3.33){};
        \node[Node](1)at(1.00,0.00){};
        \node[Node](3)at(3.00,-1.67){};
        \draw[Edge](0)--(1);
        \draw[Edge](2)--(3);
        \draw[Edge](3)--(1);
        \draw[Edge](4)--(3);
        \node(r)at(1.00,2){};
        \draw[Edge](r)--(1);
    \end{tikzpicture}\,, \quad
    \begin{tikzpicture}[xscale=.15,yscale=.12,Centering]
        \node[Leaf](0)at(0.00,-4.67){};
        \node[Leaf](2)at(2.00,-4.67){};
        \node[Leaf](4)at(4.00,-4.67){};
        \node[Leaf](6)at(6.00,-4.67){};
        \node[Node](1)at(1.00,-2.33){};
        \node[Node](3)at(3.00,0.00){};
        \node[Node](5)at(5.00,-2.33){};
        \draw[Edge](0)--(1);
        \draw[Edge](1)--(3);
        \draw[Edge](2)--(1);
        \draw[Edge](4)--(5);
        \draw[Edge](5)--(3);
        \draw[Edge](6)--(5);
        \node(r)at(3.00,2.5){};
        \draw[Edge](r)--(3);
    \end{tikzpicture}\,.
\end{equation}
By definition of the sum and the Cartesian product operations over
graded collections, the size of a binary tree $\Tfr$ satisfies
\begin{equation} \label{equ:size_binary_tree_leaves}
    |\Tfr| =
    \begin{cases}
        0 & \mbox{if } \Tfr = \Leaf, \\
        1 + |\Tfr_1| + |\Tfr_2|
            & \mbox{otherwise (} \Tfr = (\Node, (\Tfr_1, \Tfr_2))
                \mbox{)}.
    \end{cases}
\end{equation}
In other words, the size of $\Tfr$ is the number of occurrences of
$\Node$ it contains. Since $\GeneratingSeries_{\{\Leaf\}}(t) = 1$ and
$\GeneratingSeries_{\{\Node\}}(t) = t$, it follows
from~\eqref{equ:generating_series_sum}
and~\eqref{equ:generating_series_multiplication} that the generating
series of $\ColBT_\Node$ satisfies the quadratic algebraic equation
\begin{equation} \label{equ:algebraic_relation_binary_trees_leaves}
    1 - \GeneratingSeries_{\ColBT_\Node}(t)
    + t \GeneratingSeries_{\ColBT_\Node}(t)^2 = 0.
\end{equation}
The unique solution having a combinatorial meaning
of~\eqref{equ:algebraic_relation_binary_trees_leaves} is
\begin{equation}
    \GeneratingSeries_{\ColBT_\Node}(t) =
    \frac{1 - \sqrt{1 - 4t}}{2t} =
    \sum_{n \in \N}
    \frac{1}{n + 1} \binom{2n}{n} \, t^n
\end{equation}
The integer sequence of $\ColBT_\Node$ begins by
\begin{equation}
    1, 1, 2, 5, 14, 42, 132, 429, 1430
\end{equation}
and is Sequence~\OEIS{A000108} of~\cite{Slo}. These numbers are known
as \Def{Catalan numbers}.
\medbreak

\section{Posets} \label{sec:posets}
We consider now collections endowed with partial order relations
compatible with their indexations. Such structures are important in
combinatorics since they lead, for instance, to the construction of
alternative bases of combinatorial spaces (see forthcoming
Section~\ref{subsec:change_basis_posets} of Chapter~\ref{chap:algebra}).
We provide general definitions about posets and consider as examples
three important ones: the cube, Tamari, and right weak order posets.
\medbreak

\subsection{Posets on collections} \label{subsec:posets}
Let us provide the main definitions about collections endowed with the
structure of a poset.
\medbreak

\subsubsection{Elementary definitions}
An \Def{$I$-poset} is a pair $\Par{\Qca, \Ord_\Qca}$ where $\Qca$ is an
$I$-collection and $\Ord_\Qca$ is both a relation on $\Qca$ (recall that
relations on collections preserve the indexes) and a partial order
relation. The \Def{strict order relation} of $\Ord$ is the relation
$\OrdStrict$ on $\Qca$ satisfying, for all $x, y \in \Qca$,
$x \OrdStrict y$ if $x \Ord y$ and $x \ne y$.
\medbreak

The \Def{interval} between two objects $x$ and $z$ of $\Qca$ is the set
\begin{math}
    [x, z] := \left\{y \in \Qca : x \Ord_\Qca y \Ord_\Qca z\right\}.
\end{math}
When all intervals of $\Qca$ are finite, $\Qca$ is \Def{locally finite}.
Observe that when $\Qca$ is combinatorial, $\Qca$ is locally finite. For
any $i \in I$, an object $x$ of $\Qca(i)$ is a \Def{greatest} (resp.
\Def{least}) \Def{element} if for all $y \in \Qca(i)$, $y \Ord_\Qca x$
(resp. $x \Ord_\Qca y$). Moreover, for any $i \in I$, an object $x$ of
$\Qca(i)$ is a \Def{maximal} (resp. \Def{minimal}) \Def{element} if for
all $y \in \Qca(i)$, $x \Ord_\Qca y$ (resp. $y \Ord_\Qca x$) implies
$x = y$. If $x$ and $y$ are two different objects of $\Qca$, $y$
\Def{covers} $x$ if $[x, y] = \{x, y\}$. Two objects $x$ and $y$ are
\Def{comparable} (resp. \Def{incomparable}) in $\Qca$ if $x \Ord_\Qca y$
or $y \Ord_\Qca x$ (resp. neither $x \Ord_\Qca y$ nor $y \Ord_\Qca x$
holds). If for any $i \in I$ and any $i$-objects $x$ and $y$ of $\Qca$,
$x$ and $y$ are comparable, $\Qca$ is a \Def{total order}. A \Def{chain}
of $\Qca$ is a sequence $\Par{x_1, \dots, x_k}$ such that
$x_j \Ord_\Qca x_{j + 1}$ for all $j \in [k - 1]$. An \Def{antichain} of
$\Qca$ is a subset of pairwise incomparable elements of $\Qca$. The
\Def{Hasse diagram} of $\Par{\Qca, \Ord_\Qca}$ is the directed graph
having $\Qca$ as set of vertices and all the pairs $(x, y)$ where $y$
covers $x$ as set of arcs.
\medbreak

Besides, if $\Par{\Qca_1, \Ord_{\Qca_1}}$ and
$\Par{\Qca_2, \Ord_{\Qca_2}}$ are two $I$-posets, a map
$\phi : \Qca_1 \to \Qca_2$ is a \Def{poset morphism} if $\phi$ is a
collection morphism and for all $x, y \in \Qca_1$ such that
$x \Ord_{\Qca_1} y$, $\phi(x) \Ord_{\Qca_2} \phi(y)$. Besides, $\Qca_2$
is a \Def{subposet} of $\Qca_1$ if $\Qca_2$ is a subcollection of
$\Qca_1$ and $\Ord_{\Qca_2}$ is the restriction of $\Ord_{\Qca_1}$
on~$\Qca_2$. For any $i \in I$, we call \Def{$i$-subposet} of $\Qca$ the
subposet of $\Qca$  obtained by restricting $\Ord_\Qca$ on $\Qca(i)$.
\medbreak

We shall define posets $\Qca$ by drawing Hasse diagrams, where minimal
elements are drawn uppermost and vertices are labeled by the elements of
$\Qca$. For instance, the Hasse diagram
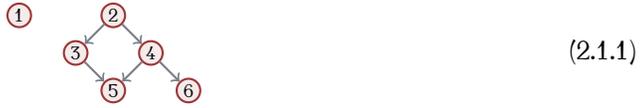
\begin{equation}
    \begin{tikzpicture}[xscale=.5,yscale=.5,Centering]
        \node[PosetVertex](1)at(-1.5,1){\begin{math}1\end{math}};
        \node[PosetVertex](2)at(1,1){\begin{math}2\end{math}};
        \node[PosetVertex](3)at(0,0){\begin{math}3\end{math}};
        \node[PosetVertex](4)at(2,0){\begin{math}4\end{math}};
        \node[PosetVertex](5)at(1,-1){\begin{math}5\end{math}};
        \node[PosetVertex](6)at(3,-1){\begin{math}6\end{math}};
        \draw[Arc](2)--(3);
        \draw[Arc](2)--(4);
        \draw[Arc](3)--(5);
        \draw[Arc](4)--(5);
        \draw[Arc](4)--(6);
    \end{tikzpicture}
\end{equation}
denotes the simple (``simple'' here means the property of collections
stated in Section~\ref{subsubsec:collections_definition}) poset
$([6], \Ord)$ satisfying among others $3 \Ord 5$ and~$2 \Ord 6$.
\medbreak

\subsubsection{Operations over posets}
\label{subsubsec:operations_posets}
If $\Par{\Qca_1, \Ord_{\Qca_1}}$ and $\Par{\Qca_2, \Ord_{\Qca_2}}$ are
two $I$-posets, the sum $\Qca_1 + \Qca_2$ is endowed with the relation
$\Ord$ satisfying, $x \Ord y$ whenever $x, y \in \Qca_1$ and
$x \Ord_{\Qca_1} y$, or $x, y \in \Qca_2$ and $x \Ord_{\Qca_2} y$. Since
$\Ord$ is an order relation, $\Qca_1 + \Qca_2$ is an $I$-poset, called
\Def{sum} of $\Par{\Qca_1, \Ord_{\Qca_1}}$ and
$\Par{\Qca_2, \Ord_{\Qca_2}}$. For any $p \in \N$ and any $I$-posets
$\Par{\Qca_1, \Ord_{\Qca_1}}$, \dots, $\Par{\Qca_p, \Ord_{\Qca_p}}$, the
Hadamard product $\BBrack{\Qca_1, \dots, \Qca_p}_\Hadamard$ is endowed
with the relation $\Ord$ satisfying
$\Par{x_1, \dots, x_p} \Ord \Par{y_1, \dots, y_p}$ for any
\begin{math}
    \Par{x_1, \dots, x_p}, \Par{y_1, \dots, y_p}
    \in \BBrack{\Qca_1, \dots, \Qca_p}_\Hadamard
\end{math}
such that $x_k \Ord_{\Qca_k} y_k$ for all $k \in [p]$. Since $\Ord$ is
an order relation, $\BBrack{\Qca_1, \dots, \Qca_p}_\Hadamard$ is an
$I$-poset, called \Def{Hadamard product} of
$\Par{\Qca_1, \Ord_{\Qca_1}}$, \dots, $\Par{\Qca_p, \Ord_{\Qca_p}}$. Let
$\Par{\Qca, \Ord_\Qca}$ be an $I$-poset. Let $\Comparable(\Qca)$ be the
subcollection of $\BBrack{\Qca, \Qca}_\Hadamard$ restrained on the
ordered pairs $(x, y)$ such that $x \Ord_\Qca y$, called
\Def{pairs of comparable objects}. By definition, $\Comparable(\Qca)$ is
endowed with the restriction of the order relation of
$\BBrack{\Qca, \Qca}_\Hadamard$ on $\Comparable(\Qca)$. We call
$\Comparable(\Qca)$ the \Def{poset of pairs of comparable objects} of
$\Qca$. Finally, the \Def{dual} of $\Qca$ is the $I$-poset
$\Par{\Qca, \bar{\Ord}_\Qca}$ such that $x \, \bar{\Ord}_\Qca \, y$
holds whenever $y \Ord_\Qca x$ for any~$x, y \in \Qca$.
\medbreak

\subsection{Examples} \label{subsec:examples_posets}
We consider here three well-known combinatorial posets.
\medbreak

\subsubsection{The cube poset} \label{subsubsec:cube_poset}
\begin{figure}[ht]
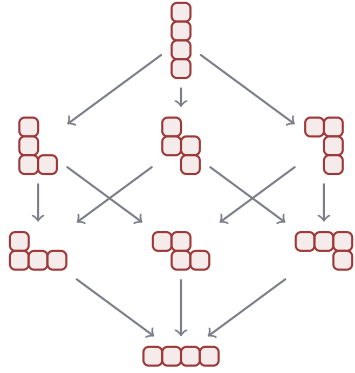

    \centering
    \begin{equation*}

            };
            \draw[Arc](1111)--(112);
            \draw[Arc](1111)--(121);
            \draw[Arc](1111)--(211);
            \draw[Arc](112)--(13);
            \draw[Arc](121)--(13);
            \draw[Arc](112)--(22);
            \draw[Arc](211)--(22);
            \draw[Arc](211)--(31);
            \draw[Arc](121)--(31);
            \draw[Arc](13)--(4);
            \draw[Arc](22)--(4);
            \draw[Arc](31)--(4);
        \end{tikzpicture}
    \end{equation*}
    \caption[The cube poset of order $4$.]
    {\footnotesize The Hasse diagram of the refinement order of
    compositions of size $4$, where each composition is represented
    through its ribbon diagram.}
    \label{fig:cube_order_4}
\end{figure}

Let $\Ord$ be the partial order relation on the combinatorial collection
$\ColComp$ of compositions generated by the covering relation
$\CoveringRel$ defined, for any composition $\LambdaB$ of length $k$, by
\begin{equation}
    \Par{\LambdaB_1, \dots, \LambdaB_{i - 1},
    \LambdaB_i, \LambdaB_{i + 1},
    \LambdaB_{i + 2}, \dots, \LambdaB_k}
    \CoveringRel
    \Par{\LambdaB_1, \dots, \LambdaB_{i - 1},
    \LambdaB_i + \LambdaB_{i + 1},
    \LambdaB_{i + 2}, \dots, \LambdaB_k}.
\end{equation}
For instance, $2123 \CoveringRel 215$ and $2123 \Ord 8$. This order is
the \Def{refinement order} of compositions. The Hasse diagram of
$(\ColComp, \Ord)$ restricted on $\ColComp(4)$ is shown in
Figure~\ref{fig:cube_order_4}.
\medbreak

Observe that for all compositions $\LambdaB$ and $\MuB$,
$\LambdaB \Ord \MuB$ if and only if
$\Des(\MuB) \subseteq \Des(\LambdaB)$. Each $n$-subposet of the
refinement order of compositions is known as the \Def{cube poset} of
dimension~$n - 1$. Moreover, the cube poset of dimension $n - 1$ is
isomorphic to the dual of the poset of all subsets of $[n - 1]$ ordered
by set inclusion. An isomorphism is provided by the map $\Des$ sending a
composition of size $n$ to a subset of~$[n - 1]$.
\medbreak

\subsubsection{The Tamari order on binary trees}
\label{subsubsec:tamari_order}
\begin{figure}[ht]
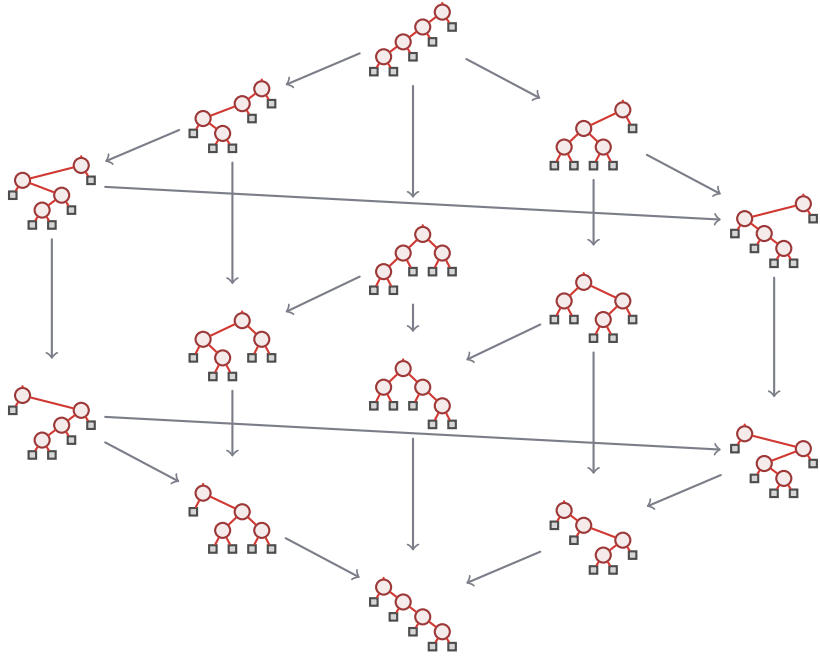

    \centering
    \begin{equation*}

            };
            \draw[Arc](222200000)--(222020000);
            \draw[Arc](222200000)--(222002000);
            \draw[Arc](222020000)--(220220000);
            \draw[Arc](222002000)--(220202000);
            \draw[Arc](220220000)--(220202000);
            \draw[Arc](222200000)--(222000200);
            \draw[Arc](222020000)--(220200200);
            \draw[Arc](222000200)--(220200200);
            \draw[Arc](222002000)--(220022000);
            \draw[Arc](222000200)--(220020200);
            \draw[Arc](220022000)--(220020200);
            \draw[Arc](220220000)--(202220000);
            \draw[Arc](220202000)--(202202000);
            \draw[Arc](202220000)--(202202000);
            \draw[Arc](202220000)--(202200200);
            \draw[Arc](220200200)--(202200200);
            \draw[Arc](220022000)--(202022000);
            \draw[Arc](202202000)--(202022000);
            \draw[Arc](202200200)--(202020200);
            \draw[Arc](220020200)--(202020200);
            \draw[Arc](202022000)--(202020200);
        \end{tikzpicture}
    \end{equation*}
    \caption[The Tamari poset of order $4$.]
    {\footnotesize
    The Hasse diagram of the Tamari poset of binary trees of size~$4$.}
    \label{fig:tamari_order_4}
\end{figure}

Let $\Ord$ be the partial order relation on the combinatorial collection
$\ColBT_\Node$ of binary trees generated by the covering relation
$\CoveringRel$ defined by
\begin{equation} \label{equ:right_rotation_tamari}
    \Par{\dots
        \Par{\Node, \Par{\Par{\Node, \Par{\Rfr_1, \Rfr_2}}, \Rfr_3}}
        \dots}
    \CoveringRel
    \Par{\dots
        \Par{\Node, \Par{\Rfr_1, \Par{\Node, \Par{\Rfr_2, \Rfr_3}}}}
        \dots},
\end{equation}
where $\Rfr_1$, $\Rfr_2$, and $\Rfr_3$ are any binary trees. We call
$\CoveringRel$ the \Def{right rotation} relation. At this moment, the
definition of this relation on binary trees is informal but, in
Section~\ref{subsec:rewrite_rules_syntax_trees} of
Chapter~\ref{chap:trees}, we shall develop precise tools to define and
handle such operations on binary trees and more generally on syntax
trees. The order $\Ord$ is the \Def{Tamari order} on binary trees. The
Hasse diagram of $\Par{\ColBT_\Node, \Ord}$ restricted on
$\ColBT_\Node(4)$ is shown in Figure~\ref{fig:tamari_order_4}.
\medbreak

\subsubsection{The right weak order on permutations}
\begin{figure}[ht]
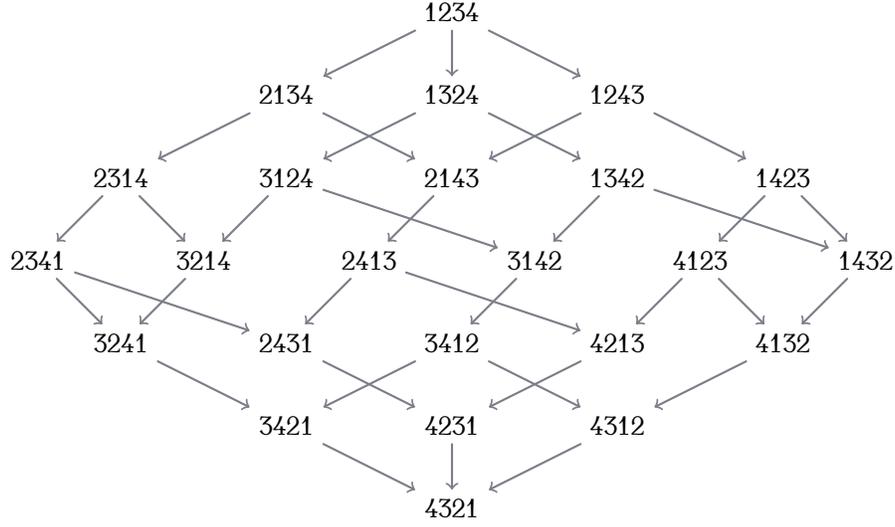

    \centering
    \begin{equation*}

    \end{equation*}
    \caption[The right weak poset of order $4$.]
    {\footnotesize
    The Hasse diagram of the right weak poset of permutations of
    size~$4$.}
    \label{fig:right_weak_order_order_4}
\end{figure}

Let $\Ord$ be the partial order relation on the combinatorial collection
$\SymmetricGroup$ of permutations generated by the covering relation
$\CoveringRel$ defined by
\begin{equation}
    u \, \Asf \Bsf \, v \CoveringRel u \, \Bsf \Asf \, v,
\end{equation}
where $u$ and $v$ are words on $\N_{\geq 1}$, and $\Asf$ and $\Bsf$ are
letters such that $\Asf < \Bsf$. This order is the
\Def{right weak order} of permutations. The Hasse diagram of
$(\SymmetricGroup, \Ord)$ restricted on $\SymmetricGroup(4)$ is shown in
Figure~\ref{fig:right_weak_order_order_4}.
\medbreak

\section{Rewrite systems} \label{sec:rewrite_systems}
A rewrite system describes a process whose goal is to transform
iteratively an object into another one. We consider rewrite systems on
$I$-collections, so that an $i$-object, $i \in I$, can be transformed
only into $i$-objects. As we shall see, rewrite systems and posets have
some close connections because it is possible, in some  cases, to
construct posets from rewrite systems.
\medbreak

\subsection{Rewrite systems on collections}
Let us provide the main definitions about collections endowed with the
structure of a rewrite system. Two properties of rewrite systems are
fundamental: the termination and the confluence. We provide strategies
to prove that a given rewrite system satisfies one or the other.
\medbreak

\subsubsection{Elementary definitions}
Let $C$ be an $I$-collection. An \Def{$I$-rewrite system} is a pair
$(C, \Rew)$ where $C$ is an $I$-collection and $\Rew$ is a relation on
$C$. We call $\Rew$ a \Def{rewrite rule}. When $x_0$, $x_1$, \dots,
$x_k$ are objects of $C$ such that $k \in \N$ and
\begin{equation}
    x_0 \Rew x_1 \Rew \cdots \Rew x_k,
\end{equation}
we say that $x_0$ is \Def{rewritable} by $\Rew$ into $x_k$ in $k$
\Def{steps}. The reflexive and transitive closure of $\Rew$ is denoted
by~$\RewRT$. The directed graph $(C, \Rew)$ consisting in $C$ as set of
vertices and $\Rew$ as set of arcs is the \Def{rewriting graph}
of~$(C, \Rew)$.
\medbreak

\subsubsection{Termination} \label{subsubsec:termination}
When there is no infinite chain
\begin{equation} \label{equ:infinite_rewrite_chain}
    x_0 \Rew x_1 \Rew x_2 \Rew \cdots
\end{equation}
where all $x_j \in C$, $j \in \N$, $(C, \Rew)$ is
\Def{terminating}. Observe that, if $C$ is combinatorial, due to the
fact that for any $i \in I$, each set $C(i)$ is finite and the fact that
the rewriting relation preserves the indexes, if such an infinite
chain~\eqref{equ:infinite_rewrite_chain} exists, then it is of the form
\begin{equation} \label{equ:infinite_rewrite_chain_repetition}
    x_0 \Rew \cdots \Rew x_r \Rew \cdots \Rew x_r \Rew \cdots,
\end{equation}
for a certain $r \in \N$. A \Def{normal form} of $(C, \Rew)$ is an
object $x$ of $C$ such that for all $x' \in C$, $x \RewRT x'$ imply
$x' = x$. In other words, a normal form of $(C, \Rew)$ is an object
which is not rewritable by~$\Rew$. This set of objects, which is a
subcollection of $C$, is denoted by $\NormalForms_{(C, \Rew)}$. The
following result provides a tool in the aim to show that a
combinatorial rewrite system is terminating.
\medbreak

\begin{Theorem} \label{thm:terminating_rewrite_rules_posets}
    Let $(C, \Rew)$ be a combinatorial rewrite system. Then, $(C, \Rew)$
    is terminating if and only if the binary relation $\RewRT$ is an
    order relation and endows~$C$ with a structure of a combinatorial
    poset.
\end{Theorem}
\medbreak

When $C$ is combinatorial and $(C, \Rew)$ is terminating, by
Theorem~\ref{thm:terminating_rewrite_rules_posets}, $\Par{C, \RewRT}$ is
a combinatorial poset and we call it the \Def{poset induced} by~$\Rew$.
\medbreak

In practice, Theorem~\ref{thm:terminating_rewrite_rules_posets} is used
as follows. To show that a combinatorial $I$-rewrite system $(C, \Rew)$
is terminating, we construct a map $\theta : C \to \Qca$ where
$(\Qca, \Ord)$ is an $I$-poset such that for any $x, x' \in C$,
$x \Rew x'$ implies $\theta(x) \OrdStrict \theta(x')$. Such a map
$\theta$ is a \Def{termination invariant}. Indeed, since each $C(i)$,
$i \in I$, is finite, this property leads to the fact that there is no
infinite chain of the
form~\eqref{equ:infinite_rewrite_chain_repetition}. In most cases,
$\Qca$ is a set of tuples of integers of a fixed length, and $\Ord$ is
the lexicographic order on these tuples.
\medbreak

\subsubsection{Confluence}
When for any objects $x$, $y_1$, and $y_2$ of $C$ such that
$x \RewRT y_1$ and $x \RewRT y_2$, there exists an object $x'$ of $C$
such that $y_1 \RewRT x'$ and $y_2 \RewRT x'$, the rewrite system
$(C, \Rew)$ is \Def{confluent}. When~$\Rew$ is both terminating and
confluent, $\Rew$ is \Def{convergent}.
\medbreak

An object $x$ of $C$ is a \Def{branching} object if there exist two
different objects $y_1$ and $y_2$ satisfying $x \Rew y_1$ and
$x \Rew y_2$. In this case, the pair $\left\{y_1, y_2\right\}$ is a
\Def{branching pair} for~$x$. We say that a branching pair
$\left\{y_1, y_2\right\}$ is \Def{joinable} if there exists an object
$z$ of $C$ such that $y_1 \RewRT z$ and $y_2 \RewRT z$. In practice,
showing that a terminating rewrite system is confluent is made simple,
thank to the following result.
\medbreak

\begin{Theorem} \label{thm:diamond_lemma}
    Let $(C, \Rew)$ be a rewrite system. If $(C, \Rew)$ is terminating
    and all its branching pairs are joinable, $(C, \Rew)$ is confluent.
\end{Theorem}
\medbreak

\subsubsection{Closures} \label{subsubsec:closures}
Let $(C, \Rew)$ be an $I$-rewrite system such that $C$ is endowed with a
set $\Pca$ of concentrated products. Then, let $\Par{C, \Rew_\Pca}$ be
the rewrite system such that $\Rew_\Pca$ contains $\Rew$ (as a binary
relation) and satisfies moreover
\begin{equation} \label{equ:closure_rewrite_rule}
    \Product\Par{x_1, \dots, x_{j - 1}, y,
    x_{j + 1}, \dots, x_p}
    \enspace \Rew_\Pca \enspace
    \Product\Par{x_1, \dots, x_{j - 1}, y',
    x_{j + 1}, \dots, x_p}
\end{equation}
for any product $\Product$ of arity $p$ of $\Pca$,
such that the left and right members of~\eqref{equ:closure_rewrite_rule}
are valid inputs for~$\Product$ and~$y \Rew y'$. The fact that all
products $\Product$ of $\Pca$ are concentrated ensures that the left
and the right members of~\eqref{equ:closure_rewrite_rule} have the
same index, so that $\Par{C, \Rew_\Pca}$ is an $I$-rewrite system. We
call $\Par{C, \Rew_\Pca}$ the \Def{$\Pca$-closure} of $(C, \Rew)$. Such
closures provide convenient and concise ways to define rewrite systems.
\medbreak

\subsection{Examples} \label{subsec:examples_rewrite_rules}
Let us review some examples of rewrite systems on various combinatorial
collections.
\medbreak

\subsubsection{A first rewrite system on words}
Let $A := \{\Asf, \Bsf\}$ be an alphabet, and consider the graded
rewrite system $(A^*, \Rew$) defined by
\begin{equation}
    u x \Rew x u
\end{equation}
for any $u \in A^*$ and $x \in A$. We have, for instance,
\begin{equation}
    \Asf \Asf \Bsf \Asf \Rew \Asf \Asf \Asf \Bsf
    \Rew \Bsf \Asf \Asf \Asf \Rew \Asf \Bsf \Asf \Asf
    \Rew \Asf \Asf \Bsf \Asf.
\end{equation}
This rewrite system is not terminating but, since for each word
$u \in A^*$ there is at most a word $v \in A^*$ satisfying $u \Rew v$,
$(A^*, \Rew)$ is confluent.
\medbreak

\subsubsection{A second rewrite system on words}
Let us now study the graded rewrite system $(A^*, \Rew)$ defined by
$\Asf \Bsf \Asf \Rew \Bsf \Asf \Bsf$ where $A$ is the alphabet of the
previous example. Consider the graded complete ternary product
\begin{math}
    \Product : A^* \times A^* \times A^* \to A^*
\end{math}
on $A^*$ defined for any $u, v, w \in A^*$ by
\begin{math}
    \Product(u, v, w) := u \Conc v \Conc w
\end{math}
where $\Conc$ is the concatenation product of words. Let
$\Pca := \{\Product\}$ and  $\Par{A^*, \Rew_\Pca}$ be the $\Pca$-closure
of $\Par{A^*, \Rew}$. By definition of closures, $\Rew_\Pca$ satisfies
\begin{equation}
    \Asf \Bsf \Asf \Rew_\Pca \Bsf \Asf \Bsf
\end{equation}
and
\begin{equation}
    \Asf \Bsf \Asf \Conc v \Conc w
    \Rew_\Pca
    \Bsf \Asf \Bsf \Conc v \Conc w,
    \qquad
    u \Conc \Asf \Bsf \Asf \Conc w
    \Rew_\Pca
    u \Conc \Bsf \Asf \Bsf \Conc w,
    \qquad
    u \Conc v \Conc \Asf \Bsf \Asf
    \Rew_\Pca
    u \Conc v \Conc \Bsf \Asf \Bsf,
\end{equation}
for any words $u$, $v$, and $w$ on $A$. All this is equivalent to the
fact that $\Rew_\Pca$ is the rewrite rule satisfying
\begin{equation}
    u \Conc \Asf \Bsf \Asf \Conc w
    \Rew_\Pca
    u \Conc \Bsf \Asf \Bsf \Conc w,
\end{equation}
for any words $u$ and $w$ on $A$. The rewrite system
$\Par{A^*, \Rew_\Pca}$ is terminating since, for any words $u$ and $v$
on $A$, if $u \Rew_\Pca v$, then $|v|_\Bsf  = |u|_\Bsf + 1$. Hence, the
map $\theta : A^n \to [0, n]$ defined for any $n \in \N$ and $u \in A^n$
by $\theta(u) := |u|_\Bsf$ is a termination invariant. The normal forms
of $\Par{A^*, \Rew_\Pca}$ are the words that do not admit
$\Asf \Bsf \Asf$ as factor. Moreover, $\Par{A^*, \Rew_\Pca}$ is not
confluent since
\begin{math}
    \Asf \Bsf \Asf \Bsf \Asf \Rew_\Pca \Bsf \Asf \Bsf \Bsf \Asf
\end{math}
and
\begin{math}
    \Asf \Bsf \Asf \Bsf \Asf \Rew_\Pca \Asf \Bsf \Bsf \Asf \Bsf
\end{math},
and $\{\Bsf \Asf \Bsf \Bsf \Asf, \Asf \Bsf \Bsf \Asf \Bsf\}$ is a
non-joinable branching pair for $\Asf \Bsf \Asf \Bsf \Asf$ (because
these two elements are normal forms).
\medbreak

\subsubsection{A rewrite system on compositions}
Let the graded rewrite system $(\ColComp, \Rew)$ defined, by seeing
compositions through their ribbon diagrams, by
\begin{equation}
    \LambdaB \Conc

    \Conc \MuB,
\end{equation}
where $\Conc$ is the concatenation of the compositions (seen as words of
integers) and $\LambdaB$ and $\MuB$ are any compositions.
We have, for instance,
\begin{equation}
    %
\,.
\end{equation}
The rewrite system $(\ColComp, \Rew)$ is terminating since, for any
compositions $\LambdaB$ and $\MuB$, if $\LambdaB \Rew \MuB$, then
$\MuB \OrdStrict \LambdaB$ where $\Ord$ is the refinement order of
compositions (see Section~\ref{subsubsec:cube_poset}).
\medbreak

\SkipTocEntry\section*{Bibliographic notes}

\SkipTocEntry\subsection*{About collections}
Our exposition about combinatorial objects through combinatorial
collections is very elementary in the sense that it requires a very
small amount of mathematical knowledge. However, combinatorial
collections form a general and powerful tool to work with algebraic
structures involving combinatorial objects. For instance, graded
collections appear in the context of operads (see
Section~\ref{subsec:operads} of Chapter~\ref{chap:operads}) or graded
associative algebras (see
Section~\ref{subsec:associative_coassociative_algebras} of
Chapter~\ref{chap:algebra}), colored collections appear in the context
of colored operads (see Section~\ref{subsec:colored_operads} of
Chapter~\ref{chap:generalizations}), cyclic collections appear in the
context of cyclic operads (see Section~\ref{subsec:cyclic_operads} of
Chapter~\ref{chap:generalizations}), symmetric collections appear in the
context of symmetric operads (see
Section~\ref{subsec:symmetric_operads} of
Chapter~\ref{chap:generalizations}), and $2$-graded collections appear
in the context of pros (see Section~\ref{sec:pros} of
Chapter~\ref{chap:generalizations}). There are other sensible tools to
encode combinatorial sets. Flajolet and Sedgewick provided a complete
description of what we call combinatorial graded collections under the
name of combinatorial classes in~\cite{FS09}, as a prelude for a
conspectus of the field of analytic combinatorics. The proofs of most of
the properties about generating series of
Section~\ref{subsec:operations_collections} can be found here. The nice
translations of most of the combinatorial operations involving
combinatorial sets as algebraic operations on their generating series,
together with its simplicity, are one of the main pros of this theory.
By shifting in the world of labeled objects, it is relevant to work with
species of structures, that are roughly speaking combinatorial graded
collections $C$ with an action of the symmetric group
$\SymmetricGroup(n)$ on each $C(n)$ which can be thought as a
relabeling action. In this context, it is more accurate to work with
exponential generating series, instead of ordinary ones when we consider
such combinatorial graded collections. This theory has been introduced
by Joyal~\cite{Joy81} and developed afterwards by the Quebec school of
combinatorics~\cite{BLL98,BLL13}. Species of structures are very good
candidates to work with symmetric operads~\cite{Men15} since the action
of the symmetric group of a symmetric operad is encapsulated into the
action of the symmetric group of an underlying species of structure. In
this book, to work with symmetric operads, we shall consider symmetric
collections. An other interesting way to describe combinatorial objects
passes through polynomial functors~\cite{Koc09}.
\medbreak

\SkipTocEntry\subsection*{About the Tamari poset}
The Tamari poset is a combinatorial poset on binary trees introduced in
the study of nonassociative operations~\cite{Tam62}. Indeed, the
covering relation generating this poset can be thought as a way to move
brackets in expressions where a nonassociative product intervenes.
Moreover, seen on binary trees, this operation translates as a right
rotation, a fundamental operation on binary search trees, used in an
algorithmic context~\cite{Knu98}. This operation is used to maintain
binary trees with a small height in order to access efficiently, from
the roots, to their internal nodes. Some of these trees are known as
balanced binary trees~\cite{AVL62} and form efficient structures to
represent dynamic sets (sets supporting the addition and the suppression
of elements). A lot of properties of the Tamari poset are known, like
the number of intervals of each of its $n$-subposets~\cite{Cha06}
(equivalently, this is the number of pairs of comparable trees
enumerated by their size), and the fact that these posets are
lattices~\cite{HT72}, for all $n \in \N$. Generalizations of this poset
have been introduced by Bergeron and Préville-Ratelle~\cite{BP12} under
the name of $m$-Tamari poset. This poset is defined on the combinatorial
collection of all $m\! +\! 1$-ary trees (see
Section~\ref{subsubsec:k_ary_trees} of Chapter~\ref{chap:trees}). The
number of intervals of each of its $n$-subposets, and the fact that
these posets are lattices are known from~\cite{BFP11}, for
all~$n \in \N$.
\medbreak

\SkipTocEntry\subsection*{About the right weak poset}
The right weak poset of permutations is, like the Tamari poset, also a
lattice~\cite{GR63,YO69}. In a surprising way, despite its apparent
simplicity, there is no known description of the number of intervals of
each $n$-subposet, $n \in \N$, of the right weak poset. Some other
combinatorial poset structures exist on $\SymmetricGroup$ like the
Bruhat order, whose generating relation is similar to the one of the
right weak poset. The definition of the Bruhat order on permutations
comes from the general notion of Bruhat order~\cite{Bjo84} in Coxeter
groups~\cite{Cox34}. As a last noteworthy fact, the cube, the Tamari,
and the right weak posets are linked through surjective morphisms of
combinatorial posets~\cite{LR02}. Indeed, a map between the right weak
poset to the Tamari poset is based upon the binary search tree insertion
algorithm~\cite{Knu98,HNT05}. This algorithm consists in inserting the
letters of a permutation to form step by step a binary tree. Moreover, a
map between the Tamari poset to the cube poset uses the
canopies~\cite{LR98} of the binary trees. The canopy of a binary tree is
a binary word encoding the orientations (to the left or to the right) of
its leaves.
\medbreak

\SkipTocEntry\subsection*{About rewrite systems}
A general reference about rewrite rules and rewrite systems
is~\cite{BN98}. In this text, a general method using maps called measure
functions to show that (not necessarily combinatorial) rewrite systems
are terminating is presented. Besides, Theorem~\ref{thm:diamond_lemma}
is a highly important result in the theory of rewrite systems, known as
the diamond lemma, and is due to Newman~\cite{New42}. There are some
additional useful tools in this theory like the Knuth-Bendix completion
algorithm~\cite{KB70}. This semi-algorithm takes as input a
non-confluent rewrite system and outputs, if possible, a confluent one
having the same reflexive, symmetric, and transitive closures. In an
algebraic context, the Knuth-Bendix completion algorithm leads to the
Buchberger algorithm~\cite{Buc76}. This algorithm computes Gröbner bases
from polynomial ideals.
\medbreak


\chapter{Treelike structures} \label{chap:trees}
This second chapter is devoted to present general notions about treelike
structures. We present more precisely the ones appearing in the
algebraic and combinatorial context of nonsymmetric operads. Rewrite
systems of syntax trees are exposed, as well as methods to prove
their termination and their confluence.
\medbreak

\section{Planar rooted trees} \label{sec:planar_rooted_trees}
Let us start with our prototypical treelike structures, the planar
rooted trees. Most of the treelike structures we shall consider in this
book are variants or enrichments of planar rooted trees.
\medbreak

\subsection{Collection of planar rooted trees}
\label{subsec:comb_collection_planar_rooted_trees}
The combinatorial graded collection of the planar rooted trees can be
defined concisely in a recursive way by using some operations over
combinatorial graded collections (see
Section~\ref{subsec:operations_collections} of
Chapter~\ref{chap:collections}). However, to define rigorously the usual
notions of internal node, leaf, child, father, path, subtree,
{\em etc.}, we need the notion of language associated with a tree.
Indeed, a planar rooted tree is in fact a finite language satisfying
some properties. Therefore, in this section, we shall adopt the point of
view of defining most of the properties of a planar rooted tree through
its language.
\medbreak

\subsubsection{First definitions}
Let $\ColPRT$ be the graded collection satisfying the
relation
\begin{equation} \label{equ:combinatorial_set_planar_rooted_trees}
    \ColPRT =
    \BBrack{\{\Node\}, \List^+\Par{\ColPRT}}_\times^+
\end{equation}
where $\Node$ is an atomic object called \Def{node}. Since
$\ColPRT(0) = \emptyset$, this collection is combinatorial. We call
\Def{planar rooted tree} each object of $\ColPRT$. By definition, a
planar rooted tree $\Tfr$ is an ordered pair
$\Par{\Node, \Par{\Tfr_1, \dots, \Tfr_k}}$, $k \in \N$, where
$\Par{\Tfr_1, \dots, \Tfr_k}$ is a (possibly empty) tuple of planar
rooted trees. This definition is recursive. By convention, the planar
rooted tree $(\Node, ())$ is denoted by $\Leaf$ and is called the
\Def{leaf}. Observe that the leaf is of size $1$. For instance,
\begin{equation} \label{equ:examples_planar_rooted_trees}
    \Leaf, \quad
    (\Node, (\Leaf)), \quad
    (\Node, (\Leaf, \Leaf)), \quad
    (\Node, (\Leaf, (\Node, (\Leaf)))), \quad
    (\Node, ((\Node, ((\Node, (\Leaf, \Leaf)))),
        \Leaf, (\Node, (\Leaf, \Leaf))))
\end{equation}
are planar rooted trees. The \Def{root arity} of a planar rooted tree
$\Tfr := \Par{\Node, \Par{\Tfr_1, \dots, \Tfr_k}}$ is $k$. If $\Tfr$ is
a planar rooted tree different from the leaf, by definition, $\Tfr$ can
be expressed as $\Tfr = \Par{\Node, \Par{\Tfr_1, \dots, \Tfr_k}}$ where
$k \in \N_{\geq 1}$ and all $\Tfr_i$, $i \in [k]$, are planar rooted
trees. In this case, for any $i \in [k]$, $\Tfr_i$ is the
\Def{$i$th suffix subtree} of $\Tfr$. Planar rooted trees are depicted
by drawing each leaf by $\LeafPic$ and each planar rooted tree different
from the leaf by a node $\NodePic$ attached below it, from left to
right, to its suffix subtrees $\Tfr_1$, \dots, $\Tfr_k$ by means of
edges $\EdgePic$. For instance, the planar rooted trees
of~\eqref{equ:examples_planar_rooted_trees} are depicted by
\begin{equation}
    \LeafPic, \quad
    \begin{tikzpicture}[scale=.3,Centering]
        \node[Leaf](1)at(0.00,-1.00){};
        \node[Node](0)at(0.00,0.00){};
        \draw[Edge](1)--(0);
        \node(r)at(0.00,1){};
        \draw[Edge](r)--(0);
    \end{tikzpicture}\,, \quad
    \begin{tikzpicture}[xscale=.2,yscale=.2,Centering]
        \node[Leaf](0)at(0.00,-1.50){};
        \node[Leaf](2)at(2.00,-1.50){};
        \node[Node](1)at(1.00,0.00){};
        \draw[Edge](0)--(1);
        \draw[Edge](2)--(1);
        \node(r)at(1.00,1.5){};
        \draw[Edge](r)--(1);
    \end{tikzpicture}\,, \quad
    \begin{tikzpicture}[xscale=.23,yscale=.2,Centering]
        \node[Leaf](0)at(0.00,-1.33){};
        \node[Leaf](3)at(2.00,-2.67){};
        \node[Node](1)at(1.00,0.00){};
        \node[Node](2)at(2.00,-1.33){};
        \draw[Edge](0)--(1);
        \draw[Edge](2)--(1);
        \draw[Edge](3)--(2);
        \node(r)at(1.00,1.5){};
        \draw[Edge](r)--(1);
    \end{tikzpicture}\,, \quad
    \begin{tikzpicture}[xscale=.18,yscale=.16,Centering]
        \node[Leaf](1)at(0.00,-6.75){};
        \node[Leaf](3)at(2.00,-6.75){};
        \node[Leaf](5)at(3.00,-2.25){};
        \node[Leaf](6)at(4.00,-4.50){};
        \node[Leaf](8)at(6.00,-4.50){};
        \node[Node](0)at(1.00,-2.25){};
        \node[Node](2)at(1.00,-4.50){};
        \node[Node](4)at(3.00,0.00){};
        \node[Node](7)at(5.00,-2.25){};
        \draw[Edge](0)--(4);
        \draw[Edge](1)--(2);
        \draw[Edge](2)--(0);
        \draw[Edge](3)--(2);
        \draw[Edge](5)--(4);
        \draw[Edge](6)--(7);
        \draw[Edge](7)--(4);
        \draw[Edge](8)--(7);
        \node(r)at(3.00,1.75){};
        \draw[Edge](r)--(4);
    \end{tikzpicture}\,.
\end{equation}
\medbreak

By definition of the Cartesian product and the list collection
operations over graded collections (see
Sections~\ref{subsubsec:cartesian_product_collections}
and~\ref{subsubsec:list_collections} of Chapter~\ref{chap:collections}),
the size of a planar rooted tree $\Tfr$ having a root arity of $k$
ssatisfies
\begin{equation} \label{equ:size_planar_rooted_trees}
    |\Tfr| = 1 + \sum\limits_{i \in [k]} |\Tfr_i|.
\end{equation}
In other words, the size of $\Tfr$ is the number of occurrences of
$\Node$ it contains.  We also deduce
from~\eqref{equ:combinatorial_set_planar_rooted_trees} that the
generating series of $\ColPRT$ satisfies
\begin{equation}
    \GeneratingSeries_\ColPRT(t) =
    \frac{t}{1 - \GeneratingSeries_\ColPRT(t)}
\end{equation}
so that it satisfies the quadratic algebraic equation
\begin{equation} \label{equ:generating_series_planar_rooted_trees}
    t - \GeneratingSeries_\ColPRT(t)
    + \GeneratingSeries_\ColPRT(t)^2 = 0.
\end{equation}
\medbreak

\subsubsection{Induction and structural induction}
\label{subsubsec:structural_induction}
One among the most obvious techniques to prove that all the planar
rooted trees of a subcollection $C$ of $\ColPRT$ satisfy a predicate
$P(\Tfr)$ (that is, a statement involving a variable $\Tfr$ taking value
in $C$) consists in performing a proof by induction on the size of the
trees of~$C$.
\medbreak

There is another method which is in some cases much more elegant than
this approach, called \Def{structural induction} on trees. A
subcollection $C$ of $\ColPRT$ is \Def{inductive} if $C$ is
nonempty and, if $\Tfr \in C$, all suffix subtrees $\Tfr_i$ of $\Tfr$
belong to~$C$. Observe in particular that $\Leaf$ belongs to any
inductive subcollection of~$\ColPRT$.
\medbreak

\begin{Theorem} \label{thm:structural_induction}
    Let $C$ be an inductive subcollection of $\ColPRT$ and $P(\Tfr)$ be
    a predicate on~$C$. If
    \begin{enumerate}[label={(\it\roman*)}]
        \item \label{item:structural_induction_1}
        the statement $P(\Leaf)$ holds;
        \item \label{item:structural_induction_2}
        for any $\Tfr_1, \dots, \Tfr_k \in C$ such that
        $\Tfr := \Par{\Node, \Par{\Tfr_1, \dots, \Tfr_k}}$ belongs to
        $C$, the fact that all $P\Par{\Tfr_i}$, $i \in [k]$, hold
        implies that $P(\Tfr)$ holds;
    \end{enumerate}
    then, all the planar rooted trees $\Sfr$ of $C$ satisfy~$P(\Sfr)$.
\end{Theorem}
\medbreak

Theorem~\ref{thm:structural_induction} provides a powerful tool to prove
properties $P(\Tfr)$ of planar rooted trees belonging to inductive
combinatorial subsets $C$. In practice, to perform a structural
induction in order to show that all the objects $\Tfr$ of $C$ satisfy
$P(\Tfr)$, we check that $C$ is inductive and that
Properties~\ref{item:structural_induction_1}
and~\ref{item:structural_induction_2} of
Theorem~\ref{thm:structural_induction} hold.
\medbreak

\subsubsection{Tree languages}
To rigorously specify nodes in planar rooted trees, we shall use a
useful interpretation of planar rooted trees as special languages on the
alphabet $\N_{\geq 1}$. Recall that a partial right monoid action of a
monoid $A^*$ of words (endowed with the concatenation product~$\Conc$)
on a set $S$ is a map $\Action : S \times A^* \to S$ satisfying
$x \Action \epsilon = x$, and for any $x \in S$, $u \in A^*$, and
$a \in A$, $x \Action u a$ is defined if and only if
$(x \Action u) \Action a$ is defined, and these two elements are the
same when they are defined. Let
\begin{equation}
    \Action : \ColPRT \times \N_{\geq 1}^* \to
    \ColPRT
\end{equation}
be the right partial monoid action defined recursively by
\begin{equation} \label{equ:action_monoid_trees}
    \Par{\Node, \Par{\Tfr_1, \dots, \Tfr_k}} \Action u :=
    \begin{cases}
        \Par{\Node, \Par{\Tfr_1, \dots, \Tfr_k}}
            & \mbox{if } u = \epsilon, \\
        \Tfr_i \Action v & \mbox{otherwise (}
            u = i v \mbox{ where } v \in \N_{\geq 1}^*
            \mbox{ and } i \in \N_{\geq 1} \mbox{)},
    \end{cases}
\end{equation}
for any $\Par{\Node, \Par{\Tfr_1, \dots, \Tfr_k}} \in \ColPRT$ and
$u \in \N_{\geq 1}^*$. Observe that this action is partial since each
$\Tfr_i$ in~\eqref{equ:action_monoid_trees} is well-defined only if $i$
is no greater than the root arity of~$\Tfr$. The \Def{tree language}
$\TreeLanguage(\Tfr)$ of $\Tfr$ is the finite language on $\N_{\geq 1}$
of all the words $u$ such that $\Tfr \Action u$ is a well-defined planar
rooted tree.
\medbreak

For instance, by setting
\begin{equation} \label{equ:example_planar_rooted_tree}
    \Tfr :=
    \begin{tikzpicture}[scale=.16,Centering]
        \node[Leaf](0)at(0.00,-2.40){};
        \node[Leaf](10)at(7.00,-7.20){};
        \node[Leaf](11)at(8.00,-2.40){};
        \node[Leaf](3)at(1.00,-9.60){};
        \node[Leaf](5)at(3.00,-9.60){};
        \node[Leaf](7)at(4.00,-4.80){};
        \node[Leaf](8)at(5.00,-7.20){};
        \node[Node](1)at(4.00,0.00){};
        \node[Node](2)at(2.00,-4.80){};
        \node[Node](4)at(2.00,-7.20){};
        \node[Node](6)at(4.00,-2.40){};
        \node[Node](9)at(6.00,-4.80){};
        \draw[Edge](0)--(1);
        \draw[Edge](10)--(9);
        \draw[Edge](11)--(1);
        \draw[Edge](2)--(6);
        \draw[Edge](3)--(4);
        \draw[Edge](4)--(2);
        \draw[Edge](5)--(4);
        \draw[Edge](6)--(1);
        \draw[Edge](7)--(6);
        \draw[Edge](8)--(9);
        \draw[Edge](9)--(6);
        \node(r)at(4.00,1.80){};
        \draw[Edge](r)--(1);
    \end{tikzpicture}\,,
\end{equation}
we have
\begin{equation}
    \Tfr \Action 1 = \LeafPic, \quad
    \Tfr \Action 231 = \LeafPic, \quad
    \Tfr \Action 3 = \LeafPic, \quad
    \Tfr \Action 21 =
    \begin{tikzpicture}[scale=.2,Centering]
        \node[Leaf](1)at(0.00,-2.67){};
        \node[Leaf](3)at(2.00,-2.67){};
        \node[Node](0)at(1.00,0.00){};
        \node[Node](2)at(1.00,-1.33){};
        \draw[Edge](1)--(2);
        \draw[Edge](2)--(0);
        \draw[Edge](3)--(2);
        \node(r)at(1.00,1.5){};
        \draw[Edge](r)--(0);
    \end{tikzpicture}\,, \quad
    \Tfr \Action 23 =
    \begin{tikzpicture}[scale=.17,Centering]
        \node[Leaf](0)at(0.00,-1.50){};
        \node[Leaf](2)at(2.00,-1.50){};
        \node[Node](1)at(1.00,0.00){};
        \draw[Edge](0)--(1);
        \draw[Edge](2)--(1);
        \node(r)at(1.00,1.75){};
        \draw[Edge](r)--(1);
    \end{tikzpicture}\,,
\end{equation}
and, among others, the actions of the words $11$, $24$, and $2321$ on
$\Tfr$ are all undefined. Moreover, the tree language of $\Tfr$ is
\begin{equation} \label{equ:example_tree_language}
    \TreeLanguage(\Tfr) =
    \{\epsilon, 1, 2, 21, 211, 2111, 2112, 22, 23, 231, 232, 3\}.
\end{equation}
\medbreak

Let $\Lca_\ColPRT$ be the combinatorial graded collection of all finite
and nonempty prefix languages $\Lca$ on $\N_{\geq 1}$ such that if
$ui \in \Lca$ where $u \in \N_{\geq 1}^*$ and $i \in \N_{\geq 2}$,
$u i' \in \Lca$ where $i' := i - 1$. The size of such a language is its
cardinality. For instance, the set $\TreeLanguage(\Tfr)$
of~\eqref{equ:example_tree_language} is an object of size $12$
of~$\Lca_\ColPRT$, and $\{\epsilon, 1, 11, 12, 2\}$ is an object of
size~$5$.
\medbreak

\begin{Proposition}
\label{prop:bijection_planar_rooted_trees_prefix_languages}
    The combinatorial graded collections $\ColPRT$ and
    $\Lca_\ColPRT$ are isomorphic. Seen as a morphism of
    combinatorial collections
    \begin{math}
        \TreeLanguage :
        \ColPRT \to \Lca_\ColPRT
    \end{math},
    $\TreeLanguage$ is an isomorphism between these two collections.
\end{Proposition}
\medbreak

Proposition~\ref{prop:bijection_planar_rooted_trees_prefix_languages} is
used in practice to define planar rooted trees through their languages.
This will be useful later when operations on planar rooted trees will be
described.
\medbreak

\subsubsection{Additional definitions}
\label{subsubsec:definitions_trees}
Let $\Tfr$ be a planar rooted tree. We say that each word of
$\TreeLanguage(\Tfr)$ is a \Def{node} of $\Tfr$. A node $u$ of $\Tfr$ is
an \Def{internal node} if there is an $i \in \N_{\geq 1}$ such that $ui$
is a node of $\Tfr$. A node $u$ of $\Tfr$ which is not an internal node
is a \Def{leaf}. The set of all internal nodes (resp. leaves) of $\Tfr$
is denoted by $\TreeLanguage_\Node(\Tfr)$ (resp.
$\TreeLanguage_\Leaf(\Tfr)$). The \Def{root} of $\Tfr$ is the node
$\epsilon$ (which can be either an internal node or a leaf). The
\Def{degree} $\deg(\Tfr)$ of $\Tfr$ is $\# \TreeLanguage_\Node(\Tfr)$
and the \Def{arity} $\Arity(\Tfr)$ of $\Tfr$ is
$\# \TreeLanguage_\Leaf(\Tfr)$. A node $u$ of $\Tfr$ is an
\Def{ancestor} of a node $v$ of $\Tfr$ if $u \ne v$ and
$u \PrefixOrder v$. Moreover, for any $i \in \N_{\geq 1}$, a node $v$ is
the \Def{$i$th child} of a node $u$ if $v = u i$. In this case, $u$ is
the (unique) \Def{father} of $v$. The \Def{arity} of a node is the
number of children it has. The lexicographic order on the words of
$\TreeLanguage(\Tfr)$ induces a total order on the nodes of $\Tfr$
called \Def{depth-first order}. The \Def{$i$th leaf} of $\Tfr$ is the
$i$th leaf encountered by considering the nodes of $\Tfr$ according to
the depth-first order. A \Def{path} in $\Tfr$ is a sequence
$\Par{u_1, \dots, u_k}$ of nodes of $\Tfr$ such that for any
$j \in [k - 1]$, $u_j$ is the father of $u_{j + 1}$. Such a path is
\Def{maximal} if $u_1$ is the root of $\Tfr$ and $u_k$ is a leaf. The
\Def{length} of a path is the number of nodes it contains. The
\Def{height} $\Height(\Tfr)$ of $\Tfr$ is the maximal length of its
maximal paths minus $1$. This is also the length of a longest word of
$\TreeLanguage(\Tfr)$ minus~$1$. For any node $u$ of $\Tfr$, the planar
rooted tree $\Tfr \Action u$ is the \Def{suffix subtree} of $\Tfr$
rooted at $u$. By extension, the \Def{$i$th suffix subtree} of $u$ is
the planar rooted tree $\Tfr \Action ui$ when $i$ is no greater than the
arity of $u$. A planar rooted tree $\Sfr$ is a \Def{prefix subtree} of
$\Tfr$ if $\TreeLanguage(\Sfr) \subseteq \TreeLanguage(\Tfr)$. A planar
rooted tree $\Sfr$ is a \Def{factor subtree} of $\Tfr$ rooted at a node
$u$ if $\Sfr$ is a prefix subtree of a suffix subtree of $\Tfr$ rooted
at~$u$.
\medbreak

Let us provide some examples for these notions. Consider the planar
rooted tree $\Tfr$ of~\eqref{equ:example_planar_rooted_tree}. Then,
\begin{subequations}
\begin{equation}
    \TreeLanguage_\Node(\Tfr) =
    \{\epsilon, 2, 21, 211, 23\},
\end{equation}
\begin{equation}
    \TreeLanguage_\Leaf(\Tfr) =
    \{1, 2111, 2112, 22, 231, 232, 3\},
\end{equation}
\end{subequations}
so that $\deg(\Tfr) = 5$ and $\Arity(\Tfr) = 7$. The $3$rd leaf of
$\Tfr$ is $2112$, and the $2$nd child of the internal node $23$ of
$\Tfr$ is $232$ and is a leaf. Besides, the sequences
$(\epsilon, 2, 21)$ and $(\epsilon, 2, 23)$ are nonmaximal paths in
$\Tfr$, and on the contrary, the paths $(\epsilon, 1)$,
$(\epsilon, 2, 21, 211, 2112)$, and $(\epsilon, 2, 22)$ are maximal. The
maximal path $(\epsilon, 2, 21, 211, 2112)$ have a maximal length among
all maximal paths of $\Tfr$ and thus, the height of $\Tfr$ is~$4$.
Finally, the planar rooted tree
\begin{equation}
    \Sfr :=
    \begin{tikzpicture}[xscale=.25,yscale=.18,Centering]
        \node[Leaf](0)at(0.00,-2.00){};
        \node[Leaf](3)at(1.00,-6.00){};
        \node[Leaf](5)at(2.00,-4.00){};
        \node[Leaf](6)at(3.00,-4.00){};
        \node[Leaf](7)at(4.00,-2.00){};
        \node[Node](1)at(2.00,0.00){};
        \node[Node](2)at(1.00,-4.00){};
        \node[Node](4)at(2.00,-2.00){};
        \draw[Edge](0)--(1);
        \draw[Edge](2)--(4);
        \draw[Edge](3)--(2);
        \draw[Edge](4)--(1);
        \draw[Edge](5)--(4);
        \draw[Edge](6)--(4);
        \draw[Edge](7)--(1);
        \node(r)at(2.00,1.50){};
        \draw[Edge](r)--(1);
    \end{tikzpicture}
\end{equation}
is a prefix subtree of $\Tfr$, and, the planar rooted tree
\begin{equation}
    \Rfr :=
    \begin{tikzpicture}[xscale=.25,yscale=.18,Centering]
        \node[Leaf](1)at(0.00,-3.33){};
        \node[Leaf](3)at(1.00,-1.67){};
        \node[Leaf](4)at(2.00,-1.67){};
        \node[Node](0)at(0.00,-1.67){};
        \node[Node](2)at(1.00,0.00){};
        \draw[Edge](0)--(2);
        \draw[Edge](1)--(0);
        \draw[Edge](3)--(2);
        \draw[Edge](4)--(2);
        \node(r)at(1.00,1.5){};
        \draw[Edge](r)--(2);
    \end{tikzpicture}\,,
\end{equation}
being a suffix subtree of $\Sfr$ rooted at the node $2$, is a factor
subtree of $\Tfr$ rooted at the node~$2$.
\medbreak

\subsection{Subcollections of planar rooted trees}
\label{subsec:families_planar_rooted_trees}
By basically restraining the possible arities of the internal nodes of
planar rooted trees, we obtain several subcollections of $\ColPRT$. We
review here the families formed by ladders, corollas, $k$-ary trees, and
Schröder trees. Besides, among these families, some admit alternative
size functions (and form therefore different combinatorial graded
collections).
\medbreak

\subsubsection{Ladders and corollas}
A \Def{ladder} is a planar rooted tree of arity $1$. The first ladders
are
\begin{equation}
    \LeafPic, \quad
\,.
\end{equation}
This set of ladders forms a subcollection $\ColLad$ of $\ColPRT$.
Besides, a \Def{corolla} is a planar rooted tree of degree $1$. The
first corollas are
\begin{equation}
    %
\,.
\end{equation}
This set of corollas forms a subcollection $\ColCor$ of $\ColPRT$.
Observe that $(\Node, (\Leaf))$ is the only planar rooted that is both a
ladder and a corolla.
\medbreak

\subsubsection{$k$-ary trees} \label{subsubsec:k_ary_trees}
Let $k \in \N_{\geq 1}$. A \Def{$k$-ary} tree is a planar rooted tree
$\Tfr$ such that all internal nodes are of arity $k$. For instance, the
first $3$-ary trees are
\begin{equation}
    \LeafPic, \quad
    %
\,.
\end{equation}
This set of $k$-ary trees forms a subcollection $\ColAry^{(k)}$ of
$\ColPRT$ expressing recursively as
\begin{equation}
    \ColAry^{(k)} = \{\Leaf\}
    + \BBrack{\{\Node\}, \List^+_{\{k\}}\Par{\ColAry^{(k)}}}_\times^+,
\end{equation}
where $\Leaf$ and $\Node$ are both atomic. One can immediately observe
that $\ColAry^{(1)} = \ColLad$.
\medbreak

By structural induction (see Theorem~\ref{thm:structural_induction}) on
$\ColAry^{(k)}$ (which is an inductive subcollection of $\ColPRT$), it
follows that for any $k$-ary tree $\Tfr$, the arity and the degree of
$\Tfr$ are related by
\begin{equation} \label{equ:k-ary_trees_arity_degree}
    \Arity(\Tfr) - \deg(\Tfr) (k - 1) = 1.
\end{equation}
This implies that a $k$-ary tree of a given arity has an imposed degree
and conversely, a $k$-ary tree of a given degree has an imposed arity.
Hence, since the size of a $k$-ary tree $\Tfr$ is
\begin{math}
    \Arity(\Tfr) + \deg(\Tfr)
\end{math}
and there are finitely many planar rooted trees of a fixed size, there
are finitely many $k$-ary trees of a fixed arity, and there are finitely
many $k$-ary trees of a fixed degree. As a consequence, the graded
collections $\ColAry^{(k)}_\Leaf$ and $\ColAry^{(k)}_\Node$ of all
$k$-ary trees such that the size of a tree of $\ColAry^{(k)}_\Leaf$ is
its arity and the size of a tree of $\ColAry^{(k)}_\Node$ is its degree
are combinatorial. On the one hand, the generating series of
$\ColAry^{(k)}_\Leaf$ satisfies the algebraic equation
\begin{equation} \label{equ:generating_series_k-ary_tree_leaf}
    t - \GeneratingSeries_{\ColAry^{(k)}_\Leaf}(t)
    + {\GeneratingSeries_{\ColAry^{(k)}_\Leaf}(t)}^k
    = 0.
\end{equation}
On the other hand, the generating series of $\ColAry^{(k)}_\Node$
satisfies the algebraic equation
\begin{equation} \label{equ:generating_series_k-ary_tree_node}
    1 - \GeneratingSeries_{\ColAry^{(k)}_\Node}(t)
    + t {\GeneratingSeries_{\ColAry^{(k)}_\Node}(t)}^k
    = 0
\end{equation}
and one can deduce that
\begin{equation} \label{equ:number_k-ary_trees}
    \# \ColAry^{(k)}_\Node(n) =
    \frac{1}{(k - 1) n + 1} \binom{kn}{n}.
\end{equation}
For instance, the integer sequences of $\ColAry^{(k)}_\Node$ begin with
\begin{subequations}
\begin{equation}
    1, 1, 1, 1, 1, 1, 1, 1,
    \qquad k = 1,
\end{equation}
\begin{equation}
    1, 1, 2, 5, 14, 42, 132, 429, 1430,
    \qquad k = 2,
\end{equation}
\begin{equation}
    1, 1, 3, 12, 55, 273, 1428, 7752, 43263,
    \qquad k = 3,
\end{equation}
\begin{equation}
    1, 1, 4, 22, 140, 969, 7084, 53820, 420732,
    \qquad k = 4.
\end{equation}
\end{subequations}
The second, third, and fourth sequences above are, respectively,
Sequences~\OEIS{A000108}, \OEIS{A001764}, and~\OEIS{A002293}
of~\cite{Slo}. These are known as the \Def{Fuss-Catalan numbers}.
\medbreak

From now on, we call \Def{binary tree} any $2$-ary tree. Recall that
these objects have been introduced in
Section~\ref{subsubsec:binary_trees} of Chapter~\ref{chap:collections}.
If $\Tfr$ is a binary tree and $u$ is an internal node of $\Tfr$, $u1$
and $u2$ are nodes of $\Tfr$. We call $u1$ (resp. $u2$) the \Def{left}
(resp. \Def{right}) \Def{child} of $u$, and $\Tfr \Action u1$ (resp.
$\Tfr \Action u2$) the \Def{left} (resp. \Def{right}) \Def{subtree} of
$u$ in $\Tfr$. The left (resp. right) subtree of $\Tfr$ is the
\Def{left} (resp. \Def{right}) \Def{subtree} of the root of $\Tfr$.
Besides, a \Def{left} (resp. \Def{right}) \Def{comb tree} is a binary
tree $\Tfr$ such that for all internal nodes $u$ of $\Tfr$, all right
(resp. left) subtrees of $u$ are leaves. The \Def{infix order} induced
by $\Tfr$ is the total order on the set of its internal nodes defined
recursively by setting that all the internal nodes of $\Tfr \Action 1$
are smaller than the root of $\Tfr$, and that the root of $\Tfr$ is
smaller than all the internal nodes of~$\Tfr \Action 2$.
\medbreak

Let us denote by $\ColBT_\Leaf$ the combinatorial graded collection
$\ColAry^{(2)}_\Leaf$ of binary trees where the size of a tree is its
arity. As a consequence
of~\eqref{equ:generating_series_planar_rooted_trees}
and~\eqref{equ:generating_series_k-ary_tree_leaf}, we observe that the
generating series of $\ColPRT$ satisfies the same algebraic relation as
the one of $\ColBT_\Leaf$. Therefore, $\ColPRT$ and $\ColBT_\Leaf$ are
isomorphic as graded collections. Let us describe an explicit
isomorphism between these two collections. Let
$\phi : \ColPRT \to \ColBT_\Leaf$ be the map recursively defined, for
any planar rooted tree $\Tfr$, by
\begin{equation} \label{equ:bijection_planar_rooted_trees_binary_trees}
    \phi(\Tfr) :=
    \begin{cases}
        \Leaf \in \ColBT_\Leaf & \mbox{if } \Tfr = \Leaf, \\
        \Par{\Node,
            \Par{\phi(\Tfr_1),
                \phi\Par{\Par{\Node,\Par{\Tfr_2, \dots, \Tfr_k}}}}}
            & \mbox{otherwise (}
            \Tfr = \Par{\Node, \Par{\Tfr_1, \Tfr_2, \dots, \Tfr_k}}
            \mbox{ with } k \in \N_{\geq 1} \mbox{)} .
    \end{cases}
\end{equation}
One has, for instance,
\begin{subequations}
\begin{equation}
    \phi\Par{
\,.
\end{equation}
\end{subequations}
\medbreak

\begin{Proposition}
\label{prop:bijection_planar_rooted_trees_binary_trees}
    The combinatorial graded collections $\ColPRT$ and
    $\ColBT_\Leaf$ are isomorphic. The map~$\phi$ defined
    by~\eqref{equ:bijection_planar_rooted_trees_binary_trees} is an
    isomorphism between these two collections.
\end{Proposition}
\medbreak

\subsubsection{Schröder trees} \label{subsubsec:schroder_trees}
A \Def{Schröder tree} is a planar rooted tree such that all internal
nodes are of arities $2$ or more. Some among the first Schröder trees
are
\begin{equation}
    \LeafPic, \quad
\,.
\end{equation}
This set of Schröder trees forms a subcollection $\ColSch$ of
$\ColPRT$ expressing recursively as
\begin{equation} \label{equ:combinatorial_set_schroder_trees}
    \ColSch = \{\Leaf\} +
    \BBrack{\{\Node\}, \List_{\N_{\geq 2}}^+\Par{\ColSch}}_\times^+,
\end{equation}
where $\Leaf$ and $\Node$ are both atomic.
\medbreak

By structural induction on $\ColSch$ (which is an inductive
subcollection of $\ColPRT$), it follows that there are finitely many
Schröder trees of a given arity $n$. For this reason, the graded
collection $\ColSch_\Leaf$ of all the Schröder trees such that the size
of a tree of $\ColSch_\Leaf$ is its arity is combinatorial. Conversely,
considering the degrees of the trees for their sizes does not form a
combinatorial graded collection since there are infinitely many Schröder
trees of degree $1$ (the corollas). The generating series of
$\ColSch_\Leaf$ satisfies the algebraic quadratic equation
\begin{equation}
    t - (1 + t) \GeneratingSeries_{\ColSch_\Leaf}(t)
    + 2 {\GeneratingSeries_{\ColSch_\Leaf}(t)}^2 = 0.
\end{equation}
Let $\Narayana(n, k)$ be the number of binary trees of arity $n$ having
exactly $k$ internal nodes having an internal node as a left child.
Then, for all $0 \leq k \leq n - 2$, it is known that
\begin{equation} \label{equ:narayana_numbers}
    \Narayana(n, k) =
    \frac{1}{k + 1} \binom{n - 2}{k} \binom{n - 1}{k}.
\end{equation}
These are \Def{Narayana numbers}. The cardinalities of the sets
$\ColSch_\Leaf(n)$ express by
\begin{equation} \label{equ:enumeration_schroder_trees}
    \# \ColSch_\Leaf(n) =
    \sum_{k \in [0, n - 2]}
    2^k \, \Narayana(n, k),
\end{equation}
for all $n \in \N_{\geq 2}$. The integer sequence of $\ColSch_\Leaf$
begins by
\begin{equation}
    1, 1, 3, 11, 45, 197, 903, 4279, 20793
\end{equation}
and forms Sequence~\OEIS{A001003} of~\cite{Slo}.
\medbreak

\section{Syntax trees} \label{sec:syntax_trees}
We are now in position to introduce syntax trees and rewrite systems
on syntax trees. These objects are central in the theory of operads
since the elements of free nonsymmetric operads can be seen as syntax
trees. Rewrite systems on syntax trees provide tools to establish
presentations by generators and relations of operads.
\medbreak

\subsection{Collections of syntax trees}
\label{subsec:graded_collection_syntax_trees}
Syntax trees are, roughly speaking, planar rooted trees where internal
nodes are labeled by objects of a fixed graded collection. These trees
can be endowed with two size functions (where the size is the degree or
the arity), leading to the definition of two distinct graded collections
of syntax trees.
\medbreak

\subsubsection{Main definitions}
Let $C$ be an augmented graded collection. A \Def{syntax tree} on $C$
(or, for short, a \Def{$C$-syntax tree}) is a planar rooted tree $\Tfr$
endowed with a map
\begin{equation}
    \omega_\Tfr : \TreeLanguage_\Node(\Tfr) \to C
\end{equation}
sending each internal node $u$ of $\Tfr$ of arity $k$ to an object of
size $k$ of $C$. This map $\omega_\Tfr$ is the \Def{labeling map} of
$\Tfr$. We say that an internal node $u$ of $\Tfr$ is \Def{labeled} by
$x \in C$ if $\omega_\Tfr(u) = x$. The collection $C$ is the
\Def{labeling collection} of $\Tfr$. The \Def{underlying planar rooted
tree} of $\Tfr$ is the planar rooted tree obtained by forgetting the map
$\omega_\Tfr$. For any $x \in C$, the \Def{corolla} labeled by $x$ is
the $C$-syntax tree $\Corolla(x)$ having exactly one internal node
labeled by $x$ and with $|x|$ leaves as children. All the notions about
planar rooted trees defined in Section~\ref{sec:planar_rooted_trees}
apply to $C$-syntax trees as well. More precisely, for any property
$P(\Sfr)$ on a planar rooted tree $\Sfr$, we say that $P(\Tfr)$ holds
if the underlying planar rooted tree $\Sfr$ of $\Tfr$ is such that
$P(\Sfr)$ holds. Moreover, the notions of suffix, prefix, and factor
subtrees of planar rooted trees naturally extend on $C$-syntax trees by
taking into account the labeling maps. In graphical representations of a
$C$-syntax tree $\Tfr$, instead of drawing each internal node $u$ of
$\Tfr$ by $\NodePic$, we draw $u$ by its label~$\omega_\Tfr(u)$.
\medbreak

For instance, consider the labeling collection
$C := C(1) \sqcup C(2) \sqcup C(3)$ where $C(1) := \{\Asf, \Bsf\}$,
$C(2) := \{\Csf\}$, and $C(3) := \{\Dsf, \Esf\}$, and the planar rooted
tree
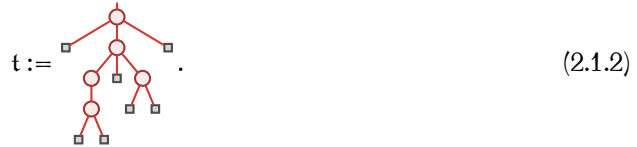
\begin{equation} \label{equ:example_syntax_tree}
    \Tfr :=
\,.
\end{equation}
By endowing $\Tfr$ with the labeling map defined by
$\omega_\Tfr(\epsilon) := \Esf$, $\omega_\Tfr(2) := \Dsf$,
$\omega_\Tfr(21) := \Asf$, $\omega_\Tfr(211) := \Csf$, and
$\omega_\Tfr(23) := \Csf$, $\Tfr$ is a $C$-syntax tree. This $C$-syntax
tree is depicted more concisely as
\begin{equation}
    %
\,.
\end{equation}
\medbreak

We denote by $\ColST^C$ the graded collection of all the
$C$-syntax trees, where the size of a $C$-syntax tree $\Tfr$ is the size
of its underlying planar rooted tree in $\ColPRT$. When $C$
is additionally combinatorial, by structural induction on planar rooted
trees, it follows that for any $\Tfr \in \ColPRT$, there are
finitely many labeling maps $\omega_\Tfr$ for $\Tfr$. For this reason,
$\ColST^C$ is in this case combinatorial. Besides, let
$\ColLad^C$, $\ColCor^C$, $\ColAry^{(k), C}$, and $\ColSch^C$ be,
respectively, the subcollections of $\ColST^C$ consisting in the
$C$-syntax trees whose underlying planar rooted trees are ladders,
corollas, $k$-ary trees, and Schröder trees. The concepts of inductive
subcollections of $\ColST^C$ and of structural induction presented
in Section~\ref{subsubsec:structural_induction} extend obviously on
$C$-syntax trees.
\medbreak

\subsubsection{Alternative definition and generating series}
The graded collection $\ColST^C$ can be described as
follows. Let $\Sca^C$ be the graded collection satisfying the relation
\begin{equation} \label{equ:expression_combinatorial_set_syntax_trees}
    \Sca^C = \{\Leaf\}
    + \BBrack{\{\Node\}, C\Compo \Sca^C}_\times^+
\end{equation}
where both $\Leaf$ and $\Node$ are atomic, and $\Compo$ is the
composition product over graded collections defined in
Section~\ref{subsubsec:composition_collections} of
Chapter~\ref{chap:collections}. Then, the combinatorial collections
$\ColST^C$ and $\Sca^C$ are isomorphic through the morphism
$\phi : \ColST^C \to \Sca^C$ of combinatorial collections
recursively defined, for any $\Tfr \in \ColST^C$ of root
arity $k$, by
\begin{equation}
    \phi(\Tfr) :=
    \begin{cases}
        \Leaf \in \Sca^C & \mbox{if } \Tfr = \Leaf, \\
        \Par{\Node,
            \Par{\omega_\Tfr(\epsilon),
                \Par{\phi\Par{\Tfr_1}, \dots, \Par{\Tfr_k}}}}
        & \mbox{otherwise}.
    \end{cases}
\end{equation}
From this equivalence
and~\eqref{equ:expression_combinatorial_set_syntax_trees}, we obtain,
when $C$ is combinatorial, that the generating series of
$\ColST^C$ satisfies
\begin{equation}
    \GeneratingSeries_{\ColST^C}(t) =
    t +
    t \, \GeneratingSeries_C
    \Par{\GeneratingSeries_{\ColST^C}(t)},
\end{equation}
where $\GeneratingSeries_C(t)$ is the generating series of $C$. For
instance, by considering the combinatorial collection $C$ defined above,
we have $\GeneratingSeries_C(t) = 2t + t^2 + 2t^3$, so that
\begin{equation}
    t + (2t - 1) \, \GeneratingSeries_{\ColST^C}(t)
    + t \, {\GeneratingSeries_{\ColST^C}(t)}^2
    + 2t \, {\GeneratingSeries_{\ColST^C}(t)}^3 = 0.
\end{equation}
\medbreak

\subsubsection{Subcollections of syntax trees}
For well-chosen combinatorial augmented graded collections $C$, it is
possible to recover a large part of the families of planar rooted
trees described in Section~\ref{subsec:families_planar_rooted_trees}.
Indeed, one has
\begin{math}
    \ColST^{\N_{\geq 1}}
    \simeq \ColPRT
\end{math},
\begin{math}
    \ColST^{\N_{\geq 2}}
    \simeq \ColSch
\end{math},
and, when $\Node_k$ is an object of size $k \in \N_{\geq 1}$,
\begin{math}
    \ColST^{\left\{\Node_k\right\}}
    \simeq \ColAry^{(k)}
\end{math}.
\medbreak

\subsubsection{Alternative sizes}
Let $\ColST^C_\Leaf$ be the graded collection of all the $C$-syntax
trees such that the size of a tree is its arity. One has
$\ColST_\Leaf^C \simeq \Sca^C$ where $\Sca^C$ is the graded collection
defined in~\eqref{equ:expression_combinatorial_set_syntax_trees} wherein
$\Leaf$ is atomic and $\Node$ is of size $0$. When $C$ is graded,
combinatorial, augmented, and has no object of size $1$, we can show by
structural induction on $\ColST_\Leaf^C$ that there are finitely many
$C$-syntax trees of a given arity $n \in \N_{\geq 1}$. For this reason,
$\ColST^C_\Leaf$ is combinatorial. In this case, the generating series
of $\ColST_\Leaf^C$ satisfies
\begin{equation}
    \GeneratingSeries_{\ColST_\Leaf^C}(t) =
    t +
    \GeneratingSeries_C
    \Par{\GeneratingSeries_{\ColST_\Leaf^C}(t)}.
\end{equation}
\medbreak

Let also $\ColST_\Node^C$ be the graded collection of all the
$C$-syntax trees such that the size of a tree is its degree. One has
$\ColST_\Node^C \simeq \Sca^C$ where $\Sca^C$ is the graded collection
defined in~\eqref{equ:expression_combinatorial_set_syntax_trees}
wherein $\Node$ is atomic and $\Leaf$ is of size $0$. When $C$ is
graded, augmented, and finite, we can show by structural induction on
$\ColST_\Node^C$ that there are finitely many $C$-syntax trees of a
given degree $n$. For this reason, $\ColST_\Node^C$ is combinatorial. In
this case, the generating series of $\ColST_\Node^C$ satisfies
\begin{equation}
    \GeneratingSeries_{\ColST_\Node^C}(t) =
    1 +
    t \, \GeneratingSeries_C
    \Par{\GeneratingSeries_{\ColST_\Node^C}(t)}.
\end{equation}
Observe that $\ColST_\Node^C$ is not an augmented graded collection.
\medbreak

\subsection{Grafting operations} \label{subsec:grafting_operations}
Three fundamental grafting operations on syntax trees are presented
here. These operations turn $\ColST^C_\Leaf$ into a collection
with concentrated products in the sense of
Section~\ref{subsubsec:collections_with_products} of
Chapter~\ref{chap:collections}.
\medbreak

\subsubsection{Partial grafting}
\label{subsubsec:partial_grafting_syntax_trees}
Let for any $n, m \in \N_{\geq 1}$ and $i \in [n]$ the product
\begin{equation}
    \circ_i^{(n, m)} :
    \ColST^C_\Leaf(n) \times \ColST^C_\Leaf(m)
    \to \ColST^C_\Leaf
\end{equation}
where for any $\Tfr \in \ColST^C_\Leaf(n)$,
$\Sfr \in \ColST^C_\Leaf(m)$, and $i \in [n]$, the syntax tree
$\Rfr := \Tfr \circ_i^{(n, m)} \Sfr$ is defined as follows. The
underlying planar rooted tree of $\Rfr$ admits the tree language
\begin{equation} \label{equ:tree_language_grafting}
    \TreeLanguage(\Rfr) s:=
    \Par{\TreeLanguage(\Tfr) \setminus \{u\}}
    \cup
    \left\{u v : v \in \TreeLanguage(\Sfr)\right\},
\end{equation}
and the labeling map of $\Rfr$ satisfies, for any
$w \in \TreeLanguage_\Node(\Rfr)$,
\begin{equation}
    \omega_\Rfr(w) :=
    \begin{cases}
            \omega_\Tfr(w)
                & \mbox{if } w \in \TreeLanguage_\Node(\Tfr), \\
            \omega_\Sfr(v)
                & \mbox{otherwise (}
                w = u v \mbox{ and } v \in \TreeLanguage_\Node(\Sfr)
                \mbox{)}.
    \end{cases}
\end{equation}
Observe that by
Proposition~\ref{prop:bijection_planar_rooted_trees_prefix_languages},
$\Rfr$ is wholly specified by its tree language $\TreeLanguage(\Rfr)$
defined in~\eqref{equ:tree_language_grafting}. In more intuitive terms,
the tree $\Rfr$ is obtained by connecting the root of $\Sfr$ onto the
$i$th leaf of $\Tfr$. For instance, by considering the same labeling
collection $C$ as above,
\begin{equation}
\,.
\end{equation}
We call each $\circ_i^{(n, m)}$ a \Def{partial grafting operation}.
\medbreak

Observe that since
\begin{equation}
    \Arity\Par{\Tfr \circ_i^{(n, m)} \Sfr} =
    \Arity(\Tfr) + \Arity(\Sfr) - 1,
\end{equation}
the product $\circ_i^{(n, m)}$ is concentrated and is not graded.
Besides, by a slight abuse of notation, we shall sometimes omit the
$(n, m)$  in the notation of $\circ_i^{(n, m)}$ in order to denote it in
a more concise way by~$\circ_i$.
\medbreak

\subsubsection{Complete grafting}
\label{subsubsec:complete_grafting_syntax_trees}
Let for any $n, m_1, \dots, m_n \in \N_{\geq 1}$ the product
\begin{equation}
    \circ^{\Par{m_1, \dots, m_n}} :
    \ColST^C_\Leaf(n) \times \ColST^C_\Leaf\Par{m_1}
    \times \dots \times \ColST^C_\Leaf\Par{m_n}
    \to \ColST^C_\Leaf
\end{equation}
where for any $\Tfr \in \ColST^C_\Leaf(n)$,
$\Sfr_1 \in \ColST^C_\Leaf\Par{m_1}$, \dots,
$\Sfr_n \in \ColST^C_\Leaf\Par{m_n}$,
\begin{equation} \label{equ:complete_grafting_syntax_trees}
    \circ^{\Par{m_1, \dots, m_n}}
    \Par{\Tfr, \Sfr_1, \dots, \Sfr_n}
    :=
    \Par{\dots \Par{\Par{\Tfr \circ_n \Sfr_n} \circ_{n - 1}
    \Sfr_{n - 1}}\dots}
    \circ_1 \Sfr_1.
\end{equation}
In more intuitive terms, the syntax tree expressed
by~\eqref{equ:complete_grafting_syntax_trees} is obtained by connecting
the root of each $\Sfr_i$ onto the $i$th leaf of~$\Tfr$. For instance,
by considering the same labeling collection $C$ as before,
\begin{equation}
    \circ^{(2, 1, 4, 2)}
    \Par{
\,.
\end{equation}
We call each $\circ^{\Par{m_1, \dots, m_n}}$ a
\Def{complete grafting operation}.
\medbreak

Observe that since
\begin{equation}
    \Arity\Par{\circ^{(m_1, \dots, m_n)}
    \Par{\Tfr, \Sfr_1, \dots, \Sfr_n}} =
    \Arity(\Sfr_1) + \dots + \Arity(\Sfr_n),
\end{equation}
the product $\circ^{\Par{m_1, \dots, m_n}}$ is concentrated and is not
graded (because the size of the first operand $\Tfr$ does not intervene
in the size of the result). Besides, by a slight abuse of notation, we
shall sometimes omit the $\Par{m_1, \dots, m_n}$  in the notation of
$\circ^{\Par{m_1, \dots, m_n}}$ in order to denote it in a more concise
way by~$\circ$. Moreover, we shall denote by
\begin{math}
    \Tfr \circ \left[\Sfr_1, \dots, \Sfr_n\right]
\end{math}
the $C$-syntax tree
\begin{math}
    \circ\Par{\Tfr, \Sfr_1, \dots, \Sfr_n}
\end{math}.
\medbreak

\subsubsection{Context grafting}
\label{subsubsec:context_grafting_syntax_trees}
Let for any $n, m, k_1, \dots, k_m \in \N_{\geq 1}$ and $i \in [n]$ the
product
\begin{equation}
    \CompoContext_i^{\Par{n, k_1, \dots, k_m}} :
    \ColST^C_\Leaf(n) \times \ColST^C_\Leaf(m)
    \times \ColST^C_\Leaf\Par{k_1}
    \times \dots \times \ColST^C_\Leaf\Par{k_m}
    \to \ColST^C_\Leaf
\end{equation}
where for any $\Tfr \in \ColST^C_\Leaf(n)$,
$\Sfr \in \ColST^C_\Leaf(m)$,
$\Rfr_1 \in \ColST^C_\Leaf(k_1)$, \dots,
$\Rfr_m \in \ColST^C_\Leaf(k_m)$,
\begin{equation} \label{equ:context_grafting_syntax_trees}
    \CompoContext_i\Par{\Tfr, \Sfr, \Rfr_1, \dots, \Rfr_m}
    :=
    \Tfr \circ_i
    \Par{\Sfr \circ \left[\Rfr_1, \dots, \Rfr_m\right]}.
\end{equation}
For instance, by considering the same labeling collection $C$ as before,
\begin{equation}
    \CompoContext_3^{(3, 1, 2, 2)}\Par{
\,.
\end{equation}
We call each $\CompoContext_i^{\Par{n, k_1, \dots, k_m}}$ a
\Def{context grafting operation}.
\medbreak

Observe that since
\begin{equation}
    \Arity\Par{\CompoContext_i^{\Par{n, k_1, \dots, k_m}}
    \Par{\Tfr, \Sfr, \Rfr_1, \dots, \Rfr_m}}
    =
    \Arity(\Tfr) - 1 +
    \Arity(\Rfr_1) + \dots + \Arity(\Rfr_m),
\end{equation}
the product $\CompoContext_i^{\Par{n, k_1, \dots, k_m}}$ is concentrated
and is not graded. Besides, by a slight abuse of notation, we shall
sometimes omit the $\Par{n, k_1, \dots, k_m}$  in the notation of
$\CompoContext_i^{\Par{n, k_1, \dots, k_m}}$ in order to denote it in a
more concise way by~$\CompoContext_i$.
\medbreak

\subsection{Patterns and rewrite systems}
\label{subsec:rewrite_rules_syntax_trees}
We focus now on the theory of rewrite systems exposed in
Section~\ref{sec:rewrite_systems} of Chapter~\ref{chap:collections} on
the particular case of syntax trees. Intuitively, a rewrite rule on
syntax trees works by replacing factor subtrees in a syntax tree by
other ones. We explain techniques to prove termination and confluence of
these particular rewrite systems.
\medbreak

\subsubsection{Occurrence and avoidance of patterns}
Let $C$ be an augmented graded collection, and $\Sfr$ and $\Tfr$ be two
$C$-syntax trees. For any node $u$ of $\Tfr$, $\Sfr$ \Def{occurs} at
position $u$ in $\Tfr$ if $\Sfr$ is a factor subtree of $\Tfr$ rooted at
$u$. In this case, we say that $\Tfr$ \Def{admits an occurrence} of the
\Def{pattern} $\Sfr$. Conversely,  $\Tfr$ \Def{avoids} $\Sfr$ if there
is no occurrence of $\Sfr$ in~$\Tfr$.
\medbreak

This property can be rephrased as follows by using the context grafting
operations. A syntax tree $\Tfr$ admits an occurrence of $\Sfr$ if there
exists syntax trees $\Tfr'$, $\Rfr_1$, \dots, $\Rfr_{|\Sfr|}$ and
$i \in \left[\left|\Tfr'\right|\right]$ such that
\begin{equation}
    \Tfr =
    \CompoContext_i\Par{\Tfr', \Sfr, \Rfr_1, \dots, \Rfr_{|\Sfr|}}.
\end{equation}
\medbreak

By extension, $\Tfr$ avoids a set $P$ of $C$-syntax trees if $\Tfr$
avoids all the patterns of $P$. For instance, consider the graded
collection $C := C(2) \sqcup C(3)$ where $C(2) := \{\Asf, \Bsf\}$ and
$C(3) := \{\Csf\}$, and the $C$-syntax tree
\begin{equation}
    \Tfr :=
\,.
\end{equation}
Then, $\Tfr$ admits an occurrence of
\begin{equation}
    %
\end{equation}
at position $1$ and two occurrences of
\begin{equation}
    %
\end{equation}
at positions $11$ and~$21$.
\medbreak

\subsubsection{Rewrite systems} \label{subsubsec:rewrite_systems}
Let $\Par{\ColST^C_\Leaf, \Rew}$ be a rewrite system on syntax trees
and
\begin{equation}
    \Pca :=
    \left\{\CompoContext^{\Par{n, k_1, \dots, k_m}}_i :
    n, m, k_1, \dots, k_m \in \N_{\geq 1}, i \in [n]
    \right\}
\end{equation}
the set of all the context grafting operations. Since, as we observed
in Section~\ref{subsubsec:context_grafting_syntax_trees}, all products
of $\Pca$ are concentrated, we can consider the $\Pca$-closure
of $\Par{\ColST^C_\Leaf, \Rew}$. Therefore, let us denote by
$\Par{\ColST^C_\Leaf, \RewContext}$ the $\Pca$-closure of
$\Par{\ColST^C_\Leaf, \Rew}$, called simply \Def{closure} of
$\Par{\ColST^C_\Leaf, \Rew}$. In other terms, $\RewContext$ is the
rewrite rule satisfying
\begin{equation} \label{equ:closure_rewrite_rule_syntax_trees}
    \CompoContext_i\Par{\Tfr, \Sfr, \Rfr_1, \dots, \Rfr_m}
    \RewContext
    \CompoContext_i\Par{\Tfr, \Sfr', \Rfr_1, \dots, \Rfr_m}
\end{equation}
for any $C$-syntax trees $\Tfr$, $\Sfr$, $\Sfr'$, $\Rfr_1$, \dots,
$\Rfr_m$ where $\Tfr$ is of arity $n$, $\Sfr$ is of arity $m$,
$i \in [n]$, and $\Sfr \Rew \Sfr'$. In intuitive terms, one has
$\Qfr \RewContext \Qfr'$ for two $C$-syntax trees $\Qfr$ and $\Qfr'$ if
there are two $C$-syntax trees $\Sfr$ and $\Sfr'$ such that
$\Sfr \Rew \Sfr'$ and, by replacing an occurrence of $\Sfr$ by $\Sfr'$
in $\Qfr$, we obtain~$\Qfr'$. For instance, by considering the same
labeling set $C$ as before, let $\Par{\ColST^C_\Leaf, \Rew}$ be the
rewrite system defined by
\begin{equation} \label{equ:example_rewrite_rule_syntax_trees}
\,.
\end{equation}
One has the following chain of rewritings
\begin{equation}
    %
\,.
\end{equation}
\medbreak

Observe by the way that the right rotation operation on binary trees
considered in Section~\ref{subsubsec:tamari_order} of
Chapter~\ref{chap:collections} can be expressed as the closure of the
rewrite system $(\ColST^B_\Leaf, \Rew)$ such that
$B := B(2) := \{\Bsf\}$ defined by
\begin{equation}
    %
\,.
\end{equation}
\medbreak

In this text, we shall mainly consider rewrite systems
$\Par{\ColST^C_\Leaf, \RewContext}$ defined as closures of rewrite
systems $\Par{\ColST^C_\Leaf, \Rew}$ such that the number of pairs
$(\Tfr, \Tfr')$ satisfying $\Tfr \Rew \Tfr'$ is finite. We say in this
case that $\Par{\ColST^C_\Leaf, \RewContext}$ is of \Def{finite type}.
In this context, the \Def{degree} of $\Par{\ColST^C_\Leaf, \RewContext}$
is the maximal degree among the $C$-syntax trees appearing as left
members of $\Rew$. The \Def{arity} of
$\Par{\ColST^C_\Leaf, \RewContext}$ is the maximal arity among the
$C$-syntax trees appearing as left (or, equivalently, as right) members
of~$\Rew$.
\medbreak

\subsubsection{Proving termination}
\label{subsubsec:proving_termination}
We have observed in Section~\ref{subsubsec:termination} of
Chapter~\ref{chap:collections} that termination invariants provide tools
to show that a combinatorial rewrite system is terminating. This idea
extends on rewrite systems on syntax trees defined as closures of other
ones in the following way.
\medbreak

Let $\Par{\ColST^C_\Leaf, \Rew}$ be a combinatorial rewrite system and
$\Par{\ColST^C_\Leaf, \RewContext}$ be its closure. Assume that
$\theta : \ColST^C_\Leaf \to \Qca$ is a termination invariant for
$\Par{\ColST^C_\Leaf, \Rew}$, where $(\Qca, \Ord)$ is a poset. We say
that $\theta$ is \Def{compatible with the closure} if, for any
$C$-syntax trees $\Sfr$ and $\Sfr'$ such that $\Sfr \Rew \Sfr'$, the
inequality
\begin{equation} \label{equ:compatibility_termination_invariant}
    \theta\Par{\CompoContext_i\Par{\Tfr, \Sfr, \Rfr_1, \dots, \Rfr_m}}
    \OrdStrict
    \theta\Par{\CompoContext_i\Par{\Tfr, \Sfr', \Rfr_1, \dots, \Rfr_m}}
\end{equation}
holds for all $C$-syntax trees $\Tfr$, $\Rfr_1$, \dots, $\Rfr_m$, and
all $i \in [\Arity(\Tfr)]$ where $m := \Arity(\Sfr) = \Arity(\Sfr')$.
Now, as a consequence of~\eqref{equ:closure_rewrite_rule_syntax_trees}
and Theorem~\ref{thm:terminating_rewrite_rules_posets} of
Chapter~\ref{chap:collections}, one has the following result.
\medbreak

\begin{Proposition} \label{prop:compatible_terminating_invariant}
    Let $C$ be a combinatorial augmented graded collection without
    object of size~$1$, $\Par{\ColST^C_\Leaf, \Rew}$ be a rewrite
    system, and $\Par{\ColST^C_\Leaf, \RewContext}$ be the closure of
    $\Par{\ColST^C_\Leaf, \Rew}$. If
    \begin{enumerate}[label={(\it\roman*)}]
        \item there exists a poset $\Qca$ and a termination invariant
        $\theta : \ColST^C_\Leaf \to \Qca$ for
        $\Par{\ColST^C_\Leaf, \Rew}$;
        \item the map $\theta$ is compatible with the closure;
    \end{enumerate}
    then, $\Par{\ColST^C_\Leaf, \RewContext}$ is terminating.
\end{Proposition}
\medbreak

Consider, for instance, the rewrite system $\Par{\ColST^C_\Leaf, \Rew}$
defined by~\eqref{equ:example_rewrite_rule_syntax_trees}. By setting
$\Qca := \N^2$ and $\Ord$ as the lexicographic order on $\N^2$, let us
define the map $\theta : \ColST^C_\Leaf \to \Qca$, for any $C$-syntax
tree $\Tfr$, by
\begin{math}
    \theta(\Tfr) := \Par{\deg(\Tfr), \TamariInvariant(\Tfr)}
\end{math},
where
\begin{equation} \label{equ:tamari_invariant}
    \TamariInvariant(\Tfr) :=
    \sum_{\substack{
        u \in \TreeLanguage_\Node(\Tfr) \\
        u \mbox{ \footnotesize of arity } 2
    }}
    \deg(\Tfr \Action u2).
\end{equation}
In other words, $\TamariInvariant(\Tfr)$ is the sum, for all binary
nodes $u$ of $\Tfr$, of the number of internal nodes appearing in the
$2$nd suffix subtrees of $u$. One can check that
$\theta(\Tfr) \OrdStrict \theta(\Tfr')$ for all the $C$-syntax trees
$\Rfr$ and $\Rfr'$ such that $\Rfr \Rew \Rfr'$. Indeed,
\begin{subequations}
\begin{equation}
    \theta\Par{

    }.
\end{equation}
\end{subequations}
Moreover, the fact that $\theta$ is compatible with the closure is a
straightforward verification. Therefore, the closure
$\Par{\ColST^C_\Leaf, \RewContext}$ of $\Par{\ColST^C_\Leaf, \Rew}$ is
terminating.
\medbreak

\subsubsection{Proving confluence}
In the same way as the tool to show that a rewrite system on $C$-syntax
trees is terminating presented in
Section~\ref{subsubsec:proving_termination}, we present here a tool to
prove that rewrite systems on syntax trees defined as closures of other
ones are confluent. This criterion requires now some precise properties.
\medbreak

\begin{Proposition} \label{prop:confluence_rewrite_rule_syntax_trees}
    Let $C$ be a combinatorial augmented graded collection without
    object of size~$1$, $\Par{\ColST^C_\Leaf, \Rew}$ be a rewrite
    system, and $\Par{\ColST^C_\Leaf, \RewContext}$ be the closure of
    $\Par{\ColST^C_\Leaf, \Rew}$. If $\Par{\ColST^C_\Leaf, \RewContext}$
    is
    \begin{enumerate}[label={(\it\roman*)}]
        \item of finite type;
        \item terminating;
        \item such that all its branching pairs consisting in trees with
        $2 \ell - 1$ internal nodes or less are joinable, where $\ell$
        is its degree;
    \end{enumerate}
    then, $\Par{\ColST^C_\Leaf, \RewContext}$ is confluent.
\end{Proposition}
\medbreak

Proposition~\ref{prop:confluence_rewrite_rule_syntax_trees} yields an
algorithmic way to check if a terminating rewrite system
$\Par{\ColST^C_\Leaf, \RewContext}$ defined as the closure of an other
one $\Par{\ColST^C_\Leaf, \Rew}$ is confluent by enumerating all the
$C$-syntax trees $\Tfr$ of degrees at most $2\ell - 1$ (where $\ell$ is
the degree of $\Par{\ColST^C_\Leaf, \RewContext}$) and by computing the
parts $G_\Tfr$ of the rewriting graphs of
$\Par{\ColST^C_\Leaf, \RewContext}$ consisting in the trees reachable
from~$\Tfr$. If each $G_\Tfr$ contains exactly one normal form (which
correspond to a vertex with no outgoing edge in $G_\Tfr$),
$\Par{\ColST^C_\Leaf, \RewContext}$ is confluent.
\medbreak

For instance, by considering the same labeling set $C$ as above, let
$\Par{\ColST^C_\Leaf, \Rew}$ be the rewrite system defined by
\begin{equation}
\,.
\end{equation}
The degree of the closure $\Par{\ColST^C_\Leaf, \RewContext}$ of
$\Par{\ColST^C_\Leaf, \Rew}$ is $\ell := 2$ and it is possible to show
that $\Par{\ColST^C_\Leaf, \RewContext}$ is terminating. Consider
\begin{equation}
    \Tfr :=
    %
\,,
\end{equation}
which is a $C$-syntax tree of degree $2\ell - 1 = 3$. The graph $G_\Tfr$
associated with $\Tfr$ is of the form
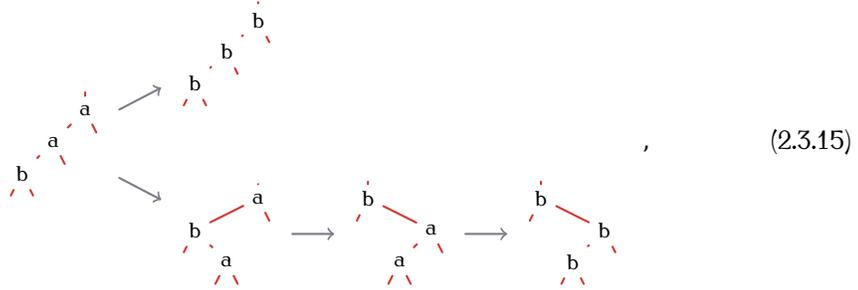
\begin{equation}
    %
\,,
\end{equation}
and shows that $\Par{\ColST^C_\Leaf, \RewContext}$ is not confluent.
Indeed, $\Tfr$ is a non-joinable branching tree. On the other hand,
consider the rewrite system $\Par{\ColST^C_\Leaf, \Rew}$ defined by
\begin{equation}
    %
\,.
\end{equation}
The degree of the closure $\Par{\ColST^C_\Leaf, \RewContext}$ of
$\Par{\ColST^C_\Leaf, \Rew}$ is $\ell := 2$ and here also, it is
possible to show that $\Par{\ColST^C_\Leaf, \RewContext}$ is
terminating. Consider
\begin{equation}
    \Tfr :=
    %
\,,
\end{equation}
a $C$-syntax tree of degree $2\ell - 1 = 3$. The graph $G_\Tfr$
associated with $\Tfr$ is of the form
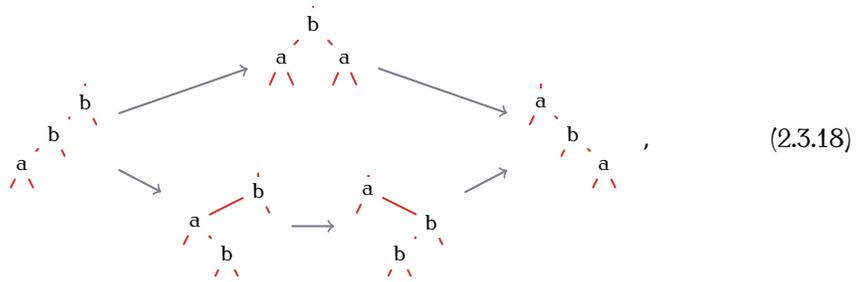
\begin{equation}
    %
        };
        \draw[Arc](2220000_bba)--(2200200_baa);
        \draw[Arc](2220000_bba)--(2202000_bab);
        \draw[Arc](2202000_bab)--(2022000_abb);
        \draw[Arc](2022000_abb)--(2020200_aba);
        \draw[Arc](2200200_baa)--(2020200_aba);
    \end{tikzpicture}\,,
\end{equation}
This graph satisfies the required property stated above, and, as a
systematic study of cases shows, all other graphs $G_\Sfr$ where $\Sfr$
is a $C$-syntax tree of degree $3$ or less, also. For this reason,
$\Par{\ColST^C_\Leaf, \RewContext}$ is confluent.
\medbreak

\section{Treelike structures} \label{sec:treelike}
We expose here two additional variants of trees. The first one are
rooted trees and are structures intervening in the study of free pre-Lie
algebras (see forthcoming Section~\ref{subsec:pre_lie_algebras} of
Chapter~\ref{chap:algebra}). The second ones are colored syntax trees
and can be seen as objects of free colored operads
(see forthcoming Section~\ref{subsec:colored_operads} of
Chapter~\ref{chap:generalizations}).
\medbreak

\subsection{Rooted trees} \label{subsec:rooted_trees}
Let $\ColRT$ be the graded collection satisfying the relation
\begin{equation}
    \ColRT = \BBrack{\{\Node\}, \Multiset^+(\ColRT)}_\times^+.
\end{equation}
where $\Node$ is an atomic object called \Def{node} and $\Multiset$ is
the multiset collection operation over graded collections defined in
Section~\ref{subsubsec:multiset_operation} of
Chapter~\ref{chap:collections}. We call \Def{rooted tree} each object of
$\ColRT$. By definition, a rooted tree $\Tfr$ is an ordered pair
$\Par{\Node, \lbag \Tfr_1, \dots, \Tfr_k \rbag}$ where
$\lbag \Tfr_1, \dots, \Tfr_k \rbag$ is a multiset of rooted trees. Like
the case of planar rooted trees, this definition is  recursive. For
instance,
\begin{equation} \label{equ:examples_rooted_trees}
    (\Node, \emptyset), \quad
    (\Node, \lbag (\Node, \emptyset) \rbag), \quad
    (\Node, \lbag (\Node, \emptyset), (\Node, \emptyset) \rbag), \quad
    (\Node, \lbag (\Node, \emptyset), (\Node, \emptyset),
        (\Node, \emptyset) \rbag), \quad
    (\Node, \lbag (\Node, \lbag (\Node, \emptyset), (\Node, \emptyset)
        \rbag) \rbag),
\end{equation}
are rooted trees. If
$\Tfr = \Par{\Node, \lbag \Tfr_1, \dots, \Tfr_k \rbag}$ is a rooted
tree, each $\Tfr_i$, $i \in [k]$, is a \Def{suffix subtree} of~$\Tfr$.
\medbreak

Rooted trees are different kinds of trees than planar rooted trees
presented in Section~\ref{sec:planar_rooted_trees}. The difference is
due to the fact that rooted trees are defined by using multisets of
rooted trees, while planar rooted trees are defined by using tuples of
planar rooted trees. Hence, the order of the suffix subtrees of a rooted
tree is not significant.
\medbreak

By drawing each rooted tree by a node $\NodePic$ attached below it to
its subtrees by means of edges~$\EdgePic$, the rooted trees
of~\eqref{equ:examples_rooted_trees} are depicted by
\begin{equation}
    \NodePic \quad,
    \begin{tikzpicture}[scale=.3,Centering]
        \node[Node](1)at(0.00,-1.00){};
        \node[Node](0)at(0.00,0.00){};
        \draw[Edge](1)--(0);
        \node(r)at(0.00,1){};
        \draw[Edge](r)--(0);
    \end{tikzpicture}\,, \quad
    \begin{tikzpicture}[scale=.2,Centering]
        \node[Node](0)at(0.00,-1.50){};
        \node[Node](2)at(2.00,-1.50){};
        \node[Node](1)at(1.00,0.00){};
        \draw[Edge](0)--(1);
        \draw[Edge](2)--(1);
        \node(r)at(1.00,1.5){};
        \draw[Edge](r)--(1);
    \end{tikzpicture}\,, \quad
    \begin{tikzpicture}[xscale=.3,yscale=.2,Centering]
        \node[Node](0)at(0.00,-2.00){};
        \node[Node](2)at(1.00,-2.00){};
        \node[Node](3)at(2.00,-2.00){};
        \node[Node](1)at(1.00,0.00){};
        \draw[Edge](0)--(1);
        \draw[Edge](2)--(1);
        \draw[Edge](3)--(1);
        \node(r)at(1.00,1.50){};
        \draw[Edge](r)--(1);
    \end{tikzpicture}\,, \quad
    \begin{tikzpicture}[scale=.2,Centering]
        \node[Node](1)at(0.00,-2.67){};
        \node[Node](3)at(2.00,-2.67){};
        \node[Node](0)at(1.00,0.00){};
        \node[Node](2)at(1.00,-1.33){};
        \draw[Edge](1)--(2);
        \draw[Edge](2)--(0);
        \draw[Edge](3)--(2);
        \node(r)at(1.00,1.5){};
        \draw[Edge](r)--(0);
    \end{tikzpicture}\,.
\end{equation}
By definition of the product and multiset operations over combinatorial
collections, the size of a rooted tree $\Tfr$ satisfies
\begin{equation} \label{equ:size_rooted_tree}
    |\Tfr| := 1 + \sum_{i \in [k]} |\Tfr_i|.
\end{equation}
The integer sequence of $\ColRT$ begins by
\begin{equation}
    1, 1, 2, 4, 9, 20, 48, 115, 286
\end{equation}
and forms Sequence~\OEIS{A000081} of~\cite{Slo}.
\medbreak

\subsection{Colored syntax trees} \label{subsec:colored_syntax_trees}
Let $\CFr$ be a set of colors and $C$ be a $\CFr$-colored collection
(see Section~\ref{subsubsec:colored_collections} of
Chapter~\ref{chap:collections}) such that the graduation of $C$ is
augmented. A \Def{$\CFr$-colored $C$-syntax tree} is a triple
$(a, \Tfr, u)$ where $\Tfr$ is a $C$-syntax tree of arity
$n \in \N_{\geq 1}$, $a \in \CFr$, $u \in \CFr^n$, and for any internal
nodes $u$ and $v$ of $\Tfr$ such that $v$ is the $i$th child of $u$,
$\Out(y) = \In_i(x)$ where $x$ (resp. $y$) is the label of $u$ (resp.
$v$). The set of all $\CFr$-colored $C$-syntax trees is denoted by
$\ColCST^C$. This set is a $\CFr$-colored collection by setting that
$\Out((a, \Tfr, u)) := a$ and $\In((a, \Tfr, u)) := u$ for all
$(a, \Tfr, u) \in \ColCST^C$. By a slight abuse of notation, if $u$ is
an internal node of $\Tfr$, we denote by $\Out(u)$ (resp. $\In(u)$) the
color $\Out(x)$ (resp. word of colors $\In(x)$) where $x$ is the label
of $u$. We say that a $\CFr$-colored $C$-syntax tree $\Tfr$ is
\Def{monochrome} if $C$ is a monochrome colored collection. In graphical
representations of a $\CFr$-colored $C$-syntax tree $(a, \Tfr, u)$, we
draw $\Tfr$ together with its output color above its root and its input
color $u(i)$ below its $i$th leaf for any~$i \in [|u|]$.
\medbreak

For instance, consider the set of colors $\CFr := \{1, 2\}$ and the
$\CFr$-colored collection $C$ defined by
$C := C(2) \sqcup C(3)$ with $C(2) := \{\Asf, \Bsf\}$,
$C(3) := \{\Csf\}$, $\Out(\Asf) := 1$, $\Out(\Bsf) := 2$,
$\Out(\Csf) := 1$, $\In(\Asf) := 11$, $\In(\Bsf) := 21$, and
$\In(\Csf) := 221$. The tree
\begin{equation}
    \begin{tikzpicture}[xscale=.35,yscale=.15,Centering]
        \node(0)at(0.00,-6.50){};
        \node(10)at(8.00,-9.75){};
        \node(12)at(10.00,-9.75){};
        \node(2)at(2.00,-6.50){};
        \node(4)at(3.00,-3.25){};
        \node(5)at(4.00,-9.75){};
        \node(7)at(5.00,-9.75){};
        \node(8)at(6.00,-9.75){};
        \node[NodeST](1)at(1.00,-3.25){\begin{math}\Bsf\end{math}};
        \node[NodeST](11)at(9.00,-6.50){\begin{math}\Asf\end{math}};
        \node[NodeST](3)at(3.00,0.00){\begin{math}\Csf\end{math}};
        \node[NodeST](6)at(5.00,-6.50){\begin{math}\Csf\end{math}};
        \node[NodeST](9)at(7.00,-3.25){\begin{math}\Asf\end{math}};
        \node(r)at(3.00,2.75){};
        \draw[Edge](0)--(1);
        \draw[Edge](1)--(3);
        \draw[Edge](10)--(11);
        \draw[Edge](11)--(9);
        \draw[Edge](12)--(11);
        \draw[Edge](2)--(1);
        \draw[Edge](4)--(3);
        \draw[Edge](5)--(6);
        \draw[Edge](6)--(9);
        \draw[Edge](7)--(6);
        \draw[Edge](8)--(6);
        \draw[Edge](9)--(3);
        \draw[Edge](r)--(3);
        \node[LeafLabel,above of=r]{\begin{math}1\end{math}};
        \node[LeafLabel,below of=0]{\begin{math}2\end{math}};
        \node[LeafLabel,below of=2]{\begin{math}1\end{math}};
        \node[LeafLabel,below of=4]{\begin{math}2\end{math}};
        \node[LeafLabel,below of=5]{\begin{math}2\end{math}};
        \node[LeafLabel,below of=7]{\begin{math}2\end{math}};
        \node[LeafLabel,below of=8]{\begin{math}1\end{math}};
        \node[LeafLabel,below of=10]{\begin{math}1\end{math}};
        \node[LeafLabel,below of=12]{\begin{math}1\end{math}};
    \end{tikzpicture}
\end{equation}
is a $\CFr$-colored $C$-syntax tree. Its output color is $1$ and its
word of input colors is $21222111$. Besides, $(1, \Leaf, 1)$ and
$(1, \Leaf, 2)$ are two $\CFr$-colored $C$-syntax trees of degree $0$
and arity~$1$.
\medbreak

The partial grafting operation of syntax trees (see
Section~\ref{subsubsec:partial_grafting_syntax_trees}) admits a
generalization on colored syntax trees. Let for any
$(a, u) \in \CFr \times \CFr^n$, $(b, v) \in \CFr \times \CFr^m$, and
$i \in[n]$ such that $b = u(i)$ the product
\begin{equation}
    \circ_i^{((a, u), (b, v))} :
    \ColCST^C(a, u) \times \ColCST^C(b, v) \to \ColCST^C
\end{equation}
defined, for any $(a, \Tfr, u) \in \ColCST^C(a, u)$ and
$(b, \Sfr, v) \in \ColCST^C(b, v)$ by
\begin{equation} \label{equ:grafting_colored_syntax_trees}
    (a, \Tfr, u) \circ_i^{((a, u), (b, v))} (b, \Sfr, v) :=
    \Par{a, \Tfr \circ_i^{(n, m)} \Sfr, u \mapsfrom_i v},
\end{equation}
where $u \mapsfrom_i v$ is the word obtained by replacing the $i$th
letter of $u$ by $v$, and $\circ_i^{(n, m)}$ is the partial grafting of
syntax trees. For instance, by considering the same labeling
$\CFr$-colored collection as above,
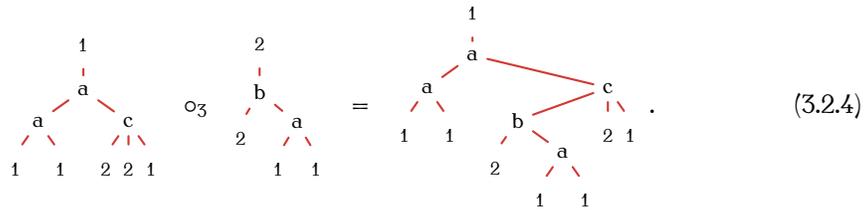
\begin{equation}
\,.
\end{equation}
We call each $\circ_i^{(n, m)}$ a \Def{partial grafting operation}.
\medbreak

Observe that since
\begin{equation}
    \Index\Par{(a, \Tfr, u) \circ_i^{(n, m)} (b, \Sfr, v)}
    = (a, u \mapsfrom_i v)
\end{equation}
where $\Index$ denotes the index of an object of a colored collection
(see Section~\ref{subsubsec:collections_definition} of
Chapter~\ref{chap:collections}), the product $\circ_i^{(n, m)}$ is
concentrated. Besides, by a slight abuse of notation, we shall sometimes
omit the $(n, m)$ in the notation of $\circ_i^{(n, m)}$ in order to
denote it in a more concise way by~$\circ_i$.
\medbreak

The generalizations of the complete and context grafting products of
syntax trees (see
Sections~\ref{subsubsec:complete_grafting_syntax_trees}
and~\ref{subsubsec:context_grafting_syntax_trees}) on colored syntax
trees follow from the definition of the partial grafting operation of
colored syntax trees just given. These two products are also
concentrated.
\medbreak

Most of the notions exposed in
Section~\ref{subsec:rewrite_rules_syntax_trees} about syntax trees and
rewrite systems on syntax trees naturally extend on colored syntax trees
like, among others, the notions of occurrences of patterns, the complete
grafting operations, and the criteria offered by
Propositions~\ref{prop:compatible_terminating_invariant}
and~\ref{prop:confluence_rewrite_rule_syntax_trees} to, respectively,
prove the termination and the confluence of rewrite system on syntax
trees.
\medbreak

\SkipTocEntry\section*{Bibliographic notes}

\SkipTocEntry\subsection*{About trees}
The concept of tree encompasses a large range of quite different
combinatorial objects. For instance, in graph theory, trees are
connected acyclic graphs while in combinatorics, one encounters mostly
rooted trees. Among rooted trees, some of these can be planar (the order
of the children of a node is relevant) or not. In addition to this, the
internal nodes, the leaves, or the edges of the trees can be labeled,
and some conditions for the arities of their nodes can be imposed. One
of the first occurrences of the concept of tree came from the work of
Cayley~\cite{Cay1857}. Nowadays, trees appear among others in computer
science as data structures~\cite{Knu98,CLRS09}, in combinatorics in
relation with enumerating questions and Lagrange
inversion~\cite{Lab81,FS09}, and in algebraic combinatorics, where
several families of trees are endowed with algebraic
structures~\cite{LR98,HNT05,Cha08}. Besides, the bijection between the
combinatorial collections of the planar rooted trees and the one of
binary trees appearing in
Proposition~\ref{prop:bijection_planar_rooted_trees_binary_trees} is
known as the rotation correspondence and is due to Knuth~\cite{Knu97}.
This bijection, offering a means of encoding a planar rooted tree by a
binary tree, admits applications in algebraic
combinatorics~\cite{NT13,EM14}.
\medbreak

\SkipTocEntry\subsection*{About enumerating properties}
Formula~\eqref{equ:number_k-ary_trees} for the Fuss-Catalan numbers,
enumerating the combinatorial collection of the $k$-ary trees with
respect to their number of internal nodes has been established
in~\cite{DM47}. Besides, Formula~\eqref{equ:enumeration_schroder_trees}
enumerating the combinatorial set of the Schröder trees with respect to
their number of leaves uses the Narayana numbers~\cite{Nar55}. These
numbers admit the following combinatorial interpretation: the $2$-graded
collection $C$ of binary trees, where the index of a binary tree $\Tfr$
is the pair $(n, k)$ where $n$ is the arity of $\Tfr$ and $k$ is the
number of internal nodes of $\Tfr$ having an internal node as a left
child satisfies~$\# C(n, k) = \Narayana(n, k)$.
\medbreak

\SkipTocEntry\subsection*{About rewrite rules on trees}
The Buchberger algorithm, which is a completion algorithm (see the end
of Chapter~\ref{chap:collections}), admits adaptations in the context of
rewrite systems of trees and operads~\cite{DK10,BD16}. Several works use
rewrite systems on trees to provide presentations of operads (see, for
instance,~\cite{Hof10,LV12,Gir16c,CCG18}).
\medbreak


\chapter{Algebraic structures} \label{chap:algebra}
This chapter deals with vector spaces obtained from graded collections.
A general framework for algebraic structures having products and
coproducts is presented. Most of the algebraic structures encountered
in algebraic combinatorics like associative, dendriform, pre-Lie
algebras, and Hopf bialgebras fit into this framework. This chapter
contains classical examples of such structures.
\medbreak

\section{Polynomials spaces} \label{sec:polynomials}
We introduce here the notion of polynomial spaces. All the algebraic
structures considered in this book are polynomial spaces endowed with
some operations or co-operations. A set of operations, analogous to
the operations on graded collections of
Section~\ref{subsec:operations_collections} of
Chapter~\ref{chap:collections}, over graded polynomial spaces are
considered. We also review some links between changes of bases of
polynomial spaces, posets, and incidence algebras.
\medbreak

\subsection{Polynomials on collections}
Intuitively, a polynomial on an $I$-collection $C$ is a finite formal
sum of objects of $C$ with coefficients in a field $\K$. In what
follows, $\K$ can be any field of characteristic~$0$.
\medbreak

\subsubsection{Polynomials} \label{subsubsec:c_polynomials}
Let $C$ be an $I$-collection. A \Def{polynomial on $C$} (or, for short,
a \Def{$C$-polynomial}) is a map
\begin{equation}
    f : C \to \K
\end{equation}
such that the set
\begin{equation} \label{equ:support_polynomial}
    \Support(f) := \left\{x \in C : f(x) \ne 0\right\}
\end{equation}
is finite, where the symbol $0$ appearing
in~\eqref{equ:support_polynomial} is the zero of $\K$. We call
$\Support(f)$ the \Def{support} of $f$. The \Def{coefficient} $f(x)$ of
$x \in C$ in $f$ is denoted by $\Angle{x, f}$. An object $x$ of $C$
\Def{appears} in $f$ if $\Angle{x, f} \ne 0$. A $C$-polynomial $f$ is a
\Def{$C$-monomial} if $\Support(f)$ is a singleton. We say that $f$ is
\Def{homogeneous} if there is an index $i \in I$ such that
$\Support(f) \subseteq C(i)$. For any finite subcollection $X$ of $C$,
the \Def{characteristic polynomial} of $X$ is the $C$-polynomial
$\Charac(X)$ defined, for any $x \in C$, by
\begin{equation}
    \Angle{x, \Charac(X)} :=
    \begin{cases}
        1 \in \K & \mbox{if } x \in X, \\
        0 \in \K & \mbox{otherwise}.
    \end{cases}
\end{equation}
Given two $C$-polynomials $f_1$ and $f_2$, the \Def{scalar product} of
$f_1$ and $f_2$ is the scalar
\begin{equation} \label{equ:scalar_product_c_polynomials}
    \Angle{f_1, f_2}
    := \sum_{x \in C} \Angle{x, f_1} \Angle{x, f_2}
\end{equation}
of $\K$. This notation for the scalar product of $C$-polynomials is
consistent with the notation $\Angle{x, f}$ for the coefficient of $x$
in $f$ because by~\eqref{equ:scalar_product_c_polynomials}, the
coefficient $\Angle{x, f}$ and the scalar product
$\Angle{\Charac(\{x\}), f}$ are equal.
\medbreak

In the particular case where $C$ is a graded collection, the
\Def{degree} $\deg(f)$ of $f$ is undefined if $\Support(f) = \emptyset$
and is otherwise the greatest size of an object appearing
in~$\Support(f)$.
\medbreak

\subsubsection{Polynomial spaces} \label{subsubsec:polynomial_spaces}
The set of all $C$-polynomials is denoted by $\K \Angle{C}$. The
\Def{underlying collection} of $\K \Angle{C}$ is $C$. For any property
$P$ of collections (see Section~\ref{sec:collections} of
Chapter~\ref{chap:collections}), we say by extension that $\K \Angle{C}$
\Def{satisfies the property $P$} if $C$ satisfies $P$.
\medbreak

This set $\K \Angle{C}$ is endowed with the following two operations.
First, the \Def{addition}
\begin{equation}
    + : \K \Angle{C} \times \K \Angle{C} \to \K \Angle{C}
\end{equation}
is defined, for any $f_1, f_2 \in \K \Angle{C}$ and $x \in C$, by
\begin{equation}
    \Angle{x, f_1 + f_2} := \Angle{x, f_1} + \Angle{x, f_2}.
\end{equation}
Second, the \Def{scalar multiplication}
\begin{equation}
    \ExtProd : \K \times \K \Angle{C} \to \K \Angle{C}
\end{equation}
is defined, for any $f \in \K \Angle{C}$, $\lambda \in \K$, and
$x \in C$, by
\begin{equation}
    \Angle{x, \lambda \ExtProd f} := \lambda \Angle{x, f}.
\end{equation}
Endowed with these two operations, $\K \Angle{C}$ is a $\K$-vector
space, named \Def{polynomial space on $C$} (or, for short,
\Def{$C$-polynomial space}). Moreover, $\K \Angle{C}$ decomposes as a
direct sum
\begin{equation} \label{equ:decomposition_polynomial_spaces}
    \K \Angle{C} = \bigoplus_{i \in I} \K \Angle{C(i)}.
\end{equation}
We call each $\K \Angle{C(i)}$ the \Def{$i$-homogeneous component} of
$\K \Angle{C}$. In the sequel, we shall also write $\K \Angle{C}(i)$ for
$\K \Angle{C(i)}$.
\medbreak

By using now the linear structure of $\K \Angle{C}$, any $C$-polynomial
$f$ can be expressed as the finite sum of $C$-monomials
\begin{equation}
    f = \sum_{x \in C} \Angle{x, f} \ExtProd \Charac(\{x\}),
\end{equation}
which is denoted, by a slight abuse of notation, by
\begin{equation} \label{equ:sum_notation_polynomials}
    f = \sum_{x \in C} \Angle{x, f} x.
\end{equation}
The notation~\eqref{equ:sum_notation_polynomials} for $f$ as a linear
combination of objects of $C$ is the \Def{sum notation} of
$C$-polynomials.
\medbreak

Since for any $C$-polynomial $f$, there is unique way to express $f$ as
a finite sum of the form~\eqref{equ:sum_notation_polynomials}, the set
\begin{equation}
    \left\{\Charac(\{x\}) : x \in C\right\}
\end{equation}
forms a basis of $\K \Angle{C}$. This basis is called
\Def{fundamental basis} of $\K \Angle{C}$, and, by a slight but
convenient abuse of notation, each basis element $\Charac(\{x\})$,
$x \in C$, is simply denoted by~$x$. Observe that each basis of
$\K \Angle{C}$ is indexed by $C$. Moreover, Let us emphasize the fact
any polynomial space $\K \Angle{C}$ is always seen through its explicit
basis $C$ (contrarily when working abstractly with a vector space $\Vca$
without explicit basis). In the sequel, we shall define products on $C$
which extend by linearity on $\K \Angle{C}$. Properties of such products
(like associativity or commutativity) can be defined and checked only
on~$C$.
\medbreak

Besides, we are sometimes led to consider several bases of
$\K \Angle{C}$ and work with many of them at the same time. In this
case, to distinguish elements expressed on different bases, we denote
them by putting elements of $C$ as indexes of a letter naming the basis.
For instance, the elements of the $\BasisB$-basis of $\K \Angle{C}$ are
denoted by $\BasisB_x$, $x \in C$.
\medbreak

Let $C_1$ and $C_2$ be two $I$-collections. A \Def{morphism} between
$\K \Angle{C_1}$ and $\K \Angle{C_2}$ is a linear map
\begin{equation}
    \phi : \K \Angle{C_1} \to \K \Angle{C_2}
\end{equation}
such that for any $x \in C_1$, $\phi(x) \in \K \Angle{C_2}(\Index(x))$.
Observe that any combinatorial collection morphism $\psi : C_1 \to C_2$
gives rise to a polynomial space morphism
\begin{math}
    \bar{\psi} : \K \Angle{C_1} \to \K \Angle{C_2}
\end{math}
obtained by extending $\psi$ linearly. Besides, $\K \Angle{C_2}$ is a
\Def{subspace} of $\K \Angle{C_1}$ if there exists an injective morphism
from $\K \Angle{C_2}$ to $\K \Angle{C_1}$. For any subset $J$ of $I$, we
denote by $\K \Angle{C}(J)$ the polynomial space $\K \Angle{C(J)}$.
Since $C(J)$ is by definition a subcollection of $C$, $\K \Angle{C}(J)$
is a subspace of~$\K \Angle{C}$.
\medbreak

\subsubsection{Combinatorial graded polynomial spaces}
When $C$ is a combinatorial graded collection, as a particular case
of~\eqref{equ:decomposition_polynomial_spaces}, $\K \Angle{C}$
decomposes as a direct sum
\begin{equation}
    \K \Angle{C} = \bigoplus_{n \in \N} \K \Angle{C}(n).
\end{equation}
Moreover, since $C$ is combinatorial, each $\K \Angle{C(n)}$,
$n \in \N$, is finite dimensional. For this reason, the
\Def{Hilbert series} of $\K \Angle{C}$, defined by
\begin{equation}
    \HilbertSeries_{\K \Angle{C}}(t)
    = \sum_{n \in \N} \dim \K \Angle{C}(n) \, t^n,
\end{equation}
is a well-defined series. We can observe that the Hilbert series
$\HilbertSeries_{\K \Angle{C}}(t)$ of $\K \Angle{C}$ and the generating
series $\GeneratingSeries_C(t)$ of $C$ are the same power series.
\medbreak

\subsubsection{Rewrite systems and quotient spaces}
\label{subsubsec:rew_quotient_space}
For any $I$-collection $C$, any rewrite system $(C, \Rew)$ gives rise
to a subspace $\RelationSpaceRewriteRule_{(C, \Rew)}$ of $\K \Angle{C}$
generated by all the homogeneous $C$-polynomials $x' - x$ whenever $x$
and $x'$ are two objects of $C$ such that $x \Rew x'$. We call
$\RelationSpaceRewriteRule_{(C, \Rew)}$ the \Def{space induced} by
$(C, \Rew)$. Conversely, when $\RelationSpace$ is a subspace of
$\K \Angle{C}$ such that there exists a rewrite system $(C, \Rew)$ such
that $\RelationSpace$ and $\RelationSpaceRewriteRule_{(C, \Rew)}$ are
isomorphic, we say that $(C, \Rew)$ is an \Def{orientation} of
$\RelationSpace$. When $(C, \Rew)$ is convergent, one has a concrete
description of the quotient space
$\K \Angle{C}/_{\RelationSpaceRewriteRule_{(C, \Rew)}}$
involving the normal forms $\NormalForms_{(C, \Rew)}$ of $(C, \Rew)$
provided by the following result.
\medbreak

\begin{Proposition}
\label{prop:space_induced_space_quotient_rewrite_rule}
    Let $(C, \Rew)$ be a convergent rewrite system. Then, as spaces
    \begin{equation}
        \K \Angle{C} /_{\RelationSpaceRewriteRule_{(C, \Rew)}}
        \simeq
        \K \Angle{\NormalForms_{(C, \Rew)}}.
    \end{equation}
\end{Proposition}
\medbreak

\subsection{Operations over polynomial spaces}
\label{subsec:operations_polynomial_spaces}
In the same way as operations over collections allow to create new
collections from already existing ones (see
Section~\ref{subsec:operations_collections} of
Chapter~\ref{chap:collections}), there exist analogous operations over
polynomial spaces. Some of these are consequences of the definitions of
operations over collections. We present here the main ones.
Alternatively, one of the aims of this section is to show that the usual
operations over spaces (direct sum, quotient, and tensor product)
produce polynomial spaces.
\medbreak

\subsubsection{Direct sum}
The sum of two collections translates as the direct sum of the
associated polynomial spaces. Indeed, for any $I$-collections $C_1$
and $C_2$,
\begin{equation} \label{equ:direct_sum_polynomial_spaces}
    \K \Angle{C_1 + C_2} \simeq \K \Angle{C_1} \oplus \K \Angle{C_2}.
\end{equation}
An isomorphism between the two spaces
of~\eqref{equ:direct_sum_polynomial_spaces} is provided by the map
\begin{equation}
    \phi : \K \Angle{C_1 + C_2}
    \to \K \Angle{C_1} \oplus \K \Angle{C_2},
\end{equation}
linearly defined for any $x \in C_1 + C_2$ by
\begin{equation}
    \phi(x) :=
    \begin{cases}
        x \in \K \Angle{C_1} & \mbox{if } x \in C_1, \\
        x \in \K \Angle{C_2} & \mbox{otherwise (} x \in C_2 \mbox{)}.
    \end{cases}
\end{equation}
For this reason, we shall identify the two
spaces of~\eqref{equ:direct_sum_polynomial_spaces}.
\medbreak

\subsubsection{Quotient space}
If $C$ is an $I$-collection and $\Vca$ is a space included in
$\K \Angle{C}$, the \Def{quotient space} of $\K \Angle{C}$ by $\Vca$ is
the space $\K \Angle{C}/_\Vca$ of all the equivalence classes
\begin{equation}
    [f]_\Vca := \left\{f + g : g \in \Vca\right\},
\end{equation}
for all $f \in \K \Angle{C}$, endowed with its natural vector space
structure. We call \Def{canonical surjection map} the linear map
$\theta : \K \Angle{C} \to \K \Angle{C}/_\Vca$ defined linearly by
$\theta(x) := [x]_\Vca$ for any object $x$ of $C$. It is always possible
to see $\K \Angle{C}/_\Vca$ as a $D$-polynomial space by providing an
adequate $I$-collection $D$ so that
\begin{equation}
    \K \Angle{C}/_\Vca \simeq \K \Angle{D}.
\end{equation}
For this reason, we shall identify any quotient space with a polynomial
space.
\medbreak

\subsubsection{Tensor product}
The Cartesian product of collections translates as the tensor product of
the associated polynomial spaces. Indeed, for any $p \in \N$, any index
sets $I_1$, \dots, $I_p$, and any $I_1$-collection $C_1$, \dots, any
$I_p$-collection $C_p$,
\begin{equation} \label{equ:tensor_product_polynomial_spaces}
    \K \Angle{\BBrack{C_1, \dots, C_p}_\times}
    \simeq \K \Angle{C_1} \otimes \dots \otimes \K \Angle{C_p}.
\end{equation}
An isomorphism between the two spaces
of~\eqref{equ:tensor_product_polynomial_spaces} is provided by the map
\begin{equation}
    \phi : \K \Angle{\BBrack{C_1, \dots, C_p}_\times}
    \to \K \Angle{C_1} \otimes \dots \otimes \K \Angle{C_p},
\end{equation}
linearly defined for any
\begin{math}
    \Par{x_1, \dots, x_p} \in \BBrack{C_1, \dots, C_p}_\times
\end{math}
by
\begin{equation}
    \phi\Par{\Par{x_1, \dots, x_p}} := x_1 \otimes \dots \otimes x_p.
\end{equation}
For this reason, we shall identify the two spaces
of~\eqref{equ:tensor_product_polynomial_spaces}. Moreover, as a
consequence, the tuple notation for tensors is linear. That is, for any
$f_1 \in \K \Angle{C_1}$, \dots, $f_p \in \K \Angle{C_p}$,
\begin{equation}
    \Par{f_1, \dots, f_p} =
    \sum_{\Par{x_1, \dots, x_p} \in \BBrack{C_1, \dots, C_p}_\times}
    \, \left(\prod_{k \in [p]} \Angle{x_k, f_k}\right)
    \, \Par{x_1, \dots, x_p}.
\end{equation}
\medbreak

\subsubsection{Tensor algebras} \label{subsubsec:tensor_algebras}
If $\Vca$ is a $\K$-vector space, the \Def{tensor algebra} of $\Vca$ is
the space $\Tensor \Vca$ defined by
\begin{equation}
    \Tensor \Vca := \bigoplus_{p \in \N} \Vca^{\otimes p}
\end{equation}
where $\Vca^{\otimes p}$, $p \in \N$, denotes the space of all
tensors on $\Vca$ of order $p \in \N$. A basis of $\Tensor \Vca$ is
formed by all tensors on any basis of $\Vca$. As a special case of the
one of tensor products discussed in the above section, the list
collection operation applied to a graded collection translates as the
tensor algebra of the associated graded polynomial space. Indeed, for
any~$p \in \N$ and $I$-collection $C$,
\begin{equation} \label{equ:fixed_tensor_polynomial_spaces}
    \K \Angle{\List_{\{p\}}(C)} \simeq {\K \Angle{C}}^{\otimes p}.
\end{equation}
so that, by using the correspondence between direct sums of spaces
and sums of collections, we obtain
\begin{equation} \label{equ:tensor_algebra_polynomial_spaces}
    \K \Angle{\List(C)} \simeq \Tensor \K \Angle{C}.
\end{equation}
An isomorphism between the two spaces
of~\eqref{equ:tensor_algebra_polynomial_spaces} is provided by the map
\begin{equation}
    \phi : \K \Angle{\List(C)} \to \Tensor \K \Angle{C},
\end{equation}
linearly defined for any $\Par{x_1, \dots, x_p} \in \List(C)$ by
\begin{equation}
    \phi\Par{\Par{x_1, \dots, x_p}} :=
    x_1 \otimes \dots \otimes x_p.
\end{equation}
For this reason, we shall identify the two spaces
of~\eqref{equ:tensor_algebra_polynomial_spaces}.
\medbreak

\subsubsection{Symmetric algebras} \label{subsubsec:symmetric_algebras}
If $\Vca$ is a $\K$-vector space, the \Def{symmetric algebra} of $\Vca$
is the space $\Symmetric \Vca$ defined by
\begin{equation}
    \Symmetric \Vca :=
    \Tensor \Vca/_{\Vca_\Symmetric},
\end{equation}
where $\Vca_\Symmetric$ is the subspace of $\Tensor \K \Angle{C}$
consisting in all the tensors
\begin{equation}
    u \otimes x_1 \otimes x_2 \otimes v
    -
    u \otimes x_2 \otimes x_1 \otimes v,
\end{equation}
where $u, v \in \Tensor \Vca$ and $x_1, x_2 \in \Vca$. A basis of
$\Symmetric \Vca$ is formed by all monomials on any basis of $\Vca$.
The multiset collection operation applied to a graded collection
translates as the symmetric algebra of the associated graded polynomial
space. Indeed, for any $I$-collection $C$,
\begin{equation} \label{equ:symmetric_algebra_polynomial_spaces}
    \K \Angle{\Multiset(C)} \simeq \Symmetric \K \Angle{C}.
\end{equation}
An isomorphism between the two spaces
of~\eqref{equ:symmetric_algebra_polynomial_spaces} is provided by the
map
\begin{equation}
    \phi : \K \Angle{\Multiset(C)} \to \Symmetric \K \Angle{C},
\end{equation}
linearly defined for any $\lbag x_1, \dots, x_p \rbag \in \Multiset(C)$
by
\begin{equation}
    \phi\Par{\lbag x_1, \dots, x_p \rbag}
    := y_1^{\alpha_1} \dots y_\ell^{\alpha_\ell},
\end{equation}
where $\ell$ is the number of distinct elements of
$\lbag x_1, \dots, x_p \rbag$ and each $\alpha_i$, $i \in [\ell]$,
denotes the multiplicity of $y_i$ in $\lbag x_1, \dots, x_p \rbag$. For
this reason, we shall identify the two spaces
of~\eqref{equ:symmetric_algebra_polynomial_spaces}.
\medbreak

\subsubsection{Exterior algebras} \label{subsubsec:exterior_algebras}
If $\Vca$ is a $\K$-vector space, the \Def{exterior algebra} of $\Vca$
is the space $\Exterior \Vca$ defined by
\begin{equation}
    \Exterior \Vca := \Tensor \Vca/_{\Vca_\Exterior},
\end{equation}
where $\Vca_\Exterior$ is the subspace of $\Tensor \Vca$ consisting in
all the tensors
\begin{equation}
    u \otimes x_1 \otimes x_2 \otimes v
    +
    u \otimes x_2 \otimes x_1 \otimes v,
\end{equation}
where $u, v \in \Tensor \Vca$ and $x_1, x_2 \in \Vca$. A basis of
$\Exterior \Vca$ is formed by all monomials on a basis of $\Vca$ without
repeated letters. The set collection operation applied to a graded
collection translates as the exterior algebra of the associated graded
polynomial space. Indeed, for any $I$-collection $C$,
\begin{equation} \label{equ:exterior_algebra_polynomial_spaces}
    \K \Angle{\Set(C)} \simeq \Exterior \K \Angle{C}.
\end{equation}
An isomorphism between the two spaces
of~\eqref{equ:exterior_algebra_polynomial_spaces} is provided by the map
\begin{equation}
    \phi : \K \Angle{\Set(C)} \to \Exterior \K \Angle{C},
\end{equation}
linearly defined for any $\left\{x_1, \dots, x_p\right\} \in \Set(C)$ by
\begin{equation}
    \phi\Par{\left\{x_1, \dots, x_p\right\}} := x_1 \dots x_p.
\end{equation}
For this reason, we shall identify the two spaces
of~\eqref{equ:exterior_algebra_polynomial_spaces}.
\medbreak

\subsubsection{Duality for combinatorial polynomial spaces}
\label{subsubsec:duality_polynomial_spaces}
Assume in this section that $C$ is combinatorial. The \Def{dual} of
$\K \Angle{C}$ is the $\K$-vector space $\K \Angle{C}^\Dual$ defined by
\begin{equation}
    \K \Angle{C}^\Dual :=
    \bigoplus_{i \in I} \K \Angle{C}(i)^\Dual,
\end{equation}
where for any $i \in I$, $\K \Angle{C}(i)^\Dual$ is the dual space of
$\K \Angle{C}(i)$. Since $C$ is combinatorial, all the
$\K \Angle{C}(i)$ are finite dimensional spaces, so that
$\K \Angle{C}(i)^\Dual \simeq \K \Angle{C}(i)$, and thus,
\begin{equation}
    \K \Angle{C}^\Dual \simeq \K \Angle{C}.
\end{equation}
For this reason, we shall identify $\K \Angle{C}$ and
$\K \Angle{C}^\Dual$ in this book once $C$ is combinatorial.
\medbreak

The \Def{duality bracket} between $\K \Angle{C}$ and
$\K \Angle{C}^\Dual$ is the linear map
\begin{equation} \label{equ:duality_bracket}
    \Angle{-} :
    \K \Angle{C} \otimes \K \Angle{C}^\Dual \simeq
    \K \Angle{\BBrack{C, C}_\times} \to \K
\end{equation}
defined linearly, for all $\Par{x, x'} \in \BBrack{C, C}_\times$, by
\begin{equation}
    \Angle{\Par{x, x'}} :=
    \begin{cases}
        1 & \mbox{if } x = x', \\
        0 & \mbox{otherwise}.
    \end{cases}
\end{equation}
To not overload the notation, we write $\Angle{x, x'}$ instead of
$\Angle{\Par{x, x'}}$. Observe that for any $f_1 \in \K \Angle{C}$ and
$f_2 \in \K \Angle{C}^\Dual$, $\Angle{f_1, f_2}$ is equal to the scalar
product~\eqref{equ:scalar_product_c_polynomials} of $f_1$ and $f_2$.
Moreover, the duality bracket extends for any $p \in \N$ on
\begin{math}
    \K \Angle{\BBrack{\List_{\{p\}}(C), \List_{\{p\}}(C)}_\times}
\end{math}
linearly by
\begin{equation}
    \Angle{\Par{x_1, \dots, x_p}, \Par{x'_1, \dots, x'_p}}
    :=
    \prod_{k \in [p]} \Angle{x_k, x'_k}
\end{equation}
for any objects $\Par{x_1, \dots, x_p}$ and $\Par{x'_1, \dots, x'_p}$
of $\List_{\{p\}}(C)$.
\medbreak

\subsection{Changes of basis and posets}
\label{subsec:change_basis_posets}
It is very usual, given a polynomial space $\K \Angle{C}$, to consider
a poset structure on $C$ to define new bases of $\K \Angle{C}$. Indeed,
such new bases are defined by considering sums of elements greater (or
smaller) than other ones. In this context, incidence algebras of posets
and their Möbius functions play an important role. We expose here these
concepts.
\medbreak

\subsubsection{Incidence algebras}
Let $(\Qca, \Ord)$ be a locally finite $I$-poset. The
\Def{incidence algebra} of $(\Qca, \Ord)$ is the polynomial space
$\K \Angle{\Comparable(\Qca)}$ ($\Comparable(\Qca)$ is defined in
Section~\ref{subsubsec:operations_posets} of
Chapter~\ref{chap:collections}) endowed with the linear binary product
$\Product$ (the notion of products in polynomial spaces is presented in
the following Section~\ref{sec:bialgebras} but here, only elementary
notions about these are needed) defined, for any objects $(x, y)$ and
$(x', y')$ of $\Comparable(\Qca)$ by
\begin{equation}
    (x, y) \Product \Par{x', y'} :=
    \begin{cases}
        \Par{x, y'} & \mbox{if } y = x', \\
        0 & \mbox{otherwise}.
    \end{cases}
\end{equation}
This product is obviously associative. Moreover, for each $i \in I$, on
the $i$-homogeneous component of $\K \Angle{\Comparable(\Qca)}$, the
$\Comparable(\Qca)$-polynomial
\begin{equation}
    \Unit_i := \sum_{x \in C(i)} (x, x)
\end{equation}
plays the role of a unit, that is,
\begin{math}
    f \Product \Unit_i = f = \Unit_i \Product f
\end{math}
for all $f \in \K \Angle{\Comparable(C)}(i)$. Let for any $i \in I$ the
$\Comparable(\Qca)$-polynomial $\zeta_i$, called
\Def{$i$-zeta polynomial} of $(\Qca, \Ord)$, defined by
\begin{equation}
    \zeta_i :=
    \sum_{\substack{
        x, y \in C(i) \\
        x \Ord y
    }}
    (x, y).
\end{equation}
This $\Comparable(\Qca)$-polynomial encodes some properties of the order
$\Ord$. For instance, the coefficient in $\zeta_i \Product \zeta_i$ of
each $(x, y) \in \Comparable(\Qca(i))$ is the cardinality of the
interval $[x, y]$ in $(\Qca, \Ord)$. The \Def{$i$-Möbius polynomial} of
$(\Qca, \Ord)$ is the $\Comparable(\Qca)$-polynomial $\mu_i$ satisfying
\begin{equation}
    \mu_i \Product \zeta_i
    = \Unit_i
    = \zeta_i \Product \mu_i.
\end{equation}
In other words, $\mu_i$ is the inverse of $\zeta_i$ with respect to the
product~$\Product$. Recall that, as exposed in
Section~\ref{subsubsec:c_polynomials}, polynomials on collections are
functions associating a coefficient with any object. For this reason,
$\zeta_i$ and $\mu_i$ are functions associating a coefficient with any
pair of comparable objects of~$\Qca$.
\medbreak

\begin{Theorem} \label{thm:mobius_polynomial_coefficients}
    Let $(\Qca, \Ord)$ be a locally finite $I$-poset. Then, the
    $i$-Möbius polynomial~$\mu_i$, $i \in I$, of $(\Qca, \Ord)$ is a
    well-defined element of $\K \Angle{\Comparable(\Qca)}$ and its
    coefficients satisfy $\Angle{(x, x), \mu_i} = 1$ for all
    $x \in \Qca(i)$, and
    \begin{equation} \label{equ:mobius_polynomial_coefficients}
        \Angle{(x, z), \mu_i}
        =
        -
        \sum_{\substack{
            y \in \Qca(i) \\
            x \Ord y \OrdStrict z
        }}
        \Angle{(x, y), \mu_i}
    \end{equation}
    for all $x, z \in C(i)$ such that $x \ne z$.
\end{Theorem}
\medbreak

Theorem~\ref{thm:mobius_polynomial_coefficients} provides a recursive
way to compute the coefficients of $\mu_i$, $i \in I$, as a consequence
of the finiteness of each interval of~$\Qca(i)$.
\medbreak

\subsubsection{Changes of basis}
\label{subsubsec:polynomial_spaces_change_basis}
Let $C$ be a combinatorial $I$-collection and $\Ord$ be a partial order
relation on $C$ such that $(C, \Ord)$ is an $I$-poset. Consider the
family
\begin{equation} \label{equ:family_basis_change_order}
    \left\{\BasisB^{\Ord}_x, x \in C\right\}
\end{equation}
of elements of $\K \Angle{C}$ defined, from the fundamental basis of
$\K \Angle{C}$, by
\begin{equation} \label{equ:partial_order_basis_definition}
    \BasisB^{\Ord}_x :=
    \sum_{\substack{
        y \in C \\
        x \Ord y
    }}
    y.
\end{equation}
Observe that since $C$ is combinatorial and $\Ord$ preserves the indexes
of the objects of $C$, each $\BasisB^{\Ord}_x$ is a homogeneous
$C$-polynomial. We call the family~\eqref{equ:family_basis_change_order}
the \Def{$\BasisB^{\Ord}$-family} of~$\K \Angle{C}$.
\medbreak

\begin{Proposition} \label{prop:partial_order_bases}
    Let $(C, \Ord)$ be a combinatorial $I$-poset. The
    $\BasisB^{\Ord}$-family forms a basis of $\K \Angle{C}$ and
    \begin{equation} \label{equ:partial_order_bases}
        x =
        \sum_{\substack{
            y \in C \\
            x \Ord y
        }}
        \Angle{(x, y), \mu_i}
        \BasisB^{\Ord}_y
    \end{equation}
    for all $x \in C(i)$, $i \in I$, where $\mu_i$ is the $i$-Möbius
    polynomial of~$(C, \Ord)$.
\end{Proposition}
\medbreak

\section{Bialgebras} \label{sec:bialgebras}
Bialgebras are polynomial spaces endowed with operations. These
operations are very general in the sense that they can have several
inputs and outputs. These structures encompass all the algebraic
structures seen in this work.
\medbreak

\subsection{Biproducts on polynomial spaces}
Polynomial spaces are rather poor algebraic structures. It is usual in
combinatorics to handle spaces endowed with several products. When a
polynomial space is graded and its products are compatible with the
sizes of the underlying combinatorial objects, all this form a graded
algebra. This notion is detailed here, as well as the concepts of
coproduct, duality, and coalgebras and bialgebras.
\medbreak

\subsubsection{Biproducts} \label{subsubsec:biproducts}
Let $C$ be an $I$-collection and $\K \Angle{C}$ be a polynomial space.
A \Def{biproduct} on $\K \Angle{C}$ is a linear map
\begin{equation} \label{equ:biproduct_on_space}
    \Biproduct :
    \K \Angle{\BBrack{C\Par{J_1}, \dots, C\Par{J_p}}_\times}
    \to \K \Angle{\List_{\{q\}}(C)}
\end{equation}
where $p, q \in \N$, and $J_1$, \dots, $J_p$ are nonempty subsets of
$I$. Equivalently, by using the interpretation of the tensor product and
of tensor algebras shown in
Section~\ref{subsec:operations_polynomial_spaces},
\eqref{equ:biproduct_on_space} is equivalent to
\begin{equation}
    \Biproduct :
    \K \Angle{C}\Par{J_1} \otimes \dots \otimes \K \Angle{C}\Par{J_p}
    \to \K \Angle{C}^{\otimes q}.
\end{equation}
The \Def{arity} (resp. \Def{coarity}) of $\Biproduct$ is $p$ (resp. $q$)
and the \Def{index domain} of $\Biproduct$ is the set
$J_1 \times \dots \times J_p$. A tuple
$\Par{x_1, \dots, x_p}$ is a \Def{valid input} for $\Biproduct$ if
$\Biproduct\Par{\Par{x_1, \dots, x_p}}$ is defined, that is,
\begin{math}
    \Par{\Index\Par{x_1}, \dots, \Index\Par{x_p}}
\end{math}
belongs to the index domain of $\Biproduct$. The \Def{image}
$\Image(\Biproduct)$ of $\Biproduct$ is the usual image of $\Biproduct$
as a linear map. To not overload the notation, we shall write
$\Biproduct\Par{x_1, \dots, x_p}$ instead of
$\Biproduct\Par{\Par{x_1, \dots, x_p}}$ for any valid input
$\Par{x_1, \dots, x_p}$ for $\Biproduct$.
\medbreak

The biproduct $\Biproduct$ can be seen as an operation taking a valid
input consisting in a bunch of $p$ objects of $C$ and outputting
bunches of $q$ objects of $C$. This biproduct is depicted by a rectangle
labeled by its name, with $p$ incoming edges (below the rectangle) and
$q$ outgoing edges (above the rectangle) as
\begin{equation}
    \begin{tikzpicture}[xscale=.45,yscale=.2,Centering,font=\scriptsize]
        \node(S1)at(0,0){};
        \node(Sq)at(2,0){};
        \node[Operator](N1)at(1,-3.5){\begin{math}\Biproduct\end{math}};
        \node(E1)at(0,-7){};
        \node(Ep)at(2,-7){};
        \draw[Edge](S1)--(N1);
        \draw[Edge](Sq)--(N1);
        \draw[Edge](N1)--(E1);
        \draw[Edge](N1)--(Ep);
        \node[above of=N1,node distance=6mm]
            {\begin{math}\dots\end{math}};
        \node[below of=N1,node distance=6mm]
            {\begin{math}\dots\end{math}};
        \node[below of=N1]
            {\begin{math}
                f \in
                \K \Angle{\BBrack{C\Par{J_1}, \dots, C\Par{J_p}}_\times}
            \end{math}};
        \node[above of=N1]
            {\begin{math}
                \Biproduct(f) \in \K \Angle{\List_{\{q\}}(C)}
            \end{math}};
    \end{tikzpicture}\,.
\end{equation}
\medbreak

\subsubsection{Completion}
Let $\Biproduct$ a biproduct on $\K \Angle{C}$ of the
form~\eqref{equ:biproduct_on_space}. In the case where the index domain
of $\Biproduct$ is $I^p$, we say that $\Biproduct$ is \Def{complete}.
Otherwise, the \Def{completion} of $\Biproduct$ is the complete
biproduct
\begin{equation}
    \dot{\Biproduct} :
    \K \Angle{\List_{\{p\}}(C)} \to \K \Angle{\List_{\{q\}}(C)}
\end{equation}
defined linearly, for any object $\Par{x_1, \dots, x_p}$ of
$\List_{\{p\}}(C)$, by
\begin{equation}
    \dot{\Biproduct}\Par{x_1, \dots, x_p} :=
    \begin{cases}
        \Biproduct\Par{x_1, \dots, x_p}
            & \mbox{if } \Par{x_1, \dots, x_p}
            \mbox{ is a valid input for } \Biproduct, \\
        0 & \mbox{otherwise}.
    \end{cases}
\end{equation}
In the sequel, we shall provide properties and constructions involving
complete biproducts. Nevertheless, all these apply also on general
biproducts since one can always work with the completion of a
noncomplete biproduct.
\medbreak

\subsubsection{Spaces of complete biproducts}
\label{subsubsec:spaces_biproducts}
The set of all the complete biproducts of arity $p$ and coarity $q$ on
$\K \Angle{C}$ has a structure of a
$\K$-vector space. Indeed, if $\Biproduct_1$ and $\Biproduct_2$ are
two such biproducts, the \Def{addition} of $\Biproduct_1$ and
$\Biproduct_2$ is the biproduct $\Biproduct_1 + \Biproduct_2$ defined by
\begin{equation}
    \Par{\Biproduct_1 + \Biproduct_2}
        \Par{x_1, \dots, x_p}
    :=
    \Biproduct_1\Par{x_1, \dots, x_p}
    +
    \Biproduct_2\Par{x_1, \dots, x_p}
\end{equation}
for any object $\Par{x_1, \dots, x_p}$ of $\List_{\{p\}}(C)$. Moreover,
for any coefficient $\lambda \in \K$, if $\Biproduct$ is such a
biproduct, the \Def{scalar multiplication} of $\Biproduct$ by $\lambda$
is the biproduct $\lambda \Biproduct$ defined by
\begin{equation}
    (\lambda \Biproduct)\Par{x_1, \dots, x_p}
    :=
    \lambda \Biproduct\Par{x_1, \dots, x_p}
\end{equation}
for any object $\Par{x_1, \dots, x_p}$ of $\List_{\{p\}}(C)$.
\medbreak

\subsubsection{Structure coefficient maps}
Let
\begin{equation} \label{equ:structure_coefficient_map}
    \xi : \List_{\{p\}}(C) \times \List_{\{q\}}(C) \to \K
\end{equation}
be a map such for any object $\Par{x_1, \dots, x_p}$ of
$\List_{\{p\}}(C)$, there are finitely many objects
$\Par{y_1, \dots y_q}$ of $\List_{\{q\}}(C)$ such that
\begin{math}
    \xi\Par{\Par{x_1, \dots, x_p}, \Par{y_1, \dots, y_q}} \ne 0.
\end{math}
From this map $\xi$, let the complete biproduct
\begin{equation}
    \Biproduct :
    \K \Angle{\List_{\{p\}}(C)} \to \K \Angle{\List_{\{q\}}(C)}
\end{equation}
satisfying, for any objects $x_1$, \dots, $x_p$ of $C$,
\begin{equation} \label{equ:biproduct_from_map}
    \Biproduct\Par{x_1, \dots, x_p} =
    \sum_{\Par{y_1, \dots, y_q} \in \List_{\{q\}}(C)}
    \xi\Par{\Par{x_1, \dots, x_p}, \Par{y_1, \dots, y_q}} \,
    \Par{y_1, \dots, y_q}.
\end{equation}
Due to the condition satisfied by $\xi$, there is a finite number of
tuples $\Par{y_1, \dots, y_q}$ appearing in the right member
of~\eqref{equ:biproduct_from_map}. Hence, all the
$\Biproduct\Par{x_1, \dots, x_p}$ are $C$-polynomials, so that
$\Biproduct$ is a well-defined complete biproduct on~$\K \Angle{C}$.
\medbreak

Conversely, from any complete biproduct $\Biproduct$ on $\K \Angle{C}$ 
of arity $p$ and coarity $q$, one can recover a map $\xi$ of the 
form~\eqref{equ:structure_coefficient_map} such
that~\eqref{equ:biproduct_from_map} holds. We call $\xi$ the
\Def{structure coefficient map} of $\Biproduct$. Beside, we say that
$\Biproduct$ is \Def{degenerate} if all its structure coefficients are
zero.
\medbreak

\subsubsection{Dual biproducts} \label{subsubsec:dual_biproducts}
Assume here that $C$ is combinatorial so that we can identify
$\K \Angle{C}$ with its dual $\K \Angle{C}^\Dual$ as mentioned in
Section~\ref{subsubsec:duality_polynomial_spaces}.
Given a complete biproduct
$\Biproduct$ on $\K \Angle{C}$ of arity $p$ and coarity $q$, let
\begin{equation}
    \Biproduct^\Dual :
    \K \Angle{\List_{\{q\}}(C)}^\Dual
    \to \K \Angle{\List_{\{p\}}(C)}^\Dual
\end{equation}
be the map linearly defined, for all objects $\Par{y_1, \dots, y_q}$ of
$\List_{\{q\}}(C)$, by
\begin{equation} \label{equ:dual_biproduct}
    \Biproduct^\Dual\Par{y_1, \dots, y_q}
    :=
    \sum_{\Par{x_1, \dots, x_p} \in \List_{\{p\}}(C)}
    \Angle{\Biproduct\Par{x_1, \dots, x_p}, \Par{y_1, \dots, y_q}}
    \, \Par{x_1, \dots, x_p}.
\end{equation}
In the case where~\eqref{equ:dual_biproduct} is a finite sum for any
object $\Par{y_1, \dots, y_q}$ of $\List_{\{q\}}(C)$, its right member
is a $C$-polynomial so that $\Biproduct^\Dual$ is a biproduct on
$\K \Angle{C}^\Dual$, called \Def{dual biproduct} of $\Biproduct$.
\medbreak

Observe that $\Biproduct^\Dual$ is of arity $q$ and coarity $p$, and is
complete. Observe also that in~\eqref{equ:dual_biproduct}, the
coefficient
\begin{math}
    \Angle{\Biproduct\Par{x_1, \dots, x_p}, \Par{y_1, \dots, y_q}}
\end{math}
is in fact equal to
\begin{math}
    \xi\Par{\Par{x_1, \dots, x_p}, \Par{y_1, \dots, y_q}}
\end{math}
where $\xi$ is the structure coefficient map of $\Biproduct$. Hence, if
one sees the map $\xi$ as a matrix whose rows are indexed by the
$\Par{x_1, \dots, x_p}$ and the columns by the $\Par{y_1, \dots, y_q}$,
the structure coefficient map of $\Biproduct^\Dual$ is the transpose of
this matrix.
\medbreak

\subsection{Products on polynomial spaces} \label{subsec:products}
We focus here on products, that are particular biproducts on polynomial
spaces. In all this section, $\K \Angle{C}$ is a polynomial space.
\medbreak

\subsubsection{Products} \label{subsubsec:products}
A \Def{product} is a biproduct of coarity $1$. Let
\begin{equation} \label{equ:product_on_space}
    \Product :
    \K \Angle{\BBrack{C\Par{J_1}, \dots, C\Par{J_p}}_\times}
    \to \K \Angle{C}
\end{equation}
be a product of arity $p \in \N$, where $J_1$, \dots, $J_p$ are
nonempty subsets of $I$. When there is a map
$\omega : J_1 \times \dots \times J_p \to I$ satisfying, for any valid
input $\Par{x_1, \dots, x_p}$ for $\Product$,
\begin{equation} \label{equ:product_on_space_condition}
    \Product\Par{x_1, \dots, x_p}
    \in
    \K \Angle{C} \Par{\omega\Par{\Index\Par{x_1},
    \dots, \Index\Par{x_p}}},
\end{equation}
we say that $\Product$ is \Def{$\omega$-concentrated} (or simply
\Def{concentrated} when it is not useful specify $\omega$). In
intuitive terms, this means that the indexes of the monomials appearing
in a product depend only on the indexes of their operands.
\medbreak

There is a close connection between products on collections (see
Section~\ref{subsubsec:collections_with_products} of
Chapter~\ref{chap:collections}) and products on polynomial spaces.
Indeed, when $C$ is an $I$-collection with a product
\begin{equation}
    \Product : C\Par{J_1} \times \dots \times C\Par{J_p}
    \to C
\end{equation}
of arity $p$, $\Product$ gives rise to a product
\begin{equation}
    \bar{\Product} :
    \K \Angle{\BBrack{C\Par{J_1}, \dots, C\Par{J_p}}_\times}
    \to \K \Angle{C}
\end{equation}
on $\K \Angle{C}$ defined by extending $\Product$ by linearity. This
product $\bar{\Product}$ is called the \Def{linearization} of~$\Product$.
When $\Product$ is a $\omega$-concentrated product on $C$
(see the aforementioned section), its linearization $\bar{\Product}$
is an $\omega$-concentrated product on $\K \Angle{C}$.
\medbreak

\subsubsection{Tensor powers} \label{subsubsec:tensor_powers}
By considering that $\Product$ is a product on $\K \Angle{C}$ of the
form~\eqref{equ:product_on_space}, let us introduce for any
$\ell \in \N_{\geq 1}$, the biproduct
\begin{equation}
    \Tensor_\ell(\Product) :
    \K \Angle{\BBrack{\List_{\{\ell\}}(C\Par{J_1}), \dots,
    \List_{\{\ell\}}(C\Par{J_p})}_\times}
    \to
    \K \Angle{\List_{\{\ell\}}(C)}
\end{equation}
defined linearly by
\begin{multline} \label{equ:product_on_tensors}
    \Tensor_\ell(\Product)\Par{\Par{x_{1, 1}, \dots, x_{\ell, 1}},
    \Par{x_{1, 2}, \dots, x_{\ell, 2}}, \dots,
    \Par{x_{1, p}, \dots, x_{\ell, p}}} \\
    :=
    \Par{\Product\Par{x_{1, 1}, \dots, x_{1, p}},
    \Product\Par{x_{2, 1}, \dots, x_{2, p}}, \dots,
    \Product\Par{x_{\ell, 1}, \dots, x_{\ell, p}}},
\end{multline}
for all
\begin{math}
    \Par{x_{1, k}, \dots, x_{\ell, k}}
    \in \List_{\{\ell\}}\Par{C\Par{J_k}},
\end{math}
$k \in [p]$.
Graphically, $\Tensor_\ell(\Product)$ is the biproduct
\begin{equation}
\,.
\end{equation}
This product $\Tensor_\ell(\Product)$ can be seen as the $\ell$th-tensor
power of $\Product$ seen as a linear map. For this reason,
$\Tensor_\ell(\Product)$ is called the \Def{$\ell$th tensor power} of
$\Product$.
\medbreak

Let us provide an example. When $\Product$ is a complete binary product
on $\K \Angle{C}$, $\Tensor_2(\Product)$ is of the form
\begin{equation}
    \Tensor_2(\Product) :
    \K \Angle{\BBrack{\List_{\{2\}}(C), \List_{\{2\}}(C)}_\times}
    \to \K \Angle{\List_{\{2\}}(C)}
\end{equation}
and it satisfies
\begin{equation} \label{equ:tensor_power_binary_product}
    \Par{x_{1, 1}, x_{2, 1}}
    \, \Tensor_2(\Product) \, \Par{x_{1, 2}, x_{2, 2}} =
    \Par{x_{1, 1} \Product x_{1, 2}, x_{2, 1} \Product x_{2, 2}}
\end{equation}
for all objects $\Par{x_{1, 1}, x_{2, 1}}$ and
$\Par{x_{1, 2}, x_{2, 2}}$ of $\List_{\{2\}}(C)$.
In~\eqref{equ:tensor_power_binary_product}, since $\Product$ and
$\Tensor_2(\Product)$ are binary products, we denote them in infix way.
We follow this convention in all this text. Graphically,
$\Tensor_2(\Product)$ is the biproduct
\begin{equation}
    \begin{tikzpicture}[xscale=.32,yscale=.15,Centering,
        font=\scriptsize]
        \node[Operator](N1)at(1,-3.5){\begin{math}\Product\end{math}};
        \node[Operator](N2)at(5,-3.5){\begin{math}\Product\end{math}};
        \node(S1)at(1,0){};
        \node(S2)at(5,0){};
        \node(E11)at(0,-7){};
        \node(E12)at(2,-7){};
        \node(E21)at(4,-7){};
        \node(E22)at(6,-7){};
        \node(I11)at(0,-12){};
        \node(I12)at(2,-12){};
        \node(I21)at(4,-12){};
        \node(I22)at(6,-12){};
        \draw[Edge](S1)--(N1);
        \draw[Edge](S2)--(N2);
        \draw[Edge](N1)--(E11);
        \draw[Edge](N1)--(E12);
        \draw[Edge](N2)--(E21);
        \draw[Edge](N2)--(E22);
        \draw[Edge](I11)--(E11);
        \draw[Edge](I22)--(E22);
        \draw[Edge](I12)--(2,-10)--(4,-9)--(E21);
        \draw[Edge](I21)--(4,-10)--(2,-9)--(E12);
        \node[below of=I11,node distance=3mm]
            {\begin{math}x_{1, 1}\end{math}};
        \node[below of=I12,node distance=3mm]
            {\begin{math}x_{2, 1}\end{math}};
        \node[below of=I21,node distance=3mm]
            {\begin{math}x_{1, 2}\end{math}};
        \node[below of=I22,node distance=3mm]
            {\begin{math}x_{2, 2}\end{math}};
        \node[above of=S1,node distance=3mm]
            {\begin{math}x_{1, 1} \Product x_{1, 2}\end{math}};
        \node[above of=S2,node distance=3mm]
            {\begin{math}x_{2, 1} \Product x_{2, 2}\end{math}};
    \end{tikzpicture}
\end{equation}
\medbreak

\subsubsection{Products of arity zero}
\label{subsubsec:products_arity_zero}
A product $\eta$ of arity $0$ on $\K \Angle{C}$ is of the form
\begin{equation}
    \eta : \K \Angle{\BBrack{}_\times} \simeq \K \to \K \Angle{C}
\end{equation}
where, as explained in
Section~\ref{subsubsec:cartesian_product_collections} of
Chapter~\ref{chap:collections}, the empty Cartesian product
$\BBrack{}_\times$ of collections contains exactly one element, namely
the empty tuple. Hence, $\eta$ is totally determined by the image
$\eta(1) \in \K \Angle{C}$ where $1 \in \K$. In this way, there is a
correspondence between products of arity zero and elements of
$\K \Angle{C}$. By a slight abuse of notation, we shall write sometimes
$\eta$ instead of $\eta(1)$. In this way, $\eta$ is no longer a map but
an element of~$\K \Angle{C}$.
\medbreak

\subsubsection{Product properties} \label{subsubsec:products_properties}
We now list some properties a product $\Product$ on $\K \Angle{C}$ of
the form~\eqref{equ:product_on_space} can satisfy.
\medbreak

In the particular case where $\K \Angle{C}$ is a graded polynomial
space, $\Product$ is \Def{graded} if $\Product$ is $\omega$-concentrated
for the map $\omega : \N^p \to \N$ defined by
$\omega\Par{\Par{n_1, \dots, n_p}} := n_1 + \dots + n_p$. This notion is
analogous to the one of the same name for collections with products
exposed in Section~\ref{subsubsec:collections_with_products} of
Chapter~\ref{chap:collections}. Observe that when $C$ is a graded
collection with a graded product~$\Product$, its linearization
$\bar{\Product}$ is a graded product on $\K \Angle{C}$.
\medbreak

We now assume that $\K \Angle{C}$ is any polynomial space. If
$\left\{\BasisB_x : x \in C\right\}$ is a basis of $\K \Angle{C}$ such
that, for any valid input $\Par{x_1, \dots, x_p}$ for $\Product$
there is an object $x$ of $C$ satisfying
\begin{equation}
    \Product\Par{\BasisB_{x_1}, \dots, \BasisB_{x_p}} = \BasisB_x,
\end{equation}
we say that the $\BasisB$-basis of $\K \Angle{C}$ is a \Def{set-basis}
with respect to~$\Product$.
\medbreak

Assume now that $\Product$ is of arity $2$ so that $\Product$ is
of the form
\begin{equation} \label{equ:binary_product_on_space}
    \Product :
    \K \Angle{\BBrack{C\Par{J_1}, C\Par{J_2}}_\times}
    \to \K \Angle{C}
\end{equation}
where $J_1$ and $J_2$ are two nonempty subsets of $I$. In the case where
$\Image(\Product)$ is contained in
$\K \Angle{C\Par{J_1 \cap J_2}}$, the \Def{associator} of
$\Product$ is the ternary product
\begin{equation}
    (-, -, -)_\Product :
    \K \Angle{
    \BBrack{C\Par{J_1}, C\Par{J_1 \cap J_2}, C\Par{J_2}}_\times}
    \to \K \Angle{C}
\end{equation}
defined linearly for all valid inputs $\Par{x_1, x_2, x_3}$
for $(-, -, -)_\Product$ by
\begin{equation}
    \Par{x_1, x_2, x_3}_\Product
    := \Par{x_1 \Product x_2} \Product x_3
    - x_1 \Product \Par{x_2 \Product x_3}.
\end{equation}
When, for all valid inputs $\Par{x_1, x_2, x_3}$ for
$(-, -, -)_\Product$, one has
\begin{equation}
    \Par{x_1, x_2, x_3}_\Product = 0,
\end{equation}
we say that $\Product$ is \Def{associative}. The \Def{commutator} of
$\Product$ is the binary product
\begin{equation}
    [-, -]_\Product :
    \K \Angle{\BBrack{C\Par{J_1 \cap J_2}, C\Par{J_2 \cap J_1}}_\times}
    \to \K \Angle{C}
\end{equation}
defined linearly for all valid inputs $\Par{x_1, x_2}$ for
$[-, -]_\Product$ by
\begin{equation}
    \left[x_1, x_2\right]_\Product
    := x_1 \Product x_2 - x_2 \Product x_1.
\end{equation}
When, for all valid inputs $\Par{x_1, x_2}$ for
$\left[x_1, x_2\right]_\Product$, one has
\begin{equation}
    \left[x_1, x_2\right]_\Product = 0,
\end{equation}
the product $\Product$ is \Def{commutative}. When there is a product
$\Unit_\Product$ of arity $0$ such that, for all
$x \in C\Par{J_1 \cap J_2}$,
\begin{equation}
    x \Product \Unit_\Product(1)
    = x
    = \Unit_\Product(1) \Product x,
\end{equation}
we say that $\Product$ is \Def{unitary} and that $\Unit_\Product$ is the
\Def{unit} of~$\Product$. Observe that if $\K \Angle{C}$ is graded and
$\Product$ is a graded product, $\Unit_\Product(1)$ is necessarily of
degree~$0$.
\medbreak

\subsubsection{Coproducts}
A \Def{coproduct} is a biproduct of arity $1$. Observe that when
$\K \Angle{C}$ is combinatorial and that $\Product$ is a concentrated
complete product, its dual $\Product^\Dual$ is a coproduct. This is not
true in general when $\K \Angle{C}$ is not combinatorial or not
concentrated since the conditions exposed in
Section~\ref{subsubsec:dual_biproducts} for the well definition of
$\Product^\Dual$ could not be satisfied.
\medbreak

A coproduct $\upsilon$ of coarity $0$ on $\K \Angle{C}$ is of the form
\begin{equation}
    \upsilon : \K \Angle{C} \to \K \Angle{\BBrack{}_\times} \simeq \K
\end{equation}
and can therefore be seen as a linear form on $\K \Angle{C}$.
\medbreak

All the properties of products defined in
Sections~\ref{subsubsec:products}
and~\ref{subsubsec:products_properties} hold for coproducts which admit
dual products in the following way. For any property $P$ on products, we
say that a coproduct $\Coproduct$ admitting a product $\Coproduct^\Dual$
as dual \Def{satisfies the property ``co$P$''} if $\Coproduct^\Dual$
satisfies $P$. For instance, $\Coproduct$ is \Def{cograded} if
$\Coproduct^\Dual$ is graded, and $\Coproduct$ is \Def{coassociative} if
$\Coproduct^\Dual$ is associative. Moreover, $\Coproduct$ is
\Def{counitary} if there exists a normal form $\Unit_\Coproduct$ on
$\K \Angle{C}$ called \Def{counit} such that its dual
$\Unit_\Coproduct^\Dual$ is the unit of~$\Coproduct^\Dual$.
\medbreak

\subsection{Polynomial bialgebras} \label{subsec:polynomial_bialgebras}
We now consider polynomial spaces endowed with a set of biproducts. The
main definitions and properties of these structures are listed.
\medbreak

\subsubsection{Elementary definitions}
\label{subsubsec:polynomial_bialgebras}
A \Def{polynomial bialgebra} is a pair $(\K \Angle{C}, \Bca)$ where
$\K \Angle{C}$ is a polynomial space endowed with a (possibly infinite)
set $\Bca$ of biproducts. When $\Bca$ contains only products (resp.
coproducts), $(\K \Angle{C}, \Bca)$ is a \Def{polynomial algebra} (resp.
\Def{polynomial coalgebra}). To not overload the notation, we shall
simply write $\K \Angle{C}$ instead of $(\K \Angle{C}, \Bca)$ when
the context is clear.
\medbreak

Let $\Par{\K \Angle{C_1}, \Bca_1}$ and $\Par{\K \Angle{C_2}, \Bca_2}$ be
polynomial bialgebras. These bialgebras are \Def{$\mu$-compatible} if
there exists a bijective map $\mu : \Bca_1 \to \Bca_2$ that sends any
biproduct of $\Bca_1$ to a biproduct of $\Bca_2$ of the same arity,
the same coarity, and the same index domain. When
$\Par{\K \Angle{C_1}, \Bca_1}$ and $\Par{\K \Angle{C_2}, \Bca_2}$ are
$\mu$-compatible, a \Def{$\mu$-polynomial bialgebra morphism} (or simply
a \Def{polynomial bialgebra morphism} when there is no ambiguity) from
$\K \Angle{C_1}$ to $\K \Angle{C_2}$ is a polynomial space morphism
\begin{math}
    \phi : \K \Angle{C_1} \to \K \Angle{C_2}
\end{math}
such that
\begin{equation} \label{equ:morphism_polynomial_bialgebras}
    \Par{\phi^{\otimes q}}
    \Par{\Biproduct\Par{x_1, \dots, x_p}}
    =
    \Par{\mu\Par{\Biproduct}}\Par{\phi\Par{x_1}, \dots, \phi\Par{x_p}}
\end{equation}
for all biproducts $\Biproduct$ or arity $p$ and coarity $q$ of
$\Bca_1$, and valid inputs $\Par{x_1, \dots, x_p}$ for $\Biproduct$,
where $\phi^{\otimes q}$ is the $q$th tensor power $\Tensor_q(\phi)$ of
$\phi$. Graphically,
\eqref{equ:morphism_polynomial_bialgebras} reads as
\begin{equation}
\,.
\end{equation}
In the special case where $\Biproduct$ is a product $\eta$ of arity $0$
(see Section~\ref{subsubsec:products_arity_zero}) of $\Bca_1$,
\eqref{equ:morphism_polynomial_bialgebras} implies that the morphism
$\phi$ satisfies $\phi\Par{\eta} = \mu(\eta)$.
\medbreak

Besides, when $\Par{\K \Angle{C_1}, \Bca_1}$ and
$\Par{\K \Angle{C_2}, \Bca_2}$ are $\mu$-compatible,
$\Par{\K \Angle{C_2}, \Bca_2}$ is a \Def{sub-bialgebra} of
$\Par{\K \Angle{C_1}, \Bca_1}$ if there is an injective
$\mu$-polynomial bialgebra morphism from $\K \Angle{C_2}$ to
$\K \Angle{C_1}$. Let $(\K \Angle{C}, \Bca)$ be a polynomial bialgebra.
For any subset $\GeneratingSet$ of $\K \Angle{C}$, the
\Def{bialgebra generated} by $\GeneratingSet$ is the smallest
sub-bialgebra $\K \Angle{C}^\GeneratingSet$ of $\K \Angle{C}$ containing
$\GeneratingSet$. When $\K \Angle{C}^\GeneratingSet = \K \Angle{C}$ and
$\GeneratingSet$ is minimal with respect to the inclusion among the
subsets of $\GeneratingSet$ satisfying this property, $\GeneratingSet$
is a \Def{minimal generating set} of $\K \Angle{C}$. A bialgebra
$\K \Angle{C}$ can have several minimal generating sets.
\medbreak

A \Def{polynomial bialgebra ideal} of $\K \Angle{C}$ is a space
$\Vca$ included in $\K \Angle{C}$ such that
\begin{equation} \label{equ:bialgebra_ideal_condition}
    \Biproduct\Par{x_1, \dots, x_{i - 1}, f, x_{i + 1}, \dots, x_p}
    \in
    \bigoplus_{j \in [q]}
    \K \Angle{\List_{\{j - 1\}}(C)}
    \otimes \Vca \otimes
    \K \Angle{\List_{\{q - j\}}(C)}
\end{equation}
for all biproducts $\Biproduct$ on $\K \Angle{C}$ of $\Bca$ of arity $p$
and coarity $q$, and all $i \in [p]$, $f \in \Vca$,  $x_k \in C$,
$k \in [p] \setminus \{i\}$ such that the left member
of~\eqref{equ:bialgebra_ideal_condition} is well-defined. Note that
in~\eqref{equ:bialgebra_ideal_condition}, we use the
identification~\eqref{equ:fixed_tensor_polynomial_spaces} between the
space of all the tensors on $\K \Angle{C}$ of a given order
$\ell \in \N$ and the polynomial space on $\List_{\{\ell\}}(C)$. Given a
polynomial bialgebra ideal $\Vca$ of $\K \Angle{C}$, the
\Def{quotient bialgebra} $\K \Angle{C}/_\Vca$ of $\K \Angle{C}$ by
$\Vca$ is defined as follows. Let
$\theta : \K \Angle{C} \to \K \Angle{C}/_\Vca$ be the canonical
surjection map from $\K \Angle{C}$ to its quotient space
$\K \Angle{C}/_\Vca$. Any biproduct $\Biproduct$ on $\K \Angle{C}$
of $\Bca$ of arity $p$ and coarity $q$ gives rise to a biproduct
$\Biproduct_\Vca$ on $\K \Angle{C}/_\Vca$ of arity $p$ and coarity $q$,
and with the same index domain as the one of $\Biproduct$. This
biproduct $\Biproduct_\Vca$ is defined linearly, for any valid input
$\Par{x_1, \dots, x_p}$ for $\Biproduct$, by
\begin{equation} \label{equ:bialgebra_quotient_biproduct}
    \Biproduct_\Vca\Par{\theta\Par{x_1}, \dots, \theta\Par{x_p}}
    :=
    \Tensor_q(\theta)\Par{\Biproduct\Par{x_1, \dots, x_p}}.
\end{equation}
Condition~\eqref{equ:bialgebra_ideal_condition} ensures
that~\eqref{equ:bialgebra_quotient_biproduct} is well-defined.
\medbreak

\subsubsection{Combinatorial polynomial bialgebras}
\label{subsubsec:combinatorial_polynomial_bialgebras}
In practice, and even more so in this book, most of the encountered
polynomial bialgebras are of the form $(\K \Angle{C}, \Bca)$ where
$\Bca$ contains only products and coproducts. We say in this case
that $\K \Angle{C}$ is \Def{uniform}. Moreover, in most practical cases,
$C$ is a graded, a bigraded, or a colored combinatorial collection.
When $(\K \Angle{C}, \Bca)$ is combinatorial, uniform, and all the
products and coproducts of $\Bca$ are complete and concentrated, the
\Def{dual bialgebra} of $\K \Angle{C}$ is the bialgebra
$\Par{\K \Angle{C}^\Dual, \Bca^\Dual}$ where $\Bca^\Dual$ is the set of
the dual biproducts of the biproducts of~$\Bca$. Note that this
bialgebra is also uniform.
\medbreak

It is very common, given a uniform combinatorial bialgebra
$(\K \Angle{C}, \Bca)$, to endow $C$ with a structure of a combinatorial
poset $(C, \Ord)$ in order to construct $\BasisB^{\Ord}$-families (see
Section~\ref{subsubsec:polynomial_spaces_change_basis}). For instance,
when a biproduct $\Biproduct$ has a complicated structure coefficient
map, considering an adequate partial order relation $\Ord$ on $C$
such that the $\BasisB^{\Ord}$-family is a set-basis with respect
to~$\Biproduct$ allows to infer properties of $\Biproduct$ (such as
minimal generating sets of $\K \Angle{C}$, a description of the
nontrivial relations satisfied by these generators, or even freeness
properties).
\medbreak

\subsubsection{Set-theoretic algebras} \label{subsubsec:set_algebras}
Let $(\K \Angle{C}, \Pca)$ be a polynomial algebra and
$\left\{\BasisB_x : x \in C\right\}$ be one of its bases which is
additionally a set-basis for all products of $\Pca$ at the same time. In
this case, it is possible to forget the linear structure
of $\K \Angle{C}$. Indeed, each product $\bar{\Product}$ of arity $p$ on
$\K \Angle{C}$ gives rise to a product $\Product$ on $C$ defined, for
any valid input $\Par{x_1, \dots, x_p}$ for $\bar{\Product}$, by
\begin{equation}
    \Product\Par{x_1, \dots, x_p} := y
\end{equation}
whenever
\begin{equation}
    \bar{\Product}\Par{\BasisB_{x_1}, \dots, \BasisB_{x_p}} = \BasisB_y
\end{equation}
for an $y \in C$. This endows the collection $C$ with products in the
sense of Section~\ref{subsubsec:collections_with_products} of
Chapter~\ref{chap:collections}. We say in this case that $C$ is a
\Def{set-theoretic algebra}.
\medbreak

A large part of the concepts presented above about bialgebras work for
the particular case of set-theoretic algebras with some adjustments. For
instance, to define quotients of a set-theoretic algebra
$\Par{C, \Pca'}$, we do not work with polynomial algebra ideals but with
congruences of set-theoretic algebras. To be a little more precise, a
\Def{set-theoretic algebra congruence} is a relation $\Congr$ on $C$
which is an equivalence relation satisfying
\begin{equation} \label{equ:set_algebra_congruence_condition}
    \Product
        \Par{x_1, \dots, x_{i - 1}, x_i, x_{i + 1}, \dots, x_p}
    \Congr
    \Product
        \Par{x_1, \dots, x_{i - 1}, x_i', x_{i + 1}, \dots, x_p}
\end{equation}
for all products $\Product$ of arities $p$, all $i \in [p]$ such that
$x_i$ and $x_i'$ are objects of $C$ satisfying $x_i \Congr x_i'$,
and all valid inputs
\begin{math}
    \Par{x_1, \dots, x_{i - 1}, x_i, x_{i + 1}, \dots, x_p}
\end{math}
and
\begin{math}
    \Par{x_1, \dots, x_{i - 1}, x_i', x_{i + 1}, \dots, x_p}
\end{math}
for $\Product$.
\medbreak

In the sequel, if ``$N$'' is the name of an algebraic structure, we call
\Def{``set-$N$''} the corresponding set-theoretic structure. For
instance, a set-theoretic unitary associative algebra is a monoid. We
shall further encounter in this way set-operads, colored set-operads,
and set-pros.
\medbreak

\section{Types of polynomial bialgebras} \label{sec:types_bialgebras}
A \Def{type of polynomial bialgebra} is specified by biproduct symbols
together with their arities and coarities, and the possible relations
between them (like, for instance, associativity, commutativity,
cocommutativity, or distributivity). In this section, we list some of
the very ordinary types of uniform polynomial bialgebras in
combinatorics like associative, dendriform and pre-Lie algebras, and
Hopf bialgebras. We give concrete examples for each of these.
\medbreak

\subsection{Associative and coassociative algebras}
\label{subsec:associative_coassociative_algebras}
An \Def{associative algebra} is a polynomial space endowed with a
complete  associative binary product. An associative algebra is
\Def{unitary} if its product is unitary. Besides, an associative
algebra is \Def{commutative} if its product is commutative. To
perfectly fit to the definition of types of bialgebras given above,
the type of unitary associative and commutative algebras is made of a
product symbol $\Product$ of arity $2$ and a product symbol $\Unit$
of arity $0$ together with the relations $(f_1, f_2, f_3)_\Product = 0$,
$[f_1, f_2]_\Product = 0$,
\begin{math}
    f \Product \Unit(1) = f = \Unit(1) \Product f
\end{math},
where and $f_1$, $f_2$, and $f_3$ are any elements of the space.
Moreover, by following the definitions of
Section~\ref{subsubsec:polynomial_bialgebras}, a polynomial space
morphism $\phi : \K \Angle{C_1} \to \K \Angle{C_2}$ is a unitary
commutative algebra morphism between two unitary commutative algebras
$\K \Angle{C_1}$ and $\K \Angle{C_2}$ if for any
$f_1, f_2 \in \K \Angle{C_1}$,
\begin{math}
    \phi\Par{f_1 \Product_1 f_2} =
    \phi\Par{f_1} \Product_2 \phi\Par{f_2}
\end{math}
and
\begin{math}
    \phi\Par{\Unit_1} = \Unit_2,
\end{math}
where $\Product_1$ (resp. $\Product_2$) is the binary product of
$\K \Angle{C_1}$ (resp. $\K \Angle{C_2}$) and $\Unit_1$ (resp.
$\Unit_2$) is the unit of $\K \Angle{C_1}$ (resp.~$\K \Angle{C_2}$).
\medbreak

A \Def{coassociative coalgebra} is a polynomial space endowed with a
coassociative coproduct. A coassociative coalgebra is \Def{counitary} if
its coproduct is counitary. Besides, a coassociative coalgebra is
\Def{cocommutative} if its coproduct is cocommutative.
\medbreak

In all this section,  $A := \{\Asf_1, \dots, \Asf_\ell\}$ is a finite
alphabet.
\medbreak

\subsubsection{Concatenation algebra}
\label{subsubsec:concatenation_algebra}
The \Def{concatenation product} is the complete binary product $\Conc$
on $\K \Angle{A^*}$ defined as the concatenation product of $A^*$
extended by linearity. Since $\Conc$ is graded and all $\K \Angle{A^n}$
are finite dimensional for all $n \in \N$,
$\Par{\K \Angle{A^*}, \Conc}$ is a combinatorial graded algebra.
Moreover, $\Conc$ is associative, noncommutative, and admits the product
of arity zero $\epsilon$, where $\epsilon$ is the empty word, as unit so
that $\Par{\K \Angle{A^*}, \Conc, \epsilon}$ is a unitary noncommutative
associative algebra called \Def{concatenation algebra} on~$A$.
\medbreak

\subsubsection{Shuffle algebra} \label{subsubsec:shuffle_algebra}
The \Def{shuffle product} is the binary product $\shuffle$ on
$\K \Angle{A^*}$ linearly and recursively defined by
\begin{subequations}
\begin{equation} \label{equ:shuffle_words_1}
    u \shuffle \epsilon := u =: \epsilon \shuffle u,
\end{equation}
\begin{equation} \label{equ:shuffle_words_2}
    ua \shuffle vb :=
    (u \shuffle vb) \Conc a + (ua \shuffle v) \Conc b
\end{equation}
\end{subequations}
for any $u, v \in A^*$ and $a, b \in A$, where $\Conc$ is the
concatenation product of the concatenation algebra on $A$. Intuitively,
$\shuffle$ consists in summing in all the ways of interlacing the
letters of the two words as input. For instance,
\begin{equation}\begin{split}
    \ColA{\Asf_1 \Asf_2}
    \shuffle \ColD{\Asf_2 \Asf_1 \Asf_1} & =
    \ColA{\Asf_1} \ColA{\Asf_2}
    \ColD{\Asf_2} \ColD{\Asf_1}
    \ColD{\Asf_1}
    +
    \ColA{\Asf_1} \ColD{\Asf_2}
    \ColA{\Asf_2} \ColD{\Asf_1}
    \ColD{\Asf_1}
    +
    \ColA{\Asf_1} \ColD{\Asf_2}
    \ColD{\Asf_1} \ColA{\Asf_2}
    \ColD{\Asf_1}
    +
    \ColA{\Asf_1} \ColD{\Asf_2}
    \ColD{\Asf_1} \ColD{\Asf_1}
    \ColA{\Asf_2} \\
    & \quad +
    \ColD{\Asf_2} \ColA{\Asf_1}
    \ColA{\Asf_2} \ColD{\Asf_1}
    \ColD{\Asf_1}
    +
    \ColD{\Asf_2} \ColA{\Asf_1}
    \ColD{\Asf_1} \ColA{\Asf_2}
    \ColD{\Asf_1}
    +
    \ColD{\Asf_2} \ColA{\Asf_1}
    \ColD{\Asf_1} \ColD{\Asf_1}
    \ColA{\Asf_2}
    +
    \ColD{\Asf_2} \ColD{\Asf_1}
    \ColA{\Asf_1} \ColA{\Asf_2}
    \ColD{\Asf_1} \\
    & \quad +
    \ColD{\Asf_2} \ColD{\Asf_1}
    \ColA{\Asf_1} \ColD{\Asf_1}
    \ColA{\Asf_2}
    +
    \ColD{\Asf_2} \ColD{\Asf_1}
    \ColD{\Asf_1} \ColA{\Asf_1}
    \ColA{\Asf_2} \\
    & =
    2 \Asf_1 \Asf_2 \Asf_2 \Asf_1 \Asf_1
    + \Asf_1 \Asf_2 \Asf_1 \Asf_2 \Asf_1
    + \Asf_1 \Asf_2 \Asf_1 \Asf_1 \Asf_2
    + \Asf_2 \Asf_1 \Asf_2 \Asf_1 \Asf_1 \\
    & \quad
    + 2 \Asf_2 \Asf_1 \Asf_1 \Asf_2 \Asf_1
    + 3 \Asf_2 \Asf_1 \Asf_1 \Asf_1 \Asf_2.
\end{split}\end{equation}
Since $\shuffle$ is graded and all $\K \Angle{A^n}$ are finite
dimensional for all $n \in \N$, $\Par{\K \Angle{A^*}, \shuffle}$ is a
combinatorial graded algebra. Moreover, $\shuffle$ is associative,
commutative, and admits $\epsilon$ as unit so that
$\Par{\K \Angle{A^*}, \shuffle, \epsilon}$ is a unitary commutative
associative algebra called \Def{shuffle algebra} on~$A$.
\medbreak

\subsubsection{Deconcatenation coalgebra}
\label{subsubsec:deconcatenation_coalgebra}
Let $\Coproduct_\Conc$ be the dual coproduct of the concatenation
product~$\Conc$ of $\K \Angle{A^*}$ considered in
Section~\ref{subsubsec:concatenation_algebra}.
By~\eqref{equ:dual_biproduct}, for all $u \in A^*$,
\begin{equation}
    \Coproduct_\Conc(u) =
    \sum_{v, w \in A^*} \Angle{v \Conc w, u} \, (v, w)
    =
    \sum_{\substack{
        v, w \in A^* \\
        v \Conc w = u
    }}
    (v, w).
\end{equation}
This coproduct is known as the \Def{deconcatenation coproduct}. For
instance,
\begin{equation}
    \Coproduct_\Conc(\Asf_1 \Asf_1 \Asf_2)
    = \Par{\epsilon, \Asf_1 \Asf_1 \Asf_2}
    + \Par{\Asf_1, \Asf_1 \Asf_2}
    + \Par{\Asf_1 \Asf_1, \Asf_2}
    + \Par{\Asf_1 \Asf_1 \Asf_2, \epsilon}.
\end{equation}
Let also $\Counit$ be the dual coproduct of the unit $\epsilon$ for
the concatenation product considered in
Section~\ref{subsubsec:concatenation_algebra}. This coproduct $\Counit$
satisfies $\Counit(\epsilon) = 1$ and $\Counit(u) = 0$ for all nonempty
words $u$. The coalgebra $\Par{\K \Angle{A^*}, \Coproduct, \Counit}$ is
a counitary noncocommutative coassociative coalgebra called
\Def{deconcatenation coalgebra} on~$A$.
\medbreak

\subsubsection{Unshuffle coalgebra}
\label{subsubsec:unshuffle_coalgebra}
Let $\Coproduct_\shuffle$ be the dual coproduct of the shuffle product
$\shuffle$ of $\K \Angle{A^*}$. Again by~\eqref{equ:dual_biproduct}, for
all $u \in A^*$,
\begin{equation}
    \Coproduct_\shuffle(u) =
    \sum_{v, w \in A^*} \Angle{v \shuffle w, u} \, (v, w).
\end{equation}
The coefficient $\Angle{v \shuffle w, u}$ counts the number of ways to
decompose $u$ as two disjoint subwords $v$ and $w$, and thus,
\begin{equation} \label{equ:coproduct_dual_shuffle}
    \Coproduct_\shuffle(u) =
    \sum_{\substack{
        J_1, J_2 \subseteq [|u|] \\
        J_1 \sqcup J_2 = [|u|]
    }}
    \Par{u_{|J_1}, u_{|J_2}}.
\end{equation}
This coproduct can alternatively be expressed by
\begin{equation}
    \Coproduct_\shuffle(\Asf)
    = (\epsilon, \Asf) + (\Asf, \epsilon)
\end{equation}
for any $\Asf \in A$, and
\begin{equation} \label{equ:coproduct_dual_shuffle_iterative}
    \Coproduct_\shuffle(u) =
    \prod_{i \in [|u|]} \Coproduct\Par{u(i)}
\end{equation}
for any $u \in A^*$, where the product
of~\eqref{equ:coproduct_dual_shuffle_iterative} denotes the iterated
version of the $2$nd tensor power $\Tensor_2(\Conc)$ of the
concatenation product $\Conc$ (see
Section~\ref{subsubsec:tensor_powers}). This product $\Tensor_2(\Conc)$
is associative due to the fact that $\Conc$ is associative, and thus,
its iterated version is well-defined. This coproduct is known as the
\Def{unshuffling coproduct}. For instance,
\begin{equation}\begin{split}
    \Coproduct_\shuffle\Par{\Asf_1 \Asf_1 \Asf_2} & =
    \Par{\Par{\epsilon, \Asf_1} + \Par{\Asf_1, \epsilon}}
    \, \Tensor_2(\Conc) \,
    \Par{\Par{\epsilon, \Asf_1} + \Par{\Asf_1, \epsilon}}
    \, \Tensor_2(\Conc) \,
    \Par{\Par{\epsilon, \Asf_2} + \Par{\Asf_2, \epsilon}} \\
    & =
    \Par{\epsilon, \Asf_1 \Asf_1 \Asf_2}
    + \Par{\Asf_2, \Asf_1 \Asf_1}
    + \Par{\Asf_1, \Asf_1 \Asf_2}
    + \Par{\Asf_1 \Asf_2, \Asf_1} \\
    & \quad +
    \Par{\Asf_1, \Asf_1 \Asf_2}
    + \Par{\Asf_1 \Asf_2, \Asf_1}
    + \Par{\Asf_1 \Asf_1, \Asf_2}
    + \Par{\Asf_1 \Asf_1 \Asf_2, \epsilon} \\
    & =
    \Par{\epsilon, \Asf_1 \Asf_1 \Asf_2}
    + \Par{\Asf_2, \Asf_1 \Asf_1}
    + 2 \Par{\Asf_1, \Asf_1 \Asf_2} \\
    & \quad +
    2 \Par{\Asf_1 \Asf_2, \Asf_1}
    + \Par{\Asf_1 \Asf_1, \Asf_2}
    + \Par{\Asf_1 \Asf_1 \Asf_2, \epsilon}.
\end{split}\end{equation}
The coalgebra $\Par{\K \Angle{A^*}, \shuffle, \Counit}$, where $\Counit$
is the counit considered in
Section~\ref{subsubsec:deconcatenation_coalgebra}, is a counitary
cocommutative coassociative coalgebra called \Def{unshuffle coalgebra}
on~$A$.
\medbreak

\subsection{Dendriform algebras} \label{subsec:dendriform_algebras}
A \Def{dendriform algebra} is a polynomial space $\K \Angle{C}$ endowed
with two complete binary products $\LDendr$ and $\RDendr$ satisfying
\begin{subequations}
    \begin{equation} \label{equ:relation_dendr_1}
        \Par{x_1 \LDendr x_2} \LDendr x_3
        = x_1 \LDendr \Par{x_2 \LDendr x_3}
        + x_1 \LDendr \Par{x_2 \RDendr x_3},
    \end{equation}
    \begin{equation} \label{equ:relation_dendr_2}
        \Par{x_1 \RDendr x_2} \LDendr x_3
        = x_1 \RDendr \Par{x_2 \LDendr x_3},
    \end{equation}
    \begin{equation} \label{equ:relation_dendr_3}
        \Par{x_1 \LDendr x_2} \RDendr x_3
        + \Par{x_1 \RDendr x_2} \RDendr x_3
        = x_1 \RDendr \Par{x_2 \RDendr x_3},
    \end{equation}
\end{subequations}
for all objects $x_1$, $x_2$, and $x_3$ of $C$.
\medbreak

A polynomial algebra $(\K \Angle{C}, \Product)$, where $\Product$ is a
binary product, admits a \Def{dendriform algebra structure} if its
product can be split into two operations
\begin{equation} \label{equ:dendriform_split}
    \Product = \LDendr + \RDendr,
\end{equation}
where $\LDendr$ and $\RDendr$ are two non-degenerate binary products
such that $(\K \Angle{C}, \LDendr, \RDendr)$ is a dendriform algebra.
Expression~\eqref{equ:dendriform_split} uses the addition of biproducts
exposed in Section~\ref{subsubsec:spaces_biproducts}. Observe that if
$(\K \Angle{C}, \Product)$ admits a dendriform algebra structure,
$\Product$ is associative. The associativity of $\LDendr + \RDendr$ is
a consequence of Relations~\eqref{equ:relation_dendr_1},
\eqref{equ:relation_dendr_2}, and~\eqref{equ:relation_dendr_3} of
dendriform algebras.
\medbreak

In all this section,  $A := \left\{\Asf_1, \dots, \Asf_\ell\right\}$ is
a finite alphabet.
\medbreak

\subsubsection{Shuffle dendriform algebra}
\label{subsubsec:shuffle_dendriform_algebra}
Consider on $\K \Angle{A^*}$ the binary products $\LDendr$ and $\RDendr$
defined linearly and recursively by
\begin{subequations}
\begin{equation} \label{equ:dendr_words_epsilson_1}
    u \LDendr \epsilon := u =: \epsilon \RDendr u,
\end{equation}
\begin{equation} \label{equ:dendr_words_epsilson_2}
    w \RDendr \epsilon =: 0 := \epsilon \LDendr w,
\end{equation}
\begin{multicols}{2}
\begin{equation} \label{equ:left_dendr_words}
    ua \LDendr v := (u \shuffle v) \Conc a,
\end{equation}

\begin{equation} \label{equ:right_dendr_words}
    u \RDendr vb := (u \shuffle v) \Conc b
\end{equation}
\end{multicols}
\end{subequations}
\noindent
for any $u, v \in A^*$, $w \in A^+$, and $a, b \in A$, where $\Conc$ is
the concatenation product of words. In other words, $u \LDendr v$ (resp.
$u \RDendr v$) is the sum of all the words $w$ obtained by shuffling $u$
and $v$ such that the last letter of $w$ comes from $u$ (resp. $v$). For
example,
\begin{subequations}
\begin{equation}\begin{split}
    \ColA{\Asf_1 \Asf_2}
    \LDendr \ColD{\Asf_2 \Asf_1 \Asf_1} & =
    \ColA{\Asf_1} \ColD{\Asf_2}
    \ColD{\Asf_1} \ColD{\Asf_1}
    \ColA{\Asf_2}
    +
    \ColD{\Asf_2} \ColA{\Asf_1}
    \ColD{\Asf_1} \ColD{\Asf_1}
    \ColA{\Asf_2}
    +
    \ColD{\Asf_2} \ColD{\Asf_1}
    \ColA{\Asf_1} \ColD{\Asf_1}
    \ColA{\Asf_2}
    +
    \ColD{\Asf_2} \ColD{\Asf_1}
    \ColD{\Asf_1} \ColA{\Asf_1}
    \ColA{\Asf_2} \\
    & =
    \Asf_1 \Asf_2 \Asf_1 \Asf_1 \Asf_2
    + 3 \Asf_2 \Asf_1 \Asf_1 \Asf_1 \Asf_2,
\end{split}\end{equation}
\begin{equation}\begin{split}
    \ColA{\Asf_1 \Asf_2}
    \RDendr \ColD{\Asf_2 \Asf_1 \Asf_1} & =
    \ColA{\Asf_1} \ColA{\Asf_2}
    \ColD{\Asf_2} \ColD{\Asf_1}
    \ColD{\Asf_1}
    +
    \ColA{\Asf_1} \ColD{\Asf_2}
    \ColA{\Asf_2} \ColD{\Asf_1}
    \ColD{\Asf_1}
    +
    \ColA{\Asf_1} \ColD{\Asf_2}
    \ColD{\Asf_1} \ColA{\Asf_2}
    \ColD{\Asf_1}
    +
    \ColD{\Asf_2} \ColA{\Asf_1}
    \ColA{\Asf_2} \ColD{\Asf_1}
    \ColD{\Asf_1} \\
    & \quad +
    \ColD{\Asf_2} \ColA{\Asf_1}
    \ColD{\Asf_1} \ColA{\Asf_2}
    \ColD{\Asf_1}
    +
    \ColD{\Asf_2} \ColD{\Asf_1}
    \ColA{\Asf_1} \ColA{\Asf_2}
    \ColD{\Asf_1} \\
    & =
    2 \Asf_1 \Asf_2 \Asf_2 \Asf_1 \Asf_1
    + \Asf_1 \Asf_2 \Asf_1 \Asf_2 \Asf_1
    + \Asf_2 \Asf_1 \Asf_2 \Asf_1 \Asf_1
    + 2 \Asf_2 \Asf_1 \Asf_1 \Asf_2 \Asf_1.
\end{split}\end{equation}
\end{subequations}
These two products endow $\K \Angle{A^*}$ with a structure of a
dendriform algebra called \Def{shuffle dendriform algebra} on $A$. This
shows moreover that the shuffle algebra $(\K \Angle{A^*}, \shuffle)$
admits a dendriform algebra structure since
\begin{equation} \label{equ:splitting_shuffle_product}
    u \shuffle v = u \LDendr v + u \RDendr v
\end{equation}
for all $u, v \in A^*$. This offers also a way to recover the
recursive definition (see~\eqref{equ:shuffle_words_1}
and~\eqref{equ:shuffle_words_2}) of~$\shuffle$.
\medbreak

\subsubsection{Max dendriform algebra}
Assume here that $A$ is a totally ordered alphabet by
$\Asf_i \leq \Asf_j$ if $i \leq j$. Consider on $\K \Angle{A^+}$ the
binary products $\LDendr$ and $\RDendr$ defined linearly by
\begin{subequations}
\begin{equation}
    u \LDendr v :=
    \begin{cases}
        u \Conc v & \mbox{if } \max_\leq(u) \geq \max_\leq(v) \\
        0 & \mbox{otherwise},
    \end{cases}
\end{equation}
\begin{equation}
    u \RDendr v :=
    \begin{cases}
        u \Conc v & \mbox{if } \max_\leq(u) < \max_\leq(v) \\
        0 & \mbox{otherwise},
    \end{cases}
\end{equation}
\end{subequations}
for all $u, v \in A^+$, where $\Conc$ is the concatenation product of
words. These two products endow $\K \Angle{A^+}$ with a structure of a
dendriform algebra called \Def{max dendriform algebra}. Moreover, we
have here $\Conc = \LDendr + \RDendr$ where~$\Conc$ is the associative
algebra product of concatenation of~$\K \Angle{A^+}$ so that
$(\K \Angle{A^+}, \Conc)$ admits a dendriform algebra structure.
\medbreak

\subsection{Pre-Lie algebras} \label{subsec:pre_lie_algebras}
A \Def{pre-Lie algebra} is a polynomial space $\K \Angle{C}$ endowed
with a binary product $\PreLieProduct$ satisfying
\begin{equation} \label{equ:relation_pre_Lie}
    \Par{x_1 \PreLieProduct x_2} \PreLieProduct x_3
    - x_1 \PreLieProduct \Par{x_2 \PreLieProduct x_3}
    = \Par{x_1 \PreLieProduct x_3} \PreLieProduct x_2
    - x_1 \PreLieProduct \Par{x_3 \PreLieProduct x_2}
\end{equation}
for all objects $x_1$, $x_2$, and $x_3$ of $C$. This
relation~\eqref{equ:relation_pre_Lie} of pre-Lie algebras says that the
associator $(-, -, -)_\PreLieProduct$ is symmetric in its two last
entries, that is
\begin{math}
    \Par{x_1, x_2, x_3}_\PreLieProduct
    =
    \Par{x_1, x_3, x_2}_\PreLieProduct.
\end{math}
\medbreak

\subsubsection{Pre-Lie algebras from associative algebras}
When $(\K \Angle{C}, \Product)$ is an associative algebra, $\Product$
satisfies in particular~\eqref{equ:relation_pre_Lie} since both left and
right members are equal to zero. For this reason,
$(\K \Angle{C}, \Product)$ is a pre-Lie algebra.
\medbreak

\subsubsection{Pre-Lie algebra of rooted trees}
\label{subsubsec:pre_lie_algebra_rooted_trees}
Recall that $\ColRT$ is the combinatorial graded collection of all
rooted trees (see Section~\ref{subsec:rooted_trees} of
Chapter~\ref{chap:trees}). Consider now on $\K \Angle{\ColRT}$ the
products
\begin{equation}
    \Graft^{(p)} :
    \K \Angle{\List_{\{p\}}\Par{\ColRT}}
    \to \K \Angle{\ColRT}
\end{equation}
defined linearly for all $p \in \N_{\geq 1}$ and all rooted trees
$\Tfr_1$, \dots, $\Tfr_p$ by
\begin{equation}
    \Graft^{(p)}\Par{\Tfr_1, \dots, \Tfr_p} :=
    \Par{\Node, \lbag \Tfr_1, \dots, \Tfr_p\rbag}.
\end{equation}
Intuitively, $\Graft^{(p)}$ consists in grafting all the trees $\Tfr_1$,
\dots, $\Tfr_p$ onto a common root. This product is symmetric with
respect to all its inputs. Now, let $\PreLieProduct$ be the binary
product on $\K \Angle{\ColRT}$ defined linearly and recursively by
\begin{equation}
    \Sfr \PreLieProduct \Tfr :=
    \Graft^{(p + 1)}
        \Par{\Sfr_1, \dots, \Sfr_p, \Tfr}
    +
    \sum_{i \in [p]}
    \Graft^{(p)}(\Sfr_1, \dots, \Sfr_{i - 1},
    \Par{\Sfr_i \PreLieProduct \Tfr},
    \Sfr_{i + 1}, \dots, \Sfr_p)
\end{equation}
for any $\Sfr, \Tfr \in \ColRT$ where
$\Sfr = \Par{\Node, \lbag \Sfr_1, \dots, \Sfr_p \rbag}$. Intuitively,
$\PreLieProduct$ consists in summing all the ways of connecting the root
of the second operand on a node of the first. For example,
\begin{equation}
\,.
\end{equation}
This product endows $\K \Angle{\ColRT}$ with a structure of a pre-Lie
algebra, called the \Def{pre-Lie algebra of rooted trees}.
\medbreak

\subsection{Hopf bialgebras}
A \Def{Hopf bialgebra} is a polynomial space $\K \Angle{C}$ endowed with
a complete binary product $\Product$ and a complete binary coproduct
$\Coproduct$ such that $\Par{\K \Angle{C}, \Product, \Unit}$ is a
unitary associative algebra, $\Par{\K \Angle{C}, \Coproduct, \Counit}$
is a counitary coassociative coalgebra, and
\begin{subequations}
\begin{equation} \label{equ:hopf_compatibility}
    \Coproduct\Par{x_1 \Product x_2} =
    \Coproduct\Par{x_1} \, \Tensor_2(\Product) \, \Coproduct\Par{x_2},
\end{equation}
\begin{equation}
    \Counit\Par{x_1 \Product x_2}
    = \Counit\Par{x_1} \Counit\Par{x_2},
\end{equation}
\begin{equation}
    \Delta(\Unit) = (\Unit, \Unit),
\end{equation}
\begin{equation}
    \Counit(\Unit) = 1,
\end{equation}
\end{subequations}
for all objects $x_1$ and $x_2$ of $C$. When
$\Par{\K \Angle{C}, \Product, \Unit, \Coproduct, \Counit}$ is
combinatorial and all its (co)-products are concentrated, its dual is
well-defined and is still a Hopf bialgebra.
\medbreak

Let us now provide some classical definitions about Hopf bialgebras.
\medbreak

\subsubsection{Primitive and group-like elements}
Let
\begin{math}
    \Par{\K \Angle{C}, \Product, \Unit, \Coproduct, \Counit}
\end{math}
be a Hopf bialgebra. An element $f$ of $\K \Angle{C}$ is \Def{primitive}
if $\Coproduct(f) = \Par{\Unit, f} + \Par{f, \Unit}$. The set
$\Pca_{\K \Angle{C}}$ of all primitive elements of $\K \Angle{C}$ forms
a subspace of $\K \Angle{C}$ and the commutator $[-, -]_\Product$ endows
$\Pca_{\K \Angle{C}}$ with a structure of a Lie algebra. Besides, an
element $f$ of $\K \Angle{C}$ is \Def{group-like} if
$\Coproduct(f) = (f, f)$.
\medbreak

\subsubsection{Convolution product and antipode}
Given two Hopf bialgebras
$\Par{\K \Angle{C_1}, \Product_1, \Unit_1, \Coproduct_1, \Counit_1}$
and
$\Par{\K \Angle{C_2}, \Product_2, \Unit_2, \Coproduct_2, \Counit_2}$, if
$\omega$ and $\omega'$ are two Hopf bialgebra morphisms from
$\K \Angle{C_1}$ to $\K \Angle{C_2}$, the \Def{convolution} of $\omega$
and $\omega'$ is the map
\begin{equation}
    \omega \Convolution \omega' : \K \Angle{C_1} \to \K \Angle{C_2}
\end{equation}
defined linearly, for any object $x$ of $C_1$, by
\begin{equation}
    \Par{\omega \Convolution \omega'}(x)
    :=
    \sum_{y_1, y_2 \in C_1}
    \xi_{\Coproduct_1}\Par{(x), \Par{y_1, y_2}} \,
    \omega\Par{y_1} \Product_2 \omega'\Par{y_2},
\end{equation}
where $\xi_{\Coproduct_1}$ is the structure coefficient map of
$\Coproduct_1$. This convolution product is associative, as a
consequence of the fact that $\Coproduct_1$ is coassociative and
$\Product_2$ is associative.
\medbreak

Now, let $\Par{\K \Angle{C}, \Product, \Unit, \Coproduct, \Counit}$ be a
Hopf bialgebra. Let $\nu : \K \Angle{C} \to \K \Angle{C}$ be the linear
map defined as the inverse of the identity map
$\Identity_{\K \Angle{C}}$ on $\K \Angle{C}$. This map $\nu$ is the
\Def{antipode} of $\K \Angle{C}$ and it can be undefined in certain
cases.
\medbreak

\subsubsection{Combinatorial connected graded Hopf bialgebras}
In algebraic combinatorics, one encounters very particular Hopf
bialgebras. Most of these are combinatorial connected graded Hopf
bialgebras $\Par{\K \Angle{C}, \Product, \Unit, \Coproduct, \Counit}$.
These structures have hence a graded product (that is,
$x_1 \Product x_2$ is homogeneous and of degree $|x_1| + |x_2|$ for all
objects $x_1$ and $x_2$ of $C$), a cograded coproduct (that is, the sum
of the sizes of each $\Par{y_1, y_2}$ appearing in $\Coproduct\Par{x}$
is $|x|$ for all objects $x$ of $C$). Moreover, each $\K \Angle{C}(n)$,
$n \in \N$, is finite dimensional and, since $\# C(0) = 1$,
$\K \Angle{C}(0)$ can be identified with $\K$. Additionally, since
$\Product$ is graded, this implies that the unit $\Unit$ is the unique
element of $C(0)$. Finally, the counit $\Counit$ is the linear map
behaving as the identity map on $\K \Angle{C}(0)$ and sending all the
elements of $\K \Angle{C}(\N_{\geq 1})$ to~$0$.
\medbreak

\begin{Proposition} \label{prop:antipode_combinatorial_Hopf_bialgebra}
    Let $\Par{\K \Angle{C}, \Product, \Unit, \Coproduct, \Counit}$ be
    a combinatorial connected graded Hopf bialgebra. Then,
    $\K \Angle{C}$ admits a unique antipode and it satisfies, for any
    $x \in C$, the recurrence
    \begin{equation} \label{equ:antipode_formula}
        \nu(x)
        =
        \delta_{\Unit, x}
        -\sum_{\substack{
            y_1, y_2 \in C \\
            y_2 \ne \Unit
        }}
        \xi_{\Coproduct}
        \Par{\Par{x}, \Par{y_1, y_2}} \,
        \nu(y_1) \Product y_2,
    \end{equation}
    where $\delta_{-, -}$ is the Kronecker symbol and
    $\xi_{\Coproduct}$ is the structure coefficient map of $\Coproduct$.
\end{Proposition}
\medbreak

Therefore, \eqref{equ:antipode_formula} implies that the antipode of
$\K \Angle{C}$ can be computed by induction on the sizes of the
objects.
\medbreak

\subsubsection{Shuffle deconcatenation Hopf bialgebra}
Let $A := \{\Asf_1, \dots, \Asf_\ell\}$ be a finite alphabet. The
concatenation product $\Conc$ (see
Section~\ref{subsubsec:concatenation_algebra}), the unit $\epsilon$ (see
Section~\ref{subsubsec:concatenation_algebra}), the unshuffling
coproduct $\Coproduct_\shuffle$ (see
Section~\ref{subsubsec:unshuffle_coalgebra}), and the counit $\Counit$
(see Section~\ref{subsubsec:deconcatenation_coalgebra}) endow
$\K \Angle{A^*}$ with a structure of a combinatorial connected graded
Hopf bialgebra
$\Par{\K \Angle{A^*}, \Conc, \epsilon, \Coproduct_\shuffle, \Counit}$.
Its dual bialgebra is the Hopf bialgebra
$\Par{\K \Angle{A^*}, \shuffle, \epsilon, \Coproduct_\Conc, \Counit}$
where $\shuffle$ is the shuffle product (see
Section~\ref{subsubsec:shuffle_algebra}) and $\Coproduct_\Conc$ is the
deconcatenation coproduct (see again
Section~\ref{subsubsec:deconcatenation_coalgebra}).
\medbreak

\subsubsection{Noncommutative symmetric functions}
Consider the combinatorial graded polynomial space
$\Sym := \K \Angle{\ColComp}$ of the compositions (see
Section~\ref{subsubsec:integer_compositions} of
Chapter~\ref{chap:collections}). Let
\begin{math}
    \left\{\BasisS_{\bm{\lambda}} : \bm{\lambda} \in \ColComp\right\}
\end{math}
be the basis of the \Def{complete noncommutative symmetric functions} of
$\Sym$ and $\Product$ be the binary product defined linearly, for any
$\bm{\lambda}, \bm{\mu} \in \ColComp$, by
\begin{equation}
    \BasisS_{\bm{\lambda}} \Product \BasisS_{\bm{\mu}} :=
    \BasisS_{\bm{\lambda} \Conc \bm{\mu}},
\end{equation}
where $\bm{\lambda} \Conc \bm{\mu}$ is the concatenation of the
compositions (seen as words of integers). Moreover, let $\Coproduct$ be
the binary coproduct defined linearly, for any
$\bm{\lambda} \in \ColComp$, by
\begin{equation} \label{equ:coproduct_sym_s_basis}
    \Coproduct\Par{\BasisS_{\bm{\lambda}}} :=
    \prod_{j \in [\Length(\bm{\lambda})]}
    \Par{
    \sum_{\substack{
        n, m \in \N \\
        n + m = \bm{\lambda}_j}
    }
    \Par{\BasisS_{(n)}, \BasisS_{(m)}}},
\end{equation}
where the product of~\eqref{equ:coproduct_sym_s_basis} denotes the
iterated version of $2$nd tensor power $\Tensor_2(\Product)$
of~$\Product$, and for any $n \in \N_{\geq 1}$, $\BasisS_{(n)}$ is the
basis element indexed by the composition of length $1$ whose only part
is $n$, and $\BasisS_{(0)}$ is identified with the unit $1$ of $\K$. For
instance,
\begin{equation}\begin{split}
    \Coproduct\Par{\BasisS_{121}} & =
    \Par{\Par{1, \BasisS_1} + \Par{\BasisS_1, 1}}
    \, \Tensor_2(\Product) \,
    \Par{\Par{1, \BasisS_2} + \Par{\BasisS_1, \BasisS_1}
        + \Par{\BasisS_2, 1}}
    \, \Tensor_2(\Product) \,
    \Par{\Par{1, \BasisS_1} + \Par{\BasisS_1, 1}} \\
    & =
    \Par{1, \BasisS_{121}} +
    \Par{\BasisS_1, \BasisS_{111}} +
    \Par{\BasisS_1, \BasisS_{12}} +
    \Par{\BasisS_1, \BasisS_{21}} +
    2 \Par{\BasisS_{11}, \BasisS_{11}} +
    \Par{\BasisS_{11}, \BasisS_2} \\
    & \quad +
    \Par{\BasisS_2, \BasisS_{11}} +
    \Par{\BasisS_{111}, \BasisS_1} +
    \Par{\BasisS_{12}, \BasisS_1} +
    \Par{\BasisS_{21}, \BasisS_1} +
    \Par{\BasisS_{121}, 1}.
\end{split}\end{equation}
The product $\Product$ and the coproduct $\Coproduct$ endow $\Sym$ with
a structure of a combinatorial connected graded Hopf bialgebra.
\medbreak

Moreover, let
\begin{math}
    \left\{\BasisR_{\bm{\lambda}} : \bm{\lambda} \in \ColComp\right\}
\end{math}
be the family defined by
\begin{equation}
    \BasisR_{\bm{\lambda}} :=
    \sum_{\substack{
        \bm{\mu} \in \ColComp \\
        \bm{\lambda} \Ord \bm{\mu}
    }}
    (-1)^{\Length(\bm{\lambda}) - \Length(\bm{\mu})} \,
    \BasisS_{\bm{\mu}},
\end{equation}
where $\Ord$ is the refinement order of compositions (see
Section~\ref{subsubsec:cube_poset} of Chapter~\ref{chap:collections}).
For instance,
\begin{equation}
    \BasisR_{212} = \BasisS_{212} - \BasisS_{23} - \BasisS_{32}
    + \BasisS_{5}.
\end{equation}
By triangularity, this family forms a basis of $\Sym$ and is known as
the basis of \Def{ribbon noncommutative symmetric functions}. On this
basis, one has, for any $\bm{\lambda}, \bm{\mu} \in \ColComp$,
\begin{equation}
    \BasisR_{\bm{\lambda}} \Product \BasisR_{\bm{\mu}} :=
    \BasisR_{\bm{\lambda} \Conc \bm{\mu}} +
    \BasisR_{\bm{\lambda} \ProductTriangleRight \bm{\mu}},
\end{equation}
for any $\bm{\lambda}, \bm{\mu} \in \ColComp$, where
$\bm{\lambda} \Conc \bm{\mu}$ is the concatenation of the compositions
and
\begin{equation}
    \bm{\lambda} \ProductTriangleRight \bm{\mu} :=
    \Par{\bm{\lambda}_1, \dots,
    \bm{\lambda}_{\Length(\bm{\lambda}) - 1},
    \bm{\lambda}_{\Length(\bm{\lambda})} + \bm{\mu}_1,
    \bm{\mu}_2, \dots, \bm{\mu}_{\Length(\bm{\mu})}}.
\end{equation}
For instance,
\begin{equation}
    \BasisR_{\ColA{3112}} \Product
    \BasisR_{\ColD{142}} =
    \BasisR_{\ColA{3112}\ColD{142}} +
    \BasisR_{\ColA{311}
        \ColB{3}\ColD{42}}.
\end{equation}
\medbreak

This Hopf bialgebra $\Sym$ is usually known as the \Def{Hopf bialgebra
of noncommutative symmetric functions}. To explain this name, consider
an alphabet $A := \left\{\Asf_1, \Asf_2, \dots\right\}$ equipped with a
total order $\Ord$ where $1 \leq i \leq j$ implies $\Asf_i \Ord \Asf_j$.
Now, let the noncommutative series
\begin{equation}
    \BasisR_{\bm{\lambda}}(A) :=
    \sum_{\substack{
        u \in A^* \\
        \Cmp(u) = \bm{\lambda}
    }}
    u,
\end{equation}
of $\K \AAngle{A^*}$ defined for all $\bm{\lambda} \in \ColComp$,
where $\Cmp$ is defined in Section~\ref{subsubsec:integer_compositions}
of Chapter~\ref{chap:collections}. Observe that all
$\BasisR_{\bm{\lambda}}(A)$ are polynomials when $A$ is finite, but are
series in the other case. For instance,
\begin{subequations}
\begin{equation}
    \BasisR_{\ColA{3}\ColD{1}}
        (\{\Asf_1, \Asf_2\}) =
    \ColA{\Asf_1 \Asf_1 \Asf_2} \ColD{\Asf_1} +
    \ColA{\Asf_1 \Asf_2 \Asf_2} \ColD{\Asf_1} +
    \ColA{\Asf_2 \Asf_2 \Asf_2} \ColD{\Asf_1},
\end{equation}
\begin{equation}\begin{split} \label{equ:example_realization_ribbon_sym}
    \BasisR_{\ColA{2}\ColD{1}}
        (\{\Asf_1, \Asf_2, \Asf_3\}) & =
    \ColA{\Asf_1 \Asf_2} \ColD{\Asf_1} +
    \ColA{\Asf_1 \Asf_3} \ColD{\Asf_1} +
    \ColA{\Asf_1 \Asf_3} \ColD{\Asf_2} +
    \ColA{\Asf_2 \Asf_2} \ColD{\Asf_1} \\
    & \quad +
    \ColA{\Asf_2 \Asf_3} \ColD{\Asf_1} +
    \ColA{\Asf_2 \Asf_3} \ColD{\Asf_2} +
    \ColA{\Asf_3 \Asf_3} \ColD{\Asf_1} +
    \ColA{\Asf_3 \Asf_3} \ColD{\Asf_2},
\end{split}\end{equation}
 \begin{equation}\begin{split}
    \BasisR_{\ColA{1}\ColD{2}
        \ColB{1}}(\{\Asf_1, \Asf_2, \Asf_3\}) & =
    \ColA{\Asf_2} \ColD{\Asf_1 \Asf_2}
        \ColB{\Asf_1} +
    \ColA{\Asf_2} \ColD{\Asf_1 \Asf_3}
        \ColB{\Asf_1} +
    \ColA{\Asf_2} \ColD{\Asf_1 \Asf_3}
        \ColB{\Asf_2} +
    \ColA{\Asf_3} \ColD{\Asf_1 \Asf_2}
        \ColB{\Asf_1} +
    \ColA{\Asf_3} \ColD{\Asf_1 \Asf_3}
        \ColB{\Asf_1} \\
    & \quad +
    \ColA{\Asf_3} \ColD{\Asf_1 \Asf_3}
        \ColB{\Asf_2} +
    \ColA{\Asf_3} \ColD{\Asf_2 \Asf_2}
        \ColB{\Asf_1} +
    \ColA{\Asf_3} \ColD{\Asf_2 \Asf_3}
        \ColB{\Asf_1} +
    \ColA{\Asf_3} \ColD{\Asf_2 \Asf_3}
        \ColB{\Asf_2}.
\end{split}\end{equation}
\end{subequations}
The linear span of all the $\BasisR_{\bm{\lambda}}(A)$,
$\bm{\lambda} \in \ColComp$, is the space of the noncommutative
symmetric functions on $A$. The associative algebra structure of
$\Sym$ is compatible with these series in the sense that
\begin{equation} \label{equ:product_sym_series}
    \BasisR_{\bm{\lambda}}(A) \Conc \BasisR_{\bm{\mu}}(A) =
    \Par{\BasisR_{\bm{\lambda}} \Product \BasisR_{\bm{\mu}}}(A)
\end{equation}
for all $\bm{\lambda}, \bm{\mu} \in \ColComp$, where the product $\Conc$
of the left member of~\eqref{equ:product_sym_series} is the usual
product of noncommutative series of~$\K \AAngle{A^*}$.
\medbreak

\subsubsection{Free quasi-symmetric noncommutative symmetric functions}
\label{subsubsec:fqsym}
Let the gra\-ded combinatorial polynomial space
$\FQSym := \K \Angle{\SymmetricGroup}$ of the permutations. Let
$\left\{\BasisF_\sigma : \sigma \in \SymmetricGroup\right\}$ be the
basis of the \Def{fundamental free quasi-symmetric functions} of
$\FQSym$ and $\Product$ be the binary product defined linearly, for any
$\sigma, \nu \in \SymmetricGroup$, by
\begin{equation} \label{equ:product_fqsym}
    \BasisF_\sigma \Product \BasisF_\nu :=
    \sum_{\pi \in \SymmetricGroup}
    \Angle{\pi, \sigma \shuffle \bar{\nu}} \BasisF_\pi,
\end{equation}
where $\bar{\nu}$ is the word obtained by incrementing each letter of
$\nu$ by $|\sigma|$, and $\shuffle$ is the shuffle product of words
defined in Section~\ref{subsubsec:shuffle_algebra}. For instance
\begin{equation}
    \BasisF_{\ColA{21}}
    \Product
    \BasisF_{\ColD{12}} =
    \BasisF_{\ColA{21}\ColD{34}}
    +
    \BasisF_{\ColA{2}\ColD{3}
    \ColA{1}\ColD{4}}
    +
    \BasisF_{\ColA{2}\ColD{3}
    \ColD{4}\ColA{1}}
    +
    \BasisF_{\ColD{3}\ColA{2}
    \ColA{1}\ColD{4}}
    +
    \BasisF_{\ColD{3}\ColA{2}
    \ColD{4}\ColA{1}}
    +
    \BasisF_{\ColD{34}\ColA{21}}.
\end{equation}
This product is known as the \Def{shifted shuffle product} and is
sometimes denoted also as $\cshuffle$. Let moreover $\Coproduct$ be the
binary coproduct defined linearly, for any $\pi \in \SymmetricGroup$, by
\begin{equation} \label{equ:coproduct_fqsym}
    \Coproduct\Par{\BasisF_\pi} :=
    \sum_{0 \leq i \leq |\pi|}
    \Par{\BasisF_{\Std\Par{\pi(1) \dots \pi(i)}},
    \BasisF_{\Std\Par{\pi(i + 1) \dots \pi\Par{|\pi|}}}},
\end{equation}
where $\Std$ is defined in Section~\ref{subsubsec:permutations} of
Chapter~\ref{chap:collections}. For instance
\begin{equation}
    \Coproduct\Par{\BasisF_{42513}} =
    \Par{1, \BasisF_{42513}} +
    \Par{\BasisF_{1}, \BasisF_{2413}} +
    \Par{\BasisF_{21}, \BasisF_{312}} +
    \Par{\BasisF_{213}, \BasisF_{12}} +
    \Par{\BasisF_{3241}, \BasisF_{1}} +
    \Par{\BasisF_{42513}, 1}.
\end{equation}
The product $\Product$ and the coproduct $\Coproduct$ endow $\FQSym$
with a structure of a combinatorial connected graded Hopf bialgebra.
\medbreak

This Hopf bialgebra $\FQSym$ is usually known as the
\Def{Hopf bialgebra of free quasi-symmetric functions}. Indeed, as for
$\Sym$, there is a way to see the elements of $\FQSym$ as noncommutative
series. For this, consider an alphabet
$A := \left\{\Asf_1, \Asf_2, \dots\right\}$ equipped with a total order
$\Ord$ where $1 \leq i \leq j$ implies $\Asf_i \Ord \Asf_j$. Let the
noncommutative series
\begin{equation}
    \BasisF_\sigma(A) :=
    \sum_{\substack{
        u \in A^* \\
        \Std(u) = \sigma^{-1}
    }}
    u,
\end{equation}
of $\K \AAngle{A^*}$ defined  for all $\sigma \in \SymmetricGroup$. For
instance
\begin{subequations}
\begin{equation} \label{equ:example_realization_fqsym_1}
    \BasisF_{312}(\{\Asf_1, \Asf_2, \Asf_3\}) =
    \Asf_2 \Asf_2 \Asf_1 + \Asf_2 \Asf_3 \Asf_1 +
    \Asf_3 \Asf_3 \Asf_1 + \Asf_3 \Asf_3 \Asf_2,
\end{equation}
\begin{equation} \label{equ:example_realization_fqsym_2}
    \BasisF_{132}(\{\Asf_1, \Asf_2, \Asf_3\}) =
    \Asf_1 \Asf_2 \Asf_1 + \Asf_1 \Asf_3 \Asf_1 +
    \Asf_1 \Asf_3 \Asf_2 + \Asf_2 \Asf_3 \Asf_2.
\end{equation}
\end{subequations}
The linear span of all the $\BasisF_\sigma(A)$,
$\sigma \in \SymmetricGroup$, is the space of the free quasi-symmetric
functions on $A$. The associative algebra structure on $\FQSym$ is
compatible with these series in the sense that
\begin{equation} \label{equ:product_fqsym_series}
    \BasisF_\sigma(A) \Conc \BasisF_\nu(A) =
    \Par{\BasisF_\sigma \Product \BasisF_\nu}(A)
\end{equation}
for all $\sigma, \nu \in \SymmetricGroup$, where the product $\Conc$
of the left member of~\eqref{equ:product_fqsym_series} is the usual
product of noncommutative series of $\K \AAngle{A^*}$.
\medbreak

Furthermore, the Hopf bialgebras $\FQSym$ and $\Sym$ are related through
the injective morphism of Hopf bialgebras $\phi : \Sym \to \FQSym$
defined linearly by
\begin{equation}
    \phi\Par{\BasisR_{\bm{\lambda}}} :=
    \sum_{\substack{
        \sigma \in \SymmetricGroup \\
        \Des\Par{\sigma^{-1}} = \Des(\bm{\lambda})
    }}
    \BasisF_\sigma
\end{equation}
for all $\bm{\lambda} \in \ColComp$. For instance,
\begin{equation} \label{equ:example_morphism_from_sym_to_fqsym}
    \phi(\BasisR_{21}) = \BasisF_{312} + \BasisF_{132}.
\end{equation}
Observe, with the help of~\eqref{equ:example_realization_ribbon_sym},
\eqref{equ:example_realization_fqsym_1},
and~\eqref{equ:example_realization_fqsym_2}, in particular
that~\eqref{equ:example_morphism_from_sym_to_fqsym} holds on the
noncommutative series associated with the elements of $\Sym$ and
$\FQSym$, that is,
$\BasisR_{21}(A) = \BasisF_{312}(A) + \BasisF_{132}(A)$.
\medbreak

\SkipTocEntry\section*{Bibliographic notes}

\SkipTocEntry\subsection*{About Incidence algebras}
One of the first apparitions of incidence algebras in combinatorics is
due to Rota~\cite{Rot64}. These structures, associated with any locally
finite poset, provide an abstraction of the principle of
inclusion-exclusion~\cite{Sta11} through their Möbius functions. Indeed,
the usual inclusion-exclusion principle comes from the Möbius function
of the cube poset. Besides, in our exposition, we have presented the
elements of incidence algebras as polynomials of pairs of comparable
elements, but in the literature~\cite{Sta11}, it is most common to see
these elements as maps associating a coefficient with each pair of
comparable elements. These two points of view are therefore equivalent
but the definition of the product of incidence algebras in terms of
polynomials is simpler.
\medbreak

\SkipTocEntry\subsection*{Dendriform algebras}
Dendriform algebras are types of polynomial algebras introduced by
Loday~\cite{Lod01}. These structures can be used as devices to split the
product of an associative algebra into two parts by putting a dendriform
algebra structure onto it. For instance, the dendriform algebra
structure put onto the shuffle algebra (see
Section~\ref{subsubsec:shuffle_algebra}
and~\ref{subsubsec:shuffle_dendriform_algebra}) leads to a recursive
expression for the shuffle product (see~\eqref{equ:shuffle_words_1}
and~\eqref{equ:shuffle_words_2}) known since Ree~\cite{Ree58}.
We invite the reader to take a look
at~\cite{LR98,Agu00,Lod02,Foi07,EMP08,EM09,LV12} for a supplementary
review of properties of dendriform algebras. Besides, in the recent
years, a lot of generalizations of dendriform algebras and their dual
notions were introduced, each of them splitting an associative product
in different ways and in more than two pieces. Tridendriform
algebras~\cite{LR04}, quadri-algebras~\cite{AL04},
ennea-algebras~\cite{Ler04}, $m$-dendriform algebras of
Leroux~\cite{Ler07}, $m$-dendriform algebras of Novelli~\cite{Nov14},
and polydendriform algebras~\cite{Gir16a,Gir16b} are examples of such
structures.
\medbreak

\SkipTocEntry\subsection*{About pre-Lie algebras}
Pre-Lie algebras were introduced by Vinberg~\cite{Vin63} and
Gerstenhaber~\cite{Ger63} independently. These structures appear under
different names in the literature, for instance as Vinberg algebras,
left-symmetric algebras, or chronological algebras. The appellation
pre-Lie algebra is now very natural since, given a pre-Lie algebra
$(\K \Angle{C}, \PreLieProduct)$, the commutator of $\PreLieProduct$
endows $\K \Angle{C}$ with a structure of a Lie algebra. In the context
of combinatorics, several pre-Lie products are defined on combinatorial
spaces by summing over all the ways to compose (in a certain sense) two
combinatorial objects. For this reason, in an intuitive way, pre-Lie
algebras encode the combinatorics of the composition of combinatorial
objects in all possible ways~\cite{Cha08}. Besides, the free objects in
the category of pre-Lie algebras have been described by Chapoton and
Livernet~\cite{CL01}. They have shown that the free pre-Lie algebra
generated by a set $\GeneratingSet$ is the combinatorial space of all
rooted trees whose nodes are labeled on $\GeneratingSet$, and the
product of two such rooted trees is the sum of all the ways to connect
the root of the second tree to a node of the first. Thereby, the pre-Lie
algebra $(\K \Angle{\ColRT}, \PreLieProduct)$ of rooted trees (see
Section~\ref{subsubsec:pre_lie_algebra_rooted_trees}) is the free
pre-Lie algebra generated by a singleton. For more details on pre-Lie
algebras, see~\cite{Man11}.
\medbreak

\SkipTocEntry\subsection*{About bialgebras}
In the field of algebraic combinatorics, many types of bialgebras have
emerged recently. In~\cite{Lod08}, Loday defined the notion of triples
of operads, leading to the constructions of various kinds of bialgebras.
This leads also to the discovery of analogs of the
Poincaré-Birkhoff-Witt and Cartier-Milnor-Moore theorems and rigidity
theorems (see as well~\cite{Cha02} and~\cite{BD18}). Loday defined among
others infinitesimal bialgebras, forming an example of bialgebras having
an associative binary product and a coassociative binary coproduct
satisfying a compatibility relation. Let us describe some other types of
bialgebras that play a role in combinatorics. Bidendriform bialgebras,
introduced by Foissy~\cite{Foi07} are one of these. These bialgebras
have two products satisfying the dendriform relations and two coproducts
satisfying the dual relations of dendriform products, and all of these
together satisfy some compatibility relations. There is a notion of
bidendriform bialgebra structure onto a Hopf bialgebra which leads to a
rigidity theorem in the sense that a Hopf bialgebra admitting a
bidendriform bialgebra structure is self-dual, free as an associative
algebra, and free as a coassociative coalgebra. Moreover,
in~\cite{Foi12}, Foissy considered algebraic structures, named
$\Dup$-$\Dendr$ bialgebras, having two binary products satisfying the
duplicial relations~\cite{BF03,Lod08}, two binary coproducts such that
their dual products satisfy the dendriform relations, and such that
these four (co)products satisfy several compatibility relations. These
structures lead to rigidity theorems in the sense that any
$\Dup$-$\Dendr$ bialgebra is free as a duplicial algebra. In the same
way, Foissy introduced also in~\cite{Foi15} structures named
$\Com$-$\PreLie$ bialgebras, that are spaces with an associative and
commutative binary product, a pre-Lie product, and a binary coproduct
that satisfy compatibility relations. Another interesting example has
been brought by Livernet~\cite{Liv06} wherein bialgebra structures
having a pre-Lie product and a coproduct satisfying the dual relation of
the so-called nonassociative permutative relation have been considered
to construct here again a rigidity theorem.
\medbreak

\SkipTocEntry\subsection*{About Hopf bialgebras}
The Hopf bialgebra $\Sym$ of noncommutative symmetric functions has
been introduced in~\cite{GKLLRT94} as a generalization of the usual
symmetric functions~\cite{Mac15}. This generalization is a consequence
of the fact that there is a surjective morphism from $\Sym$ to the
algebra of symmetric functions. The Hopf bialgebra $\FQSym$ of free
quasi-symmetric functions has been introduced by Malvenuto and
Reutenauer~\cite{MR95} and is sometimes called the Malvenuto-Reutenauer
algebra. Due to its interpretation~\cite{DHT02} as an algebra of
noncommutative series $\BasisF_\sigma(A)$, each element of $\FQSym$ can
be seen as a particular function, whence its name. Other classical
examples of Hopf bialgebras include the Poirier-Reutenauer Hopf
bialgebra of tableaux~\cite{PR95}, also known as the Hopf bialgebra of
free symmetric functions $\FSym$~\cite{DHT02,HNT05}. This Hopf bialgebra
is defined on the combinatorial space of all standard Young tableaux.
The Loday-Ronco Hopf bialgebra~\cite{LR98}, also known as the Hopf
bialgebra of binary search trees $\PBT$~\cite{HNT05} is defined on the
combinatorial space of all binary trees. As other modern examples of
combinatorial spaces endowed with a Hopf bialgebra structure, one can
cite $\WQSym$~\cite{Hiv99} involving packed words, $\PQSym$~\cite{NT07}
involving parking functions, $\Bell$~\cite{Rey07} involving set
partitions, $\Baxter$~\cite{LR12,Gir12a} involving ordered pairs of twin
binary trees, and $\Camb$~\cite{CP17} involving Cambrian trees. The
study of all these structures uses a large set of tools. Indeed, it
relies on algorithms transforming words into combinatorial objects,
congruences of free monoids, partials orders structures and lattices,
and polytopes and their geometric realizations. Besides, a polynomial
realization of a combinatorial Hopf bialgebra $\K \Angle{C}$ consists in
seeing $\K \Angle{C}$ as an algebra of noncommutative series so that its
product is the usual product of series and its coproduct is obtained by
alphabet doubling (see, for instance,~\cite{Hiv03}). In this text, only
the polynomial realizations of $\Sym$ and $\FQSym$ have been detailed,
but all the Hopf bialgebras discussed here have polynomial realizations.
\medbreak


\chapter{Nonsymmetric operads} \label{chap:operads}
This chapter introduces nonsymmetric operads. Our presentation relies on
the framework of graded collections and graded spaces introduced in the
previous chapters. We consider here also set-operads, algebras over
operads, free operads, presentations by generators and relations, Koszul
duality and Koszulity of operads. At the end of the chapter, several
examples of operads on a large family of combinatorial collections are
provided.
\medbreak

\section{Operads as polynomial spaces} \label{sec:operads}
Let us start by posing the main definitions about operads. We first
present the notion of partial composition maps on abstract operators and
then focus on operads and algebras over operads.
\medbreak

\subsection{Composition maps}
Intuitively, the elements of an operad are operations with several
inputs and one output that can be composed. We introduce here the
notion of abstract operator and two ways to compose them by the
so-called partial or full compositions. To consider operads, we shall
expect that partial and full compositions satisfy some relations. One of
the aims of this section is to give an intuition about these relations.
In all this section, $\K \Angle{C}$ is an augmented graded polynomial
space.
\medbreak

\subsubsection{Abstract operators} \label{subsubsec:abstract_operators}
From now, we shall see any homogeneous element $f$ of $\K \Angle{C}$ of
degree $n$ as an operator with $n$ inputs and a single output, called
\Def{abstract operator}. We will use in the context the term \Def{arity}
instead of size or degree, so that the arity of $f$ is $n$. Any abstract
operator is depicted by following the drawing conventions of biproducts
exposed in Section~\ref{subsubsec:biproducts} of
Chapter~\ref{chap:algebra}. Therefore, $f$ is depicted by
\begin{equation}
    \begin{tikzpicture}
        [xscale=.25,yscale=.3,Centering,font=\scriptsize]
        \node[Operator](x)at(0,0){\begin{math}f\end{math}};
        \node(r)at(0,2){};
        \node(x1)at(-3,-2){};
        \node(xn)at(3,-2){};
        \node[below of=x1,node distance=1mm](ex1)
            {\begin{math}1\end{math}};
        \node[below of=xn,node distance=1mm](exn)
            {\begin{math}n\end{math}};
        \draw[Edge](r)--(x);
        \draw[Edge](x)--(x1);
        \draw[Edge](x)--(xn);
        \node[below of=x,node distance=7mm]
            {\begin{math}\dots\end{math}};
    \end{tikzpicture}\,.
\end{equation}
The reason to see the elements of $\K \Angle{C}$ in this way is based on
the fact that we shall consider compositions operations on
$\K \Angle{C}$ subjected to relations that are easy to understand
through this formalism.
\medbreak

\subsubsection{Partial composition maps}
\label{subsubsec:partial_compo_maps}
Let for all $n, m \in \N_{\geq 1}$ and $i \in [n]$ binary products of
the form
\begin{equation}
    \circ_i^{(n, m)} :
    \K \Angle{\BBrack{C(n), C(m)}_\times} \to
    \K \Angle{C}(n + m - 1).
\end{equation}
On abstract operators, these products $\circ_i^{(n, m)}$ behave in the
following way. For any $f \in \K \Angle{C}(n)$ and
$g \in \K \Angle{C}(m)$, $f \circ_i^{(n, m)} g$ is the abstract operator
\begin{equation} \label{equ:partial_compostion_on_operators}
\,.
\end{equation}
In words, $f \circ_i^{(n, m)} g$ is obtained by plugging the output of
$g$ onto the $i$th input of $f$. Observe that since one input of $f$ is
used to make the connection with the output of $g$, the right member
of~\eqref{equ:partial_compostion_on_operators} is of arity $n + m - 1$.
Moreover, observe also that the products $\circ_i^{(n, m)}$ are
concentrated. By a slight abuse of notation, we shall sometimes omit the
$(n, m)$ in the notation of $\circ_i^{(n, m)}$ in order to denote it in
a more concise way by $\circ_i$.
\medbreak

When for any objects $x \in C(n)$, $y \in C(m)$, $z \in C(k)$, and any
integers $i \in [n]$ and $j \in [m]$, the relations
\begin{equation} \label{equ:operad_axiom_assoc_series}
    \Par{x \circ_i y} \circ_{i + j - 1} z = x \circ_i \Par{y \circ_j z}
\end{equation}
hold, we say that the products $\circ_i$ are \Def{series associative}.
To understand this relation, let us consider the abstract operators
expressed by the left and right members
of~\eqref{equ:operad_axiom_assoc_series}. On the one hand, we have
\begin{multline}
    \Par{
\,,
\end{multline}
and on the other,
\begin{multline}
    %
\,.
\end{multline}
We observe that the two obtained abstract operators are the same, as
expressed by~\eqref{equ:operad_axiom_assoc_series}.
\medbreak

Besides, when for any objects $x \in C(n)$, $y \in C(m)$, $z \in C(k)$,
and any integers $i, j \in [n]$ such that $i < j$, the relations
\begin{equation} \label{equ:operad_axiom_assoc_parallel}
    \Par{x \circ_i y} \circ_{j + m - 1} z = \Par{x \circ_j z} \circ_i y
\end{equation}
hold, we say that the products $\circ_i$ are \Def{parallel associative}.
To understand this relation, let us consider the abstract operators
expressed by the left and right members
of~\eqref{equ:operad_axiom_assoc_parallel}. On the one hand, we have
\begin{multline}
    \Par{
\,,
\end{multline}
and on the other,
\begin{multline}
    \Par{
    %
\,.
\end{multline}
We observe that the two obtained abstract operators are the same, as
expressed by~\eqref{equ:operad_axiom_assoc_parallel}.
\medbreak

Finally, when there exists a product of arity $0$ (see
Section~\ref{subsubsec:products_arity_zero} of
Chapter~\ref{chap:algebra}) of the form
\begin{equation}
    \Unit : \K \Angle{\BBrack{}_\times} \to \K \Angle{C}(1),
\end{equation}
such that for any object $x \in C(n)$ and any integer $i \in [n]$ the
relations
\begin{equation} \label{equ:operad_axiom_unit}
    \Unit \circ_1 x = x = x \circ_i \Unit
\end{equation}
hold, we say that the products $\circ_i$ are \Def{unital} and that
$\Unit$ is the \Def{unit}. To understand this relation, let us consider
the abstract operators associated with each member
of~\eqref{equ:operad_axiom_unit}. This leads to the relation
\begin{equation}

\end{equation}
for the left part of~\eqref{equ:operad_axiom_unit} and
\begin{equation}
    %
\end{equation}
for its right part, saying that $\Unit$ is an operator of arity $1$
behaving as the identity map.
\medbreak

When the products $\circ_i$ are series associative, parallel
associative, and unital, the $\circ_i$ are called
\Def{partial composition maps}.
\medbreak

\subsubsection{Full composition maps}
Let for all $n, m_1, \dots, m_n \in \N_{\geq 1}$ products of the form
\begin{equation}
    \circ^{\Par{m_1, \dots, m_n}} :
    \K \Angle{\BBrack{C(n), C\Par{m_1}, \dots, C\Par{m_n}}_\times}
    \to \K \Angle{C}\Par{m_1 + \dots + m_n}.
\end{equation}
On abstract operators, these products $\circ^{(m_1, \dots, m_n)}$ behave
in the following way. For any $f \in \K \Angle{C}(n)$ and
$g_i \in \K \Angle{C}(m_i)$, $i \in [n]$,
\begin{math}
    \circ^{\Par{m_1, \dots, m_n}}\Par{f, g_1, \dots, g_n}
\end{math}
is the abstract operator
\begin{multline} \label{equ:full_compostion_on_operators}
    \circ^{\Par{m_1, \dots, m_n}}\Par{
\,.
\end{multline}
In words,
\begin{math}
    \circ^{\Par{m_1, \dots, m_n}}\Par{f, g_1, \dots, g_n}
\end{math}
is obtained by plugging the outputs of the $g_i$, $i \in [n]$, onto the
$i$th inputs of $f$ simultaneously. Observe that since each input of
$f$ is connected to a $g_i$, the inputs of the right member
of~\eqref{equ:full_compostion_on_operators} are the ones of the
$g_i$, $i \in [n]$, so that its arity is $m_1 + \dots + m_n$.
Moreover, observe also that the products $\circ^{\Par{m_1, \dots, m_n}}$
are concentrated. By a slight abuse of notation, we shall sometimes omit
the $\Par{m_1, \dots, m_n}$ in the notation of
$\circ^{\Par{m_1, \dots, m_n}}$ in order to denote it in a more
concise way by $\circ$. Moreover, we shall write
$f \circ \left[g_1, \dots, g_n\right]$ instead of
$\circ \Par{f, g_1, \dots, g_n}$.
\medbreak

When for any objects $x \in C(n)$, $y_i \in C(m_i)$, $i \in [n]$,
$z_{i, j} \in C(k_{i, j})$, $i \in [n]$, $j \in [m_i]$, the relations
\begin{multline} \label{equ:operad_axiom_full_associativity}
    \Par{x \circ\left[y_1, \dots, y_n\right]} \circ
    \left[z_{1, 1}, \dots, z_{1, m_1},
    \dots,
    z_{n, 1}, \dots, z_{n, m_n}\right] \\
    =
    x \circ \left[
    y_1 \circ \left[z_{1, 1}, \dots, z_{1, m_1}\right],
    \dots
    y_n \circ \left[z_{n, 1}, \dots, z_{n, m_n}\right]
    \right]
\end{multline}
hold, we say that the product $\circ$ is \Def{associative}. To
understand this relation, let us consider the abstract operators
expressed by the left and right members
of~\eqref{equ:operad_axiom_full_associativity}. On the one hand, we
have
\begin{multline}
    \Par{
\,,
\end{multline}
and on the other,
\begin{multline}
    %
\,.
\end{multline}
We observe that the two obtained abstract operators are the same, as
expressed by~\eqref{equ:operad_axiom_full_associativity}.
\medbreak

Besides, when there exists a product $\Unit$ of arity $0$ on
$\K \Angle{C}$ satisfying $\Unit \in \K \Angle{C}(1)$ and such that
for any object $x \in C(n)$ the relations
\begin{equation} \label{equ:operad_axiom_full_unit}
    \Unit \circ [x]
    = x
    = x \circ \left[
    \underbrace{\Unit, \dots, \Unit}_{n \mbox{ \footnotesize terms}}
    \right]
\end{equation}
hold, we say that the products $\circ$ are \Def{unital} and that $\Unit$
is the \Def{unit}. To understand this relation, let us consider the
abstract operators associated with each member
of~\eqref{equ:operad_axiom_full_unit}. This leads to the relation
\begin{equation}

\end{equation}
for the left part of~\eqref{equ:operad_axiom_full_unit} and
\begin{equation}
    %
\end{equation}
for its right part, saying that $\Unit$ is an operator of arity $1$
behaving as the identity map.
\medbreak

When the products $\circ$ are associative and unital, the $\circ$ are
called \Def{full composition maps}.
\medbreak

\subsubsection{Equivalence between partial and full composition maps}
Let $\circ_i$ be partial composition maps on $\K \Angle{C}$. We
construct from the $\circ_i$ the products
$\circ^{\Par{m_1, \dots, m_n}}$,
$n, m_1, \dots, m_n \in \N_{\geq 1}$, on $\K \Angle{C}$ defined
linearly in the following way. For any $x \in C(n)$, $y_i \in C(m_i)$,
$i \in [n]$, let us set
\begin{equation} \label{equ:partial_to_full_composition}
    \circ^{\Par{m_1, \dots, m_n}}
    \Par{x, y_1, \dots, y_n}
    :=
    \Par{\dots \Par{\Par{x \circ_n y_n}
    \circ_{n - 1} y_{n - 1}} \dots } \circ_1 y_1.
\end{equation}
\medbreak

\begin{Proposition} \label{prop:partial_to_full_composition}
    Let $\K \Angle{C}$ be an augmented graded polynomial space endowed
    with partial composition maps $\circ_i$. Then, the products $\circ$
    on $\K \Angle{C}$ defined by~\eqref{equ:partial_to_full_composition}
    are full composition maps.
\end{Proposition}
\medbreak

Conversely, let $\circ$ be full composition maps on $\K \Angle{C}$
and $\Unit$ their unit. We construct from the $\circ$ and $\Unit$ the
products $\circ_i^{(n, m)}$, $n, m \in \N_{\geq 1}$, $i \in [n]$, on
$\K \Angle{C}$ defined linearly in the following way. For any
$x \in C(n)$ and $y \in C(m)$, let us set
\begin{equation} \label{equ:full_to_partial_composition}
    x \circ_i^{(n, m)} y :=
    x \circ \left[
    \underbrace{\Unit, \dots, \Unit}_
    {i - 1 \mbox{ \footnotesize terms}},
    y,
    \underbrace{\Unit, \dots, \Unit}_
    {n - i \mbox{ \footnotesize terms}}
    \right],
\end{equation}
where $\Unit$ is the unit of the $\circ$.
\medbreak

\begin{Proposition} \label{prop:full_to_partial_composition}
    Let $\K \Angle{C}$ be an augmented graded polynomial space endowed
    with full composition maps $\circ$ and their unit $\Unit$. Then,
    the products $\circ_i$ on $\K \Angle{C}$ defined
    by~\eqref{equ:full_to_partial_composition} are partial composition
    maps.
\end{Proposition}
\medbreak

\subsection{Operads} \label{subsec:operads}
Operads are algebraic structures furnishing a formalization of the
notion of abstract operators and their partial and full compositions.
They allow, for instance, to mimic the composition of abstract operators
for various collections of combinatorial objects (words, trees, graphs,
{\em etc.}) and make them behave like operators. We provide here
definitions about these algebraic structures and present set-operads.
\medbreak

\subsubsection{First definitions}
A \Def{nonsymmetric operad} (or an \Def{operad} for short) is a triple
\begin{equation}
    \Par{\K \Angle{C},
    \left\{\circ^{(n, m)}_i : n, m \in \N_{\geq 1}, i \in [n]\right\},
    \Unit}
\end{equation}
where $\K \Angle{C}$ is an augmented graded polynomial space, the
$\circ_i$ are partial composition maps, and $\Unit$ is their unit.
Equivalently, by Propositions~\ref{prop:partial_to_full_composition}
and~\ref{prop:full_to_partial_composition}, an operad is a triple
\begin{equation}
    \Par{\K \Angle{C},
    \left\{\circ^{\Par{m_1, \dots, m_n}} : n \in \N_{\geq 1},
    m_i \in \N_{\geq 1}, i \in [n]\right\},
    \Unit}
\end{equation}
where $\K \Angle{C}$ is an augmented graded polynomial space, the
$\circ$ are full composition maps, and $\Unit$ is their unit. For this
reason, in the sequel, we shall consider operads through partial or full
composition maps indifferently. Moreover, given an operad $\K \Angle{C}$
defined through partial composition maps $\circ_i$, we call
\Def{full composition maps} of $\K \Angle{C}$ the full composition maps
$\circ$ defined in~\eqref{equ:partial_to_full_composition}. Conversely,
if $\K \Angle{C}$ is defined through full composition maps $\circ$, we
call \Def{partial composition maps} of $\K \Angle{C}$ the partial
composition maps $\circ_i$ defined
in~\eqref{equ:full_to_partial_composition}.
\medbreak

Since an operad is a particular polynomial algebra, all the
properties, definitions, and notations about polynomial algebras
exposed in Section~\ref{sec:bialgebras} of Chapter~\ref{chap:algebra}
remain valid for operads (like operad morphisms, suboperads, generating
sets, operad ideals and quotients, {\em etc.}). In particular, to be
more precise, if $\K \Angle{C_1}$ and $\K \Angle{C_2}$ are two operads,
a map $\phi : \K \Angle{C_1} \to \K \Angle{C_2}$ is an
\Def{operad morphism} if $\phi$ is a graded polynomial space morphism,
it sends the unit of $\K \Angle{C_1}$ to the unit of $\K \Angle{C_2}$,
and
\begin{equation}
    \phi\Par{x \circ_i y} = \phi(x) \circ_i \phi(y)
\end{equation}
for all $x \in C_1(n)$, $y \in C_1$, and $i \in [n]$. If $\K \Angle{C}$
is an operad and $\GeneratingSet$ is a subset of $\K \Angle{C}$, the
\Def{operad generated} by $\GeneratingSet$ is the smallest suboperad
$\K \Angle{C}^\GeneratingSet$ of $\K \Angle{C}$ containing
$\GeneratingSet$. A space $\Vca$ included in $\K \Angle{C}$ is an
\Def{operad ideal} of $\K \Angle{C}$ if $x \circ_i f \in \Vca$ and
$f \circ_j y \in \Vca$ for all homogeneous element $f$ of $\Vca$ of
degree $m$, $x \in C(n)$, $y \in C$, $i \in [n]$, and $j \in [m]$. The
\Def{quotient operad} $\K \Angle{C}/_\Vca$ of $\K \Angle{C}$ by $\Vca$
is defined as follows. Let
\begin{math}
    \theta : \K \Angle{C} \to \K \Angle{C}/_\Vca
\end{math}
be the canonical surjection map from $\K \Angle{C}$ to
$\K \Angle{C}/_\Vca$. The space $\K \Angle{C}/_\Vca$ is endowed with the
structure of an operad  through the partial composition maps defined by
\begin{equation} \label{equ:partial_composition_quotient}
    \theta(x) \circ_i \theta(y) := \theta\Par{x \circ_i y}
\end{equation}
for any $x \in C(n)$, $y \in C$, and $i \in [n]$, where the second
occurrence of $\circ_i$ in~\eqref{equ:partial_composition_quotient} is
the partial composition map of~$\K \Angle{C}$.
\medbreak

\subsubsection{Additional definitions}
Let $\K \Angle{C}$ be an operad. An element $f$ of arity $2$ of
$\K \Angle{C}$ is \Def{associative} if $f \circ_1 f - f \circ_2 f = 0$.
If $\K \Angle{C_1}$ and $\K \Angle{C_2}$ are two operads, an
\Def{operad antimorphism} is a graded polynomial space morphism
\begin{math}
    \phi : \K \Angle{C_1} \to \K \Angle{C_2}
\end{math}
sending the unit of $\K \Angle{C_1}$ to the unit of $\K \Angle{C_2}$ and
satisfying
\begin{equation}
    \phi\Par{x \circ_i y} = \phi(x) \circ_{n - i + 1} \phi(y)
\end{equation}
for any $x \in C_1(n)$, $y \in C_1$, and $i \in [n]$.
A \Def{symmetry} of $\K \Angle{C}$ is either an operad
automorphism or an operad antiautomorphism of $\K \Angle{C}$. The set of
all symmetries of $\K \Angle{C}$ forms a group for the map composition,
called \Def{group of symmetries} of~$\K \Angle{C}$.
\medbreak

The \Def{Hadamard product} of a sequence $\K \Angle{C_1}$, \dots,
$\K \Angle{C_p}$, $p \in \N$, of operads is the operad on the graded
polynomial space $\K \Angle{\BBrack{C_1, \dots, C_p}_\Hadamard}$ where
$\BBrack{C_1, \dots, C_p}_\Hadamard$ is the Hadamard product on
collections (see Section~\ref{subsubsec:hadamard_product_collections} of
Chapter~\ref{chap:collections}). The partial composition maps $\circ_i$
of this operad are defined linearly by
\begin{equation} \label{equ:partial_composition_map_Hadamard}
    \Par{x_1, \dots, x_p} \circ_i \Par{y_1, \dots, y_p} :=
    \Par{x_1 \circ_i y_1, \dots, y_p \circ_i y_p},
\end{equation}
for any objects $\Par{x_1, \dots, x_p}$ and $\Par{y_1, \dots, y_p}$ of
$\BBrack{C_1, \dots, C_p}_\Hadamard$, where the occurrences of $\circ_i$
in~\eqref{equ:partial_composition_map_Hadamard} are, from left to right,
the partial composition maps of $\K \Angle{C_1}$, \dots,
$\K \Angle{C_p}$. The unit of the operad
$\K \Angle{\BBrack{C_1, \dots, C_p}_\Hadamard}$ is
$\Par{\Unit_1, \dots, \Unit_p}$ where $\Unit_k$ is the unit of
$\K \Angle{C_k}$ for any $k \in [p]$.
\medbreak

\subsubsection{Set-operads} \label{subsubsec:set_operads}
An operad $\K \Angle{C}$ is a \Def{set-operad} if $C$ is a set-basis
(see Section~\ref{subsubsec:products_properties} of
Chapter~\ref{chap:algebra}) with respect to all the partial composition
maps of $\K \Angle{C}$, and the unit $\Unit$ is an object of $C$. This
implies in particular that for all $x \in C(n)$, $y \in C(m)$, and
$i \in [n]$, $x \circ_i y$ is an object of $C$. To study a set-operad
$\K \Angle{C}$, it is in some cases convenient to forget about its
linear structure and see its partial composition maps $\circ_i$ as
set-theoretic maps (see Section~\ref{subsubsec:set_algebras} of
Chapter~\ref{chap:algebra}). Let us consider now that $\K \Angle{C}$ is
a set-operad and let us review some properties that operads of this kind
of operad can satisfy.
\medbreak

A collection of maps
\begin{equation}
    \tau_n : C(n) \to [n]
\end{equation}
where $n \in \N_{\geq 1}$ are \Def{root maps} of $\K \Angle{C}$ if, for
any $x \in C(n)$, $y \in C(m)$, and $i \in [n]$,
\begin{equation}
    \tau_{n + m - 1}\Par{x \circ_i y} =
    \begin{cases}
        \tau_n(x) + m - 1 & \mbox{if } i \leq \tau_n(x) - 1, \\
        \tau_n(x) + \tau_m(y) - 1 & \mbox{if } i = \tau_n(x), \\
        \tau_n(x) & \mbox{otherwise (} i \geq \tau_n(x)\mbox{)}.
    \end{cases}
\end{equation}
In this case, we say that $\K \Angle{C}$ is a \Def{rooted operad} with
respect to the maps $\tau_n$, $n \in \N_{\geq 1}$. More intuitively,
this property says that in a rooted operad, each object $x$ of $C(n)$
has a particular input $\tau_n(x)$ which is preserved by the partial
composition maps.
\medbreak

Besides, let for any $y \in C(m)$, $n \in \N_{\geq 1}$, and $i \in [n]$
the maps
\begin{equation}
    \Basic^{(n, y)}_i : C(n) \to C(n + m - 1),
\end{equation}
defined for any $x \in C(n)$ by
\begin{equation}
    \Basic^{(n, y)}_i(x) := x \circ_i y.
\end{equation}
When for all $y \in C(m)$, $n \in \N_{\geq 1}$, and $i \in [n]$, all the
maps $\Basic^{(n, y)}_i$ are injective, $\K \Angle{C}$ is a
\Def{basic operad}. More intuitively, this property says that in a basic
operad, one can recover the object $x$ from $x \circ_i y$ with the
knowledge of $i$ and~$y$.
\medbreak

\subsection{Algebras over operads} \label{subsec:algebras_over_operads}
One of the main interests of the theory of operads is that each operad
encodes a category of type of algebras. In this way, by studying a
single operad, it is possible to get general results about all the
algebras of the encoded category. Moreover, morphisms between operads
offer general constructions to, given an algebra of one type, obtain
an algebra of another type. We explain here all these notions and expose
also the concept of free algebras over operads.
\medbreak

\subsubsection{From operads to types of algebras}
Any operad $\K \Angle{C}$ encodes a type of polynomial algebras (see
Section~\ref{sec:types_bialgebras} of Chapter~\ref{chap:algebra}) called
\Def{algebras over $\K \Angle{C}$} (or, for short,
\Def{$\K \Angle{C}$-algebras}). A $\K \Angle{C}$-algebra is a (not
necessarily graded) polynomial space $\K \Angle{D}$, where $D$ is a
collection, which is endowed for all $n \in \N_{\geq 1}$ with linear
maps
\begin{equation} \label{equ:action_algebra_over_operad}
    \Action_n : \K \Angle{\BBrack{C(n), \List_{\{n\}}(D)}_\times}
    \to \K \Angle{D}
\end{equation}
satisfying the relations imposed by the operad structure of
$\K \Angle{C}$, that are,
\begin{subequations}
for all $x \in C(n)$, $y \in C(m)$, $i \in [n]$, and
\begin{math}
    \Par{a_1, \dots, a_{n + m - 1}} \in \List_{\{n + m - 1\}}(D),
\end{math}
\begin{multline} \label{equ:algebra_over_operad}
    \Action_{n + m - 1}\Par{x \circ_i y,
    \Par{a_1, \dots, a_{n + m - 1}}} = \\
    \Action_n\Par{x, \Par{a_1, \dots,
        a_{i - 1},
        \Action_m\Par{y, \Par{a_i, \dots, a_{i + m - 1}}},
        a_{i + m}, \dots, a_{n + m - 1}}},
\end{multline}
and for all $a_1 \in D$,
\begin{equation} \label{equ:algebra_over_operad_unit}
    \Action_1\Par{\Unit, \Par{a_1}} = a_1.
\end{equation}
\end{subequations}
In other words, any object $x$ of $C$ of arity $n$ plays the role of a
complete product (in the sense of Section~\ref{subsubsec:biproducts} of
Chapter~\ref{chap:algebra}) of the form
\begin{equation}
    x : \K \Angle{\List_{\{n\}}(D)} \to \K \Angle{D},
\end{equation}
defined, for any
\begin{math}
    \Par{a_1, \dots, a_n} \in \List_{\{n\}}(D)
\end{math}
by
\begin{equation}
    x\Par{a_1, \dots, a_n} :=
    \Action_n\Par{x, \Par{a_1, \dots, a_n}}.
\end{equation}
Under this point of view, Relation~\eqref{equ:algebra_over_operad} reads
as
\begin{equation}
    \begin{tikzpicture}
        [xscale=.25,yscale=.3,Centering,font=\scriptsize]
        \node[OperatorColorB](x)at(0,0)
            {\begin{math}x \circ_i y\end{math}};
        \node(r)at(0,2){};
        \node(x1)at(-3,-2){};
        \node(xn)at(3,-2){};
        \node[below of=x1,node distance=1mm](ex1)
            {\begin{math}a_1\end{math}};
        \node[below of=xn,node distance=1mm](exn)
            {\begin{math}a_{n + m - 1}\end{math}};
        \draw[Edge](r)--(x);
        \draw[Edge](x)--(x1);
        \draw[Edge](x)--(xn);
        \node[below of=x,node distance=7mm]
            {\begin{math}\dots\end{math}};
    \end{tikzpicture}
    \enspace = \enspace
    \begin{tikzpicture}
        [xscale=.5,yscale=.4,Centering,font=\scriptsize]
        \node[Operator](x)at(0,0){\begin{math}x\end{math}};
        \node(r)at(0,1.5){};
        \node(x1)at(-3,-2){};
        \node(xn)at(3,-2){};
        \node[below of=x1,node distance=1mm](ex1)
            {\begin{math}a_1\end{math}};
        \node[below of=xn,node distance=1mm](exn)
            {\begin{math}a_{n + m - 1}\end{math}};
        \node[right of=ex1,node distance=8mm]
            {\begin{math}\dots\end{math}};
        \node[left of=exn,node distance=9mm]
            {\begin{math}\dots\end{math}};
        \draw[Edge](r)--(x);
        \draw[Edge](x)--(x1);
        \draw[Edge](x)--(xn);
        \node[OperatorColorA](y)at(0,-2.5){\begin{math}y\end{math}};
        \node(y1)at(-1.6,-4.5){};
        \node(yn)at(1.6,-4.5){};
        \node[below of=y1,node distance=1mm](ey1)
            {\begin{math}a_i\end{math}};
        \node[below of=yn,node distance=1mm](eyn)
            {\begin{math}a_{i + m - 1}\end{math}};
        \draw[Edge](y)--(y1);
        \draw[Edge](y)--(yn);
        \node[below of=y,node distance=9mm]
            {\begin{math}\dots\end{math}};
        \draw[Edge](x)--(y);
    \end{tikzpicture}\,,
\end{equation}
and Relation~\eqref{equ:algebra_over_operad_unit} says that $\Unit$ is
the identity map on $\K \Angle{D}$. From now, to define an algebra
$\K \Angle{D}$ over an operad $\K \Angle{C}$, we shall simply describe
how the objects $x$ of $C$ behave as linear products on~$\K \Angle{D}$.
\medbreak

Observe that by~\eqref{equ:algebra_over_operad}, any associative
element of $\K \Angle{C}$ gives rise to an associative operation
on~$\K \Angle{D}$ (details are given in further
Section~\ref{subsubsec:associative_operad}).
\medbreak

\subsubsection{Categories of algebras}
The class of all the $\K \Angle{C}$-algebras forms a category, called
\Def{category of $\K \Angle{C}$-algebras}, wherein morphisms between
$\K \Angle{C}$-algebras are polynomial algebra morphisms (see
Section~\ref{subsubsec:polynomial_bialgebras} of
Chapter~\ref{chap:algebra}). More concretely, if $\K \Angle{D_1}$ and
$\K \Angle{D_2}$ are two $\K \Angle{C}$-algebras such that $D_1$
and $D_2$ are two collections on a same set of indexes, a map
\begin{equation}
    \phi : \K \Angle{D_1} \to \K \Angle{D_2}
\end{equation}
is a $\K \Angle{C}$-algebra morphism if $\phi$ is a polynomial space
morphism and satisfies
\begin{equation}
    \phi\Par{x\Par{a_1, \dots, a_n}}
    =
    x\Par{\phi\Par{a_1}, \dots, \phi\Par{a_n}}
\end{equation}
for all $n \in \N_{\geq 1}$, $x \in C(n)$, and
$a_1, \dots, a_n \in D_1$.
\medbreak

\begin{Proposition} \label{prop:operad_morphisms_and_algebras}
    Let $\K \Angle{C_1}$ and $\K \Angle{C_2}$ be two operads and
    $\phi : \K \Angle{C_1} \to \K \Angle{C_2}$ be an operad morphism.
    Then, if $\K \Angle{D}$ is a $\K \Angle{C_2}$-algebra, by setting
    for any $n \in \N_{\geq 1}$, $x \in C_1(n)$, and
    \begin{math}
        a_1, \dots, a_n \in D,
    \end{math}
    \begin{equation}
        x\Par{a_1, \dots, a_n} :=
        \Par{\phi(x)}\Par{a_1, \dots, a_n},
    \end{equation}
    the space $\K \Angle{D}$ becomes a $\K \Angle{C_1}$-algebra.
\end{Proposition}
\medbreak

Proposition~\ref{prop:operad_morphisms_and_algebras} brings a way to
construct $\K \Angle{C_1}$-algebras from both a $\K \Angle{C_2}$-algebra
and an operad morphism between $\K \Angle{C_1}$ and $\K \Angle{C_2}$.
Some classical constructions of algebras come within this framework. For
instance, it is well-known that any dendriform algebra leads to an
associative algebra by considering the product obtained by summing the
two dendriform products (see Section~\ref{subsec:dendriform_algebras}
of Chapter~\ref{chap:algebra}). This construction is in fact the
consequence of an operad morphism from the associative operad to the
dendriform operad (see the forthcoming
Sections~\ref{subsubsec:associative_operad}
and~\ref{subsubsec:dendriform_operad}).
\medbreak

\subsubsection{Free algebras over operads}
Let us now describe particular algebras over operads. Let $S$ be a
graded collection and let us consider the graded space
\begin{equation}
    \K \Angle{C}^{(S)} := \K \Angle{C \Compo S}
\end{equation}
where $\Compo$ is the composition product of graded collections (see
Section~\ref{subsubsec:composition_collections} of
Chapter~\ref{chap:collections}). Let us endow $\K \Angle{C}^{(S)}$ with
the products $x \in C(n)$, $n \in \N_{\geq 1}$, defined linearly,
for all objects
\begin{math}
    \Par{y_i, \Par{s_{i, 1}, \dots, s_{i, |y_i|}}}
\end{math}
of $\Par{C \Compo S}\Par{m_i}$, $m_i \in \N_{\geq 1}$, $i \in [n]$, by
\begin{multline} \label{equ:products_free_algebra_over_operads}
    x\Par{
    \Par{y_1, \Par{s_{1, 1}, \dots, s_{1, |y_1|}}},
    \dots,
    \Par{y_n, \Par{s_{n, 1}, \dots, s_{n, |y_n|}}}
    } \\
    :=
    \Par{x \circ \left[y_1, \dots, y_n\right],
    \Par{
    s_{1, 1}, \dots, s_{1, |y_1|},
    \dots,
    s_{n, 1}, \dots, s_{n, |y_n|}
    }}.
\end{multline}
\medbreak

\begin{Proposition} \label{prop:free_algebra_over_operads}
    Let $\K \Angle{C}$ be an operad and $S$ be a graded collection.
    Then, the space $\K \Angle{C}^{(S)}$ endowed with the linear
    products $x \in C$ defined
    by~\eqref{equ:products_free_algebra_over_operads} is a
    $\K \Angle{C}$-algebra.
\end{Proposition}
\medbreak

Let now
\begin{equation}
    \iota : S \to \K \Angle{C}^{(S)}
\end{equation}
be the map defined for any $s \in S$ by $\iota(s) := (\Unit, (s))$,
where as usual, $\Unit$ denotes the unit of $\K \Angle{C}$. This map can
be seen as an inclusion of $S$ into $\K \Angle{C}^{(S)}$.
\medbreak

\begin{Theorem} \label{thm:free_algebra_over_operads}
    Let $\K \Angle{C}$ be an operad and $S$ be a graded collection.
    Then, $\K \Angle{C}^{(S)}$ is the unique $\K \Angle{C}$-algebra
    (up to isomorphism) such that for any graded $\K \Angle{C}$-algebra
    $\K \Angle{D}$ and any map $f : S \to \K \Angle{D}$ respecting the
    sizes, there exists a unique $\K \Angle{C}$-algebra morphism
    $\phi : \K \Angle{C}^{(S)} \to \K \Angle{D}$ such that
    $f = \phi \circ \iota$.
\end{Theorem}
\medbreak

Theorem~\ref{thm:free_algebra_over_operads} provides the fact that
$\K \Angle{C}^{(S)}$ satisfies a universality property saying (with
the notations of the statement of the theorem) that the diagram
\begin{equation}
    \begin{tikzpicture}[xscale=1.4,yscale=1.2,Centering]
        \node(S)at(0,0){\begin{math}S\end{math}};
        \node(A)at(2,0){\begin{math}\K \Angle{D}\end{math}};
        \node(FA)at(0,-2)
            {\begin{math}\K \Angle{C}^{(S)}\end{math}};
        \draw[Map](S)--(A)node[midway,above]{\begin{math}f\end{math}};
        \draw[Injection](S)--(FA)node[midway,left]
            {\begin{math}\iota\end{math}};
        \draw[Map,dashed](FA)--(A)node[midway,right]
            {\begin{math}\phi\end{math}};
    \end{tikzpicture}
\end{equation}
commutes and therefore, that $\K \Angle{C}^{(S)}$ is a free object in
the category of the $\K \Angle{C}$-algebras. For this reason, we call
$\K \Angle{C}^{(S)}$ the \Def{free $\K \Angle{C}$-algebra over $S$}.
\medbreak

When $\AtomElement$ is an atom, the free $\K \Angle{C}$-algebra
$\K \Angle{C}^{(\{\AtomElement\})}$ admits the following description.
First, since $\AtomElement$ is of size $1$, by
Relation~\eqref{equ:relation_unit_composition_collection} of
Chapter~\ref{chap:collections}, $\K \Angle{C}^{(\{\AtomElement\})}$ is
isomorphic (as a polynomial space) to $\K \Angle{C}$ and each basis
element $\Par{x, \AtomElement^{|x|}}$ of
$\K \Angle{C}^{(\{\AtomElement\})}$ can be identified with the basis
element $x \in C$ of $\K \Angle{C}$. Moreover,
by~\eqref{equ:products_free_algebra_over_operads}, the operations
$x \in C(n)$, $n \in \N_{\geq 1}$, of $\K \Angle{C}$ satisfy, for any
$y_1, \dots, y_n \in C$,
\begin{equation}
    x\Par{y_1, \dots y_n} := x \circ \left[y_1, \dots y_n\right].
\end{equation}
\medbreak

\section{Free operads, presentations, and Koszulity}
\label{sec:free_operads}
Free operads are intuitively operads wherein partial composition maps
satisfy only the required relations. These operads can be realized as
spaces of syntax trees. We present here some general notions for operads
related to free operads: presentations by generators and relations,
Koszul duality, and Koszulity for binary and quadratic operads.
\medbreak

\subsection{Free operads} \label{subsec:free_operads}
Let us start by defining free operads, exposing the universality
property they satisfy, and a notion of factorization of the elements of
a operad relying on free operads.
\medbreak

\subsubsection{Operads of syntax trees}
Let $\GeneratingSet$ be an augmented graded collection. The
\Def{free operad} over $\GeneratingSet$ is the operad
\begin{equation}
    \FreeOperad(\GeneratingSet)
    :=
    \K \Angle{\ColST_\Leaf^\GeneratingSet},
\end{equation}
where $\ColST_\Leaf^\GeneratingSet$ is the graded collection of all the
$\GeneratingSet$-syntax trees (see
Section~\ref{subsec:graded_collection_syntax_trees} of
Chapter~\ref{chap:trees}). The space $\FreeOperad(\GeneratingSet)$ is
endowed with the linearizations of the partial grafting operations
$\circ_i$, $i \in \N_{\geq 1}$, defined in
Section~\ref{subsubsec:partial_grafting_syntax_trees} of
Chapter~\ref{chap:trees}. The unit of
$\FreeOperad(\GeneratingSet)$ is the only $\GeneratingSet$-syntax tree
$\Leaf$ of arity~$1$ and degree~$0$.
\medbreak

Recall, as defined in
Section~\ref{subsec:graded_collection_syntax_trees} of
Chapter~\ref{chap:trees}, that for any $x \in \GeneratingSet$,
$\Corolla(x)$ is the corolla labeled by~$x$. We shall from now see
$\Corolla$ as a map
\begin{equation}
    \Corolla : \GeneratingSet \to \FreeOperad(\GeneratingSet)
\end{equation}
called \Def{inclusion map}. In the sequel, if required by the context,
we shall implicitly see any element $x$ of $\GeneratingSet$ as the
corolla $\Corolla(x)$ of $\FreeOperad(\GeneratingSet)$. For instance,
when $x \in \GeneratingSet(n)$ and $y \in \GeneratingSet$, we shall
simply denote by $x \circ_i y$ the syntax tree
$\Corolla(x) \circ_i \Corolla(y)$ for any $i \in [n]$.
\medbreak

Free operads satisfy the following universality property. The free
operad $\FreeOperad(\GeneratingSet)$ is the unique operad (up to
isomorphism) such that for any operad $\K \Angle{C}$ and any map
$f : \GeneratingSet \to \K \Angle{C}$ respecting the arities, there
exists a unique operad morphism
$\phi : \FreeOperad(\GeneratingSet) \to \K \Angle{C}$ such that
$f = \phi \circ \Corolla$. In other terms, the diagram
\begin{equation}
    \begin{tikzpicture}[xscale=1.4,yscale=1.2,Centering]
        \node(G)at(0,0){\begin{math}\GeneratingSet\end{math}};
        \node(O)at(2,0){\begin{math}\K \Angle{C}\end{math}};
        \node(AG)at(0,-2)
            {\begin{math}\FreeOperad(\GeneratingSet)\end{math}};
        \draw[Map](G)--(O)node[midway,above]{\begin{math}f\end{math}};
        \draw[Injection](G)--(AG)node[midway,left]
            {\begin{math}\Corolla\end{math}};
        \draw[Map,dashed](AG)--(O)node[midway,right]
            {\begin{math}\phi\end{math}};
    \end{tikzpicture}
\end{equation}
commutes.
\medbreak

\subsubsection{Evaluations and treelike expressions}
Let $\K \Angle{C}$ be an operad. Since $C$ is an augmented graded
collection, one can consider the free operad $\FreeOperad(C)$ of the
$C$-syntax trees. By definition, the fundamental basis of
$\FreeOperad(C)$ is the set of the syntax trees on $C$. The
\Def{evaluation map} of $\K \Angle{C}$ is the map
\begin{equation}
    \Eval : \FreeOperad(C) \to \K \Angle{C}
\end{equation}
defined linearly by induction, for any $C$-syntax tree $\Tfr$, by
\begin{equation}
    \Eval(\Tfr) :=
    \begin{cases}
        \Unit \in \K \Angle{C}
            & \mbox{if } \Tfr = \Leaf, \\
        \omega_\Tfr(\epsilon) \circ
        \left[\Eval\Par{\Tfr_1}, \dots,
            \Eval\Par{\Tfr_k}\right]
            & \mbox{otherwise},
    \end{cases}
\end{equation}
where the $\circ$ are the full composition maps of $\K \Angle{C}$,
$\omega_\Tfr(\epsilon)$ is the label of the root of $\Tfr$, and $k$ is
the root arity of $\Tfr$. This map $\Eval$ is the unique surjective
operad morphism from $\FreeOperad(C)$ to $\K \Angle{C}$ satisfying
$\Eval(\Corolla(x)) = x$ for all $x \in C$.
\medbreak

For any element $f$ of $\K \Angle{C}$, a \Def{treelike expression} of
$f$ is an element $f'$ of $\FreeOperad(C)$ such that $\Eval(f') = f$. A
treelike expression can therefore be thought as a factorization in an
operad.
\medbreak

\subsection{Presentations by generators and relations}
\label{subsec:presentations}
To understand the structure of an operad, it is in most of the cases
fruitful to see it as a quotient of a free operad, leading to the notion
of presentation by generators and relations. Indeed, by comparing the
presentations of two operads, it is most of the time easy to construct
injective or surjective morphisms between them. Moreover, knowing a
presentation of an operad facilitates the description of the category of
the algebras it encodes. We present here a tool coming from the theory
of rewrite systems on syntax trees to establish presentations.
\medbreak

\subsubsection{Presentations}
A \Def{presentation} of an operad $\K \Angle{C}$ consists in a pair
$(\GeneratingSet, \RelationSpace)$ such that $\GeneratingSet$ is an
augmented graded collection, $\RelationSpace$ is a subspace of
$\FreeOperad(\GeneratingSet)$ and
\begin{equation}
    \K \Angle{C} \simeq
    \FreeOperad(\GeneratingSet)/_{\Angle{\RelationSpace}}
\end{equation}
where $\Angle{\RelationSpace}$ is the operad ideal of
$\FreeOperad(\GeneratingSet)$ generated by $\RelationSpace$. We call
$\GeneratingSet$ the \Def{set of generators} and $\RelationSpace$ the
\Def{space of relations} of~$\K \Angle{C}$.
\medbreak

We say that a presentation $(\GeneratingSet, \RelationSpace)$ of
$\K \Angle{C}$ is \Def{quadratic} if $\RelationSpace$ is a homogeneous
subspace of $\FreeOperad(\GeneratingSet)$ consisting in syntax trees of
degree $2$. Besides, we say that $(\GeneratingSet, \RelationSpace)$  is
\Def{binary} if $\GeneratingSet$ has only elements of size (arity)~$2$.
By extension, we say also that $\K \Angle{C}$ is \Def{quadratic} (resp.
\Def{binary}) if it admits a quadratic (resp. binary) presentation.
\medbreak

There is a close link between operad ideals, closures of rewrite rules
of syntax trees (see Section~\ref{subsubsec:rewrite_systems} of
Chapter~\ref{chap:trees}), and spaces induced by rewrite rules (see
Section~\ref{subsubsec:rew_quotient_space} of
Chapter~\ref{chap:algebra}) brought by the following statement.
\medbreak

\begin{Proposition} \label{prop:closure_rewrite_rule_operad_ideal}
    Let $\GeneratingSet$ be an augmented graded collection and
    $\Par{\ColST^\GeneratingSet_\Leaf, \Rew}$ be a rewrite system. Then,
    \begin{equation}
        \Angle{\RelationSpaceRewriteRule
            _{\Par{\ColST^\GeneratingSet_\Leaf, \Rew}}}
        =
        \RelationSpaceRewriteRule
            _{\Par{\ColST^\GeneratingSet_\Leaf, \RewContext}}.
    \end{equation}
\end{Proposition}
\medbreak

In the statement of
Proposition~\ref{prop:closure_rewrite_rule_operad_ideal}, recall that
\begin{math}
    \RelationSpaceRewriteRule
        _{\Par{\ColST^\GeneratingSet_\Leaf, \Rew}}
\end{math}
denote the space induced by $\Par{\ColST^\GeneratingSet_\Leaf, \Rew}$
and
\begin{math}
    \RelationSpaceRewriteRule
        _{\Par{\ColST^\GeneratingSet_\Leaf, \RewContext}}
\end{math}
denotes the space induced by the closure of
$\Par{\ColST^\GeneratingSet_\Leaf, \Rew}$.
\medbreak

\subsubsection{Proving presentations through rewrite systems}
\label{subsubsec:presentation_rewrite_systems}
Rewrite systems on syntax trees (see Section~\ref{sec:rewrite_systems}
of Chapter~\ref{chap:collections} and
Section~\ref{subsec:rewrite_rules_syntax_trees} of
Chapter~\ref{chap:trees}) are powerful tools to prove that a given
operad admits a conjectured presentation. The following result provides
a way to establish presentations of operads.
\medbreak

\begin{Theorem} \label{thm:presentation_operads}
    Let $\K \Angle{C}$ be an operad, $\GeneratingSet$ be a subcollection
    of $C$, and $\RelationSpace$ be a subspace of
    $\FreeOperad(\GeneratingSet)$ of syntax trees of degrees $2$ or
    more. If
    \begin{enumerate}[label={(\it\roman*)}]
        \item \label{item:presentation_operads_1}
        the collection $\GeneratingSet$ is a generating set of
        $\K \Angle{C}$ as an operad;
        \item \label{item:presentation_operads_2}
        for any $f \in \RelationSpace$, $\Eval(f) = 0$;
        \item \label{item:presentation_operads_3}
        there exists a rewrite system
        $\Par{\ColST^\GeneratingSet_\Leaf, \Rew}$ being an
        orientation of $\RelationSpace$, such that its closure
        $\Par{\ColST^\GeneratingSet_\Leaf, \RewContext}$ is
        convergent, and its set of normal forms
        \begin{math}
            \NormalForms
                _{\Par{\ColST^\GeneratingSet_\Leaf, \RewContext}}
        \end{math}
        is isomorphic (as a graded collection) to $C$,
    \end{enumerate}
    then $(\GeneratingSet, \RelationSpace)$ is a presentation of
    $\K \Angle{C}$.
\end{Theorem}
\medbreak

In practice, there are at least two ways to use
Theorem~\ref{thm:presentation_operads} to establish a presentation of an
operad $\K \Angle{C}$. The first one is the most obvious: it consists
first in finding a generating set $\GeneratingSet$ of $\K \Angle{C}$,
then conjecturing (likely with the help of the computer) a space of
relations $\RelationSpace$ and a rewrite system
$\Par{\ColST^\GeneratingSet_\Leaf, \Rew}$ such that all
conditions~\ref{item:presentation_operads_1},
\ref{item:presentation_operads_2}, and~\ref{item:presentation_operads_3}
are satisfied. This can be technical (especially to prove that the
closure $\Par{\ColST^\GeneratingSet_\Leaf, \RewContext}$ is convergent),
and relies heavily on computer exploration. The second way requires as a
prerequisite that $\K \Angle{C}$ is combinatorial (and thus, all its
homogeneous components are finite dimensional). In this case, we need
here also to find a generating set $\GeneratingSet$ of $\K \Angle{C}$, a
space of relations $\RelationSpace$ and a rewrite system
$\Par{\ColST^\GeneratingSet_\Leaf, \Rew}$ such
that~\ref{item:presentation_operads_1},
and~\ref{item:presentation_operads_2} hold, and that $C$ and
\begin{math}
    \NormalForms
        _{\Par{\ColST^\GeneratingSet_\Leaf, \RewContext}}
\end{math}
are isomorphic as graded combinatorial collections. The difference with
the first way occurs for~\ref{item:presentation_operads_3}: it is now
sufficient to prove that
$\Par{\ColST^\GeneratingSet_\Leaf, \RewContext}$ is
terminating (and not necessarily convergent). Indeed, if
$\Par{\ColST^\GeneratingSet_\Leaf, \RewContext}$ is terminating, since
$\K \Angle{C}$ is combinatorial,
\begin{equation} \label{equ:inequality_dim_presentation_operad}
    \dim \K \Angle{C}(n) =
    \# \NormalForms
        _{\Par{\ColST^\GeneratingSet_\Leaf, \RewContext}}(n)
    \geq
    \dim \FreeOperad(\GeneratingSet)/_{\Angle{\RelationSpace}}(n)
\end{equation}
for all $n \in \N_{\geq 1}$. The inequality
of~\eqref{equ:inequality_dim_presentation_operad} comes from the fact
that, since we do not know if
$\Par{\ColST^\GeneratingSet_\Leaf, \RewContext}$ is confluent, it can
have more normal forms of arity $n$ than the dimension of
$\FreeOperad(\GeneratingSet)/_{\Angle{\RelationSpace}}$ in arity $n$. It
follows from~\eqref{equ:inequality_dim_presentation_operad}, by using
straightforward arguments, that there is an operad isomorphism from
$\FreeOperad(\GeneratingSet)/_{\Angle{\RelationSpace}}$
to~$\K \Angle{C}$.
\medbreak

\subsubsection{Realizations and presentations}
Defining an operad can be done in at least two different ways. The first
way consists in describing explicitly an augmented graded polynomial
space $\K \Angle{C}$ together with algorithms for the computation of the
partial composition maps $\circ_i$ involving objects of $C$. This
concrete manner provides a \Def{realization} of an operad. The second
way consists in defining an operad through its presentation
$(\GeneratingSet, \RelationSpace)$, that is, an operad which is by
definition isomorphic to
$\FreeOperad(\GeneratingSet)/_{\Angle{\RelationSpace}}$. This manner
provides only an abstract definition of an operad since nor the
underlying space neither the partial composition maps of the operad are
known at this stage. In practice, to fully understand an operad, it is
most of the time useful to know one of its realizations and one of its
presentations.
\medbreak

\subsubsection{From presentations to types of algebras}
The knowledge of a presentation $(\GeneratingSet, \RelationSpace)$ of an
operad $\K \Angle{C}$ leads to a simple description of the category of
$\K \Angle{C}$-algebras. Indeed, the symbols of $\GeneratingSet$ specify
the products of the algebras of the category, and the relations of
$\RelationSpace$ specify the relations between these products. This
relies on the fact that since $\GeneratingSet$ is a generating set of
$\K \Angle{C}$, any $f \in \K \Angle{C}(n)$ writes as an expression
involving the linear structure of $\K \Angle{C}$, its partial
composition maps, and elements of $\GeneratingSet$. Now, for any
$\K \Angle{C}$-algebra $\K \Angle{D}$,
Relation~\eqref{equ:algebra_over_operad} implies that one can write any
$f\Par{a_1, \dots, a_n}$, $a_1, \dots, a_n \in \K \Angle{D}$, in
terms of a linear combination of compositions of products of
$\GeneratingSet$. Hence, the knowledge of the behavior of each product
$x \in \GeneratingSet$ on $\K \Angle{D}$ is enough to know
the behavior of any product $f \in \K \Angle{C}(n)$,
$n \in \N_{\geq 1}$, on $\K \Angle{D}$. Moreover, the relations between
the products of $\GeneratingSet$ satisfied by any
$\K \Angle{C}$-algebra are encoded by the elements of $\RelationSpace$.
Indeed, each element $f$ of $\RelationSpace$ is a formal sum of
$\GeneratingSet$-syntax trees which is, by definition, equated with~$0$
(that is, $\Eval(f) = 0$).
\medbreak

\subsection{Koszulity} \label{subsec:koszulity}
Given a presentation of a quadratic and binary operad, one can compute
a presentation of another operad, namely of its Koszul dual. This kind
of duality has a close connection with the concept of Koszulity of
operads which is defined originally in an algebraic way. This property
on operads can be rephrased in terms of properties of orientations of
spaces of relations and rewrite systems. As a concrete consequence of
Koszulity, given a combinatorial Koszul operad, its Hilbert series and
the one of its Koszul dual are inverse (in a certain sense) one of the
other.
\medbreak

\subsubsection{Koszul duality} \label{subsubsec:koszul_duality}
Let $\K \Angle{C}$ be an operad admitting a binary and quadratic
presentation $(\GeneratingSet, \RelationSpace)$ where $\GeneratingSet$
is finite, the \Def{Koszul dual} of $\K \Angle{C}$ is the operad
$\K \Angle{C}^!$, isomorphic to the operad admitting the presentation
$\Par{\GeneratingSet, \RelationSpace^\perp}$ where
$\RelationSpace^\perp$ is the annihilator of
$\RelationSpace$ in $\FreeOperad(\GeneratingSet)$ with respect to the
linear map
\begin{equation}
    \Angle{-} :
    \K \Angle{
    \BBrack{\ColST_\Leaf^\GeneratingSet(3),
    \ColST_\Leaf^\GeneratingSet(3)}_\times}
    \to \K
\end{equation}
linearly defined, for all $x, x', y, y' \in \GeneratingSet(2)$, by
\begin{equation} \label{equ:bilinear_map_koszul}
    \Angle{\Par{x \circ_i y, x' \circ_{i'} y'}} :=
    \begin{cases}
        1 & \mbox{if }
            x = x', y = y', \mbox{ and } i = i' = 1, \\
        -1 & \mbox{if }
            x = x', y = y', \mbox{ and } i = i' = 2, \\
        0 & \mbox{otherwise}.
    \end{cases}
\end{equation}
To not overload the notation, we write $\Angle{\Tfr, \Sfr}$ instead of
$\Angle{(\Tfr, \Sfr)}$ for any pair $(\Tfr, \Sfr)$ of
$\GeneratingSet$-syntax trees of arity $3$ and degree $2$.
\medbreak

Then, with knowledge of a presentation of $\K \Angle{C}$, one can
compute a presentation of~$\K \Angle{C}^!$.
\medbreak

\subsubsection{Koszulity} \label{subsubsec:koszulity}
An operad $\K \Angle{C}$ admitting a quadratic presentation is
\Def{Koszul} if its Koszul complex is acyclic. Furthermore, when
$\K \Angle{C}$ is Koszul, combinatorial, and admits a binary and
quadratic presentation, the Hilbert series of $\K \Angle{C}$ and of its
Koszul dual $\K \Angle{C}^!$ are related by
\begin{equation} \label{equ:hilbert_series_koszul_operads}
    \HilbertSeries_{\K \Angle{C}}
        \Par{-\HilbertSeries_{\K \Angle{C}^!}(-t)}
    = t
    = \HilbertSeries_{\K \Angle{C}^!}
        \Par{-\HilbertSeries_{\K \Angle{C}}(-t)}.
\end{equation}
Relation~\eqref{equ:hilbert_series_koszul_operads} can be used either to
prove that an operad is not Koszul (it is the case when the coefficients
of the hypothetical Hilbert series of the Koszul dual admits
coefficients that are not nonnegative integers) or to compute the
Hilbert series of the Koszul dual of a Koszul operad.
\medbreak

The Koszulity of an operad $\K \Angle{C}$ can be proved by using rewrite
systems on syntax trees, in the following way.
\medbreak

\begin{Proposition} \label{prop:koszulity_criterion_pbw}
    Let $\K \Angle{C}$ be an operad admitting a quadratic presentation
    $(\GeneratingSet, \RelationSpace)$. If there exists an orientation
    $\Par{\ColST^\GeneratingSet_\Leaf, \Rew}$ of $\RelationSpace$ such
    that its closure $\Par{\ColST^\GeneratingSet_\Leaf, \RewContext}$ is
    a convergent rewrite system, then $\K \Angle{C}$ is Koszul.
\end{Proposition}
\medbreak

When $\Par{\ColST^\GeneratingSet_\Leaf, \RewContext}$ satisfies the
conditions contained in the statement of
Proposition~\ref{prop:koszulity_criterion_pbw}, the set of
$\GeneratingSet$-syntax trees that are normal forms
\begin{math}
    \NormalForms_{\Par{\ColST^\GeneratingSet_\Leaf, \RewContext}}
\end{math}
forms a basis of
\begin{math}
    \FreeOperad(\GeneratingSet)/_{\Angle{\RelationSpace}}
\end{math},
called \Def{Poincaré-Birkhoff-Witt basis}.
\medbreak

One of the main merits of Koszul operads with Poincaré-Birkhoff-Witt
bases is that they come with a generic way to build an associated
realization. Assume that $(\GeneratingSet, \RelationSpace)$ is a
quadratic presentation and let us set the goal to find a realization of
$\FreeOperad(\GeneratingSet)/_{\Angle{\RelationSpace}}$. If one can
construct an orientation $\Par{\ColST^\GeneratingSet_\Leaf, \Rew}$ of
$\RelationSpace$ such that
\begin{math}
    \Par{\ColST^\GeneratingSet_\Leaf, \RewContext}
\end{math}
is a convergent rewrite system, by
Proposition~\ref{prop:koszulity_criterion_pbw}, the set of all normal
forms of $\Par{\ColST^\GeneratingSet_\Leaf, \RewContext}$ forms a basis
of $\FreeOperad(\GeneratingSet)/_{\Angle{\RelationSpace}}$. Moreover, to
compute the partial composition $\Tfr \circ_i \Sfr$ of two such normal
forms $\Tfr$ and $\Sfr$, start with the syntax tree
$\Rfr := \Tfr \circ_i \Sfr$ obtained by using the partial composition
map $\circ_i$ of $\FreeOperad(\GeneratingSet)$, and then rewrite $\Rfr$
using $\RewContext$ as much as possible in order to obtain a normal form
$\Rfr'$. This process is well-defined since
$\Par{\ColST^\GeneratingSet_\Leaf, \RewContext}$ is convergent. We have
established the fact that the space
\begin{math}
    \K \Angle{
    \NormalForms_{\Par{\ColST^\GeneratingSet_\Leaf, \RewContext}}}
\end{math}
is isomorphic to $\FreeOperad(\GeneratingSet)/_{\Angle{\RelationSpace}}$
and that the partial composition maps just described endow this first
space with an operad structure having $(\GeneratingSet, \RelationSpace)$
as presentation.
\medbreak

\section{Main operads} \label{sec:main_operads}
We provide here classical examples of operads. These examples are
divided into three categories depending on the general families of the
involved combinatorial objects: words, trees, or graphs. We also present
two general constructions to obtain, respectively, operads on words and
operads on graphs. Table~\ref{tab:examples_operads} contains an
overview of these.
\begin{table}[ht]
    \centering
    \begin{tabular}{|c|c|c|c|c|c|} \hline
        Operad & Objects & Arity & Set-operad & Binary & Quadratic \\
        \hline \hline
        $\As$ & Integers & Value & Yes & Yes & Yes \\
        $\Per$ & Permutations & Length
            & Yes & No & Yes \\
        $\Dias$ & Word on $\{0, 1\}$ with one $0$ & Length
            & Yes & Yes & Yes \\
        $\T \Mca$ & Words on $\Mca$ & Length & Yes & No & No \\
        $\Motz$ & Motzkin paths & Points & Yes & No & Yes \\
        \hline
        $\Mag$ & Binary trees & Leaves & Yes & Yes & Yes \\
        $\Dup$ & Binary trees & Int. nodes & Yes & Yes & Yes \\
        $\Dendr$ & Binary trees & Int. nodes & No & Yes & Yes \\
        $\BS$ & Schröder trees & Leaves & Yes & Yes & Yes \\
        $\PreLie$ & Standard rooted trees & Nodes & No & No & Yes \\
        $\NAP$ & Standard rooted trees & Nodes & Yes & No & ? \\
        \hline
        $\NCT$ & Noncrossing trees & Sides & Yes & Yes & Yes \\
        $\BNC$ & Bicolored noncross. config. & Sides
            & Yes & Yes & Yes\\
        $\Grav$ & Gravity chord config. & Sides & Yes & No & No \\
        $\CliMonoid \Mca$ & $\bar{\Mca}$-config. & Sides &
            Yes & No & No \\
        $\NCliMonoid \Mca$ & Noncross. $\bar{\Mca}$-config.
            & Sides & Yes & No & No \\ \hline
    \end{tabular}
    \bigbreak

    \caption[Properties of some operads]
    {\footnotesize Main properties of some operads. Here, $\Mca$
    is a monoid.}
    \label{tab:examples_operads}
\end{table}
\medbreak

\subsection{Operads of words}
Five examples of operads are provided here. Their common point is that
they are defined on graded spaces of families of words. The associative
and diassociative operads seem not, at first glance, operads of words.
We shall explain how to provide a realization of these two operads as
operads of words through a general construction of operads from monoids.
\medbreak

\subsubsection{Associative operad} \label{subsubsec:associative_operad}
Let $A := \left\{\Afr_n : n \in \N_{\geq 1}\right\}$ be the graded
collection where $\left|\Afr_n\right| := n$ for any $n \in \N_{\geq 1}$.
The \Def{associative operad} $\As$ is the space $\K \Angle{A}$ endowed
with the partial composition maps $\circ_i$ defined linearly, for any
$\Afr_n \in A(n)$, $\Afr_m \in A(m)$, and $i \in [n]$, by
\begin{equation}
    \Afr_n \circ_i \Afr_m := \Afr_{n + m - 1}.
\end{equation}
The unit of $\As$ is $\Afr_1$. This operad is a set-operad, is
combinatorial, and its Hilbert series satisfies
\begin{equation}
    \HilbertSeries_\As(t) = \frac{t}{1 - t}.
\end{equation}
Moreover, $\As$ admits the presentation
$(\GeneratingSet, \RelationSpace)$ where $\GeneratingSet := \{\Afr_2\}$
and $\RelationSpace$ is the space generated by
\begin{equation} \label{equ:relation_space_As}
    \Corolla\Par{\Afr_2} \circ_1 \Corolla\Par{\Afr_2} -
    \Corolla\Par{\Afr_2} \circ_2 \Corolla\Par{\Afr_2}.
\end{equation}
\medbreak

Since $\GeneratingSet$ contains only $\Afr_2$, any algebra over $\As$
is a space $\K \Angle{D}$ endowed with a binary product~$\Afr_2$.
Moreover, since $\RelationSpace$ contains the
element~\eqref{equ:relation_space_As}, we have for any
$f_1, f_2, f_3 \in \K \Angle{D}$,
\begin{equation}\begin{split}
    0 & = \Par{\Afr_2 \circ_1 \Afr_2 - \Afr_2 \circ_2 \Afr_2}
    \Par{f_1, f_2, f_3} \\
    & =
    \Par{\Afr_2 \circ_1 \Afr_2}\Par{f_1, f_2, f_3}
    - \Par{\Afr_2 \circ_2 \Afr_2}\Par{f_1, f_2, f_3} \\
    & = \Afr_2\Par{\Afr_2\Par{f_1, f_2}, f_3}
    - \Afr_2\Par{f_1, \Afr_2\Par{f_2, f_3}}.
\end{split}\end{equation}
This is equivalent to the relation
\begin{equation}
    \Par{f_1 \, \Afr_2 \, f_2} \, \Afr_2 \, f_3
    - f_1 \, \Afr_2 \, \Par{f_2 \, \Afr_2 \, f_3} = 0
\end{equation}
written in infix way, implying that $\Afr_2$ is associative. Hence, any
$\As$-algebra is an associative algebra.
\medbreak

\subsubsection{Operad of permutations} \label{subsubsec:operad_per}
For any permutation $\sigma$ of $\SymmetricGroup(n)$, $i \in [n]$, and
$k \in \N$, let $\uparrow_i^k(\sigma)$ be the word on $\N$ obtained by
incrementing by $k$ the letters of $\sigma$ greater than $i$. The
\Def{operad of permutations} $\Per$ is the space
$\K \Angle{\SymmetricGroup}$ endowed with the partial composition maps
$\circ_i$ defined linearly, for any $\sigma \in \SymmetricGroup(n)$,
$\nu \in \SymmetricGroup(m)$, and $i \in [n]$ in the following way.
First, let $\sigma' := \uparrow_{\sigma(i)}^{m - 1}(\sigma)$ and
$\nu' := \uparrow_0^{\sigma(i) - 1}(\nu)$. The partial composition of
$\sigma$ and $\nu$ is defined as
\begin{equation}
    \sigma \circ_i \nu :=
    \sigma'_{|[1, i - 1]} \, \nu' \, \sigma'_{|[i + 1, n]}.
\end{equation}
For instance,
\begin{subequations}
\begin{equation}
    \ColA{1}{\bf 2}\ColA{3}
    \circ_2
    \ColD{12}
    = \ColA{1}\ColD{23}\ColA{4},
\end{equation}
\begin{equation}
    \ColA{741}{\bf 5}\ColA{623}
    \circ_4
    \ColD{231}
    = \ColA{941}\ColD{675}\ColA{823}
\end{equation}
\end{subequations}
are two partial compositions in $\Per$. The unit of $\Per$ is the
permutation $1 \in \SymmetricGroup(1)$. This operad is a set-operad, is
combinatorial, and its Hilbert series satisfies
\begin{equation}
    \HilbertSeries_{\Per}(t) = \sum_{n \in \N_{\geq 1}} n! \, t^n.
\end{equation}
\medbreak

A \Def{simple permutation} is a permutation $\sigma$ such that for all
factors $u$ of $\sigma$, if the letters of $u$ form an interval of $\N$
then $|u| = 1$ or $|u| = |\sigma|$. For instance, the permutation
$6{\bf 2413}57$ is not simple since the letters of the factor
$u := 2413$ form an interval of $\N$. On the other hand, the permutation
$5137462$ is simple.
\medbreak

The operad $\Per$ admits the presentation
$(\GeneratingSet, \RelationSpace)$ where $\GeneratingSet$ is the set of
all simple permutations of sizes $2$ or more and $\RelationSpace$ is the
space generated by
\begin{subequations}
\begin{equation}
    \Corolla(12) \circ_1 \Corolla(12)
    - \Corolla(12) \circ_2 \Corolla(12),
\end{equation}
\begin{equation}
    \Corolla(21) \circ_1 \Corolla(21)
    - \Corolla(21) \circ_2 \Corolla(21).
\end{equation}
\end{subequations}
\medbreak

\subsubsection{Diassociative operad} \label{subsubsec:dias_operad}
Let $E := \left\{\Efr_{n, k} : n \in \N_{\geq 1}, k \in [n]\right\}$ be
the graded collection where $\left|\Efr_{n, k}\right| := n$ for any
$n \in \N_{\geq 1}$ and $k \in [n]$. The \Def{diassociative operad}
$\Dias$ is the space $\K \Angle{E}$ endowed with the partial composition
maps $\circ_i$ defined linearly, for any $\Efr_{n, k} \in E(n)$,
$\Efr_{m, \ell} \in E(m)$, and $i \in [n]$, by
\begin{equation} \label{equ:partial_composition_Dias}
    \Efr_{n, k} \circ_i \Efr_{m, \ell} =
    \begin{cases}
        \Efr_{n + m - 1, k + m - 1}
            & \mbox{if } i < k, \\
        \Efr_{n + m - 1, k + \ell - 1}
            & \mbox{if } i = k, \\
        \Efr_{n + m - 1, k}
            & \mbox{otherwise (} i > k \mbox{)}.
    \end{cases}
\end{equation}
The unit of $\Dias$ is $\Efr_{1, 1}$. This operad is a set-operad, is
combinatorial, and its Hilbert series satisfies
\begin{equation}
    \HilbertSeries_{\Dias}(t) = \frac{t}{(1 - t)^2} =
    \sum_{n \in \N_{\geq 1}} n \, t^n.
\end{equation}
Moreover, $\Dias$ admits the presentation
$(\GeneratingSet, \RelationSpace)$ where
$\GeneratingSet := \left\{\Efr_{2, 1}, \Efr_{2, 2}\right\}$ and
$\RelationSpace$ is the space generated by, by denoting by $\LDias$
(resp. $\RDias$) the elements $\Efr_{2, 1}$ (resp. $\Efr_{2, 2}$),
\begin{subequations}
\begin{equation} \label{equ:operad_dias_relation_1}
    \Corolla(\LDias) \circ_1 \Corolla(\LDias) -
    \Corolla(\LDias) \circ_2 \Corolla(\LDias),
    \quad
    \Corolla(\LDias) \circ_1 \Corolla(\LDias) -
    \Corolla(\LDias) \circ_2 \Corolla(\RDias),
\end{equation}
\begin{equation} \label{equ:operad_dias_relation_2}
    \Corolla(\LDias) \circ_1 \Corolla(\RDias) -
    \Corolla(\RDias) \circ_2 \Corolla(\LDias),
\end{equation}
\begin{equation} \label{equ:operad_dias_relation_3}
    \Corolla(\RDias) \circ_1 \Corolla(\LDias) -
    \Corolla(\RDias) \circ_2 \Corolla(\RDias),
    \quad
    \Corolla(\RDias) \circ_1 \Corolla(\RDias) -
    \Corolla(\RDias) \circ_2 \Corolla(\RDias).
\end{equation}
\end{subequations}
It is possible to show that the closure
$\Par{\ColST_\Leaf^\GeneratingSet, \RewContext}$ of the orientation
$\Par{\ColST_\Leaf^\GeneratingSet, \Rew}$ of $\RelationSpace$ defined by
\begin{subequations}
\begin{equation} \label{equ:dias_rewrite_1}
    \Corolla(\LDias) \circ_2 \Corolla(\LDias)
    \Rew
    \Corolla(\LDias) \circ_1 \Corolla(\LDias),
    \quad
    \Corolla(\LDias) \circ_2 \Corolla(\RDias)
    \Rew
    \Corolla(\LDias) \circ_1 \Corolla(\LDias),
\end{equation}
\begin{equation} \label{equ:dias_rewrite_2}
    \Corolla(\RDias) \circ_2 \Corolla(\LDias)
    \Rew
    \Corolla(\LDias) \circ_1 \Corolla(\RDias),
\end{equation}
\begin{equation} \label{equ:dias_rewrite_3}
    \Corolla(\RDias) \circ_1 \Corolla(\LDias)
    \Rew
    \Corolla(\RDias) \circ_2 \Corolla(\RDias),
    \quad
    \Corolla(\RDias) \circ_1 \Corolla(\RDias)
    \Rew
    \Corolla(\RDias) \circ_2 \Corolla(\RDias).
\end{equation}
\end{subequations}
is convergent. Its normal forms are the syntax trees that avoid the
trees appearing in the left members of~\eqref{equ:dias_rewrite_1},
\eqref{equ:dias_rewrite_2}, and~\eqref{equ:dias_rewrite_3}. All this
implies, by Proposition~\ref{prop:koszulity_criterion_pbw}, that $\Dias$
is Koszul.
\medbreak

Besides, any algebra over $\Dias$ is space $\K \Angle{D}$ endowed with
two binary products $\LDias$ and $\RDias$ such that both $\LDias$ and
$\RDias$ are associative (as consequences
of~\eqref{equ:operad_dias_relation_1}
and~\eqref{equ:operad_dias_relation_3}), and, for any $x, y, z \in D$,
\begin{subequations}
\begin{equation}
    x \LDias y \LDias z = x \LDias (y \RDias z),
\end{equation}
\begin{equation}
    (x \RDias y) \LDias z = x \RDias (y \LDias z),
\end{equation}
\begin{equation}
    (x \LDias y) \RDias z = x \RDias y \RDias z.
\end{equation}
\end{subequations}
These structures are called \Def{diassociative algebras}.
\medbreak

\subsubsection{From monoids to operads}
We describe here a general way for constructing operads of words. Let
$\Mca$ be a monoid with an associative product $\Product$ admitting
$\Unit$ as unit. We denote by $\T \Mca$ the space $\K \Angle{\Mca^+}$
where $\Mca^+$ is the graded collection of all nonempty words on $\Mca$
seen as an alphabet. The space $\T \Mca$ is endowed with the partial
composition maps $\circ_i$ defined linearly, for any $u \in \Mca(n)$,
$v \in \Mca(m)$, and $i \in [n]$, by
\begin{equation}
    u \circ_i v := u(1) \dots u(i - 1)
    \, (u(i) \Product v(1)) \, \dots \, (u(i) \Product v(m)) \,
    u(i + 1) \dots u(n).
\end{equation}
\medbreak

\begin{Proposition} \label{prop:monoids_to_operads}
    For any monoid $\Mca$, $\T \Mca$ is an operad.
\end{Proposition}
\medbreak

The unit of $\T \Mca$ is the unit $\Unit$ of the monoid $\Mca$, seen as
a word of length $1$. The operad $\T \Mca$ is a set-operad. Moreover,
when $\Mca$ is finite, $\T \Mca$ is combinatorial and its Hilbert series
satisfies
\begin{equation}
    \HilbertSeries_{\T \Mca}(t) =  \frac{t}{1 - mt}
    = \sum_{n \in \N_{\geq 1}} m^n \, t^n
\end{equation}
where $m := \# \Mca$.
\medbreak

Let us consider an example. Let $\Mca := \{\Asf, \Bsf\}^*$ be a free
monoid of words. Then, $\T \Mca$ is the space of all words whose letters
are words on $\{\Asf, \Bsf\}$. We call such element \Def{multiwords}.
For instance,
\begin{math}
    (\Asf \Asf, \Bsf \Asf, \Bsf, \epsilon, \Asf)
\end{math}
is a multiword of arity $5$ of $\T \Mca$ and
\begin{equation}
    (\ColA{\Asf \Asf},
    \ColA{\Bsf \Asf}, \mathbf{\Bsf},
    \ColA{\epsilon},
    \ColA{\Asf})
    \circ_3
    (\ColD{\Asf \Bsf},
    \ColD{\epsilon},
    \ColD{\Asf})
    =
    (\ColA{\Asf \Asf},
    \ColA{\Bsf \Asf},
    \mathbf{\Bsf} \ColD{\Asf \Bsf},
    \mathbf{\Bsf}, \mathbf{\Bsf} \ColD{\Asf},
    \ColA{\epsilon}, \ColA{\Asf})
\end{equation}
is a partial composition in $\T \Mca$.
\medbreak

\begin{Proposition} \label{prop:presentation_t}
    Let $\Mca$ be a monoid. Then, the operad $\T \Mca$ admits the
    presentation $(\GeneratingSet, \RelationSpace)$ where
    \begin{math}
        \GeneratingSet :=
        \Mca \sqcup \left\{\Unit\Unit\right\}
    \end{math}
    and $\RelationSpace$ is the space generated by
    \begin{subequations}
    \begin{equation} \label{equ:relation_t_1}
        \Corolla(\Unit\Unit)\circ_1 \Corolla(\Unit\Unit)
        -
        \Corolla(\Unit\Unit) \circ_2 \Corolla(\Unit\Unit),
    \end{equation}
    \begin{equation} \label{equ:relation_t_2}
        \Corolla(x) \circ_1 \Corolla(y)
        -
        \Corolla(x \Product y),
        \qquad x, y \in \Mca,
    \end{equation}
    \begin{equation} \label{equ:relation_t_3}
        \Corolla(\Unit\Unit) \circ
        \left[\Corolla(x), \Corolla(x)\right]
        -
        \Corolla(x) \circ_1 \Corolla(\Unit\Unit),
        \qquad x \in \Mca.
    \end{equation}
    \end{subequations}
\end{Proposition}
\medbreak

Observe that the presentation of $\T \Mca$ provided by
Proposition~\ref{prop:presentation_t} is not minimal in the sense
that the exhibited generating set $\GeneratingSet$ may be not minimal.
\medbreak

The operads $\As$ and $\Dias$ can be obtained through this construction
$\T$. First, one can check that $\As \simeq \T \{\Unit\}$ where
$\{\Unit\}$ is the trivial monoid. An isomorphism between $\As$ and
$\T \{\Unit\}$ is provided by the linear map
$\phi : \As \to \T \{\Unit\}$ satisfying $\phi\Par{\Afr_n} = \Unit^n$
for all $n \in \N_{\geq 1}$. For instance,
\begin{equation}
    \ColA{\Unit\Unit} {\bf \ColF{\Unit}} \ColA{\Unit\Unit\Unit}
    \circ_3
    \ColD{\Unit\Unit}
    =
    \ColA{\Unit\Unit} \ColD{\Unit\Unit} \ColA{\Unit\Unit\Unit}
\end{equation}
is a partial composition in this realization of $\As$. Besides, $\Dias$
is isomorphic to the suboperad of $\T (\N, \max)$ generated by the words
$01$ and $10$. An isomorphism between $\Dias$ and
$\T (\N, \max)^{\{01, 10\}}$ is provided by the linear map
$\phi : \Dias \to \T (\N, \max)^{\{01, 10\}}$ satisfying
$\phi\Par{1^k 0 1^\ell} = \Efr_{k + 1 + \ell, k + 1}$ for all
$k, \ell \in \N$. For instance,
\begin{subequations}
\begin{equation}
    \ColA{11} {\bf \ColF{0}} \ColA{11} \circ_3 \ColD{01}
    =
    \ColA{11} \ColD{01} \ColA{11},
\end{equation}
\begin{equation}
    \ColA{110} {\bf \ColF{1}} \ColA{1} \circ_4 \ColD{01}
    =
    \ColA{110} \ColD{11} \ColA{1}
\end{equation}
\end{subequations}
are two partial compositions in this realization of~$\Dias$.
\medbreak

\subsubsection{Operad of Motzkin words} \label{subsubsec:motzkin_operad}
A \Def{Motzkin word} is a nonempty word $u$ on $\N$ starting and
finishing by $0$ and such that $|u(i) - u(i + 1)| \leq 1$ for all
$i \in [|u| - 1]$. We denote here by $M$ the graded collection of all
the Motzkin words where the size of a word is its length. Let $\Motz$ be
the suboperad of $\T (\N, +)$ generated by the set $\{00, 010\}$. It is
possible to show by induction on the arities that
$\Motz = \K \Angle{M}$. From the definition of the construction $\T$,
the partial composition maps $\circ_i$ of $\Motz$ behave as follows.
Given two Motzkin words $u$ and $v$, $u \circ_i v$ is the Motzkin
word obtained by replacing the letter at position $i$ in $u$ by a copy
of $v$ wherein each of its letters is incremented by $u(i)$. For
instance,
\begin{equation} \label{equ:example_motz}
    \ColA{0112} {\bf \ColF{3}} \ColA{21010} \circ_4 \ColD{0122110}
    =
    \ColA{011} \ColD{2344332} \ColA{321010}
\end{equation}
is a partial composition in $\Motz$. By representing a Motzkin word $u$
as a path in the quarter plane (that is, by drawing points
$(i - 1, u(i))$ for all positions $i$ and by connecting all pairs of
adjacent points by lines), \eqref{equ:example_motz} becomes
\begin{equation}
\,.
\end{equation}
The unit of $\Motz$ is $\UnitPath$, the Motzkin word $0$. This operad is
a set-operad, is combinatorial, and its Hilbert series satisfies
\begin{equation}
    \HilbertSeries_\Motz(t) =
    \frac{1 - t - \sqrt{1 - 2t - 3t^2}}{2t}.
\end{equation}
The first coefficients of its Hilbert series are
\begin{equation}
    1, 1, 2, 4, 9, 21, 51, 127, 323
\end{equation}
and form Sequence~\OEIS{A001006} of~\cite{Slo}. Moreover, $\Motz$ admits
the presentation $(\GeneratingSet, \RelationSpace)$ where
\begin{equation}
    \GeneratingSet := \left\{\PathStable, \PathPeak\right\}
\end{equation}
and $\RelationSpace$ is the space generated by
\begin{subequations}
\begin{equation}
    \Corolla\left(\PathStable\right)
    \circ_1 \Corolla\left(\PathStable\right)
    -
    \Corolla\left(\PathStable\right)
    \circ_2 \Corolla\left(\PathStable\right),
\end{equation}
\begin{equation}
    \Corolla\left(\PathPeak\right)
    \circ_1 \Corolla\left(\PathStable\right)
    -
    \Corolla\left(\PathStable\right)
    \circ_2 \Corolla\left(\PathPeak\right),
\end{equation}
\begin{equation}
    \Corolla\left(\PathStable\right)
    \circ_1 \Corolla\left(\PathPeak\right)
    -
    \Corolla\left(\PathPeak\right)
    \circ_3 \Corolla\left(\PathStable\right),
\end{equation}
\begin{equation}
    \Corolla\left(\PathPeak\right)
    \circ_1 \Corolla\left(\PathPeak\right)
    -
    \Corolla\left(\PathPeak\right)
    \circ_3 \Corolla\left(\PathPeak\right).
\end{equation}
\end{subequations}
\medbreak

\subsection{Operads of trees}
Six examples of operads are provided here. Their common point is that
they are defined on augmented  graded spaces of families of trees:
binary trees (seen endowed with several size functions), bicolored
Schröder trees, and labeled rooted trees.
\medbreak

\subsubsection{Magmatic operad} \label{subsubsec:magmatic_operad}
The \Def{magmatic operad} $\Mag$ is the space $\K \Angle{\ColBT_\Leaf}$
(where $\ColBT_\Leaf$ is the combinatorial graded collection of binary
trees defined in Section~\ref{subsubsec:k_ary_trees} of
Chapter~\ref{chap:trees}) endowed with the partial composition maps
$\circ_i$ defined as the linearizations of the partial grafting defined
in Section~\ref{subsubsec:partial_grafting_syntax_trees} of
Chapter~\ref{chap:trees}. For instance,
\begin{equation}

\end{equation}
is a partial composition in $\Mag$. The unit of $\Mag$ is $\LeafPic$.
This operad is a set-operad, is combinatorial, and its Hilbert series
satisfies
\begin{equation}
    \HilbertSeries_{\Mag}(t)
    = \frac{1 - \sqrt{1 - 4t}}{2}
    = \sum_{n \in \N_{\geq 1}} \binom{2n - 1}{n - 1} \frac{1}{n} t^n.
\end{equation}
Moreover, $\Mag$ admits the presentation
$(\GeneratingSet, \RelationSpace)$ where
\begin{equation}
    \GeneratingSet := \left\{\BinaryNode\right\}
\end{equation}
and $\RelationSpace$ is the trivial space.
\medbreak

Any algebra over $\Mag$ is a space $\K \Angle{D}$ with a binary product
which does satisfy any required relation.
\medbreak

\subsubsection{Duplicial operad}
The \Def{duplicial operad} $\Dup$ is the space
$\K \Angle{\Augmentation\Par{\ColBT_\Node}}$ (where
$\ColBT_\Node$ is the combinatorial graded collection of binary trees
defined in Section~\ref{subsubsec:binary_trees} of
Chapter~\ref{chap:collections}) endowed with the partial composition
maps $\circ_i$ defined linearly, for any
$\Tfr \in \Augmentation\Par{\ColBT_\Node}(n)$,
$\Sfr \in \Augmentation\Par{\ColBT_\Node}(m)$, and $i \in [n]$,
by $\Rfr := \Tfr \circ_i \Sfr$ where $\Rfr$ is the binary tree obtained
by replacing the $i$th (with respect to the infix order) internal node
$u$ of $\Tfr$ by a copy of $\Sfr$, and by grafting the left subtree of
$u$ to the first leaf of the copy, and the right subtree of $u$ to the
last leaf of the copy. For instance,
\begin{equation}

\end{equation}
is a partial composition in $\Dup$. The unit of $\Dup$ is $\BinaryNode$.
This operad is a set-operad, is combinatorial, and its Hilbert series
satisfies
\begin{equation}
    \HilbertSeries_{\Dup}(t)
    = \frac{1 - 2t - \sqrt{1 - 4t}}{2t}
    = \sum_{n \in \N_{\geq 1}} \binom{2n}{n} \frac{1}{n + 1} \, t^n.
\end{equation}
Moreover, $\Dup$ admits the presentation
$(\GeneratingSet, \RelationSpace)$ where
\begin{equation} \label{equ:generating_set_dup}
    \GeneratingSet :=
    \left\{\BinaryTreeLeft\,, \BinaryTreeRight\right\}
\end{equation}
and $\RelationSpace$ is the space generated by, by denoting by $\LDup$
(resp. $\RDup$) the first (resp. second) tree
of~\eqref{equ:generating_set_dup},
\begin{subequations}
\begin{equation} \label{equ:operad_dup_relation_1}
    \Corolla(\LDup) \circ_1 \Corolla(\LDup) -
    \Corolla(\LDup) \circ_2 \Corolla(\LDup),
\end{equation}
\begin{equation} \label{equ:operad_dup_relation_2}
    \Corolla(\RDup) \circ_1 \Corolla(\LDup) -
    \Corolla(\LDup) \circ_2 \Corolla(\RDup),
\end{equation}
\begin{equation} \label{equ:operad_dup_relation_3}
    \Corolla(\RDup) \circ_1 \Corolla(\RDup) -
    \Corolla(\RDup) \circ_2 \Corolla(\RDup).
\end{equation}
\end{subequations}
\medbreak

Any algebra over $\Dup$ is a space $\K \Angle{D}$ endowed with two
binary products $\LDup$ and $\RDup$ such that both $\LDup$ and $\RDup$
are associative (as consequences of~\eqref{equ:operad_dup_relation_1}
and~\eqref{equ:operad_dup_relation_2}), and, for any $x, y, z \in D$,
\begin{equation}
    (x \LDup y) \RDup z = x \LDup (y \RDup z).
\end{equation}
These structures are called \Def{duplicial algebras}.
\medbreak

\subsubsection{Dendriform operad} \label{subsubsec:dendriform_operad}
The \Def{dendriform operad} $\Dendr$ is defined as the operad admitting
the presentation $(\GeneratingSet, \RelationSpace)$ where
\begin{math}
    \GeneratingSet := \GeneratingSet(2) := \{\LDendr, \RDendr\}
\end{math}
and $\RelationSpace$ is the space generated by
\begin{subequations}
\begin{equation} \label{equ:operad_dendr_relation_1}
    \Corolla(\LDendr) \circ_1 \Corolla(\LDendr) -
    \Corolla(\LDendr) \circ_2 \Corolla(\LDendr) -
    \Corolla(\LDendr) \circ_2 \Corolla(\RDendr),
\end{equation}
\begin{equation} \label{equ:operad_dendr_relation_2}
    \Corolla(\LDendr) \circ_1 \Corolla(\RDendr) -
    \Corolla(\RDendr) \circ_2 \Corolla(\LDendr),
\end{equation}
\begin{equation} \label{equ:operad_dendr_relation_3}
    \Corolla(\RDendr) \circ_1 \Corolla(\LDendr) +
    \Corolla(\RDendr) \circ_1 \Corolla(\RDendr) -
    \Corolla(\RDendr) \circ_2 \Corolla(\RDendr).
\end{equation}
\end{subequations}
\medbreak

The dendriform operad and the diassociative operad are the Koszul duals
one of the other. This can be shown by computing a basis of
$\RelationSpace^\perp$ where $\RelationSpace$ is the space of relations
of $\Dendr$, and by observing that $\RelationSpace^\perp$ and the space
of relations of $\Dias$ shown in Section~\ref{subsubsec:dias_operad} are
the same (by replacing, respectively, by $\LDias$ and $\RDias$ the
generators $\LDendr$ and $\RDendr$ appearing in it). As a consequence of
this fact and the Koszulity of $\Dias$, the Hilbert series
$\HilbertSeries_\Dendr(t)$ and $\HilbertSeries_\Dias(t)$
satisfy~\eqref{equ:hilbert_series_koszul_operads}. It is then possible
to obtain the explicit description
\begin{equation}
    \HilbertSeries_\Dendr(t)
    = \frac{1 - 2t - \sqrt{1 - 4t}}{2t}
    = \sum_{n \in \N_{\geq 1}} \binom{2n}{n} \frac{1}{n + 1} \, t^n
\end{equation}
for the Hilbert series of $\Dendr$. This shows that $\Dendr$ is, as a
combinatorial polynomial space, the space
$\K \Angle{\Augmentation\Par{\ColBT_\Node}}$.
\medbreak

From the definition of $\Dendr$ by generators and relations, one can
observe that any algebra over $\Dendr$ is a dendriform algebra (see
Section~\ref{subsec:dendriform_algebras} of
Chapter~\ref{chap:algebra}). Moreover, the free dendriform algebra over
one generator is the space $\Dendr$, that is the linear span of all
nonempty binary trees, endowed with the linear binary products
$\LDendr$ and $\RDendr$ defined recursively, for any nonempty  tree
$\Sfr$, and binary trees $\Tfr_1$ and $\Tfr_2$ by
\begin{subequations}
\begin{multicols}{2}
\begin{equation}
    \Sfr \LDendr \LeafPic
    := \Sfr =:
    \LeafPic \RDendr \Sfr,
\end{equation}

\begin{equation}
    \LeafPic \LDendr \Sfr := 0 =: \Sfr \RDendr \LeafPic,
\end{equation}
\end{multicols}
\begin{multicols}{2}
\begin{equation}
    \BinaryRoot{\Tfr_1}{\Tfr_2} \LDendr \Sfr :=
    \BinaryRoot{\Tfr_1}{\Tfr_2 \LDendr \Sfr}
    + \BinaryRoot{\Tfr_1}{\Tfr_2 \RDendr \Sfr}\,,
\end{equation}

\begin{equation}
    \begin{split} \Sfr \RDendr \end{split}
    \BinaryRoot{\Tfr_1}{\Tfr_2} :=
    \BinaryRoot{\Sfr \RDendr \Tfr_1}{\Tfr_2}
    + \BinaryRoot{\Sfr \LDendr \Tfr_1}{\Tfr_2}.
\end{equation}
\end{multicols}
\end{subequations}
\noindent
Note that neither $\LeafPic \LDendr \LeafPic$ nor
$\LeafPic \RDendr \LeafPic$ need to be defined. We have for instance,
\begin{subequations}
\begin{equation}
\,.
\end{equation}
\end{subequations}
\medbreak

Besides, one can check that the element $\LDendr + \RDendr$ of $\Dendr$
is associative. This implies that the linear map
$\phi : \As \to \Dendr$ defined by $\phi(\Afr_2) := \LDendr + \RDendr$
extends in a unique way into an operad morphism. Now, by
Proposition~\ref{prop:operad_morphisms_and_algebras}, we obtain that
if $\K \Angle{D}$ is a dendriform algebra, the binary product $\Afr_2$
defined for any $f_1, f_2 \in \K \Angle{D}$ by
\begin{equation}
    f_1 \, \Afr_2 \, f_2
    := f_1 \, \phi(\Afr_2) \, f_2
    = f_1 \LDendr f_2 + f_1 \RDendr f_2
\end{equation}
is associative and endows $\K \Angle{D}$ with the structure of an
associative algebra.
\medbreak

\subsubsection{Bicolored Schröder tree operad}
A \Def{bicolored Schröder tree} is a Schröder tree $\Tfr$ (see
Section~\ref{subsubsec:schroder_trees} of Chapter~\ref{chap:trees}) such
that each internal node is assigned with an element of the set
$\{0, 1\}$ and all internal nodes that have a father labeled by $0$
(resp. $1$) are labeled by $1$ (resp. $0$). Let $\ColBSch_\Leaf$ be the
graded collection of all bicolored Schröder trees wherein the size of
such trees is their number of leaves. The \Def{bicolored Schröder tree
operad} $\BS$ is the space $\K \Angle{\ColBSch_\Leaf}$ endowed with the
partial composition maps $\circ_i$ defined linearly, for any
$\Tfr \in \ColBSch_\Leaf(n)$, $\Sfr \in \ColBSch_\Leaf(m)$, and
$i \in [n]$ by $\Rfr := \Tfr \circ_i \Sfr$ where $\Rfr$ is the bicolored
Schröder tree obtained by grafting a copy of $\Sfr$ onto the $i$th leaf
of $\Tfr$ and, in the case where the edge connecting this leaf and the
copy of $\Sfr$ have the extremities that are internal nodes labeled by
the same element $a \in \{0, 1\}$, by contracting this edge to form a
single internal node labeled by $a$. For instance,
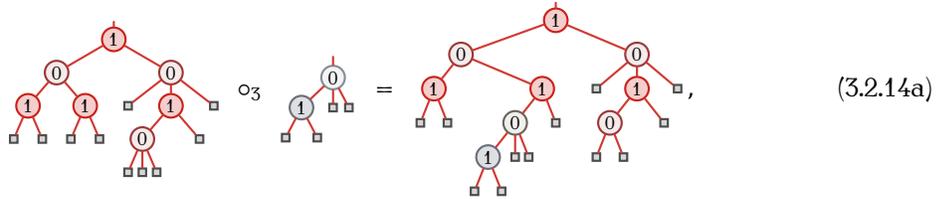
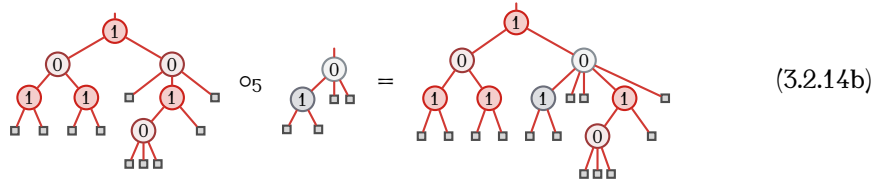
\begin{subequations}
\begin{equation}

\end{equation}
\end{subequations}
are two partial compositions in $\BS$. The unit of $\BS$ is $\LeafPic$.
This operad is a set-operad, is combinatorial, and its Hilbert series
satisfies
\begin{equation}
    \HilbertSeries_{\BS}(t) = \frac{1 - t - \sqrt{1 - 6 t + t^2}}{2}.
\end{equation}
The first coefficients of its Hilbert series are
\begin{equation}
    1, 2, 6, 22, 90, 394, 1806, 8558, 41586
\end{equation}
and form Sequence~\OEIS{A006318} of~\cite{Slo}. Moreover, $\BS$ admits
the presentation $(\GeneratingSet, \RelationSpace)$ where
\begin{equation} \label{equ:generating_set_bs}
    \GeneratingSet :=
    \left\{
    \begin{tikzpicture}[xscale=.21,yscale=.2,Centering]
        \node[Leaf](0)at(0.00,-1.50){};
        \node[Leaf](2)at(2.00,-1.50){};
        \node[Node](1)at(1.00,0.00){\begin{math}0\end{math}};
        \draw[Edge](0)--(1);
        \draw[Edge](2)--(1);
        \node(r)at(1.00,2){};
        \draw[Edge](r)--(1);
    \end{tikzpicture}\,,
    \begin{tikzpicture}[xscale=.21,yscale=.2,Centering]
        \node[Leaf](0)at(0.00,-1.50){};
        \node[Leaf](2)at(2.00,-1.50){};
        \node[NodeColorC](1)at(1.00,0.00){\begin{math}1\end{math}};
        \draw[Edge](0)--(1);
        \draw[Edge](2)--(1);
        \node(r)at(1.00,2){};
        \draw[Edge](r)--(1);
    \end{tikzpicture}
    \right\}
\end{equation}
and $\RelationSpace$ is the space generated by, by denoting by
$\Product_0$ (resp. $\Product_1$) the first (resp. second) tree
of~\eqref{equ:generating_set_bs},
\begin{subequations}
\begin{equation}
    \Corolla\Par{\Product_0}
    \circ_1 \Corolla\Par{\Product_0}
    -
    \Corolla\Par{\Product_0}
    \circ_2 \Corolla\Par{\Product_0},
\end{equation}
\begin{equation}
    \Corolla\Par{\Product_1}
    \circ_1 \Corolla\Par{\Product_1}
    -
    \Corolla\Par{\Product_1}
    \circ_2 \Corolla\Par{\Product_1}.
\end{equation}
\end{subequations}
\medbreak

Any algebra over $\BS$ is a space $\K \Angle{D}$ endowed with two binary
associative products $\Product_0$ and $\Product_1$. These structures are
called \Def{two-associative algebras}.
\medbreak

\subsubsection{Labeled rooted trees}
\label{subsubsec:labeled_rooted_trees}
A \Def{labeled rooted tree} is a rooted tree $\Tfr$ (see
Section~\ref{subsec:rooted_trees} of Chapter~\ref{chap:trees}) endowed
with an injective map sending each internal node of $\Tfr$ to an element
of $\N$ called \Def{label}. Due to the injective labeling of the nodes
of any labeled rooted tree $\Tfr$, we shall identify each node of $\Tfr$
with its label. The set of all labels appearing in $\Tfr$ is denoted by
$\Labels(\Tfr)$. For any $i \in \Labels(\Tfr)$, we denote by
$\Tfr^{(i)}$ the set of the suffix subtrees rooted at the children of
the node $i$ in $\Tfr$. Moreover, for any $k \in \N$, we denote by
$\uparrow_i^k(\Tfr)$ the labeled rooted tree obtained from $\Tfr$ by
incrementing by $k$ its nodes greater than $i$. Let $\Tfr$ and $\Sfr$ be
two labeled rooted trees such that $i \in \Labels(\Tfr)$ and
\begin{math}
    \Par{\Labels(\Tfr) \setminus \{i\}} \cap \Labels(\Sfr) = \emptyset,
\end{math}
and $\phi : \Tfr^{(i)} \to \Labels(\Sfr)$ be a map. We denote by
$\Tfr \PreLieGrafting_i^\phi \Sfr$ the labeled rooted tree obtained by
replacing the node $i$ in $\Tfr$ by the root of a copy of $\Sfr$, and by
grafting each tree $\Rfr$ of $\Tfr^{(i)}$ as a child of the node
$\phi(\Rfr)$ in the copy of~$\Sfr$.
\medbreak

A \Def{standard rooted tree} is a labeled rooted tree $\Tfr$ having all
its labels in the set $[n]$ where $n$ is the number of nodes of $\Tfr$.
We denote by $\ColSRT$ the combinatorial graded collection of the
standard rooted trees wherein the size of such trees is their number of
nodes. As usual, we draw standard rooted trees as rooted trees where the
label of each internal node is written inside it.
\medbreak

\subsubsection{Pre-Lie operad} \label{subsubsec:pre_Lie_operad}
The \Def{pre-Lie operad} $\PreLie$ is the space $\K \Angle{\ColSRT}$
endowed with the partial composition maps $\circ_i$ defined linearly,
for any $\Tfr \in \ColSRT(n)$, $\Sfr \in \ColSRT(m)$, and $i \in [n]$ in
the following way. First, let $\Tfr' := \uparrow_i^{m - 1}(\Tfr)$ and
$\Sfr' := \uparrow_0^{i - 1}(\Sfr)$. The partial composition of $\Tfr$
and $\Sfr$ is defined as the sum
\begin{equation}
    \Tfr \circ_i \Sfr :=
    \sum_{\phi \, : \, {\Tfr'}^{(i)} \to \Labels\Par{\Sfr'}}
    \Tfr' \PreLieGrafting_i^\phi \Sfr'.
\end{equation}
For instance,
\begin{subequations}
\begin{equation}

\end{equation}
\end{subequations}
are two partial compositions in $\PreLie$. The unit of $\PreLie$ is
\begin{math}
    \raisebox{.5em}{%
    %
}.
\end{math}
This operad is combinatorial and its Hilbert series satisfies
\begin{equation}
    \HilbertSeries_{\PreLie}(t)
    = \sum_{n \in \N_{\geq 1}} n^{n - 1} \, t^n.
\end{equation}
The first coefficients of its Hilbert series are
\begin{equation}
    1, 2, 9, 64, 625, 7776, 117649, 2097152, 43046721
\end{equation}
and form Sequence~\OEIS{A000169} of~\cite{Slo}.
\medbreak

\subsubsection{Nonassociative permutative operad}
\label{subsubsec:operad_nap}
The \Def{nonassociative permutative operad} $\NAP$ is the space
$\K \Angle{\ColSRT}$ endowed with the partial composition maps $\circ_i$
defined linearly, for any $\Tfr \in \ColSRT(n)$, $\Sfr \in \ColSRT(m)$,
and $i \in [n]$ in the following way. By using the notations of
Section~\ref{subsubsec:labeled_rooted_trees} about labeled rooted trees,
let $\Tfr' := \uparrow_i^{m - 1}(\Tfr)$,
$\Sfr' := \uparrow_0^{i - 1}(\Sfr)$, and
$\phi : {\Tfr'}^{(i)} \to \Labels\Par{\Sfr'}$ be the map defined
for any $\Rfr \in {\Tfr'}^{(i)}$ by $\phi(\Rfr) := j$ where $j$ is the
label of the root of $\Sfr'$. The partial composition of $\Tfr$ and
$\Sfr$ is defined as
\begin{equation}
    \Tfr \circ_i \Sfr := \Tfr' \PreLieGrafting_i^\phi \Sfr'.
\end{equation}
Observe that $\Tfr \circ_i \Sfr$ is a particular element appearing in
the partial composition $\Tfr \circ_i \Sfr$ of the operad $\PreLie$. For
instance,
\begin{subequations}
\begin{equation}

\end{equation}
\end{subequations}
are two partial compositions in $\NAP$. The unit of $\NAP$ is
\begin{math}
    \raisebox{.5em}{%
    %
}.
\end{math}
This operad is a set-operad, is combinatorial, and its Hilbert series
is the same as the one of~$\PreLie$.
\medbreak

\subsection{Operads of graphs} \label{subsec:operads_graphs}
As last examples, we expose here operads defined on graded spaces of
families of graphs. These graphs are configurations of chords in
polygons having labeled arcs. We shall also provide a general
construction of operads of graphs from unitary magmas.
\medbreak

\subsubsection{Configurations of chords}
A \Def{polygon} of \Def{size} $n \in \N_{\geq 1}$ is a directed graph
$\Pfr$ on the set of vertices $[n + 1]$. An \Def{arc} of $\Pfr$ is a
pair of integers $(x, y)$ with $1 \leq x < y \leq n + 1$, a
\Def{diagonal} is an arc $(x, y)$ different from $(x, x + 1)$ and
$(1, n + 1)$, and an \Def{edge} is an arc of the form $(x, x + 1)$ and
different from $(1, n + 1)$. We denote by $\Arcs_\Pfr$ (resp.
$\Diagonals_\Pfr$, $\Edges_\Pfr$) the set of all arcs (resp. diagonals,
edges) of $\Pfr$. For any $i \in [n]$, the \Def{$i$th edge} of $\Pfr$ is
the edge $(i, i + 1)$, and the arc $(1, n + 1)$ is the \Def{base}
of~$\Pfr$.
\medbreak

For any set $S$, an \Def{$S$-configuration of chords} (or simply an
\Def{$S$-configuration}) is a polygon $\Cfr$ endowed with a partial
function
\begin{equation}
    \phi_\Cfr : \Arcs_\Cfr \to S.
\end{equation}
When $\phi_\Cfr((x, y))$ is defined, we say that the arc $(x, y)$ is
\Def{labeled} and we denote it by $\Cfr(x, y)$, otherwise, $(x, y)$ is
\Def{unlabeled}. When the base of $\Cfr$ is labeled, we denote it by
$\Cfr_0$, and when the $i$th edge of $\Cfr$ is labeled, we denote it by
$\Cfr_i$. Two diagonals $(x, y)$ and $\Par{x', y'}$ of $\Cfr$ are
\Def{crossing} if $x < x' < y < y'$ or $x' < x < y' < y$. The
$S$-configuration $\Cfr$ is \Def{noncrossing} if it does not admit any
pair of crossing labeled diagonals. The graded collection of all
$S$-configurations (resp. noncrossing $S$-configurations) is denoted by
$\ColCC^S$ (resp. $\ColNCC^S$).
\medbreak

In our graphical representations, each polygon is depicted so that its
base is the bottommost segment, vertices are implicitly numbered from
$1$ to $n + 1$ in the clockwise direction We shall represent any
$S$-configuration $\Cfr$ by drawing a polygon of the same size as the
one of $\Cfr$ and by labeling its arcs accordingly. For instance
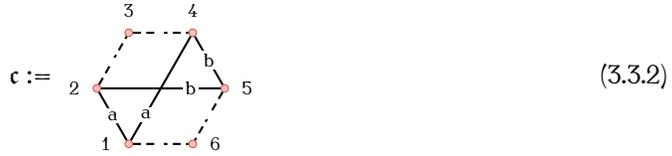
\begin{equation}
    \Cfr :=
    \begin{tikzpicture}[scale=.85,Centering]
        \node[CliquePoint](1)at(-0.50,-0.87){};
        \node[CliquePoint](2)at(-1.00,-0.00){};
        \node[CliquePoint](3)at(-0.50,0.87){};
        \node[CliquePoint](4)at(0.50,0.87){};
        \node[CliquePoint](5)at(1.00,0.00){};
        \node[CliquePoint](6)at(0.50,-0.87){};
        \draw[CliqueEdge](1)edge[]node[CliqueLabel]
            {\begin{math}\Asf\end{math}}(2);
        \draw[CliqueEmptyEdge](1)edge[]node[CliqueLabel]{}(6);
        \draw[CliqueEmptyEdge](2)edge[]node[CliqueLabel]{}(3);
        \draw[CliqueEmptyEdge](3)edge[]node[CliqueLabel]{}(4);
        \draw[CliqueEdge](4)edge[]node[CliqueLabel]
            {\begin{math}\Bsf\end{math}}(5);
        \draw[CliqueEmptyEdge](5)edge[]node[CliqueLabel]{}(6);
        \draw[CliqueEdge](1)edge[]node[CliqueLabel,near start]
            {\begin{math}\Asf\end{math}}(4);
        \draw[CliqueEdge](2)edge[]node[CliqueLabel,near end]
            {\begin{math}\Bsf\end{math}}(5);
        \node[left of=1,node distance=3mm,font=\scriptsize]
            {\begin{math}1\end{math}};
        \node[left of=2,node distance=3mm,font=\scriptsize]
            {\begin{math}2\end{math}};
        \node[above of=3,node distance=3mm,font=\scriptsize]
            {\begin{math}3\end{math}};
        \node[above of=4,node distance=3mm,font=\scriptsize]
            {\begin{math}4\end{math}};
        \node[right of=5,node distance=3mm,font=\scriptsize]
            {\begin{math}5\end{math}};
        \node[right of=6,node distance=3mm,font=\scriptsize]
            {\begin{math}6\end{math}};
    \end{tikzpicture}
\end{equation}
is an $\{\Asf, \Bsf\}$-configuration of size $5$. Its set of all
diagonals is
\begin{equation}
    \Diagonals_\Cfr =
    \{(1, 3), (1, 4), (1, 5), (2, 4), (2, 5), (2, 6),
    (3, 5), (3, 6), (4, 6)\},
\end{equation}
its set of all edges is
\begin{equation}
    \Edges_\Cfr = \{(1, 2), (2, 3), (3, 4), (4, 5), (5, 6)\},
\end{equation}
and its set of all arcs is
\begin{equation}
    \Arcs_\Cfr = \Diagonals_\Cfr \sqcup \Edges_\Cfr \sqcup \{(1, 6)\}.
\end{equation}
The arcs $(1, 2)$ and $(1, 4)$ of $\Cfr$ are labeled by $\Asf$, the arcs
$(2, 5)$ and $(4, 5)$ are labeled by $\Bsf$, and the other arcs are
unlabeled. The labeled diagonals $(1, 4)$ and $(2, 5)$ are crossing so
that $\Cfr$ is not noncrossing.
\medbreak

\subsubsection{Noncrossing tree operad} \label{subsubsec:nct_operad}
A \Def{noncrossing tree} is a $\{\star\}$-configuration $\Cfr$, where
$\star$ is any symbol, satisfying the following conditions. First,
$\Cfr$ is noncrossing and its base is labeled, and, by denoting by $n$
the size of $\Cfr$, the graph on $[n + 1]$ consisting in the edges
$\{x, y\}$ if $(x, y)$ is labeled in $\Cfr$, is connected and simply
connected. For instance,
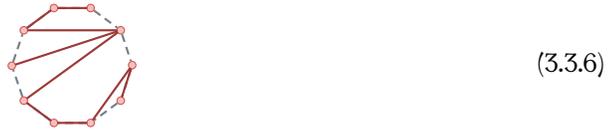
\begin{equation}
    \begin{tikzpicture}[scale=.8,Centering]
        \node[CliquePoint](0)at(-0.3,-0.95){};
        \node[CliquePoint](1)at(-0.8,-0.58){};
        \node[CliquePoint](2)at(-1.,-0.){};
        \node[CliquePoint](3)at(-0.8,0.59){};
        \node[CliquePoint](4)at(-0.3,0.96){};
        \node[CliquePoint](5)at(0.31,0.96){};
        \node[CliquePoint](6)at(0.81,0.59){};
        \node[CliquePoint](7)at(1.,0.01){};
        \node[CliquePoint](8)at(0.81,-0.58){};
        \node[CliquePoint](9)at(0.31,-0.95){};
        \draw[CliqueEmptyEdge](0)--(1);
        \draw[CliqueEmptyEdge](1)--(2);
        \draw[CliqueEmptyEdge](2)--(3);
        \draw[CliqueEmptyEdge](3)--(4);
        \draw[CliqueEmptyEdge](4)--(5);
        \draw[CliqueEmptyEdge](5)--(6);
        \draw[CliqueEmptyEdge](6)--(7);
        \draw[CliqueEmptyEdge](7)--(8);
        \draw[CliqueEmptyEdge](8)--(9);
        \draw[CliqueEmptyEdge](9)--(0);
        \draw[CliqueEdgeBlue](0)--(1);
        \draw[CliqueEdgeBlue](0)--(9);
        \draw[CliqueEdgeBlue](1)--(6);
        \draw[CliqueEdgeBlue](2)--(6);
        \draw[CliqueEdgeBlue](3)--(6);
        \draw[CliqueEdgeBlue](3)--(4);
        \draw[CliqueEdgeBlue](4)--(5);
        \draw[CliqueEdgeBlue](7)--(8);
        \draw[CliqueEdgeBlue](7)--(9);
    \end{tikzpicture}
\end{equation}
is a noncrossing tree of size $9$. The graded collection of all
noncrossing trees is denoted by $\ColNCT$. The
\Def{operad of noncrossing trees} $\NCT$ is the space
$\K \Angle{\ColNCT}$ endowed with the partial composition maps $\circ_i$
defined graphically as follows. For any $\Cfr \in \ColNCT(n)$,
$\Dfr \in \ColNCT(m)$, and $i \in [n]$, the noncrossing tree
$\Cfr \circ_i \Dfr$ is obtained by gluing the base of $\Dfr$ onto the
$i$th edge of $\Cfr$, so that the arc $(i, i + m)$ of
$\Cfr \circ_i \Dfr$ is labeled when the $i$th edge of $\Cfr$ is labeled
and is unlabeled otherwise. For example,
\begin{subequations}
\begin{equation}

\end{equation}
\end{subequations}
are two partial compositions in $\NCT$. The unit of $\NCT$ is
$\UnitSolidClique$. This operad is a set-operad, is combinatorial, and
its Hilbert series satisfies
\begin{equation}
    \HilbertSeries_{\NCT}(t) =
    \sum_{n \in \N_{\geq 1}} \binom{3n - 2}{n - 1} \frac{1}{n} \, t^n.
\end{equation}
The first coefficients of its Hilbert series are
\begin{equation}
    1, 2, 7, 30, 143, 728, 3876, 21318, 120175
\end{equation}
and form Sequence~\OEIS{A006013} of~\cite{Slo}. Moreover, $\NCT$ admits
the presentation $(\GeneratingSet, \RelationSpace)$ where
\begin{equation} \label{equ:generating_set_nct}
    \GeneratingSet
    := \left\{\TriangleLeft, \TriangleRight\right\}
\end{equation}
and $\RelationSpace$ is the space generated by, by denoting by $\LNCT$
(resp. $\RNCT$) the first (resp. second) noncrossing tree
of~\eqref{equ:generating_set_nct},
\begin{equation}
    \Corolla(\RNCT) \circ_1 \Corolla(\LNCT)
    - \Corolla(\LNCT) \circ_2 \Corolla(\RNCT).
\end{equation}
\medbreak

\subsubsection{Bicolored noncrossing configuration operad}
A \Def{bicolored noncrossing configuration} is a noncrossing
$\left\{\star, \bar{\star}\right\}$-configuration $\Cfr$, where $\star$
and $\bar{\star}$ are any symbols, such that all arcs labeled by
$\bar{\star}$ are diagonals. We draw any arc labeled by $\star$ (resp.
$\bar{\star}$) by a thick (resp. dotted) line. For instance,
\begin{equation}
    \begin{tikzpicture}[scale=.8,Centering]
        \node[CliquePoint](0)at(-0.3,-0.95){};
        \node[CliquePoint](1)at(-0.8,-0.58){};
        \node[CliquePoint](2)at(-1.,-0.){};
        \node[CliquePoint](3)at(-0.8,0.59){};
        \node[CliquePoint](4)at(-0.3,0.96){};
        \node[CliquePoint](5)at(0.31,0.96){};
        \node[CliquePoint](6)at(0.81,0.59){};
        \node[CliquePoint](7)at(1.,0.01){};
        \node[CliquePoint](8)at(0.81,-0.58){};
        \node[CliquePoint](9)at(0.31,-0.95){};
        \draw[CliqueEmptyEdge](0)--(1);
        \draw[CliqueEmptyEdge](1)--(2);
        \draw[CliqueEmptyEdge](2)--(3);
        \draw[CliqueEmptyEdge](3)--(4);
        \draw[CliqueEmptyEdge](4)--(5);
        \draw[CliqueEmptyEdge](5)--(6);
        \draw[CliqueEmptyEdge](6)--(7);
        \draw[CliqueEmptyEdge](7)--(8);
        \draw[CliqueEmptyEdge](8)--(9);
        \draw[CliqueEmptyEdge](9)--(0);
        \draw[CliqueEdgeBlue](0)--(1);
        \draw[CliqueEdgeBlue](1)--(7);
        \draw[CliqueEdgeBlue](3)--(5);
        \draw[CliqueEdgeBlue](6)--(7);
        \draw[CliqueEdgeBlue](8)--(9);
        \draw[CliqueEdgeRed](1)--(5);
        \draw[CliqueEdgeRed](1)--(9);
    \end{tikzpicture}
\end{equation}
is a bicolored noncrossing configuration of size~$9$. By definition, we
set that there is only one bicolored noncrossing configuration of size
$1$, having its only arc unlabeled. The graded collection of all
bicolored noncrossing configurations is denoted by $\ColBNC$. The
\Def{operad of bicolored noncrossing configurations} $\BNC$ is the space
$\K \Angle{\ColBNC}$ endowed with the partial composition maps $\circ_i$
defined graphically as follows. For any $\Cfr \in \ColBNC(n)$,
$\Dfr \in \ColBNC(m)$, and $i \in [n]$, the bicolored noncrossing
configuration $\Cfr \circ_i \Dfr$ is obtained by gluing the base of
$\Dfr$ onto the $i$th edge of $\Cfr$, and then, if the base of $\Dfr$
and the $i$th edge of $\Cfr$ are both unlabeled, the arc $(i, i + m)$ of
$\Cfr \circ_i \Dfr$ becomes labeled by $\bar{\star}$; if the base of
$\Dfr$ and the $i$th edge of $\Cfr$ are both labeled by $\star$, the arc
$(i, i + m)$ of $\Cfr \circ_i \Dfr$ becomes labeled by $\star$;
otherwise, the arc $(i, i + m)$ of $\Cfr \circ_i \Dfr$ becomes
unlabeled. For example,
\begin{subequations}
\begin{equation}

\end{equation}
\end{subequations}
are three partial compositions in $\BNC$. The unit of $\BNC$ is
\begin{math}
%
.
\end{math}
This operad is a set-operad, is combinatorial, and its Hilbert series
satisfies
\begin{equation}
    \HilbertSeries_{\BNC}(t) =
    \frac{1 - 4t - \sqrt{1 - 20t + 4t^2}}{6}.
\end{equation}
The first coefficients of its Hilbert series are
\begin{equation}
    1, 8, 80, 992, 13760, 204416, 3180800, 51176960, 844467200
\end{equation}
and form Sequence~\OEIS{A234596} of~\cite{Slo}.
\medbreak

\subsubsection{Gravity operad}
A \Def{gravity chord configuration} is an $\{\star\}$-configuration
$\Cfr$, where $\star$ is any symbol, satisfying the following
conditions. By denoting by $n$ the size of $\Cfr$, all the edges and the
base of $\Cfr$ are labeled, and if $(x, y)$ and $\Par{x', y'}$ are two
labeled crossing diagonals of $\Cfr$ such that $x < x'$, the arc
$\Par{x', y}$ is unlabeled. In other words, the quadrilateral formed by
the vertices $x$, $x'$, $y$, and $y'$ of $\Cfr$ is such that its side
$(x', y)$ is unlabeled. For instance,
\begin{equation}
    \begin{tikzpicture}[scale=.8,Centering]
        \node[CliquePoint](1)at(-0.38,-0.92){};
        \node[CliquePoint](2)at(-0.92,-0.38){};
        \node[CliquePoint](3)at(-0.92,0.38){};
        \node[CliquePoint](4)at(-0.38,0.92){};
        \node[CliquePoint](5)at(0.38,0.92){};
        \node[CliquePoint](6)at(0.92,0.38){};
        \node[CliquePoint](7)at(0.92,-0.38){};
        \node[CliquePoint](8)at(0.38,-0.92){};
        \draw[CliqueEdgeBlue](1)--(2);
        \draw[CliqueEdgeBlue](1)--(8);
        \draw[CliqueEdgeBlue](2)--(3);
        \draw[CliqueEdgeBlue](2)--(5);
        \draw[CliqueEdgeBlue](2)--(6);
        \draw[CliqueEdgeBlue](2)--(7);
        \draw[CliqueEdgeBlue](3)--(4);
        \draw[CliqueEdgeBlue](3)--(6);
        \draw[CliqueEdgeBlue](4)--(5);
        \draw[CliqueEdgeBlue](5)--(6);
        \draw[CliqueEdgeBlue](6)--(7);
        \draw[CliqueEdgeBlue](7)--(8);
    \end{tikzpicture}
\end{equation}
is a gravity chord configuration of size $7$ having four labeled
diagonals (observe in particular that, as required, the arc $(3, 5)$ is
not labeled). By definition, we set that there is only one  gravity
chord configuration of size $1$, having its only arc unlabeled. The
graded collection of all gravity chord configurations is denoted by
$\ColGravCC$. The \Def{operad of gravity chord configurations} $\Grav$
is the space $\K \Angle{\ColGravCC}$ endowed with the partial
composition maps $\circ_i$ defined graphically as follows. For any
$\Cfr \in \ColGravCC(n)$, $\Dfr \in \ColGravCC(m)$, and $i \in [n]$, the
gravity chord configuration $\Cfr \circ_i \Dfr$ is obtained by gluing
the base of $\Dfr$ onto the $i$th edge of $\Cfr$, so that the arc
$(i, i + m)$ of $\Cfr \circ_i \Dfr$ is labeled. For example,
\begin{equation}

\end{equation}
is a partial composition in $\Grav$. The unit of $\Grav$ is
\begin{math}
%
.
\end{math}
This operad is a set-operad, is combinatorial, and its Hilbert series
satisfies
\begin{equation}
    \HilbertSeries_{\Grav}(t)
    = t + \sum_{n \in \N_{\geq 2}} \frac{n!}{2} \, t^n.
\end{equation}
The first coefficients of its Hilbert series are
\begin{equation}
    1, 1, 3, 12, 60, 360, 2520, 20160, 181440
\end{equation}
and form Sequence~\OEIS{A001710} of~\cite{Slo}.
\medbreak

\subsubsection{From unitary magmas to graph operads}
\label{subsubsec:configuration_operads}
We describe here a general way for constructing operads of
configurations of chords. Let $\Mca$ be a unitary magma with binary
product $\Product$ admitting $\Unit$ as unit. By setting
$\bar{\Mca} := \Mca \setminus \{\Unit\}$, $\ColCC_{\bar{\Mca}}$ is the
graded collection of all $\Mca$-configurations where all labeled arcs
have labels different from $\Unit$. We denote by $\CliMonoid \Mca$ the
space $\K \Angle{\ColCC_{\bar{\Mca}}}$. The space $\CliMonoid \Mca$ is
endowed with the partial composition maps $\circ_i$ defined linearly,
for any $\Cfr \in \ColCC_{\bar{\Mca}}(n)$,
$\Dfr \in \ColCC_{\bar{\Mca}}(m)$, and $i \in [n]$ in the following way.
First, let $\Cfr'$ and $\Dfr'$ be the two $\Mca$-configurations,
respectively, obtained from $\Cfr$ and $\Dfr$ by labeling by $\Unit$ the
possible unlabeled edges and the possible unlabeled bases. The
configuration $\Cfr \circ_i \Dfr$ is obtained by gluing the base of
$\Dfr'$ onto the $i$th edge of $\Cfr'$, then by relabeling the arc
$(i, i + m)$ of $\Cfr \circ_i \Dfr$ by $\Cfr_i \Product \Dfr_0$, and
finally by making unlabeled the possible arcs labeled by $\Unit$.
\medbreak

\begin{Proposition} \label{prop:monoid_to_configuration_operads}
    For any monoid $\Mca$, $\CliMonoid \Mca$ is an operad.
\end{Proposition}
\medbreak

The unit of $\CliMonoid \Mca$ is the $\bar{\Mca}$-configuration
$\UnitClique$ of size $1$ having its only arc unlabeled. The operad
$\CliMonoid \Mca$ is a set-operad. Moreover, when $\Mca$ is finite,
$\CliMonoid \Mca$ is combinatorial and its Hilbert series satisfies
\begin{equation}
    \HilbertSeries_{\CliMonoid \Mca}(t) =
    \sum_{n \in \N_{\geq 1}} m^{\binom{n + 1}{2}} \, t^n
\end{equation}
where $m := \# \Mca$.
\medbreak

Let us consider an example. By considering the monoid $\Z$ for the usual
integer addition, $\CliMonoid \Z$ is the space of all configurations
with labels in $\Z \setminus \{0\}$. Moreover,
\begin{subequations}
\begin{equation}

\end{equation}
\end{subequations}
are two partial compositions in~$\CliMonoid \Z$.
\medbreak

Let now the space
\begin{equation}
    \CliMagma \Mca :=
    \K \Angle{
        \left\{\UnitClique\right\} +
        \bigsqcup_{n \in \N_{\geq 2}}
        \ColCC_{\bar{\Mca}}(n)}.
\end{equation}
\medbreak

\begin{Proposition} \label{prop:magma_to_configuration_operads}
    When $\Mca$ is a unitary magma, $\CliMagma \Mca$ is an operad.
    Moreover, when $\Mca$ is also a monoid, $\CliMagma \Mca$ is a
    suboperad of~$\CliMonoid \Mca$.
\end{Proposition}
\medbreak

Now, when $\Mca$ is a monoid, let us define $\NCliMonoid \Mca$ as the
space $\K \Angle{\ColNCC_{\bar{\Mca}}}$ of the noncrossing
$\bar{\Mca}$-configurations. Immediately by definition of the partial
composition of $\CliMonoid \Mca$, one can observe that if $\Cfr$ and
$\Dfr$ are noncrossing, any partial composition $\Cfr \circ_i \Dfr$ is
also noncrossing. For this reason, $\NCliMonoid \Mca$ is a suboperad of
$\CliMonoid \Mca$.
\medbreak

\begin{Proposition} \label{prop:presentation_nc}
    Let $\Mca$ be a monoid. Then, the operad $\NCliMonoid \Mca$ admits
    the presentation $(\GeneratingSet, \RelationSpace)$ where
    \begin{equation}
        \GeneratingSet :=
        \left\{
            \CliqueOne{x} : x \in \bar{\Mca}
        \right\}
        \sqcup
        \left\{\TriangleEmpty\right\}
    \end{equation}
    and $\RelationSpace$ is the space generated by
    \begin{subequations}
    \begin{equation}
        \Corolla\Par{\TriangleEmpty}
        \circ_1 \Corolla\Par{\TriangleEmpty}
        -
        \Corolla\Par{\TriangleEmpty}
        \circ_2 \Corolla\Par{\TriangleEmpty},
    \end{equation}
    \begin{equation}
        \Corolla\Par{\CliqueOne{x}}
        \circ_1 \Corolla\Par{\CliqueOne{y}}
        -
        \Corolla\Par{\CliqueOne{x \Product y}},
        \qquad x, y \in \bar{\Mca}.
    \end{equation}
    \end{subequations}
\end{Proposition}
\medbreak

From the presentation of $\NCliMonoid \Mca$ provided by
Proposition~\ref{prop:presentation_nc} and the one of the operad of
words $\T \Mca$ provided by Proposition~\ref{prop:presentation_t}, we
can observe that $\T \Mca$ is a quotient of $\NCliMonoid \Mca$. The
linear map $\phi : \NCliMonoid \Mca \to \T \Mca$ satisfying, for any
$x \in \bar{\Mca}$,
\begin{subequations}
\begin{equation}
    \phi\Par{\CliqueOne{x}} = x \in \T \Mca(1),
\end{equation}
\begin{equation}
    \phi\Par{\TriangleEmpty} = \Unit\Unit \in \T \Mca(2)
\end{equation}
\end{subequations}
extends in a unique way into an operad morphism and this morphism is
surjective.
\medbreak

All the operads of graphs presented in
Section~\ref{subsec:operads_graphs} can be constructed directly or as
suboperads of $\CliMonoid \Mca$ or $\CliMagma \Mca$ for suitable monoids
or unitary magmas~$\Mca$.
\medbreak

\SkipTocEntry\section*{Bibliographic notes}

\SkipTocEntry\subsection*{About operads, pre-Lie systems, and
combinatorics}
The theory of operads arose first in the 1970s in the field of
algebraic topology through the works of May~\cite{May72} and Boardman
and Vogt~\cite{BV73}. The first motivation was to study loop spaces. By
``operad'', most of the authors mean what we call symmetric operad, that
is a nonsymmetric operad wherein each subspace consisting in the
elements of arity $n$ is endowed with the action of the symmetric group
$\SymmetricGroup(n)$ (see forthcoming
Section~\ref{subsec:symmetric_operads} of
Chapter~\ref{chap:generalizations}). Nonsymmetric operads appeared a
little earlier in the work of Gerstenhaber under the name of pre-Lie
systems~\cite{Ger63}. To be more precise, a pre-Lie system is a space
endowed with products $\circ_i$ which are series associative
(see~\eqref{equ:operad_axiom_assoc_series}) and parallel associative
(see~\eqref{equ:operad_axiom_assoc_parallel}), but not necessarily
unital (see~\eqref{equ:operad_axiom_unit}). The theory of operads has
somewhat been neglected in the next twenty years following its
discovery but it was put back on the front of the stage in the
1990s~\cite{Lod96}. At this moment, an increasing number of
combinatorists began to take an interest in the subject and several
works relating combinatorics and operads were performed. One can cite,
for instance,~\cite{MY91} dealing with Möbius species and compositions of
trees, \cite{Lod01,CL01,Lod08,Gir15a,Gir16c,Gir17b} where operads on
many combinatorial families are defined,
and~\cite{Lod05,Cha06,Liv06,CL07,CG14} where operad structures lead to
the discovery of algebraic and combinatorial properties. The main
philosophy is here twofold: on the one hand, the structure thereby added
on combinatorial families enables to see these in a new light, and on
the other, techniques coming from combinatorics lead to establish
algebraic properties of operads and the category of algebras they
encode. Classical and complementary references about operads
are~\cite{Mar08,Cha08,LV12,Men15,Yau16}.
\medbreak

\SkipTocEntry\subsection*{About set-operads}
Due to the fact that their linear structure can be forgotten,
set-operads form a class of operads which is in some sense simpler than
the class of general ones. Despite this apparent simplicity, set-operads
remain very rich structures and, as suggested in
Section~\ref{sec:main_operads}, a lot of operads appearing in
combinatorics are set-operads. Moreover, as a consequence of the lack of
linear structure, there are simple techniques to establish presentations
by generators and relations of set-operads by using rewrite systems on
trees (as exposed in
Section~\ref{subsubsec:presentation_rewrite_systems}). Computer
exploration is a crucial tool in this context. For instance, the
works~\cite{CG14,Gir16a,Gir16c,Gir17b,CCG18} use the computer to
conjecture orientations of spaces of relations (used as a prerequisite
of Theorem~\ref{thm:presentation_operads}). Besides, we exposed in
Section~\ref{subsubsec:set_operads} two special notions about
set-operads: the one of rooted operads is, up to a slight variation,
introduced by Chapoton in~\cite{Cha14} as a tool to study series on
operads, and the one of basic operads is due to Vallette~\cite{Val07}
and intervenes as a prerequisite for a tool for showing that an operad
is Koszul.
\medbreak

\SkipTocEntry\subsection*{About Koszul duality and Koszulity}
Koszul duality for binary and quadratic operads has been introduced by
Ginzburg and Kapranov~\cite{GK94}. This duality is an extension of the
Koszul duality of quadratic associative algebras~\cite{Pri70}. The
so-called rebirth of the operads in the 1990s~\cite{Lod96} was in part
due to this duality. Note that the duality exposed in
Section~\ref{subsubsec:koszul_duality} concerns only nonsymmetric
operads but the theory also includes the case of symmetric operads.
Besides, the definition of the Koszul property for an operad consisting
in asking for the acyclicity of its Koszul complex (see, for
instance,~\cite{LV12}) admits several reformulations. A first criterion
is due to Vallette~\cite{Val07} (see also~\cite{Men15}) passing by the
construction of a family of posets from an operad and showing that they
are Cohen-Macaulay (see, for instance,~\cite{BGS82}). The criterion 
using convergent rewrite systems exhibited in
Proposition~\ref{prop:koszulity_criterion_pbw} is a consequence of the
work of Dotsenko and Khoroshkin~\cite{DK10}. The concept of
Poincaré-Birkhoff-Witt bases (which, as explained in
Section~\ref{subsubsec:koszulity}, form bases of Koszul operads) arises
in the work of Hoffbeck~\cite{Hof10}.
\medbreak

\SkipTocEntry\subsection*{About the presented examples of operads}
Let us now finally give some details about the operads reviewed in
Section~\ref{sec:main_operads}. The operad of permutations $\Per$ is in
some cases called ``associative operad''~\cite{AL07}. Indeed, $\Per$ can
be seen as the regularization of the associative operad $\As$ (see
forthcoming Section~\ref{subsec:symmetric_operads} of
Chapter~\ref{chap:generalizations}). The diassociative operad $\Dias$
has been  introduced by Loday in~\cite{Lod01} within its presentation by
generators and relations, and its realization in terms of the elements
$\Efr_{n, k}$ is due to Chapoton~\cite{Cha05}. The construction $\T$,
associating an operad with any monoid has been brought
in~\cite{Gir15a}. As explained, it provides alternative constructions of
$\As$ and $\Dias$, and also for the triassociative operad $\Trias$
(see~\cite{LR04}). As illustrated in~\cite{Gir15a}, the construction
$\T$ can be used in order to build operads on a large range of
combinatorial graded collections (words, permutations, $k$-ary trees,
integer compositions, directed animals, {\em etc.}). We have presented
here in this context only the operad $\Motz$ of Motzkin words, obtained
as a suboperad of an operad obtained from the construction $\T$.
Besides, the duplicial operad appeared in~\cite{Lod08} in the context of
the study of types of bialgebras. The dendriform operad $\Dendr$ was
defined as the Koszul dual of $\Dias$ in~\cite{Lod01}. The presentations
of $\Dup$ and $\Dendr$ are very similar and they share the same graded
space of binary trees. There are also two alternative realizations of
$\Dendr$ in terms of rational functions~\cite{Cha07,Lod10}. The
bicolored Schröder tree operad $\BS$ was considered in~\cite{LR06} under
the name $2as$. The pre-Lie operad $\PreLie$ has been defined by
Chapoton and Livernet~\cite{CL01} as the operad such that algebras of
the category it encodes are pre-Lie algebras (see
Section~\ref{subsec:pre_lie_algebras} of Chapter~\ref{chap:algebra}).
This operad is usually studied as a symmetric operad but its
nonsymmetric version has the interesting property to be
free~\cite{BL10}. The nonassociative permutative operad $\NAP$ has been
introduced in~\cite{Liv06} as a symmetric operad. Unlike $\PreLie$,
$\NAP$ is not free as a nonsymmetric operad (it is easy to find a
nontrivial relation in degree $2$ for instance). Some links between
$\PreLie$ and $\NAP$ have been exploited in~\cite{Sai14}. The operad of
noncrossing trees $\NCT$ was defined in~\cite{Cha07} as a suboperad of a
bigger operad $\Mould$, the operad of mould. The algebras over $\NCT$
are sometimes called $\mathrm{L}$-algebras and have been studied
in~\cite{Ler11}. In the same text~\cite{Cha07}, a generalization
of $\NCT$ involving noncrossing plants was brought, that are
combinatorial objects defined as configurations of chords satisfying
some conditions. The operad of bicolored noncrossing configurations
$\BNC$, introduced in~\cite{CG14}, is a further generalization of this
latter. The operad $\BNC$ contains as suboperads generated by binary
elements the operads of noncrossing plants, of noncrossing trees, the
dipterous operad~\cite{LR03}, and the operad of bicolored Schröder
trees. The gravity operad $\Grav$ is, as a symmetric operad, defined by
Getzler~\cite{Get94}. It has been studied as a nonsymmetric one
in~\cite{AP17}. The general constructions $\CliMonoid$, $\CliMagma$, and
$\NCliMonoid$, introduced in~\cite{Gir17b}, produce operads on
configurations of chords. These constructions can be used to provide
alternative realizations of all the operads presented in
Section~\ref{subsec:operads_graphs}, and also of some other, as the ones
of multi-tildes and double multi-tildes~\cite{LMN13,GLMN16} coming from
a context of formal language theory~\cite{CCM11}. Finally, let us
mention that a general reference about some of these operads (and some
others) is~\cite{Zin12}, where a large number of  morphisms between them
are referenced.
\medbreak


\chapter{Applications and generalizations} \label{chap:generalizations}
This last chapter is devoted to review some applications of the theory
of operads for enumerative prospects. To this aim, we present formal
power series on operads, generalizing usual generating series. We also
provide an overview on enrichments of operads: colored operads, cyclic
operads, symmetric operads, and pros.
\medbreak

\section{Series on operads} \label{sec:series_operads}
We consider here the notion of spaces of formal power series on
collections, forming a generalization of usual generating series. Any
product~$\Product$ (in the sense of
Section~\ref{subsubsec:collections_with_products} of
Chapter~\ref{chap:collections}) on a collection $C$ gives rise to a
product on the series on $C$, leading potentially to the discovery of
enumerative properties on the objects of $C$. We shall present how
associative and graded products on graded collections lead to
generalizations of the usual multiplication product of generating
series, and how full composition maps of set-operads lead to
generalizations of the usual composition product of generating series.
\medbreak

\subsection{Series on algebraic structures}
We introduce now series on collections, which are intuitively possibly
infinite formal sums of objects of a collection. Elementary definitions
about these series are reviewed here.
\medbreak

\subsubsection{Series spaces}
Let $C$ be an $I$-collection. A \Def{series on $C$} (or, for short,
a \Def{$C$-series}) is a map $\Fbf : C \to \K$. The \Def{coefficient}
$\Fbf(x)$ of $x \in C$ in $\Fbf$ is denoted by $\Angle{x, \Fbf}$. The
set of all $C$-series is denoted by $\K \AAngle{C}$. This set
$\K \AAngle{C}$ is endowed with the following two operations. First, the
\Def{addition} $\Fbf_1 + \Fbf_2$ of two $C$-series $\Fbf_1$ and $\Fbf_2$
is defined, for any $x \in C$, by
\begin{math}
    \Angle{x, \Fbf_1 + \Fbf_2} := \Angle{x, \Fbf_1} + \Angle{x, \Fbf_2}.
\end{math}
Second, the \Def{scalar multiplication} of a $C$-series $\Fbf$ by
$\lambda \in \K$ is defined, for any $x \in C$, by
\begin{math}
    \Angle{x, \lambda \ExtProd \Fbf} := \lambda \Angle{x, \Fbf}.
\end{math}
Endowed with these two operations, $\K \AAngle{C}$ is a $\K$-vector
space, named \Def{series space on $C$} (or, for short,
\Def{$C$-series space}).
\medbreak

Observe that $C$-polynomials (see Section~\ref{sec:polynomials} of
Chapter~\ref{chap:algebra}) are particular $C$-series and that
$\K \Angle{C}$ is a subspace of $\K \AAngle{C}$. One among the crucial
differences between $\K \AAngle{C}$ and $\K \Angle{C}$ is that this last
admits $C$ as a basis while $\K \AAngle{C}$ has no explicit basis. As a
side remark, a $C$-series can be seen as a linear form on
$\K \Angle{C}$. For this reason, $\K \AAngle{C}$ can be seen as the
(usual) dual space of $\K \Angle{C}$ (this has not to be confused with
the dual of a combinatorial polynomial space considered in
Section~\ref{subsubsec:duality_polynomial_spaces} of
Chapter~\ref{chap:algebra}).
\medbreak

For any subcollection $X$ of $C$, the \Def{characteristic series} of $X$
is the $C$-series $\Charac(X)$ defined, for any $x \in C$, by
$\Angle{x, \Charac(X)} := 1$ for all $x \in X$ and
$\Angle{y, \Charac(X)} := 0$ for all $y \in C \setminus X$. By using now
the linear structure of $\K \AAngle{C}$, any $C$-series $\Fbf$ can be
expressed as the possibly infinite sum
\begin{equation}
    \Fbf = \sum_{x \in C} \Angle{x, \Fbf} \ExtProd \Charac(\{x\}),
\end{equation}
which is denoted, by a slight abuse of notation, by
\begin{equation} \label{equ:sum_notation_series}
    \Fbf = \sum_{x \in C} \Angle{x, \Fbf} x.
\end{equation}
The notation~\eqref{equ:sum_notation_series} for $\Fbf$ as a (possibly
infinite) linear combination of objects of $C$ is the \Def{infinite
sum notation} of $C$-series.
\medbreak

\subsubsection{Generating series and products}
By setting that $\{t\}$ is a graded collection wherein $t$ is an atomic
object, $\K \AAngle{\Multiset(\{t\})}$ is the space of the usual
generating series. Indeed, by denoting by $t^n$ each object
$\lbag t, \dots, t \rbag$ of $\Multiset(\{t\})$ of size $n \in \N$,
any element $\Fbf$ of
$\K \AAngle{\Multiset(\{t\})}$ is by definition of the form
\begin{equation}
    \Fbf = \sum_{n \in \N} \Angle{t^n, \Fbf} t^n.
\end{equation}
To not overload the notation, we shall write $\K \AAngle{t}$ for
$\K \AAngle{\Multiset(\{t\})}$.
\medbreak

The usual multiplication (resp. composition) of generating series is
denoted by $\Conc$ (resp. $\circ$). In this way, by setting
$1 := \Charac\Par{\left\{t^0\right\}}$, $\Par{\K \AAngle{t}, \Conc, 1}$
is a unitary associative algebra. Moreover, by denoting by
$t\K \AAngle{t}$ the subspace of $\K \AAngle{t}$ of all the series
$\Fbf$ such that $\Angle{t^0, \Fbf} = 0$,
$\Par{t\K \AAngle{t}, \circ, t}$ is a unitary associative algebra.
\medbreak

\subsubsection{Index series}
When $C$ is combinatorial, the \Def{index series} of $C$ is
the $I$-series $\IndexSeries(C)$, where $I$ is seen as a simple
collection, defined by
\begin{equation}
    \IndexSeries(C)
    := \sum_{x \in C} \Index(x)
    = \sum_{i \in I} \# C(i) \, i.
\end{equation}
Since the coefficient $\Angle{i, \IndexSeries(C)}$ is the number of
elements of index $i \in I$ of $C$, the series $\IndexSeries(C)$ encodes
enumerating data about $C$. It is then worthwhile to provide ways of
expressing $\IndexSeries(C)$ in order to compute its coefficients. We
shall besides consider in the sequel index series of colored operads and
of pros as analogs of usual Hilbert series.
\medbreak

Moreover, let $\omega : I \to \N$ be a map. The
\Def{$\omega$-evaluation map} is the map
\begin{equation}
    \Eval^\omega : \K \AAngle{C} \to \K \AAngle{t}
\end{equation}
defined, for any $C$-series $\Fbf$, by
\begin{equation}
    \Eval^\omega(\Fbf) :=
    \sum_{x \in C} \Angle{x, \Fbf} t^{\omega(x)}.
\end{equation}
This series is well-defined if each fiber $\omega^{-1}(n)$ is finite
for any $n \in \N$ and is called the \Def{$\omega$-evaluation} of
$\Fbf$. When $C$ is combinatorial and graded, one has
\begin{equation} \label{equ:eval_graded_charac_series}
    \Eval^\Identity(\IndexSeries(C))
    = \Eval^{|-|}(\Charac(C))
    = \GeneratingSeries_C(t)
\end{equation}
where $\Identity$ is the identity map on the index set $I = \N$ of $C$,
and $|-|$ is the size function of $C$. Recall that
$\GeneratingSeries_C(t)$ denotes the generating series of~$C$.
\medbreak

\subsubsection{Products on series}
Assume that the $I$-collection $C$ is endowed with a product
\begin{equation}
   \Product : C\Par{J_1} \times \dots \times C\Par{J_p} \to C
\end{equation}
where $p \in \N$ and $J_1$, \dots, $J_p$ are nonempty subsets of $I$
(see Section~\ref{subsubsec:collections_with_products} of
Chapter~\ref{chap:collections}). Then, let the linear map
\begin{equation}
    \bar{\Product} :
    \K \AAngle{C}^{\otimes p} \to \K \AAngle{C}
\end{equation}
defined, for any $\Fbf_1, \dots, \Fbf_p \in \K \AAngle{C}$ and
$x \in C$, by
\begin{equation} \label{equ:extension_product_on_series}
    \Angle{x, \bar{\Product}\Par{\Fbf_1, \dots, \Fbf_p}} :=
    \sum_{\substack{
        y_1, \dots, y_p \in C \\
        \Product\Par{y_1, \dots, y_p} = x
    }}
    \prod_{k \in [p]} \Angle{y_k, \Fbf_k}.
\end{equation}
In other terms, by using the sum notation of series,
\begin{equation}
    \bar{\Product}\Par{\Fbf_1, \dots, \Fbf_p} =
    \sum_{x_k \in C\Par{J_k}, k \in [p]} \enspace
    \prod_{k \in [p]} \Angle{x_k, \Fbf_k} \,
    \Product\Par{x_1, \dots, x_p}.
\end{equation}
We call $\bar{\Product}$ the \Def{series extension} of $\Product$. In
this way, series extensions of products on $C$ allow to translate
set-theoretic algebraic structures on $C$ into products on series.
\medbreak

For instance, by considering the product $\Product$ on
$\Multiset(\{t\})$ defined by $t^n \Product t^m := t^{n + m}$ for any
$n, m \in \N$, the series extension of $\Product$ on $\K \AAngle{t}$ is
the multiplication $\Conc$ of generating series.
\medbreak

\subsection{Generalizing series multiplication}
We consider series extensions of binary products on graded collections
satisfying some conditions and explain how they provide generalizations
of the multiplication product of generating series.
\medbreak

\subsubsection{Series on monoids and multiplication}
Let $C$ be a graded collection endowed with a graded complete binary
associative product $\Product$. In this case, its series extension
$\bar{\Product}$ endows $\K \AAngle{C}$ with the structure of an
associative algebra. Moreover, when $\Product$ admits a unit $\Unit$,
$C$ is a monoid and the series $\UnitSeries := \Charac(\{\Unit\})$ of
$\K \AAngle{C}$ is the unit of $\bar{\Product}$.
\medbreak

\begin{Proposition} \label{prop:series_monoids_multiplication}
    Let $C$ be a graded collection endowed with a graded complete binary
    associative product $\Product$ admitting a unit. Then, the map
    $\Eval^{|-|}$ is a unitary associative algebra morphism between
    $\Par{\K \AAngle{C}, \bar{\Product}, \UnitSeries}$ and
    $\Par{\K \AAngle{t}, \Conc, 1}$. Moreover, $\Eval^{|-|}$ is
    surjective when $C(n) \ne \emptyset$ for all $n \in \N$.
\end{Proposition}
\medbreak

Proposition~\ref{prop:series_monoids_multiplication} implies in
particular that if one obtains a nontrivial expression for the
characteristic series $\Charac(C)$ of $C$ by using the sum of series and
the product $\bar{\Product}$, its $|-|$-evaluation will provide a
nontrivial expression for the generating series $\GeneratingSeries_C(t)$
of $C$. We shall present examples in the further sections.
\medbreak

\subsubsection{A monoid of paths}
Let  us consider an example of series on monoids and an application to
enumeration. We call \Def{path} any nonempty word $u$ on $\N$ and we
denote by $\ColPath$ be the graded collection of all paths, where the
size of a path is its length as a word minus $1$. This collection is
endowed with the complete binary product $\Product$ defined, for any
$u \in \ColPath(n - 1)$ and $v \in \ColPath(m - 1)$, by
\begin{equation}
    u \Product v :=
    \begin{cases}
        \uparrow_k(u(1) \dots u(n - 1)) \, \Conc \, v
        & \mbox{if } v(1) \geq u(n), \\
        u \, \Conc \, \uparrow_k(v(2) \dots v(m)) & \mbox{otherwise},
    \end{cases}
\end{equation}
where $\Conc$ is the concatenation product of words,
$k := |u(n) - v(1)|$, and for any path $w$, $\uparrow_k(w)$ is the path
obtained by incrementing all the letters of $w$ by $k$. For instance,
\begin{equation} \label{equ:example_product_paths}
    \ColA{0101121} \Product \ColD{210011}
    =
    \ColA{121223} \ColF{2} \ColD{10011}.
\end{equation}
By depicting a path through its graph in the quarter plane (that is, by
drawing points $(i - 1, u(i))$ for all positions $i$ and by connecting
all pairs of adjacent points by lines),
\eqref{equ:example_product_paths} becomes
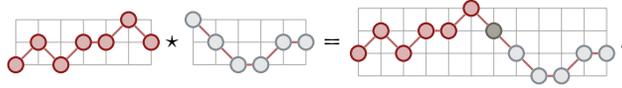
\begin{equation*}
\,.
\end{equation*}
In intuitive terms, the product $\Product$ consists in concatenating the
paths by superimposing the last letter of the first operand with the
first letter of the second. Observe that the path $0$, denoted by
$\UnitPath$, has size zero in $\ColPath$ and is the unit of $\Product$.
\medbreak

\begin{Proposition} \label{prop:graded_monoid_paths}
    The triple $\Par{\ColPath, \Product, \UnitPath}$ is a graded monoid.
\end{Proposition}
\medbreak

By Propositions~\ref{prop:series_monoids_multiplication}
and~\ref{prop:graded_monoid_paths}, we can consider series of paths and
the series extension of the product of paths. To specify particular
families of paths, let us introduce the following tool. For any
subcollection $X$ of $\ColPath$, the \Def{$\Product$-Kleene star}
$X^{\Product_\star}$ of $X$ is the subcollection defined by
\begin{equation}
    X^{\Product_\star} :=
    \bigcup_{p \in \N} \,
    \underbrace{X \Product \cdots \Product X}
    _{p \mbox{ \footnotesize terms}}.
\end{equation}
In other words, $X^{\Product_\star}$ is the submonoid of $\ColPath$
generated by $X$. As a collection, $X^{\Product_\star}$ contains all the
paths obtained by concatenating elements of~$X$.
\medbreak

\subsubsection{Enumeration of families of paths}
Let us first study the collection $\ColPathSch$ of
\Def{Schröder paths}, that are the paths of
\begin{math}
    \left\{\PathUp, \PathStableStable, \PathDown
    \right\}^{\Product_\star}
\end{math}
starting and finishing by $0$. This collection is combinatorial and its
characteristic series is
\begin{equation}
    \Charac\Par{\ColPathSch}:=
    \UnitPath
    +
    \PathStableStable
    +
    \PathPeak
    +

    +
    \cdots~.
\end{equation}
By reasoning on the non-ambiguous decomposition of Schröder paths, one
can establish that this series satisfies the nontrivial relation
\begin{equation} \label{equ:relation_series_schroder_paths}
    \Charac\Par{\ColPathSch} =
    \UnitPath
    + \PathStableStable \, \bar{\Product} \, \Charac\Par{\ColPathSch}
    + \PathUp \, \bar{\Product} \, \Charac\Par{\ColPathSch}
        \, \bar{\Product} \, \PathDown \, \bar{\Product} \,
        \Charac\Par{\ColPathSch}.
\end{equation}
The $|-|$-evaluation of the left and right members
of~\eqref{equ:relation_series_schroder_paths} leads, by
using~\eqref{equ:eval_graded_charac_series}, to the algebraic relation
\begin{equation}
    \GeneratingSeries_{\ColPathSch}(t)
    = 1
    + t^2 \, \GeneratingSeries_{\ColPathSch}(t)
    + t^2 \, {\GeneratingSeries_{\ColPathSch}(t)}^2
\end{equation}
for the generating series of $\ColPathSch$.
\medbreak

We can use similar mechanisms to obtain expressions for the generating
series of other families of paths. Let us consider three of these.
\medbreak

\paragraph{Dyck paths}
The characteristic series of the collection $\ColPathDyck$ of
\Def{Dyck paths}, that are paths of
\begin{math}
    \left\{\PathUp, \PathDown \right\}^{\Product_\star}
\end{math}
starting and finishing by $0$ satisfies
\begin{equation}
    \Charac\Par{\ColPathDyck} = \UnitPath
    + \PathUp \, \bar{\Product} \, \Charac\Par{\ColPathDyck}
    \, \bar{\Product} \, \PathDown \, \bar{\Product} \,
    \Charac\Par{\ColPathDyck}.
\end{equation}
\medbreak

\paragraph{Motzkin paths}
The characteristic series of the collection $\ColPathMotz$ of
\Def{Motzkin paths}, that are paths of
\begin{math}
    \left\{\PathUp, \PathStable, \PathDown \right\}^{\Product_\star}
\end{math}
starting and finishing by $0$ satisfies
\begin{equation} \label{equ:series_motz_monoid}
    \Charac\Par{\ColPathMotz} = \UnitPath
    + \PathStable \, \bar{\Product} \, \Charac\Par{\ColPathMotz}
    + \PathUp \, \bar{\Product} \, \Charac\Par{\ColPathMotz}
    \, \bar{\Product} \, \PathDown \, \bar{\Product} \,
    \Charac\Par{\ColPathMotz}.
\end{equation}
\medbreak

\paragraph{Fibonacci paths}
The characteristic series of the collection $\ColPathFib$ of
\Def{Fibonacci paths}, that are paths of
\begin{math}
    \left\{\PathStable, \PathPeak\right\}^{\Product_\star}
\end{math}
satisfies
\begin{equation}
    \Charac\Par{\ColPathFib} = \UnitPath
    + \PathStable \, \bar{\Product} \, \Charac\Par{\ColPathFib}
    + \PathPeak \, \bar{\Product} \, \Charac\Par{\ColPathFib}.
\end{equation}
\medbreak

\subsection{Generalizing series composition} \label{subsec:series_compo}
We consider now series on combinatorial collections endowed with the
structure of set-operads and explain how they provide generalizations of
the composition product of generating series.
\medbreak

\subsubsection{Series on operads and composition}
Let $C$ be a set-operad (see Section~\ref{subsubsec:set_operads} of
Chapter~\ref{chap:operads}). In this case, the series extensions
$\bar{\circ}$ of its full composition maps $\circ$ satisfy, for any
$n \in \N_{\geq 1}$ and any $C$-series $\Fbf$, $\Gbf_1$, \dots,
$\Gbf_n$,
\begin{equation}
    \Angle{x, \Fbf \, \bar{\circ} \, \left[\Gbf_1, \dots, \Gbf_n\right]}
    =
    \sum_{\substack{
        y \in C(n) \\
        z_1, \dots, z_n \in C \\
        y \circ \left[z_1, \dots, z_n\right] = x
    }}
    \Angle{y, \Fbf} \prod_{i \in [n]} \Angle{z_i, \Gbf_i}.
\end{equation}
Let now $\Compo$ be the binary product on $\K \AAngle{C}$ defined as the
sum of all the series extensions of the full composition maps of $C$.
More precisely, for any $C$-series $\Fbf$ and $\Gbf$,
\begin{equation}
    \Angle{x, \Fbf \Compo \Gbf} :=
    \sum_{p \in \N_{\geq 1}}
    \Angle{x,
    \Fbf \, \bar{\circ} \,
    \left[
    \underbrace{\Gbf, \dots, \Gbf}_{p \mbox{ \footnotesize terms}}
    \right]}
    =
    \sum_{\substack{
        y \in C(n), n \in \N \\
        z_1, \dots, z_n \in C \\
        y \circ \left[z_1, \dots, z_n\right] = x
    }}
    \Angle{y, \Fbf} \prod_{i \in [n]} \Angle{z_i, \Gbf}.
\end{equation}
As immediate observations, remark that $\Compo$ is linear on left, is
not linear on the right, and admits $\UnitSeries := \Charac(\{\Unit\})$
as a left and a right unit, where $\Unit$ is the operad unit of~$C$.
\medbreak

\begin{Proposition} \label{prop:series_operads_eval_composition}
    Let $C$ be a set-operad, $x \in C(n)$, $n \in \N_{\geq 1}$, and
    $\Gbf_1$, \dots, $\Gbf_n$ be $C$-series. Then,
    \begin{equation}
        \Eval^{|-|}\Par{x \, \bar{\circ} \,
        \left[\Gbf_1, \dots, \Gbf_n\right]}
        = \prod_{i \in [n]} \Eval^{|-|}\Par{\Gbf_i}.
    \end{equation}
\end{Proposition}
\medbreak

\begin{Proposition} \label{prop:series_operads_complete_composition}
    Let $C$ be a set-operad. Then, the map $\Eval^{|-|}$ is a unitary
    associative algebra morphism between
    $\Par{\K \AAngle{C}, \Compo, \UnitSeries}$ and
    $\Par{t\K \AAngle{t}, \circ, t}$. Moreover, $\Eval^{|-|}$ is
    surjective when $C(n) \ne \emptyset$ for all $n \in \N_{\geq 1}$.
\end{Proposition}
\medbreak

Propositions~\ref{prop:series_operads_eval_composition}
and~\ref{prop:series_operads_complete_composition} imply in particular
that if one obtains a nontrivial expression for the characteristic
series $\Charac(C)$ of $C$ by using the sum of series, and the products
$\bar{\circ}$ and $\Compo$, its $|-|$-evaluation will provide a
nontrivial expression for the generating series $\GeneratingSeries_C(t)$
of $C$. We shall present examples in the further sections.
\medbreak

\subsubsection{Enumeration of Motzkin paths}
Let us consider the set-operad $\Motz$ of Motzkin words introduced in
Section~\ref{subsubsec:motzkin_operad} of Chapter~\ref{chap:operads}. By
using the operad structure on the underlying graded collection
$\ColPathMotz'$ of $\Motz$, one obtains the nontrivial relation
\begin{equation} \label{equ:series_motz_operad}
    \Charac\Par{\ColPathMotz'}
    = \UnitPath
    +
    \PathStable
    \, \bar{\circ} \, \left[\UnitPath, \Charac\Par{\ColPathMotz'}\right]
    + \PathPeak
    \, \bar{\circ} \, \left[\UnitPath, \Charac\Par{\ColPathMotz'},
    \Charac\Par{\ColPathMotz'}\right]
\end{equation}
for the characteristic series of $\ColPathMotz'$. Observe
that~\eqref{equ:series_motz_operad} and~\eqref{equ:series_motz_monoid}
are two equivalent expressions for enumerating collections of Motzkin
paths (yet with different size functions). Nevertheless,
\eqref{equ:series_motz_operad} has the advantage of not requiring
the definition of a general algebraic structure of paths
(in~\eqref{equ:series_motz_operad}, all the terms are series of Motzkin
paths, while in~\eqref{equ:series_motz_monoid}, $\PathUp$
and~$\PathDown$ are not Motzkin paths).
\medbreak

\subsubsection{Enumeration of noncrossing trees}
Let us consider the set-operad $\NCT$ of noncrossing trees introduced
in Section~\ref{subsubsec:nct_operad} of Chapter~\ref{chap:operads}.
By using the operad structure on the underlying collection $\ColNCT$ of
$\NCT$, one obtains the nontrivial relation
\begin{equation} \label{equ:series_nct_operad_1}
    \Charac\Par{\ColNCT} = \UnitSolidClique
    + \TriangleLeft \, \bar{\circ} \,
    \left[\Charac\Par{\ColNCT}, \Charac\Par{\ColNCT}\right]
    + \TriangleRight \, \bar{\circ} \,
    \left[\Fbf, \Charac\Par{\ColNCT}\right],
\end{equation}
where $\Fbf$ is the $\ColNCT$-series satisfying
\begin{equation} \label{equ:series_nct_operad_2}
    \Fbf = \UnitSolidClique
    + \TriangleRight \, \bar{\circ} \,
    \left[\Fbf, \Charac\Par{\ColNCT}\right],
\end{equation}
for the characteristic series of~$\ColNCT$.
\medbreak

\section{Enriched operads} \label{sec:enriched_operads}
Three enrichments of nonsymmetric operads are presented here: colored
operads, cyclic operads, and symmetric operads.
\medbreak

\subsection{Colored operads} \label{subsec:colored_operads}
We begin by introducing colored operads. These variations of operads
involve colored collections. Intuitively, each element of a colored
operad has a color for its output and colors for each of its inputs. The
partial composition of two elements is defined if and only if the colors
of the involved input and output coincide.
\medbreak

\subsubsection{Colored polynomial spaces}
Given a $\CFr$-colored collection $C$ (see
Section~\ref{subsubsec:colored_collections} of
Chapter~\ref{chap:collections}), the polynomial space $\K \Angle{C}$ is
said \Def{$\CFr$-colored} and is endowed with the maps $\Out$ and $\In$
associating with each nonzero homogeneous element $f$ of $\K \Angle{C}$,
respectively, its output color $\Out(f)$ and its word of input colors
$\In(f)$ (see the aforementioned section). The \Def{arity} $|f|$ of an
$(a, u)$-homogeneous element $f$ of $\K \Angle{C}$, where $(a, u)$ is a
$\CFr$-colored index, is the length $|u|$ of the word $u$.
Alternatively, the arity of $f$ is the degree of $f$ in $\K \Angle{C'}$
where $C'$ is the graduation of $C$ (see
Section~\ref{subsubsec:casting_collections} of
Chapter~\ref{chap:collections}). Moreover, to not overload the notation,
we denote by $\K \Angle{C}(a, u)$ the homogeneous component
$\K \Angle{C}((a, u))$ of $\K \Angle{C}$ for any $\CFr$-colored
index~$(a, u)$.
\medbreak

\subsubsection{Colored abstract operators}
We regard any homogeneous element $f$ of an augmented $\CFr$-colored
polynomial space $\K \Angle{C}$ as a \Def{colored abstract operator},
that is, an abstract operator wherein the output and each input are
associated with an element of $\CFr$. If $f$ is of arity $n$,
$\Out(f) = a$, and $\In(f) = u$, $f$ is depicted as
\begin{equation}
    \begin{tikzpicture}
        [xscale=.35,yscale=.4,Centering,font=\scriptsize]
        \node[Operator](x)at(0,0){\begin{math}f\end{math}};
        \node(r)at(0,2){};
        \node(x1)at(-3,-2){};
        \node(xn)at(3,-2){};
        \node[below of=x1,node distance=1mm](ex1)
            {\begin{math}1\end{math}};
        \node[below of=xn,node distance=1mm](exn)
            {\begin{math}n\end{math}};
        \draw[Edge](r)edge[]node[EdgeLabel]{\begin{math}a\end{math}}(x);
        \draw[Edge](x)
            edge[]node[EdgeLabel]{\begin{math}u(1)\end{math}}(x1);
        \draw[Edge](x)
            edge[]node[EdgeLabel]{\begin{math}u(n)\end{math}}(xn);
        \node[below of=x,node distance=9mm]
            {\begin{math}\dots\end{math}};
    \end{tikzpicture}\,.
\end{equation}
The output and input colors of $f$ are written onto the output and input
edges.
\medbreak

\subsubsection{Colored operads} \label{subsubsec:colored_operads}
Let $C$ be an augmented $\CFr$-colored collection. Let for all
$\CFr$-colored indexes $(a, u)$ and $(b, v)$, and $i \in [|u|]$ such
that $b = u(i)$, binary products of the form
\begin{equation}
   \circ_i^{((a, u), (b, v))} :
   \K \Angle{\BBrack{C(a, u), C(b, v)}_\times}
   \to \K \Angle{C}(a, u \mapsfrom_i v),
\end{equation}
where $u \mapsfrom_i v$ is the word obtained by replacing the $i$th
letter of $u$ by $v$. On abstract operators, these products
$\circ_i^{((a, u), (b, v))}$ behave as the products
$\circ_i^{\Par{|u|, |v|}}$ of operads (see
Section~\ref{subsubsec:partial_compo_maps} of
Chapter~\ref{chap:operads}) but with the addition of taking into account
of the output and input colors of the colored abstract operators.
Indeed, for any $f \in \K \Angle{C}(a, u)$ and
$g \in \K \Angle{C}(b, v)$, $f \circ_i^{((a, u), (b, v))} g$ is the
abstract operator
\begin{multline} \label{equ:partial_compostion_on_colored_operators}
\,.
\end{multline}
Let us emphasize the fact that these products require that the output 
color of $g$ is equal to the $i$th input color of $f$. By a slight abuse 
of notation, we shall sometimes omit the $((a, u), (b, v))$ in the 
notation of $\circ_i^{((a, u), (b, v))}$ in order to denote it in a more 
concise way by~$\circ_i$.
\medbreak

When for any objects $x$, $y$, and $z$ of $C$, the fact that the left
and right members of Relation~\eqref{equ:operad_axiom_assoc_series}
(resp. Relation~\eqref{equ:operad_axiom_assoc_parallel}) of
Chapter~\ref{chap:operads} are well-defined implies that they are equal,
$\circ_i$ is \Def{series associative} (resp.
\Def{parallel associative}). Moreover, assume that there exists a set of
elements $\left\{\Unit_a : a \in \CFr\right\}$ of arity $1$ of
$\K \Angle{C}$ such that for any $a \in \CFr$,
\begin{math}
    \Out\Par{\Unit_a} = a = \In\Par{\Unit_a}.
\end{math}
When for any object $x$ of $C$, the fact that, by replacing by $\Unit_a$
each occurrence of $\Unit$ in Relation~\eqref{equ:operad_axiom_unit} of
Chapter~\ref{chap:operads}, the first or the last members of the
relation are well-defined implies that they are equal to $x$, $\circ_i$
is \Def{unital}. We call in this case each $\Unit_a$, $a \in \CFr$, a
\Def{unit of color $a$}.
\medbreak

When the products $\circ_i$ are series associative, parallel
associative, and unital, the $\circ_i$ are called
\Def{partial composition maps}. A \Def{$\CFr$-colored operad} is a
$\CFr$-colored polynomial space $\K \Angle{C}$ endowed with partial
composition maps. The main algebraic notions presented in
Section~\ref{subsec:operads} of Chapter~\ref{chap:operads} for operads
(like full composition maps associated with partial composition maps,
morphisms, quotients, group of symmetries, set-operads, {\em etc.}) hold
straightforwardly for colored operads. When $\K \Angle{C}$ is
combinatorial, its \Def{Hilbert series} is the
$\CFr \times \CFr^+$-series
\begin{equation}
    \HilbertSeries_{\K \Angle{C}} :=
    \IndexSeries(C) = \sum_{(a, u) \in \CFr \times \CFr^+}
    \Par{\dim \K \Angle{C}(a, u)} \, (a, u).
\end{equation}
\medbreak

\subsubsection{Categorical point of view}
\label{subsubsec:categorical_operads}
Recall that a monoid $\Mca$ can be seen as a category with exactly one
object $\AtomElement$ wherein the elements of $\Mca$ are interpreted as
morphisms $\phi : \AtomElement \to \AtomElement$. In the same way, an
operad $\K \Angle{C}$ can be seen as a multicategory with exactly one
object $\AtomElement$. In this case, the elements of $\K \Angle{C}(n)$,
$n \in \N_{\geq 1}$, are interpreted as multimorphisms
$\phi : \AtomElement^n \to \AtomElement$. The full composition maps of
$\K \Angle{C}$ translate as the composition of multimorphisms.
\medbreak

In a similar way, a $\CFr$-colored operad $\K \Angle{C}$ can be seen as
a multicategory having $\CFr$ as set of objects. In this case, the
elements of $\K \Angle{C}(a, u)$ where $(a, u)$ is a $\CFr$-colored
index, are interpreted as multimorphisms
$\phi :u(1) \times \dots \times u(n) \to a$. The full composition maps
of $\K \Angle{C}$ translate as the composition of multimorphisms, where
the constraints imposed by the colors in $\K \Angle{C}$ become
constraints imposed by the domains and codomains of multimorphisms.
\medbreak

\subsubsection{Free colored operads}
Let $\CFr$ be a set of colors and $\GeneratingSet$ be an augmented
$\CFr$-colored collection. The \Def{free colored operad} over
$\GeneratingSet$ is the operad
\begin{equation}
    \FreeColoredOperad(\GeneratingSet)
    :=
    \K \Angle{\ColCST_\Leaf^\GeneratingSet},
\end{equation}
where $\ColCST_\Leaf^\GeneratingSet$ is the graded collection of all the
$\CFr$-colored $\GeneratingSet$-syntax trees (see
Section~\ref{subsec:colored_syntax_trees} of Chapter~\ref{chap:trees}).
The space $\FreeColoredOperad(\GeneratingSet)$ is endowed with the
linearizations of the partial grafting operations $\circ_i$,
$i \in \N_{\geq 1}$, defined in
Section~\ref{subsec:colored_syntax_trees} of Chapter~\ref{chap:trees}.
The unit of color $a$, $a \in \CFr$, of
$\FreeColoredOperad(\GeneratingSet)$ is the only $\CFr$-colored
$\GeneratingSet$-syntax tree $\LeafPicCST{a}$ of arity~$1$ and
degree~$0$ and having $a$ as output and input color.
\medbreak

\subsubsection{Example: bud operads} \label{subsubsec:bud_operads}
Given an operad $\K \Angle{C}$ and a set of colors $\CFr$, there is an
easy way to construct a $\CFr$-colored operad. Let
$\Bud_\CFr(\K \Angle{C})$ be the $\CFr$-colored space
$\K \Angle{\Coloration_\CFr(\K \Angle{C})}$ where $\Coloration_\CFr(C)$
denotes the $\CFr$-coloration of $C$
(see Section~\ref{subsubsec:coloration_collections} of
Chapter~\ref{chap:collections}). Let us endow $\Bud_\CFr(\K \Angle{C})$
with the partial composition maps defined linearly, for any objects
$(a, x, u)$ and $(b, y, v)$ of $\Coloration_\CFr(C)$, and $i \in [|u|]$
such that $b = u(i)$, by
\begin{equation} \label{equ:partial_compo_bud_operads}
    (a, x, u) \circ_i (b, y, v) :=
    \Par{a, x \circ_i y, u \mapsfrom_i v}
\end{equation}
where the second occurrence of $\circ_i$
in~\eqref{equ:partial_compo_bud_operads} is the partial composition map
of the operad $\K \Angle{C}$ and $\mapsfrom_i$ is the operation on words
on $\CFr$ defined in Section~\ref{subsubsec:colored_operads}.
\medbreak

\begin{Proposition} \label{prop:bud_operads}
    For any set of colors $\CFr$ and any operad $\K \Angle{C}$,
    $\Bud_\CFr(\K \Angle{C})$ is a $\CFr$-colored operad.
\end{Proposition}
\medbreak

We call $\Bud_\CFr(\K \Angle{C})$ the \Def{$\CFr$-bud operad} of
$\K \Angle{C}$. For instance, one has in
$\Bud_{\{\Asf, \Bsf\}}\Par{\As}$ (where $\As$ is the associative operad
defined in Section~\ref{subsubsec:associative_operad} of
Chapter~\ref{chap:operads}) the partial composition
\begin{equation}
    \Par{\ColA{\Asf}, \Afr_4, \ColA{\Bsf}
    \ColF{\Bsf} \ColA{\Asf \Bsf}}
    \circ_2 \Par{\Bsf, \Afr_3, \ColD{\Asf \Asf \Bsf}}
    = \Par{\ColA{\Asf}, \Afr_6,
    \ColA{\Bsf} \ColD{\Asf \Asf \Bsf} \ColA{\Asf \Bsf}}.
\end{equation}
\medbreak

Let us observe that the bud operad of a free operad can be not free as
a colored operad. For instance, consider the $\{\Asf, \Bsf\}$-bud
operad of the magmatic operad $\Mag$ (see
Section~\ref{subsubsec:magmatic_operad} of Chapter~\ref{chap:operads}).
One has in $\Bud_{\{\Asf, \Bsf\}}(\Mag)$ among others the nontrivial
relation
\begin{equation}
    \Par{
    \Asf,
    \begin{tikzpicture}[xscale=.2,yscale=.15,Centering]
        \node[Leaf](0)at(0.00,-1.50){};
        \node[Leaf](2)at(2.00,-1.50){};
        \node[Node](1)at(1.00,0.00){};
        \draw[Edge](0)--(1);
        \draw[Edge](2)--(1);
        \node(r)at(1.00,2){};
        \draw[Edge](r)--(1);
    \end{tikzpicture}\,,
    \Asf \Asf}
    \circ_1
    \Par{
    \Asf,
    \begin{tikzpicture}[xscale=.2,yscale=.15,Centering]
        \node[Leaf](0)at(0.00,-1.50){};
        \node[Leaf](2)at(2.00,-1.50){};
        \node[Node](1)at(1.00,0.00){};
        \draw[Edge](0)--(1);
        \draw[Edge](2)--(1);
        \node(r)at(1.00,2){};
        \draw[Edge](r)--(1);
    \end{tikzpicture}\,,
    \Asf\Asf}
    =
    \Par{
    \Asf,
    \begin{tikzpicture}[xscale=.16,yscale=.18,Centering]
        \node[Leaf](0)at(0.00,-3.33){};
        \node[Leaf](2)at(2.00,-3.33){};
        \node[Leaf](4)at(4.00,-1.67){};
        \node[Node](1)at(1.00,-1.67){};
        \node[Node](3)at(3.00,0.00){};
        \draw[Edge](0)--(1);
        \draw[Edge](1)--(3);
        \draw[Edge](2)--(1);
        \draw[Edge](4)--(3);
        \node(r)at(3.00,1.75){};
        \draw[Edge](r)--(3);
    \end{tikzpicture}\,,
    \Asf \Asf \Asf}
    =
    \Par{
    \Asf,
    \begin{tikzpicture}[xscale=.2,yscale=.15,Centering]
        \node[Leaf](0)at(0.00,-1.50){};
        \node[Leaf](2)at(2.00,-1.50){};
        \node[Node](1)at(1.00,0.00){};
        \draw[Edge](0)--(1);
        \draw[Edge](2)--(1);
        \node(r)at(1.00,2){};
        \draw[Edge](r)--(1);
    \end{tikzpicture}\,,
    \Bsf \Asf}
    \circ_1
    \Par{
    \Bsf,
    \begin{tikzpicture}[xscale=.2,yscale=.15,Centering]
        \node[Leaf](0)at(0.00,-1.50){};
        \node[Leaf](2)at(2.00,-1.50){};
        \node[Node](1)at(1.00,0.00){};
        \draw[Edge](0)--(1);
        \draw[Edge](2)--(1);
        \node(r)at(1.00,2){};
        \draw[Edge](r)--(1);
    \end{tikzpicture}\,,
    \Asf\Asf}\,,
\end{equation}
implying that $\Bud_{\{\Asf, \Bsf\}}(\Mag)$ is not free as a colored
operad.
\medbreak

\subsection{Cyclic operads} \label{subsec:cyclic_operads}
We focus now on cyclic operads. These variations of operads involve
cyclic collections. Intuitively, in a cyclic operad, the output and the
inputs of the elements play an interchangeable role. This is due to
the fact that these structures are endowed with a map performing a
cyclic action on the inputs and outputs of its elements.
\medbreak

\subsubsection{Cyclic polynomial spaces}
A graded polynomial space $\K \Angle{C}$ is \Def{cyclic} if it is
endowed for all $n \in \N$ with unary products
\begin{equation}
    \CyclicAction_n : \K \Angle{C}(n) \to \K \Angle{C}(n)
\end{equation}
such that $\CyclicAction_n^{n + 1}$ is the identity map on
$\K \Angle{C}(n)$. We say that the $\CyclicAction_n$, $n \in \N$, are
\Def{cycle maps} of $\K \Angle{C}$. By a slight abuse of notation, we
shall sometimes omit the $n$ in the notation of~$\CyclicAction_n$ in
order to denote it in a more concise way by~$\CyclicAction$. As usual,
if $\phi : \K \Angle{C_1} \to \K \Angle{C_2}$ is a morphism between two
graded polynomial spaces $\K \Angle{C_1}$ and $\K \Angle{C_2}$, $\phi$
is a \Def{cyclic polynomial space morphism} if it commutes with the
cycle maps of $\K \Angle{C_1}$ and $\K \Angle{C_2}$.
\medbreak

Remark that when $C$ is a cyclic collection (see
Section~\ref{subsubsec:cyclic_collections} of
Chapter~\ref{chap:collections}), the linearizations of the cycle maps
of $C$ endow $\K \Angle{C}$ with the structure of a cyclic polynomial
space.
\medbreak

\subsubsection{Cyclic abstract operators}
We regard any homogeneous element $f$ of an augmented cyclic polynomial
space $\K \Angle{C}$ as a \Def{cyclic abstract operator}, that is, an
abstract operator wherein the output can play the role of an input and
an input can play the role of the output. In this case, $\CyclicAction$
behaves in the following way. For any $f \in \K \Angle{C}(n)$,
\begin{equation}
    \CyclicAction\Par{
\,.
\end{equation}
In words, $\CyclicAction(f)$ is obtained by transforming each input of
$f$ of index $i + 1$ into an input of index $i$ for any $i \in [n - 1]$,
by transforming the $1$st input of $f$ into an output, and by
transforming the output of $f$ into an input of index $n$. It is
straightforward to check that $\CyclicAction^{n + 1}(f) = f$ as required
due to the fact that $C$ is a cyclic collection.
\medbreak

\subsubsection{Cyclic operads}
A \Def{cyclic operad} is an operad $\K \Angle{C}$ such that
$\K \Angle{C}$ is also cyclic as a polynomial space and satisfies, for
any $x \in C(n)$, $y \in C(m)$, and $i \in [n - 1]$, the compatibility
relations,
\begin{subequations}
\begin{equation} \label{equ:cyclic_operad_axiom_cycle_first}
    \CyclicAction\Par{x \circ_1 y} =
    \CyclicAction(y) \, \circ_m \CyclicAction(x),
\end{equation}
\begin{equation} \label{equ:cyclic_operad_axiom_cycle_others}
    \CyclicAction\Par{x \circ_{i + 1} y} =
    \CyclicAction(x) \circ_i y,
\end{equation}
\begin{equation} \label{equ:cyclic_operad_axiom_unit}
    \CyclicAction(\Unit) = \Unit.
\end{equation}
\end{subequations}
To understand these relations, let us consider first the abstract
operators expressed by the left and right members
of~\eqref{equ:cyclic_operad_axiom_cycle_first}. On the one hand, we have
\begin{equation}
    \CyclicAction\Par{
\,,
\end{equation}
and on the other,
\begin{equation}
    \CyclicAction\Par{
    %
\,.
\end{equation}
We observe that the two obtained abstract operators are the same.
Indeed, for both of them, the connections between $x$ and $y$ are the
same. This is what is expressed
by~\eqref{equ:cyclic_operad_axiom_cycle_first}. Let us now consider the
abstract operators expressed by the left and right members
of~\eqref{equ:cyclic_operad_axiom_cycle_others}. On the one hand, we
have
\begin{equation}
    \CyclicAction\Par{
    %
\,,
\end{equation}
and on the other,
\begin{equation}
    \CyclicAction\Par{
    %
\,.
\end{equation}
We observe that the two obtained abstract operators are the same.
Indeed, for both of them, the connections between $x$ and $y$ are the
same. This is what is expressed
by~\eqref{equ:cyclic_operad_axiom_cycle_others}. Last
relation~\eqref{equ:cyclic_operad_axiom_unit} expresses that
\begin{equation}
  \CyclicAction\Par{
  %
\,.
\end{equation}
Since, in an operad, the unit $\Unit$ can be seen as the identity map
(see Section~\ref{subsubsec:partial_compo_maps} of
Chapter~\ref{chap:operads}), this map is also invertible. This is what
is expressed by~\eqref{equ:cyclic_operad_axiom_unit}.
\medbreak

\subsubsection{Example: operads of configurations of chords}
Let us consider the construction $\CliMonoid$ associating with any
monoid $\Mca$ the operad of configurations of chords $\CliMonoid \Mca$
exposed in Section~\ref{subsubsec:configuration_operads} of
Chapter~\ref{chap:operads}.
\medbreak

Let $\CyclicAction$ be the cyclic map on $\ColCC_{\bar{\Mca}}$ defined,
for any $\bar{\Mca}$-configuration $\Cfr$ in the following way.
The configuration $\CyclicAction(\Cfr)$ is obtained by
applying a rotation of one step of $\Cfr$ in the counterclockwise
direction. For instance, one has in $\CliMonoid \Z$,
\begin{equation}
    \CyclicAction\left(
    \begin{tikzpicture}[scale=.85,Centering]
        \node[CliquePoint](1)at(-0.50,-0.87){};
        \node[CliquePoint](2)at(-1.00,-0.00){};
        \node[CliquePoint](3)at(-0.50,0.87){};
        \node[CliquePoint](4)at(0.50,0.87){};
        \node[CliquePoint](5)at(1.00,0.00){};
        \node[CliquePoint](6)at(0.50,-0.87){};
        \draw[CliqueEdge](1)edge[]node[CliqueLabel]
            {\begin{math}1\end{math}}(2);
        \draw[CliqueEdge](1)edge[]node[CliqueLabel]
            {\begin{math}-2\end{math}}(5);
        \draw[CliqueEmptyEdge](1)edge[]node[CliqueLabel]{}(6);
        \draw[CliqueEdge](2)edge[]node[CliqueLabel]
            {\begin{math}-2\end{math}}(3);
        \draw[CliqueEmptyEdge](3)edge[]node[CliqueLabel]{}(4);
        \draw[CliqueEdge](3)edge[]node[CliqueLabel]
            {\begin{math}1\end{math}}(5);
        \draw[CliqueEmptyEdge](4)edge[]node[CliqueLabel]{}(5);
        \draw[CliqueEmptyEdge](5)edge[]node[CliqueLabel]{}(6);
    \end{tikzpicture}
    \right)
    \enspace = \enspace
    \begin{tikzpicture}[scale=.85,Centering]
        \node[CliquePoint](1)at(-0.50,-0.87){};
        \node[CliquePoint](2)at(-1.00,-0.00){};
        \node[CliquePoint](3)at(-0.50,0.87){};
        \node[CliquePoint](4)at(0.50,0.87){};
        \node[CliquePoint](5)at(1.00,0.00){};
        \node[CliquePoint](6)at(0.50,-0.87){};
        \draw[CliqueEdge](1)edge[]node[CliqueLabel]
            {\begin{math}-2\end{math}}(2);
        \draw[CliqueEdge](1)edge[]node[CliqueLabel]
            {\begin{math}1\end{math}}(6);
        \draw[CliqueEmptyEdge](2)edge[]node[CliqueLabel]{}(3);
        \draw[CliqueEdge](2)edge[]node[CliqueLabel]
            {\begin{math}1\end{math}}(4);
        \draw[CliqueEmptyEdge](3)edge[]node[CliqueLabel]{}(4);
        \draw[CliqueEmptyEdge](4)edge[]node[CliqueLabel]{}(5);
        \draw[CliqueEdge](4)edge[]node[CliqueLabel]
            {\begin{math}-2\end{math}}(6);
        \draw[CliqueEmptyEdge](5)edge[]node[CliqueLabel]{}(6);
    \end{tikzpicture}\,.
\end{equation}
\medbreak

\begin{Proposition} \label{prop:cyclic_configuration_operads}
    For any monoid $\Mca$, $\CliMonoid \Mca$ is a cyclic operad for
    the cycle maps~$\CyclicAction$.
\end{Proposition}
\medbreak

\subsection{Symmetric operads} \label{subsec:symmetric_operads}
As a last variant of operads, we consider now symmetric operads. These
variations of operads involve symmetric collections. Intuitively, in a
symmetric operad, the inputs of the elements can be permuted. This is
due to the fact that these structures are endowed with maps letting the
symmetric group of order $n$ acting on its elements of arity~$n$.
\medbreak

\subsubsection{Symmetric polynomial spaces}
A graded polynomial space $\K \Angle{C}$ is \Def{symmetric} if it is
endowed for all $n \in \N$ and $\sigma \in \SymmetricGroup(n)$ with
unary products
\begin{equation}
    \SymmetricAction_\sigma : \K \Angle{C}(n) \to \K \Angle{C}(n)
\end{equation}
such that $\SymmetricAction_{\Identity_n}$ is the identity map on
$\K \Angle{C}(n)$, where $\Identity_n$ denotes the identity map of
$\SymmetricGroup_n$, and
\begin{math}
    \SymmetricAction_{\sigma_1} \circ \SymmetricAction_{\sigma_2}
    = \SymmetricAction_{\sigma_2 \circ \sigma_1}
\end{math}
for any permutations $\sigma_1$ and $\sigma_2$ of $\SymmetricGroup(n)$.
We say that the $\SymmetricAction_\sigma$, $\sigma \in \SymmetricGroup$,
are \Def{symmetric maps} of $\K \Angle{C}$. As usual, if
$\phi : \K \Angle{C_1} \to \K \Angle{C_2}$ is a morphism between two
graded polynomial spaces $\K \Angle{C_1}$ and $\K \Angle{C_2}$, $\phi$
is a \Def{symmetric polynomial space morphism} if it commutes with the
symmetric maps of $\K \Angle{C_1}$ and $\K \Angle{C_2}$.
\medbreak

Remark that when $C$ is a symmetric collection (see
Section~\ref{subsubsec:symmetric_collections} of
Chapter~\ref{chap:collections}), the linearizations of the symmetric
maps of $C$ endow $\K \Angle{C}$ with the structure of a symmetric
polynomial space.
\medbreak

\subsubsection{Symmetric abstract operators}
We regard any homogeneous element $f$ of an augmented symmetric
polynomial space $\K \Angle{C}$ as a \Def{symmetric abstract operator},
that is, an abstract operator wherein the inputs are endowed with a
total order. More precisely, the inputs of a symmetric abstract operator
of arity $n$ are number from $1$ to $n$, but not necessarily from left
to right in the increasing order as is the case for usual abstract
operators (see Section~\ref{subsubsec:abstract_operators} of
Chapter~\ref{chap:operads}). A symmetric abstract operator $f$ of arity
$n$ is depicted as
\begin{equation}

\end{equation}
where $\pi$ is a permutation of size $n$. Moreover, each symmetric map
$\SymmetricAction_\sigma$, $\sigma \in \SymmetricGroup$, behaves in the
following way. For any $f \in \K \Angle{C}(n)$ and
$\sigma \in \SymmetricGroup(n)$,
\begin{equation}
    \SymmetricAction_\sigma\Par{
    %
\,.
\end{equation}
In words, $\SymmetricAction_\sigma(f)$ is obtained by permuting the
inputs of $f$ as specified by $\sigma$. It is straightforward to check
that
\begin{math}
    \Par{\SymmetricAction_{\sigma_1} \circ \SymmetricAction_{\sigma_2}}
    (f)
    = \SymmetricAction_{\sigma_2 \circ \sigma_1}(f)
\end{math}
for all $\sigma_1, \sigma_2 \in \SymmetricGroup(n)$ as expected since
$\K \Angle{C}$ is a symmetric polynomial space. On symmetric abstract
operators, for any $f \in \K \Angle{C}(n)$ and $g \in \K \Angle{C}(m)$,
the partial composition $f \circ_i g$ is the symmetric abstract operator
\begin{equation} \label{equ:partial_compostion_on_symmetric_operators}
\,,
\end{equation}
where $j \in [n]$ is such that $\pi(j) = i$ and
$\mu = \pi \circ_j \tau$, where the occurrence of $\circ_i$ is the
partial composition map of the operad $\Per$ of permutations (see
Section~\ref{subsubsec:operad_per} of Chapter~\ref{chap:operads}).
\medbreak

\subsubsection{Symmetric operads}
A \Def{symmetric operad} is an operad $\K \Angle{C}$ such that
$\K \Angle{C}$ is also symmetric as polynomial space and satisfies,
for any $x \in  C(n)$, $\sigma \in \SymmetricGroup(n)$, $y \in C(m)$,
$\nu \in \SymmetricGroup(m)$, and $i \in [n]$, the compatibility
relation
\begin{equation} \label{equ:symmetric_operad_axiom}
    \SymmetricAction_\sigma(x) \circ_i \SymmetricAction_\nu(y)
    = \SymmetricAction_{\sigma \circ_i \nu}\Par{x \circ_{\sigma(i)} y},
\end{equation}
where the occurrence of $\circ_i$ in the right member
of~\eqref{equ:symmetric_operad_axiom} refers to the partial composition
maps of $\Per$. To understand this relation, let us consider the
abstract operators expressed by the left and right members
of~\eqref{equ:symmetric_operad_axiom}. On the one hand, we have
\begin{multline} \label{equ:symmetric_operad_axiom_developed_1}
    \SymmetricAction_\sigma\Par{
\,,
\end{multline}
where $j \in [n]$ is such that $\Par{\sigma^{-1} \circ \pi}(j) = i$ and
$\mu = \Par{\sigma^{-1} \circ \pi} \circ_j \Par{\nu^{-1} \circ \tau}$.
On the other, we have
\begin{multline} \label{equ:symmetric_operad_axiom_developed_2}
    \SymmetricAction_{\sigma \circ_i \nu}\Par{
    %
\,,
\end{multline}
where $j' \in [n]$ is such that $\pi\Par{j'} = \sigma(i)$ and
$\mu' = \pi \circ_{j'} \tau$. Let us now explain why the last symmetric
abstract operators of~\eqref{equ:symmetric_operad_axiom_developed_1}
and~\eqref{equ:symmetric_operad_axiom_developed_2} are equal. First,
from the hypothesis $\sigma^{-1}(\pi(j)) = i$ and
$\pi\Par{j'} = \sigma(i)$, we deduce that
\begin{math}
    \sigma^{-1}(\pi\Par{j'}) = \sigma^{-1}(\sigma(i)) = i,
\end{math}
implying that $j' = j$. This provides the fact that the considered
abstract operators have the same shape: the output of $y$ is connected
to the $j$th input of $x$ in both cases. Now, consider the following
result connecting the group theoretic composition $\circ$ of
permutations and the partial composition maps $\circ_i$ of~$\Per$.
\medbreak

\begin{Lemma} \label{lem:composition_permutations_compatibility}
    Let $n, m \in \N_{\geq 1}$, $i \in [n]$, and four permutations
    $\pi, \sigma \in \SymmetricGroup(n)$,
    $\tau, \nu \in \SymmetricGroup(m)$. Then, in the operad $\Per$,
    \begin{equation}
        \Par{\sigma \circ_i \nu}^{-1}
        \circ \Par{\pi \circ_{\pi^{-1}(\sigma(i))} \tau}
        =
        \Par{\sigma^{-1} \circ \pi}
        \circ_{\pi^{-1}(\sigma(i))} \Par{\nu^{-1} \circ \tau}.
    \end{equation}
\end{Lemma}
\medbreak

By Lemma~\ref{lem:composition_permutations_compatibility}, the
permutations $\mu$ and $\Par{\sigma \circ_i \nu}^{-1} \circ \mu'$ are
equal. Therefore, the two obtained abstract operators are the same.
\medbreak

\subsubsection{Example: the symmetric operad $\NAP$}
\label{subsubsec:symmetric_operad_nap}
Let us consider the operad $\NAP$ introduced in
Section~\ref{subsubsec:operad_nap} of Chapter~\ref{chap:operads}. We
endow the underlying collection $\ColSRT$ of $\NAP$ with the symmetric
maps $\SymmetricAction_\sigma$, $\sigma \in \SymmetricGroup$, defined
for any standard rooted tree $\Tfr$ of size $n$ in the following way.
The tree $\SymmetricAction_\sigma(\Tfr)$ is obtained by relabeling all
the nodes $i \in [n]$ by $\sigma^{-1}(i)$. For instance,
\begin{subequations}
\begin{equation}
    \SymmetricAction_{3142}
    \left(
\,.
\end{equation}
\end{subequations}
Then, the linearization of these symmetric maps endows $\NAP$ with
the structure of a symmetric polynomial space. Together with the partial
composition maps of $\NAP$ and its unit, $\NAP$ is a symmetric operad.
\medbreak

It is possible to show that, as a symmetric operad, the set
\begin{equation} \label{equ:generating_set_nap}
    \GeneratingSet :=
    \left\{
    %
    \right\}
\end{equation}
is a minimal generating set of $\NAP$. Moreover, $\NAP$ contains the
single nontrivial relation
\begin{equation} \label{equ:relation_nap}
    %
    \right)
    \enspace = \enspace
    0
\end{equation}
involving its generator.
\medbreak

\subsubsection{Constructions involving operads}
If $\K \Angle{C}$ is a symmetric operad, by forgetting its symmetric
maps, the space $\K \Angle{C}$ seen as an augmented graded polynomial
space is an operad. We call this operad the \Def{symmetric oblivion} of
the symmetric operad~$\K \Angle{C}$.
\medbreak

Besides, given an operad $\K \Angle{C}$, there are at least two ways to
construct a symmetric operad. Let us present these. The first one
consists in turning the augmented graded collection $C$ into a symmetric
collection by endowing it with the symmetric maps
$\SymmetricAction_\sigma$, $\sigma \in \SymmetricGroup$, defined by
$\SymmetricAction_\sigma(x) := x$ for all $x \in C$ and
$\sigma \in \SymmetricGroup(|x|)$. The space $\K \Angle{C}$ endowed with
the linearizations of these symmetric maps on $C$, together with the
partial composition maps and unit of the operad $\K \Angle{C}$ forms a
symmetric operad, called \Def{trivial symmetric operad}
of~$\K \Angle{C}$.
\medbreak

The second one consists in turning $C$ into a symmetric collection by
considering its regularization $\Regularization(C)$ (see
Section~\ref{subsubsec:regularization} of
Chapter~\ref{chap:collections}). Let us endow the symmetric polynomial
space $\K \Angle{\Regularization(C)}$ with the partial composition maps
defined linearly for any $(x, \sigma) \in \Par{\Regularization(C)}(n)$,
$(y, \nu) \in \Regularization(C)$, and $i \in [n]$, by
\begin{equation} \label{equ:composition_regularization}
    (x, \sigma) \circ_i (y, \nu)
    :=
    \Par{x \circ_{\sigma^{-1}(i)} y, \sigma \circ_{\sigma^{-1}(i)} \nu},
\end{equation}
where the first occurrence of a partial composition map in the right
member of~\eqref{equ:composition_regularization} refers to the partial
composition map of the operad $\K \Angle{C}$, and the second one refers
to the partial composition map of the operad~$\Per$.
\medbreak

\begin{Proposition} \label{prop:regularization_operads}
    For any operad $\K \Angle{C}$, $\K \Angle{\Regularization(C)}$ is a
    symmetric operad.
\end{Proposition}
\medbreak

Proposition~\ref{prop:regularization_operads} is a consequence of
Lemma~\ref{lem:composition_permutations_compatibility}. The symmetric
operad $\K \Angle{\Regularization(C)}$, denoted by a slight abuse of
notation by $\Regularization(\K \Angle{C})$, is the \Def{regularization}
of $\K \Angle{C}$. Let us notice that in general, the regularization of
the symmetric oblivion of a symmetric operad $\K \Angle{C}$ is different
to $\K \Angle{C}$. For instance, consider the symmetric operad $\NAP$
considered in Section~\ref{subsubsec:symmetric_operad_nap}. The
regularization of the oblivion of $\NAP$ is the symmetric operad
$\Regularization(\NAP)$, and this last is not isomorphic to $\NAP$ since
$\dim \NAP(3) = 9$ and
\begin{math}
    \dim (\Regularization(\NAP))(3) = 9 \times 3! = 54.
\end{math}
\medbreak

\subsubsection{Algebras over symmetric operads}
Let $\K \Angle{C}$ be a symmetric operad. An
\Def{algebra over $\K \Angle{C}$} (or, for short, a
\Def{$\K \Angle{C}$-algebra}) is an algebra $\K \Angle{D}$ over the
oblivion of $\K \Angle{C}$ (see
Section~\ref{subsec:algebras_over_operads} of
Chapter~\ref{chap:operads}) satisfying the following additional
condition. By denoting by $\Action_n$, $n \in \N_{\geq 1}$, the linear
maps~\eqref{equ:action_algebra_over_operad} of
Chapter~\ref{chap:operads} endowing the space $\K \Angle{D}$ with the
structure of a $\K \Angle{C}$-algebra, for any $x \in C(n)$,
$\sigma \in \SymmetricGroup(n)$, and
$\Par{a_1, \dots, a_n} \in \List_{\{n\}}(D)$,
\begin{equation}
    \Action_n\Par{
    \SymmetricAction_\sigma(x), \Par{a_1, \dots, a_n}}
    =
    \Action_n\Par{
    x, \Par{a_{\sigma^{-1}(1)}, \dots, a_{\sigma^{-1}(n)}}}.
\end{equation}
As for operads, one can regard each homogeneous element of arity
$n \in \N_{\geq 1}$ of the symmetric operad $\K \Angle{C}$ as a complete
product of arity $n$ on~$\K \Angle{D}$.
\medbreak

For instance, let us consider algebras on the symmetric operad $\NAP$.
We denote by $\Tfr$ the single generator of $\NAP$ appearing
in~\eqref{equ:generating_set_nap}. Since this generator is subjected to
Relation~\eqref{equ:relation_nap}, for any
$f_1, f_2, f_3 \in \K \Angle{D}$,
\begin{equation}\begin{split}
    0 & = \Par{\Tfr \circ_1 \Tfr
    - \SymmetricAction_{132}\Par{\Tfr \circ_1 \Tfr}}
    \Par{f_1, f_2, f_3} \\
    & = \Par{\Tfr \circ_1 \Tfr}\Par{f_1, f_2, f_3}
    - \Par{\SymmetricAction_{132}\Par{\Tfr \circ_1 \Tfr}}
    \Par{f_1, f_2, f_3} \\
    & = \Par{\Tfr \circ_1 \Tfr}\Par{f_1, f_2, f_3}
        - \Par{\Tfr \circ_1 \Tfr}\Par{f_1, f_3, f_2} \\
    & = \Tfr \Par{\Tfr \Par{f_1, f_2}, f_3}
        - \Tfr \Par{\Tfr \Par{f_1, f_3}, f_2}.
\end{split}\end{equation}
This is equivalent to the relation
\begin{equation}
    \Par{f_1 \, \Tfr \, f_2} \, \Tfr \, f_3
    - \Par{f_1 \, \Tfr \, f_3} \, \Tfr \, f_2
    = 0
\end{equation}
written in infix way.
\medbreak

\section{Product categories} \label{sec:pros}
A very intuitive generalization of operads arises when one thinks of
considering elements with several outputs, instead of only one as in the
case of operads. This leads to a sort of extension of operads, called
pros. We present here these algebraic structures.
\medbreak

\subsection{Abstract bioperators}
While operads are defined over augmented graded polynomial spaces,
pros require spaces on $2$-graded collections. Let us provide elementary
definitions about these.
\medbreak

\subsubsection{Bigraded polynomial spaces}
Given a $2$-graded collection $C$ (see
Section~\ref{subsubsec:multigraded_collections} of
Chapter~\ref{chap:collections}), the polynomial space $\K \Angle{C}$ is
said \Def{bigraded}. The \Def{arity} (resp. \Def{coarity}) of a nonzero
$(p, q)$-homogeneous element $f$ of $\K \Angle{C}$ is $p$ (resp. $q$).
Moreover, to not overload the notation, we denote by
$\K \Angle{C}(p, q)$ the homogeneous component $\K \Angle{C}((p, q))$ of
$\K \Angle{C}$ for any $(p, q) \in \N^2$. When $\K \Angle{C}$ is
combinatorial, its \Def{Hilbert series} is the $\N^2$-series
\begin{equation}
    \HilbertSeries_{\K \Angle{C}}
    := \IndexSeries(C)
    = \sum_{(p, q) \in \N^2} \Par{\dim \K \Angle{C}(p, q)} \, (p, q).
\end{equation}
\medbreak

\subsubsection{Abstract bioperators}
We regard any homogeneous element $f$ of a bigraded polynomial space
$\K \Angle{C}$ as an \Def{abstract bioperator}, that is, an abstract
operator having zero or more inputs and zero or more outputs. These
abstract bioperators are depicted by following the drawing conventions
of biproducts exposed in Section~\ref{subsubsec:biproducts} of
Chapter~\ref{chap:algebra}. Therefore, if $f$ is of arity $p$ and of
coarity $q$, $f$ is depicted by
\begin{equation}
    \begin{tikzpicture}
        [xscale=.45,yscale=.35,Centering,font=\scriptsize]
        \node[Operator](0)at(0,0){\begin{math}f\end{math}};
        \node(1)at(-1,-2){};
        \node(2)at(1,-2){};
        \node(3)at(-1,2){};
        \node(4)at(1,2){};
        \draw[Edge](0)--(1);
        \draw[Edge](0)--(2);
        \draw[Edge](0)--(3);
        \draw[Edge](0)--(4);
        \node[below of=1,node distance=.05cm]
            {\begin{math}1\end{math}};
        \node[below of=2,node distance=.05cm]
            {\begin{math}p\end{math}};
        \node[above of=3,node distance=.05cm]
            {\begin{math}1\end{math}};
        \node[above of=4,node distance=.05cm]
            {\begin{math}q\end{math}};
        \node[below of=0,node distance=.75cm]
            {\begin{math}\dots\end{math}};
        \node[above of=0,node distance=.75cm]
            {\begin{math}\dots\end{math}};
    \end{tikzpicture}\,.
\end{equation}
\medbreak

\subsubsection{Composing abstract bioperators}
Let $\K \Angle{C}$ be bigraded polynomial space. A
\Def{horizontal composition map} on $\K \Angle{C}$ is a complete binary
$\omega$-concentrated product on $\K \Angle{C}$ of the form
\begin{equation}
    \ProHoriz : \K \Angle{\BBrack{C, C}_\times} \to \K \Angle{C}
\end{equation}
for the map $\omega : \N^2 \times \N^2 \to \N^2$ satisfying
\begin{math}
    \omega\Par{\Par{p_1, q_1}, \Par{p_2, q_2}}
    := \Par{p_1 + p_2, q_1 + q_2}.
\end{math}
On abstract bioperators, for any $f \in \K \Angle{C}(p, q)$ and
$g \in \K \Angle{C}\Par{p', q'}$, $f \ProHoriz g$ is the abstract
bioperator
\begin{equation}
\,.
\end{equation}
\medbreak

A \Def{vertical composition map} on $\K \Angle{C}$ is a binary product
on $\K \Angle{C}$ of the form, for any $p, q, r \in \N$,
\begin{equation}
    \ProVerti^{(p, q, r)} :
    \K \Angle{\BBrack{C(q, r), C(p, q)}_\times} \to \K \Angle{C}(p, r).
\end{equation}
By a slight abuse of notation, we shall sometimes omit the $(p, q, r)$
in the notation of $\ProVerti^{(p, q, r)}$ in order to denote it in a
more concise way by $\ProVerti$. On abstract bioperators, for any
$f \in \K \Angle{C}(q, r)$ and $g \in \K \Angle{C}(p, q)$,
$f \ProVerti g$ is the abstract bioperator
\begin{equation}
\,.
\end{equation}
\medbreak

Finally, a \Def{unit map} on $\K \Angle{C}$ is a product of arity zero
on $\K \Angle{C}$ of the form, for any $p \in \N$,
\begin{equation}
    \Unit_p : \K \Angle{\BBrack{}_\times} \to \K \Angle{C}(p, p).
\end{equation}
For any $p \in \N$, $\Unit_p$ is the abstract bioperator
\begin{equation}
    %
\,.
\end{equation}
\medbreak

\subsection{Pros}
Pros are algebraic structures furnishing a formalization of the notion
of abstract bioperators and their compositions. We provide here
definitions about these structures and about bialgebras over pros.
\medbreak

\subsubsection{Elementary definitions}
The bigraded space $\K \Angle{C}$ is a \Def{product category} (or, for
short, a \Def{pro}) if it is endowed with a horizontal composition map
$\ProHoriz$, vertical composition maps~$\ProVerti$, and unit maps
$\Unit_p$, $p \in \N$, satisfying  the following six
relations~\eqref{equ:assoc_compo_h}, \eqref{equ:assoc_compo_v},
\eqref{equ:compo_square}, \eqref{equ:relation_unit_h},
\eqref{equ:relation_unit_zero} and~\eqref{equ:relation_unit_v}.
\medbreak

First, $\ProHoriz$ is associative, that is, for any $x, y, z \in C$,
\begin{equation} \label{equ:assoc_compo_h}
    (x \ProHoriz y) \ProHoriz z = x \ProHoriz (y \ProHoriz z).
\end{equation}
Second, the products $\ProVerti$ are associative, that is, for any
$x \in C(r, s)$, $y \in C(q, r)$, $z \in C(p, q)$,
\begin{equation} \label{equ:assoc_compo_v}
    (x \ProVerti y) \ProVerti z = x \ProVerti (y \ProVerti z).
\end{equation}
Moreover, the products $\ProHoriz$ and $\ProVerti$ satisfy the
\Def{square compatibility relation}, that is for any $x \in C(q, r)$,
$y \in C(p, q)$, $x' \in C\Par{q', r'}$, $y' \in C\Par{p', q'}$,
\begin{equation} \label{equ:compo_square}
    (x \ProVerti y) \ProHoriz \Par{x' \ProVerti y'}
    = \Par{x \ProHoriz x'} \ProVerti \Par{y \ProHoriz y'}.
\end{equation}
Finally, the unit maps satisfy the following three sorts of relations.
For any $p, q \in \N$,
\begin{equation} \label{equ:relation_unit_h}
    \Unit_p \ProHoriz \Unit_q = \Unit_{p + q},
\end{equation}
for any $x \in C(p, q)$,
\begin{equation} \label{equ:relation_unit_zero}
    x \ProHoriz \Unit_0 = x = \Unit_0 \ProHoriz x,
\end{equation}
and for any $x \in C(p, q)$,
\begin{equation} \label{equ:relation_unit_v}
    x \ProVerti \Unit_p = x = \Unit_q \ProVerti x.
\end{equation}
\medbreak

Let us understand these relations with the help of abstract bioperators
and the behavior of their compositions. First, the left and
right members of~\eqref{equ:assoc_compo_h} are both equal to
\begin{equation}
\,,
\end{equation}
where $x \in C(p, q)$, $y \in C\Par{p', q'}$, and
$z \in C\Par{p'', q''}$. Second, the left and right members
of~\eqref{equ:assoc_compo_v} are both equal to
\begin{equation}
    %
\,.
\end{equation}
Moreover, the left and right members of~\eqref{equ:compo_square} are
both equal to
\begin{equation}
    %
\,.
\end{equation}
Finally, concerning the relations involving the unit maps, the left and
right members of~\eqref{equ:relation_unit_h} are both equal to
\begin{equation}
    %
\,.
\end{equation}
Relation~\eqref{equ:relation_unit_zero} expresses the fact that
$\Unit_0$ is the unit for the product $\ProHoriz$, and
Relation~\eqref{equ:relation_unit_v} says that $\Unit_p$, $p \in \N$,
is the unit of the products $\ProVerti^{(p, p, q)}$ and
$\ProVerti^{(q, p, p)}$ for any $q \in \N$.
\medbreak

Since a pro is a particular polynomial algebra, all the properties and
definitions about polynomial algebras exposed in
Section~\ref{subsec:polynomial_bialgebras} of
Chapter~\ref{chap:algebra} remain valid for pros (like pros morphisms,
sub-pros, generating  sets, pro ideals and quotients, {\em etc.}).
\medbreak

\subsubsection{Categorical definition}
In the same way as monoids, operads, and colored operads can be defined
precisely and concisely by using the language of category theory (see
Section~\ref{subsubsec:categorical_operads}), pros admit a similar
definition using this language. Indeed, a pro $\K \Angle{C}$ can be seen
as a category where the objects are the elements of $\N$ and which is
equipped with a bifunctor $\ProHoriz$ defined by
$p \ProHoriz p' := p + p'$. The elements of $\K \Angle{C}$ are
interpreted as maps $\phi : p \to q$. The horizontal composition of
$\K \Angle{C}$ translates as the bifunctor of the category, the vertical
composition translates as the composition of morphisms, and the unit map
translates as the identity maps $\Unit_p : p \to p$ of the category.
\medbreak

\subsubsection{Bialgebras over pros}
Let $\K \Angle{C}$ be a pro. A \Def{bialgebra over $\K \Angle{C}$} (or,
for short, a \Def{$\K \Angle{C}$-bialgebra}) is a polynomial space
$\K \Angle{D}$, where $D$ is a (not necessarily $2$-graded) collection,
which is endowed for all $p, q \in \N$ with linear maps
\begin{equation}
    \Action_{p, q} :
    \K \Angle{\BBrack{C(p, q), \List_{\{p\}}(D)}_\times}
    \to \K \Angle{\List_{\{q\}}(D)}
\end{equation}
satisfying the relations imposed by the pro structure of $\K \Angle{C}$,
that are,
\begin{subequations}
for all $x \in C\Par{p', q'}$, $y \in C\Par{p'', q''}$, and
\begin{math}
    \Par{a_1, \dots, a_{p' + p''}}
    \in \List_{\left\{p' + p''\right\}}(D),
\end{math}
by denoting by $\Conc$ the concatenation of tuples extended by
linearity,
\begin{multline} \label{equ_bialgebras_over_pros_horiz}
    \Action_{p' + p'', q' + q''}
    \Par{x \ProHoriz y,
    \Par{a_1, \dots, a_{p'}, a_{p' + 1}, \dots, a_{p' + p''}}} \\
    =
    \Action_{p', q'}\Par{x, \Par{a_1, \dots, a_{p'}}}
    \Conc
    \Action_{p'', q''}\Par{y, \Par{a_{p' + 1}, \dots, a_{p' + p''}}},
\end{multline}
for all $x \in C(q, r)$, $y \in C(p, q)$, and
$\Par{a_1, \dots, a_p} \in \List_{\{p\}}(D)$,
\begin{equation} \label{equ_bialgebras_over_pros_verti}
    \Action_{p, r}\Par{x \circ y, \Par{a_1, \dots, a_p}}
    =
    \Action_{q, r}\Par{x, \Action_{p, q}\Par{y, \Par{a_1, \dots, a_p}}},
\end{equation}
and for all $\Par{a_1, \dots, a_p} \in \List_{\{p\}}(D)$,
\begin{equation} \label{equ_bialgebras_over_pros_units}
    \Action_{p, p}\Par{\Unit_p, \Par{a_1, \dots, a_p}}
    =
    \Par{a_1, \dots, a_p}.
\end{equation}
\end{subequations}
In other words, any object $x$ of $C$ of arity $p$ and coarity $q$ plays
the role of a complete biproduct (in the sense of
Section~\ref{subsubsec:biproducts} of Chapter~\ref{chap:algebra}) of the
form
\begin{equation}
    x : \K \Angle{\List_{\{p\}}(D)} \to \K \Angle{\List_{\{q\}}(D)},
\end{equation}
defined, for any
\begin{math}
    \Par{a_1, \dots, a_p} \in \List_{\{p\}}(D)
\end{math}
by
\begin{equation}
    x\Par{a_1, \dots, a_p} :=
    \Action_{p, q}\Par{x, \Par{a_1, \dots, a_p}}.
\end{equation}
\medbreak

\subsection{Main pros}
We provide here classical examples of pros involving for most of these
combinatorial objects (see Table~\ref{tab:examples_pros}). The first one
is the associative pro, a sort of generalization of the associative
operad. The next one is a pro of matrices, where the usual matrix
operations (direct sum and multiplication) are revisited in the context
of pros. This structure contains a lot of interesting sub-pros, like a
pro of binary relations, a pro of maps, and a pro of permutations.
\begin{table}[ht]
    \centering
    \begin{tabular}{|c|c|c|c|c|} \hline
        Pro & Objects & Arity & Coarity \\ \hline \hline
        $\PAs$ & Pairs $(p, q)$ of nonneg. int. & $p$ & $q$
            \\ \hline
        $\Mat_{\Sbb}$ & Matrices on $\Sbb$ & num. of rows
            & num. of columns \\ \hline
        $\BRel$ & Binary relations & card. of domain
            & card. of codomain \\
        $\K \Angle{\ColMap}$ & Maps & card. of domain
            & card. of codomain \\
        $\K \Angle{\ColNDMap}$ & Nondecreasing maps & card. of domain
            & card. of codomain \\
        $\K \Angle{\ColInj}$ & Injective maps & card. of domain
            & card. of codomain \\
        $\K \Angle{\ColSur}$ & Surjective maps & card. of domain
            & card. of codomain \\
        $\K \Angle{\ColBij}$ & Bijective maps & card. of domain
            & card. of domain \\ \hline
    \end{tabular}
    \bigbreak

    \caption[Overview of some operads]
    {\footnotesize Overview of some pros. Here, $\Sbb$ is a semiring.}
    \label{tab:examples_pros}
\end{table}
\medbreak

\subsubsection{Associative pro}
Let $A := \left\{\Afr_{p, q} : p, q \in \N \right\}$ be the $2$-graded
collection where the index of each $\Afr_{p, q}$, $p, q \in \N$, is
$(p, q)$. The \Def{associative pro} $\PAs$ is the space $\K \Angle{A}$
endowed with the horizontal composition map $\ProHoriz$ defined
linearly, for any $\Afr_{p, q} \in A(p, q)$ and
$\Afr_{p', q'} \in A\Par{p', q'}$ by
\begin{math}
    \Afr_{p, q} \ProHoriz \Afr_{p', q'} := \Afr_{p + p', q + q'},
\end{math}
with the vertical composition maps $\ProVerti$ defined linearly, for any
$\Afr_{q, r} \in A(q, r)$, $\Afr_{p, q} \in A(p, q)$ by
\begin{math}
    \Afr_{q, r} \ProVerti \Afr_{p, q} := \Afr_{p, r},
\end{math}
and with the unit maps $\Unit_p$, $p \in \N$, defined
by~$\Unit_p := \Afr_{p, p}$.
\medbreak

\subsubsection{Pro of matrices}
Let $\Sbb$ be a semiring (that is, a ring without the condition
that each element have an additive inverse). Let $\ColMat_\Sbb$ be the
$2$-graded collection of all matrices on $\Sbb$ where
$\ColMat_\Sbb(p, q)$ is the set of such matrices of dimension
$p \times q$, $p, q \in \N$. The \Def{pro of matrices} $\Mat_\Sbb$ is
the space $\K \Angle{\ColMat_\Sbb}$ endowed with the horizontal
composition map $\ProHoriz$ defined linearly, for any
$\Mfr_1, \Mfr_2 \in \Mat_\Sbb$ by
\begin{equation}
    \Mfr_1 \ProHoriz \Mfr_2 := \Mfr_1 \oplus \Mfr_2,
\end{equation}
where $\oplus$ is the matrix direct sum. Moreover, $\Mat_\Sbb$ is
endowed with the vertical composition maps $\ProVerti$ defined linearly,
for any $\Mfr_1 \in \ColMat_\Sbb(q, r)$ and
$\Mfr_2 \in \ColMat_\Sbb(p, q)$, $p, q, r \in \N$, by
\begin{equation}
    \Mfr_1 \ProVerti \Mfr_2 := \Mfr_2 \Conc \Mfr_1,
\end{equation}
where $\Conc$ is the matrix multiplication. Finally, we define the unit
maps $\Unit_p$, $p \in \N$, of $\Mat_\Sbb$ as the identity matrix of
order $p$. For instance, by setting $\Sbb$ as the semiring
$(\N, +, \Conc)$,
\begin{equation}
    \Matrix{
        \ColA{2} & \ColA{0} & \ColA{0} \\
        \ColA{1} & \ColA{1} & \ColA{0}
    }
    \ProHoriz
    \Matrix{
        \ColD{1} & \ColD{1}  \\
        \ColD{3} & \ColD{0}
    }
    =
    \Matrix{
        \ColA{2} & \ColA{0} & \ColA{0} & 0 & 0 \\
        \ColA{1} & \ColA{1} & \ColA{0} & 0 & 0 \\
        0 & 0 & 0 & \ColD{1} & \ColD{1} \\
        0 & 0 & 0 & \ColD{3} & \ColD{0}
    }
\end{equation}
and
\begin{equation}
    \Matrix{
        1 & 0 \\
        0 & 1 \\
        1 & 1
    }
    \ProVerti
    \Matrix{
        0 & 1 & 1 \\
        0 & 2 & 1
    }
    =
    \Matrix{
        1 & 2 \\
        1 & 3
    }
\end{equation}
are, respectively, horizontal and vertical compositions in $\Mat_\Sbb$.
\medbreak

\subsubsection{Pro of binary relations}
Let $\Bool$ be the Boolean semiring, that is the set $\{0, 1\}$ equipped
with the addition $+$ satisfying $0 + 0 = 0$ and
\begin{math}
    0 + 1 = 1 + 0 = 1 + 1 = 1,
\end{math}
and the multiplication $\Conc$ satisfying $1 \Conc 1 = 1$ and
\begin{math}
    0 \Conc 1 = 1 \Conc 0 = 0 \Conc 0 = 0.
\end{math}
The \Def{pro of binary relations} $\BRel$ is the pro $\Mat_\Bool$. By
definition, the collection of all matrices on $\Bool$, called
\Def{Boolean matrices} forms a basis of $\BRel$. For instance,
\begin{equation} \label{equ:pro_brel_example_1}
    \Matrix{
        \ColA{1} & \ColA{1} & \ColA{0} \\
        \ColA{1} & \ColA{0} & \ColA{0} \\
        \ColA{0} & \ColA{0} & \ColA{0} \\
        \ColA{0} & \ColA{1} & \ColA{0}
    }
    \ProHoriz
    \Matrix{
        \ColD{0} & \ColD{1} \\
        \ColD{1} & \ColD{0}
    }
    =
    \Matrix{
        \ColA{1} & \ColA{1} & \ColA{0} & 0 & 0 \\
        \ColA{1} & \ColA{0} & \ColA{0} & 0 & 0 \\
        \ColA{0} & \ColA{0} & \ColA{0} & 0 & 0 \\
        \ColA{0} & \ColA{1} & \ColA{0} & 0 & 0 \\
        0 & 0 & 0 & \ColD{0} & \ColD{1} \\
        0 & 0 & 0 & \ColD{1} & \ColD{0}
    }
\end{equation}
and
\begin{equation} \label{equ:pro_brel_example_2}
    \Matrix{
        1 & 1 & 0 \\
        1 & 0 & 0 \\
        0 & 1 & 0 \\
        0 & 1 & 0
    }
    \ProVerti
    \Matrix{
        1 & 0 & 0 & 0 \\
        0 & 1 & 0 & 1 \\
        0 & 0 & 1 & 1
    }
    =
    \Matrix{
        1 & 1 & 0 \\
        1 & 1 & 0 \\
        0 & 1 & 0
    }
\end{equation}
are, respectively, horizontal and vertical compositions in $\BRel$. There
is, for any $p, q \in \N$, a one-to-one correspondence between the set
$\Mat_\Bool(p, q)$ and the set of all binary relations between $[p]$
and $[q]$. Indeed, a matrix $\Mfr \in \Mat_\Bool(p, q)$, $p, q \in \N$,
and a binary relation $\BinRel$ between $[p]$ and $[q]$ are in
correspondence if, for any $x \in [p]$ and $y \in [q]$, one has
$\Mfr_{x, y} = 1$ if and only if $x \BinRel y$. By using this
correspondence, and by drawing a binary relation $\BinRel$ through a
graph connecting $x$ to $y$ if $x \BinRel y$, where the elements of the
domain are depicted below and the elements of the codomain are depicted
above, \eqref{equ:pro_brel_example_1} and~\eqref{equ:pro_brel_example_2}
translate, respectively, as
\begin{equation}
\,.
\end{equation}
Therefore, $\ProHoriz$ is the concatenation of binary relations and
$\ProVerti$ is their composition.
\medbreak

\subsubsection{Pros of maps and variations}
The pro $\BRel$ admits many sub-pros by considering subspaces of binary
relations satisfying particular conditions. For instance, by setting
$\ColMap$ as the subcollection of $\ColMat_\Bool$ restrained on binary
relations that are maps (that is, $\Mfr \in \ColMap(p, q)$ if for any
$x \in [p]$, there is exactly one $y \in [q]$ such that
$\Mfr_{x, y} = 1$), $\K \Angle{\ColMap}$ is a sub-pro of $\BRel$.
Moreover, by setting $\ColNDMap$ as the subcollection of $\ColMap$
restrained on maps that are nondecreasing (that is,
$\Mfr \in \ColNDMap(p, q)$ if for any $x, x' \in [p]$ such that
$x \leq x'$, $\Mfr_{x, y} = \Mfr_{x', y'} = 1$ implies $y \leq y'$),
$\K \Angle{\ColNDMap}$ is a sub-pro of $\K \Angle{\ColMap}$. Besides, by
setting $\ColInj$ (resp. $\ColSur$) as the subcollection of $\ColMap$
restrained on maps that are injections (resp. surjections),
$\K \Angle{\ColInj}$ (resp. $\K \Angle{\ColSur}$) is a sub-pro of
$\K \Angle{\ColMap}$. Finally, by setting $\ColBij$ as the
subcollection of $\ColMap$ restrained on maps that are bijections,
$\K \Angle{\ColBij}$ is a sub-pro of both $\K \Angle{\ColInj}$
and~$\K \Angle{\ColSur}$.
\medbreak

\SkipTocEntry\section*{Bibliographic notes}

\SkipTocEntry\subsection*{About series on operads}
Our approach concerning formal power series through the framework of
collections encompasses the classical case of power series, as
explained in Section~\ref{sec:series_operads}. Since the introduction of
formal power series, a lot of generalizations were proposed in order to
extend the range of enumerative problems they can help to solve. The
most obvious ones are multivariate series allowing to count objects not
only with respect to their sizes but additionally  with respect to
various other statistics (see
Section~\ref{subsubsec:multigraded_collections} of
Chapter~\ref{chap:collections}). Such series are elements of
$\K \AAngle{\Multiset\Par{\left\{t_1, \dots, t_k\right\}}}$ where all
the $t_i$, $i \in [k]$, are atomic objects. Another one consists in
considering noncommutative series on words~\cite{Eil74,SS78,BR10} (and
thus, elements of $\K \AAngle{A^*}$, where $A$ is an alphabet), or
even, pushing the generalization one step further, on elements of a
monoid~\cite{Sak09} (and thus, elements of $\K \AAngle{\Mca}$, where
$\Mca$ is a monoid). Besides, as another generalization, series on trees
have been considered~\cite{BR82,Boz01}. Series on operads increase the
list of these generalizations. Chapoton was the first to have considered
such series on operads~\cite{Cha02,Cha08,Cha09}. Several authors have
contributed to this field by considering slight variations in the
definitions of these series. Among these, one can cite van der
Laan~\cite{Vdl04}, Frabetti~\cite{Fra08}, and Loday and
Nikolov~\cite{LN13}. In this text, we have presented series on
set-operads as powerful tools to provide descriptions of the generating
series of some combinatorial graded collections $C$. All this relies on
a set-operad structure on $C$ having the property to be a Koszul operad
and admitting a Poincaré-Birkhoff-Witt basis. The obtained descriptions
of the generating series of $C$ are in fact the generating series of the
syntax trees on the generators of $C$ (as an operad) that avoid some
patterns (that are the left members of an orientation of the space of
relations of $C$, see Proposition~\ref{prop:koszulity_criterion_pbw} of
Chapter~\ref{chap:operads}). Similar ideas were brought by Khoroshkin 
and Piontkovski~\cite{KP15}, focused on the theory of Gröbner bases for
symmetric operads. In~\cite{Gir18}, the emphasis was put on the
combinatorial and enumerative consequences of set-operads admitting
Poincaré-Birkhoff-Witt bases, leading to refinements of their Hilbert
series.
\medbreak

\SkipTocEntry\subsection*{About colored operads}
Classical references about colored operads are~\cite{BV73}
and~\cite{Yau16}. By looking at this theory in the shoes of a computer
scientist, one can think the colors as data types in computer
programming. In the same way as two functions can be composed only if
the output type of the one is equal to the type of an input of the
second, the (partial and full) composition maps of a colored operad
require a property on the colors of the operands. Moreover, to pursue
the analogy, the elements of arity one of a colored operad having
different input and output colors can be interpreted as casting
operators in computer programming (that are, operators taking as input
one object and changing its type).  Besides, colored operads are
interesting devices for enumerative prospects when combined with series
on operads~\cite{Gir16d}. In this cited work, a generalization of both
context-free grammars (see~\cite{Har78,HMU06}) and regular tree grammars
(see~\cite{CDGJLLTT07}) using colored operads was proposed. The colors
play here the role of terminal and nonterminal symbols of the grammars.
The bud operad construction presented in
Section~\ref{subsubsec:bud_operads} appears in this context.
\medbreak

\SkipTocEntry\subsection*{About cyclic operads}
In a cyclic operad, the distinction between inputs and outputs of the
elements is diminished due to the fact that the cycle maps change
outputs into inputs and conversely. These structures appeared first
in~\cite{GK95}. Alternative descriptions of cyclic operads have been
provided. Among them, in~\cite{CO17}, the authors proposed an
axiomatization wherein composition maps are parametrized by two
vertices (corresponding to inputs or outputs of the elements involved in
the composition). Intuitively, this amounts to create a link between the
two vertices of the abstract operators.
\medbreak

\SkipTocEntry\subsection*{About symmetric operads}
In the literature, symmetric operads are simply called ``operads'' and
are the main considered variants among all the algebraic structures of
this sort. There are several alternative ways to define symmetric
operads, and an interesting one relies on the theory of species of
structures~\cite{Men15} (see the bibliographic notes of
Chapter~\ref{chap:collections}). Besides, some of the tools and notions
presented in Chapter~\ref{chap:operads} admit generalizations for
symmetric operads. Among these, free structures can be described by
syntax trees with labeled leaves, Poincaré-Birkhoff-Witt bases can be
described by pattern avoidance in these trees, and Koszul duality
remains a well-defined notion. General references about symmetric
operads are~\cite{Mar08,LV12,Men15}. As a matter of fact, Koszul duality
is a very important topic in the theory of symmetric operads. One of the
most beautiful properties of this duality is the fact that the
symmetric operad of commutative associative algebras is the Koszul dual
of the symmetric operad of Lie algebras. An important object here is the
operadic butterfly~\cite{Lod01,Lod06}, a diagram of symmetric operads
related by symmetric operad morphisms and links established by Koszul
duality. This diagram contains, for instance, the operads of associative
algebras, commutative associative algebras, Lie algebras, and dendriform
algebras.
\medbreak

\SkipTocEntry\subsection*{About pros}
Surprisingly, even if pros are in some sense generalizations of operads,
they appeared earlier than these lasts in the work of Mac
Lane~\cite{McL65}. Equally surprising is the fact that the
axiomatization of pros is in some way simpler than the one of operads
(there is no need to do arithmetic on the indexes to describe the
relations the horizontal and vertical composition maps have to satisfy,
contrarily to partial composition maps of operads). Basic and modern
references about pros are~\cite{Lei04} and~\cite{Mar08}. Besides, pros
are related to the Hopf bialgebra theory since in~\cite{BG16} a
construction from set-pros to Hopf bialgebras was proposed. Moreover,
in~\cite{Laf03,Laf11}, several examples of pros were provided, and links
with rewrite systems on elements of free pros were presented. As a last
remark, free pros are difficult objects to describe explicitly due
mainly to the fact that they can contain elements of arity or coarity
zero. The works~\cite{Cor18,LLMN18} collected results in this direction.
\medbreak

\backmatter

\bibliographystyle{alpha}
\bibliography{Bibliography}

\end{document}